\documentclass[reqno]{amsart}
\usepackage[letterpaper,hmargin=27mm,vmargin=27mm,bmargin=27mm,tmargin=22mm]{geometry}
\usepackage[utf8]{inputenc}
\usepackage{amsthm,amsmath,amsfonts,amssymb}
\usepackage{bm,mathrsfs,stmaryrd,enumitem}
\usepackage[colorlinks,linkcolor=blue,citecolor=blue,urlcolor=blue]{hyperref}
\usepackage{graphicx}

\newtheorem{theorem}{Theorem}[section]
\newtheorem{proposition}[theorem]{Proposition}
\newtheorem{assumption}[theorem]{Assumption}
\newtheorem{corollary}[theorem]{Corollary}
\newtheorem{lemma}[theorem]{Lemma}
\newtheorem{definition}[theorem]{Definition}
\newtheorem{remark}[theorem]{Remark}
\theoremstyle{definition}

\newcommand{\rr}{\mathrm} 
\newcommand{\f}[1]{\boldsymbol{\rr{#1}}} 
\newcommand{\bb}{\mathbb} 
\newcommand{\mb}{\mathbf} 
\newcommand{\cal}{\mathcal}
\newcommand{\scr}{\mathscr}
\newcommand{\fra}{\mathfrak}
\newcommand{\wh}{\widehat}
\newcommand{\wt}{\widetilde}
\newcommand{\txt}[1]{\text{#1}}

\newcommand{\C}{\bb C}
\newcommand{\R}{\bb R}
\DeclareMathOperator{\Tr}{Tr}
\DeclareMathOperator{\supp}{supp}
\DeclareMathOperator{\spec}{spec}
\DeclareMathOperator{\dist}{dist}
\DeclareMathOperator{\diag}{diag}
\newcommand{\dd}{\rr d}
\newcommand{\ii}{\rr i}

\newcommand{\N}{\bb N}
\newcommand{\Z}{\bb Z}

\newcommand{\bbu}{\mb u}
\newcommand{\bbh}{\mb h}
\newcommand{\bbq}{\mb q}
\newcommand{\bbv}{\mb v}
\newcommand{\bbw}{\mb w}
\newcommand{\bbx}{\mb x}
\newcommand{\bbs}{\mb s}
\newcommand{\bbe}{\mb e}
\newcommand{\bhv}{\wh\bbv}

\newcommand{\btv}{\wt\bbv}
\newcommand{\btw}{\wt\bbw}
\newcommand{\bfru}{\bm{\mathfrak{u}}}
\newcommand{\bfrv}{\bm{\mathfrak{v}}}

\newcommand{\calA}{\cal A}
\newcommand{\calB}{\cal B}
\newcommand{\calC}{\cal C}
\newcommand{\calE}{\cal E}

\newcommand{\calH}{\cal H}
\newcommand{\calI}{\cal I}
\newcommand{\calJ}{\cal J}

\newcommand{\calL}{\cal L}
\newcommand{\calM}{\cal M}

\newcommand{\calP}{\cal P}

\newcommand{\calR}{\cal R}
\newcommand{\calS}{\cal S}
\newcommand{\calT}{\cal T}
\newcommand{\calU}{\cal U}
\newcommand{\calV}{\cal V}
\newcommand{\calW}{\cal W}
\newcommand{\calX}{\cal X}

\newcommand{\scrD}{\scr D}

\newcommand{\fram}{\fra m}

\newcommand{\eps}{\varepsilon}
\newcommand{\la}{\lambda}
\newcommand{\La}{\Lambda}
\newcommand{\Om}{\Omega}
\newcommand{\te}{\theta}
\newcommand{\vphi}{\varphi}
\newcommand{\bnu}{\boldsymbol{\nu}}

\newcommand{\paren}[1]{\left(#1\right)}
\newcommand{\abs}[1]{\left|#1\right|}
\newcommand{\norm}[1]{\left\|#1\right\|}

\newcommand{\scalar}[2]{\langle #1,#2\rangle}
\newcommand{\avg}[1]{\langle #1 \rangle}

\newcommand{\Sg}{\Sigma}
\newcommand{\sg}{\sigma}
\newcommand{\cals}{\cal S}

\numberwithin{equation}{section}

\begin{document}

\begin{center}

	\begin{minipage}{0.85\textwidth}
		\vspace{2.5cm}
		
		\begin{center}
			\large\bf
			Local Laws and Edge Universality for Noncentral Sample Covariance Matrices
		\end{center}
	\end{minipage}
\end{center}

\vspace{0.5cm}

\begin{center}
	\begin{minipage}[t]{0.32\textwidth}
		\centering
		{Can Hu} \\
		{\small Northeast Normal University \\ \texttt{hucan@nenu.edu.cn}}
	\end{minipage}
	\hfill
	\begin{minipage}[t]{0.32\textwidth}
		\centering
		{Jiang Hu} \\
		{\small Northeast Normal University \\ \texttt{huj156@nenu.edu.cn}}
	\end{minipage}
	\hfill
	\begin{minipage}[t]{0.32\textwidth}
		\centering
		{Zhidong Bai} \\
		{\small Northeast Normal University \\ \texttt{baizd@nenu.edu.cn}}
	\end{minipage}
	
	 \vspace{0.6cm}

\end{center}

\vspace{0.8cm}

\begin{center}
        \begin{minipage}{0.85\textwidth}
                \small
                \textbf{Abstract.}
                We consider the real noncentral sample covariance matrices $\calW=YY^\top$ with $Y=A+\Sg^{1/2}X$. Here $A\in\R^{M\times N}$ is deterministic, $\Sg\in\R^{M\times M}$ is a deterministic positive definite population covariance matrix and $X\in\R^{M\times N}$  has independent centered entries with variance $N^{-1}$. We prove local laws near regular right edges down to optimal spectral scales without requiring the commutativity of $AA^\top$ and $\Sg$. As a consequence, we obtain optimal eigenvalue rigidity at the rightmost regular edge and delocalization of the corresponding left and right singular vectors. We also show that, with high probability, there are no eigenvalues in the adjacent spectral gap beyond the optimal $N^{-2/3}$ edge scale, up to an arbitrarily small $N^\eps$ loss. Finally, we establish edge universality at the rightmost regular edge: after centering and scaling, the largest eigenvalue converges to the Tracy--Widom distribution. The main technical ingredient is a stability analysis of the matrix Dyson equation (MDE) associated with the linearization of $Y$, whose self-energy operator does not satisfy the flatness condition of the general MDE theory. Exploiting the special block structure, we reduce the stability analysis exactly to a two-dimensional operator. This reduction yields regularity of the spectral density and square-root behavior at regular right edges, together with sharp stability bounds near such edges.
        \end{minipage}
\end{center}

\vspace{0.15cm}

{\footnotesize \textit{Keywords:} noncentral sample covariance matrices, matrix Dyson equation, local laws, edge universality.}

{\footnotesize \textit{2020 Mathematics Subject Classification:} Primary 60B20; Secondary 15B52.}

\vspace{0.5cm}

\makeatletter

\renewcommand{\@secnumfont}{\bfseries}

\renewcommand{\section}{\@startsection{section}{1}{\z@}{.7\linespacing\@plus\linespacing}{.5\linespacing}{\normalfont\bfseries\centering}}

\makeatother

\pagenumbering{arabic}

\setcounter{secnumdepth}{3}
\setcounter{tocdepth}{1}
\hypersetup{bookmarksdepth=3}

\makeatletter
\newcommand{\settocline}[4]{%
  \expandafter\def\csname l@#1\endcsname{\@tocline{#2}{0pt}{#3}{#4}{}}%
}
\settocline{section}{1}{1pc}{}
\settocline{subsection}{2}{2pc}{5pc}
\settocline{subsubsection}{3}{3pc}{7pc}
\settocline{paragraph}{4}{4pc}{9pc}
\makeatother

\makeatletter
\renewcommand{\paragraph}{\@startsection{paragraph}{4}
  {\z@}
  {.5\linespacing\@plus.7\linespacing}
  {-\fontdimen2\font}
  {\normalfont\bfseries}}
\makeatother

\tableofcontents

\section{Introduction}

Sample covariance matrices are fundamental objects in multivariate statistics and high-dimensional data analysis; see \cite{bai2010spectral,yao2015large}. Their spectral properties provide important information on the underlying covariance structure and play a central role in principal component analysis and related high-dimensional inference problems; see \cite{johnstone2001distribution,paul2007asymptotics}. In the high-dimensional regime, where the dimension and the sample size grow at comparable rates, the extreme eigenvalues are of particular interest. A central question is to understand their fluctuations near the spectral edges and the extent to which these fluctuations are universal.

In this paper, we consider the real rectangular data matrix
\begin{align}
        Y:=A+\Sg^{1/2}X, \label{eq:intro_model}
\end{align}
where $X\in\R^{M\times N}$ has independent centered entries with variance $N^{-1}$, $A\in\R^{M\times N}$ is deterministic, and $\Sg\in\R^{M\times M}$ is a deterministic positive definite population covariance matrix. We work in the high-dimensional regime where $M,N\to\infty$ and the aspect ratio $M/N$ remains bounded away from zero and infinity. The deterministic matrix $A$ represents a noncentral component and is allowed to have general rank. The precise moment assumptions on $X$ and uniform bounds on the deterministic parameters will be given in Assumptions~\ref{ass:basic} and \ref{ass:rand_ent}.

The sample covariance matrix associated with \eqref{eq:intro_model} is
\begin{align}
        \calW:=YY^\top . \nonumber
\end{align}
Our aim is to develop a local spectral theory for $\calW$ near its spectral edges without assuming that $AA^\top$ and $\Sg$ commute. Particular attention is given to the rightmost edge and the fluctuations of the largest eigenvalue.

To study the spectrum of $\calW$ on local scales, we use the linearization of the rectangular matrix $Y$. We introduce the companion Gram matrix $W:=Y^\top Y$, which has the same nonzero eigenvalues as $\calW$, and the associated Green function
\begin{align}
        G(z):=\begin{pmatrix} -I_M   & Y     \\
                Y^\top & -zI_N\end{pmatrix}^{-1},\qquad z\in\C^+. \nonumber
\end{align}
Its diagonal blocks are directly related to the resolvents of $\calW$ and $W$, so that the local spectral analysis can be formulated in terms of $G$.
The deterministic approximation to $G$ is given by the unique solution $\Pi=\Pi(z)$ to the \emph{matrix Dyson equation} (MDE) \eqref{eq:intro_mde}. 

\pdfbookmark[2]{Main difficulty}{bm:main_difficulty}
\paragraph*{Main difficulty}
Matrix Dyson equations provide a general framework for deterministic resolvent approximations in random matrix models with nontrivial matrix structure; see \cite{erdos2019matrixa,ajanki2019stability,alt2020dyson}. In the present setting, this framework does not require $AA^\top$ and $\Sg$ to commute. The main difficulty, however, is that the self-energy operator $\calS$ in \eqref{eq:self_energy} does not satisfy the flatness assumptions used in the general MDE theory \cite{ajanki2019stability,alt2020correlated,alt2020dyson}. In particular, the general MDE theory does not directly yield the required bounds on the inverse of the stability operator $\calB$ (see \eqref{eq:B_def}) near the spectral edges. Such a loss of flatness has also been identified as a genuine stability obstruction in other matrix Dyson equations \cite{alt2018local,alt2019location,alt2021inhomogeneous}. For instance, \cite{alt2021inhomogeneous} encountered a non-flat self-energy operator arising from the block structure of the Hermitized model and treated the resulting instability through a nonlinear transformation that removes the unstable direction from the MDE.

\pdfbookmark[2]{Main contributions}{bm:main_contributions}
\paragraph*{Main contributions}
To overcome the stability obstruction, we exploit the special block structure of the MDE \eqref{eq:intro_mde}. Since the self-energy operator $\calS$ depends on a matrix only through the two normalized block traces in \eqref{eq:mg}, the analysis of the stability operator $\calB:\C^{(M+N)\times(M+N)}\to\C^{(M+N)\times(M+N)}$ reduces exactly to that of the two-dimensional operator $\calI-\calJ:\C^2\to\C^2$ in Lemma~\ref{lem:stab_red}. Away from the spectral edges, we track the two eigenvalues of $\calJ$ and show that only one of them can approach $1$, while the other remains uniformly separated from $1$. This yields quantitative control of the stability operator and of the derivatives of the MDE solution. At a regular right edge, $\calI-\calJ$ has a single unstable direction and is uniformly invertible on the complementary subspace. Resolving the stable component and projecting onto the unstable direction yields the scalar normal form \eqref{eq:scalar_nf}. The macroscopic gap to the right of the edge ensures that the quadratic coefficient in \eqref{eq:scalar_nf} is uniformly nondegenerate. This analysis also yields the sharp inverse bound \eqref{eq:sharp_stab} for the stability operator of the MDE \eqref{eq:intro_mde} near the regular right edge.

Using the stability analysis developed above, we first establish the deterministic spectral theory of the matrix Dyson equation \eqref{eq:intro_mde}, whose solution defines a compactly supported deterministic probability measure on $[0,\infty)$. We prove that the density of this measure is locally $1/3$-H\"older continuous on $(0,\infty)$ and real analytic wherever it is positive; see Theorem~\ref{thm:det_spec}. At a regular right edge, we further obtain the square-root behavior of the density; see Proposition~\ref{prop:sqrt_edge}. We next prove anisotropic and averaged local laws near regular right edges down to spectral scales $\Im z\gg N^{-1}$; see Theorem~\ref{thm:opt_ll}. As a consequence, with high probability, there are no eigenvalues in the adjacent spectral gap beyond the $N^{-2/3}$ edge scale, up to an arbitrarily small $N^\eps$ loss. At the rightmost regular edge, we further obtain optimal eigenvalue rigidity and delocalization of the corresponding left and right singular vectors; see Corollary~\ref{cor:edge_rig}. Finally, we establish edge universality at the rightmost regular edge. The rescaled largest eigenvalue has the same asymptotic distribution as the rescaled largest eigenvalue of the GOE, as stated in Theorem~\ref{thm:edge_univ}. Consequently, the largest eigenvalue has Tracy--Widom fluctuations by Corollary~\ref{cor:edge_tw}.

\pdfbookmark[2]{Related works}{bm:related_works}
\paragraph*{Related works}
Two benchmark regimes of \eqref{eq:intro_model} have been studied extensively. When $A=0$, one recovers the centered sample covariance model with a general population covariance matrix. Its global spectral behavior is described by the deformed Mar\v{c}enko--Pastur law \cite{marcenko1967distributiona,silverstein1995analysis}. For the standard sample covariance model, the almost sure limits of the largest and smallest eigenvalues were established in \cite{yin1988limit,bai1993limit}. Anisotropic local laws and edge universality, including Tracy--Widom fluctuations, have been established in \cite{pillai2014universality,bao2015universality,lee2016tracy,knowles2017anisotropic,ding2018necessary,fan2022tracy}. When $\Sg=I_M$, the data matrix takes the additively deformed form $Y=A+X$. The associated information-plus-noise models have an extensive deterministic spectral theory \cite{dozier2007empirical,dozier2007analysis}. Subspace estimation for deterministic signals in this high-dimensional setting was studied in \cite{vallet2012improved}, while edge statistics have more recently been investigated in \cite{ding2022edge,zhang2024tracywidom}.

Several related extensions have also been studied. One such extension allows more general variance profiles in the noise, for which deterministic equivalents and bilinear resolvent estimates for information-plus-noise matrices were developed in \cite{hachem2007deterministic,hachem2013bilinear}. Optimal entrywise and averaged local laws were established for random Gram matrices with general variance profiles in \cite{alt2017local,alt2017singularities}, and Tracy--Widom fluctuations at regular edges were obtained in \cite{ding2022tracywidom}. Edge universality for separable covariance matrices was proved in \cite{yang2019edge}, while \cite{fan2026anisotropic} recently established an anisotropic local law for a class of non-separable sample covariance matrices. Another important special regime is the finite-rank case of $A$. For white noise, outlier singular values and singular vector statistics have been studied extensively in \cite{benaych-georges2012singular,ding2020high,bao2021singular}. More recent work has considered low-rank signals in the presence of a general population covariance matrix \cite{lin2024asymptotic,liu2025asymptotica,bao2025signal}.

The analysis becomes more involved when both $A$ and $\Sg$ are general, since $AA^\top$ and $\Sg$ need not admit a common eigenbasis. By contrast, under the additional commutativity condition $AA^\top\Sg=\Sg AA^\top$, the limiting spectral distribution and its support were studied in \cite{zhou2023limitinga,zhou2024analysis}. It was further shown in \cite{bai2025no,zhou2025exact} that there are no eigenvalues outside the limiting support and that exact separation holds. This condition permits simultaneous diagonalization and reduces the associated self-consistent relations to scalar equations. More recently, \cite{zhuang2026no} proved that there are no eigenvalues outside the limiting support for generally correlated and noncentral covariance models beyond this scalar reduction. These works, however, do not provide a local spectral theory for \eqref{eq:intro_model} near its spectral edges without assuming commutativity between $AA^\top$ and $\Sg$.

Beyond the works discussed above, an important development in the local spectral theory of general random matrices is the theory of vector and matrix Dyson equations. This framework extends the scalar self-consistent equations of classical ensembles to more general random matrix models. For Wigner-type matrices with inhomogeneous variance profiles, a systematic theory of the vector Dyson equation and the corresponding local spectral analysis was developed in \cite{ajanki2017singularities,ajanki2017universality}. This framework was subsequently extended to correlated random matrices through the matrix Dyson equation. Its stability theory and local spectral analysis were developed in \cite{ajanki2019stability,alt2020dyson,alt2020correlated}. For detailed treatments and broader perspectives on these developments, we refer to \cite{ajanki2019quadratic,erdos2019matrixa}.

\pdfbookmark[2]{Strategy of proof}{bm:strategy_of_proof}
\paragraph*{Strategy of proof}
We follow the standard \emph{three-step strategy} for edge universality, consisting of a local law, universality for a Gaussian divisible model and a comparison with the original ensemble; see \cite{erdhos2017dynamical,erdos2019matrixa}. The first step contains the main model-dependent analysis. In the present setting, the local law is obtained by combining the stochastic estimate for the perturbed MDE in Theorem~\ref{thm:rand_sc_eq} with the deterministic stability analysis developed above. Under an a priori control of $G-\Pi$, the anisotropic and averaged error estimates are given by \eqref{eq:anis_rand_error} and \eqref{eq:avg_rand_error}, respectively. The proof of these estimates follows the approach of \cite{he2018isotropic,erdos2019random,bao2024random}. The deterministic input is the sharp stability estimate \eqref{eq:sharp_stab} for the MDE \eqref{eq:intro_mde}. A continuity argument in $\Im z$, combined with a bootstrap scheme, upgrades the a priori estimates and extends the local law from large imaginary parts down to the optimal local scale.

For the second step, we construct a covariance-preserving Gaussian divisible model. We reduce the covariance to $\Sg_t:=\Sg-tI_M$ and set
\begin{align}
        Y_0:=A+\Sg_t^{1/2}X,\qquad Y_t:=Y_0+\sqrt{t}\,X^{\rr G}. \nonumber
\end{align}
Since $\Sg_t+tI_M=\Sg$, the Gaussian component restores exactly the original covariance, and the deterministic MDE for $Y_t$ is the same as \eqref{eq:intro_mde}. The local laws and rigidity estimates provide the spectral input needed to verify the $\eta_*$-regularity of the initial matrix $Y_0$; see Proposition~\ref{prop:eta_star_reg}. We then identify the edge location and scaling of the Gaussian divisible model in Proposition~\ref{prop:cond_edge_par} and apply the rectangular Gaussian divisible edge results of \cite{ding2022edge,ding2022tracywidom}.

For the third step, we remove the Gaussian component by a short-time Green function comparison. We interpolate through
\begin{align}
        Y(s):=A+(\Sg-sI_M)^{1/2}X+\sqrt{s}\,X^{\rr G}, \qquad 0\leq s\leq t. \nonumber
\end{align}
The covariance remains equal to $\Sg$ along the interpolation. The covariance-preserving interpolation yields an exact cancellation of the second-order terms between the non-Gaussian and Gaussian parts \eqref{eq:comp_second}. The leading remaining contribution comes from the third cumulant.  \eqref{eq:comp_cum_rem} gives the bound $N^{1/6+C\chi}$ at the level of the comparison derivative. We choose $t=N^{-1/3+\omega}$ with $0<\omega<1/12$, so that, after integration over $0\leq s\leq t$, the third-cumulant contribution is bounded by $N^{-1/6+\omega+C\chi}$ and hence tends to zero for sufficiently small $\chi$. The fourth and higher cumulant terms are smaller. Proposition~\ref{prop:short_comp} therefore transfers the edge statistics from the Gaussian divisible model back to the original matrix; see \cite{huang2020transition,alt2020correlated,huang2026edge}.

\pdfbookmark[2]{Organization}{bm:organization}
\paragraph*{Organization}
The paper is organized as follows.
Section~\ref{sec:main_results} states the main results.
Section~\ref{sec:det_mde} develops the deterministic theory of the MDE \eqref{eq:intro_mde} and the associated measure.
Section~\ref{sec:reg_edge} studies the MDE \eqref{eq:intro_mde} near a regular right edge.
It establishes the scalar normal form, the square-root behavior of the density, and the sharp stability estimate.
Section~\ref{sec:ll} proves the local laws and their spectral consequences.
Section~\ref{sec:edge_univ} proves the edge universality theorem by constructing a Gaussian divisible model that preserves the covariance and using a short-time comparison.
The appendices collect auxiliary estimates for the deterministic MDE analysis and the technical estimates used in the local law and edge universality proofs.

\pdfbookmark[2]{Notation and conventions}{bm:notation_conventions}
\paragraph*{Notation and conventions}
Throughout the paper, $N$ is the fundamental large parameter.
All asymptotic statements are understood for sufficiently large $N$.
Unless stated otherwise, all constants are independent of $N$.
Their values may change from line to line.
We set $L:=M+N$.
For integers $a\leq b$, we write $\llbracket a,b\rrbracket:=[a,b]\cap\Z$.
We use Latin indices $i,j\in\llbracket1,M\rrbracket$ and Greek indices $\mu,\nu\in\llbracket M+1,L\rrbracket$.
Unless stated otherwise, $a,b\in\llbracket1,L\rrbracket$.
The columns of $X$ and the canonical basis of the second block are indexed by $\llbracket M+1,L\rrbracket$.
The ambient space of a canonical basis vector $\bbe_i$ or $\bbe_\mu$ is determined by the context.
For nonnegative quantities $x$ and $y$, we write $x\lesssim y$ if $x\leq Cy$ for some constant $C$.
We write $x\sim y$ if $x\lesssim y$ and $y\lesssim x$.
The notation $x\gtrsim y$ is defined analogously.
We also write $x=O(y)$ if $x\lesssim y$.
For positive semidefinite matrices, inequalities involving $\lesssim$ and $\gtrsim$ are understood in the sense of quadratic forms.
For vectors, $\|\cdot\|$ denotes the Euclidean norm, while for matrices it denotes the operator norm.
For $R\in\C^{k\times k}$, we set $\|R\|_{\rr{hs}}:=(k^{-1}\Tr R^*R)^{1/2}$ and $\Im R = (R-R^*)/(2\ii)$.
Moreover, for a square matrix $R$, we set $\avg{R}_k:=k^{-1}\Tr R$. The subscript $k$ specifies the normalization and may not equal the dimension of $R$.
If $\calT$ is a linear operator acting on matrices, then $\|\calT\|:=\sup_{\|R\|=1}\|\calT[R]\|$ denotes the operator norm induced by the matrix operator norm.
The symbol $\calI$ denotes the identity operator on the space under consideration.
Vectors are written in boldface, and $I_k$ denotes the $k\times k$ identity matrix.
We write $z=E+\ii\eta\in\C^+$ with $\eta>0$.

\section{Main results}\label{sec:main_results}

\subsection{Model and assumptions}
We consider the data matrix
\begin{align}
        Y=A+\Sg^{1/2}X, \nonumber
\end{align}
where $A\in\R^{M\times N}$ is deterministic, $\Sg\in\R^{M\times M}$ is deterministic positive definite, and $X\in\R^{M\times N}$ is random. The associated noncentral sample covariance matrix is
\begin{align}
        \calW=YY^\top. \nonumber
\end{align}
Its companion Gram matrix is
\begin{align}
        W=Y^\top Y. \nonumber
\end{align}
The matrices $\calW$ and $W$ have the same nonzero eigenvalues. For our spectral analysis, it is more convenient to work with $W$. The corresponding results for the nonzero eigenvalues of $\calW$ then follow immediately.

We impose the following standard high-dimensional assumptions.

\begin{assumption}[High-dimensional regime]\label{ass:basic}
        Set $c_N:=M/N$. There exist constants $0<c_0\leq C_0<\infty$ such that
        \begin{align}
                c_0\leq c_N\leq C_0, \qquad \|\Sg\|+\|\Sg^{-1}\|+\|AA^\top\|\leq C_0. \nonumber
        \end{align}
\end{assumption}

\begin{assumption}[Random entries]\label{ass:rand_ent}
        Let $X=(X_{i\mu})\in\R^{M\times N}$ be a real random matrix with independent entries satisfying
        \begin{align}
                \bb E X_{i\mu}=0, \qquad \bb E X_{i\mu}^2=\frac1N \nonumber
        \end{align}
        for all $i,\mu$. Moreover, for every fixed $p\in\N$, there exists a constant $C_p>0$ such that
        \begin{align}
                \max_{i,\mu}\bb E\abs{\sqrt N X_{i\mu}}^p\leq C_p. \nonumber
        \end{align}
\end{assumption}

\subsection{Deterministic approximation and regular spectral edges}
We introduce the symmetric linearizations of the random and deterministic parts of $Y$,
\begin{align}
        \calH:=\begin{pmatrix}0               & \Sg^{1/2}X \\
               X^\top\Sg^{1/2} & 0\end{pmatrix},\qquad \La:=\begin{pmatrix}0      & A \\
               A^\top & 0\end{pmatrix}. \nonumber
\end{align}
For $z\in\C^+$, set
\begin{align}
        \wh I:=\begin{pmatrix}I_M & 0 \\
               0   & 0\end{pmatrix},\qquad \wt I:=\begin{pmatrix}0 & 0   \\
               0 & I_N\end{pmatrix},\qquad I_z:=\wh I+z\wt I. \nonumber
\end{align}
The Green function of the linearized model is defined by
\begin{align}
        G(z):=(\calH+\La-I_z)^{-1} = \begin{pmatrix} -I_M     & Y     \\
                Y^{\top} & -zI_N\end{pmatrix}^{-1}. \nonumber
\end{align}
By the Schur complement formula, the two diagonal blocks are
\begin{align}
        \wh G(z):=z(\calW-zI_M)^{-1},\qquad \wt G(z):=(W-zI_N)^{-1}. \label{eq:G_schur}
\end{align}
Thus
\begin{align}
        G(z)=\begin{pmatrix}\wh G(z) & *        \\
               *        & \wt G(z)\end{pmatrix}, \nonumber
\end{align}
and $\wh G$ and $\wt G$ contain the resolvents of $\calW$ and $W$, respectively.

For $R\in\C^{L\times L}$, define the self-energy operator $\cals:\C^{L\times L}\to\C^{L\times L}$ by
\begin{align}
        \cals[R]:=\begin{pmatrix} \avg{R_2}_N\Sg & 0                  \\
                0              & \avg{R_1\Sg}_N I_N\end{pmatrix}, \label{eq:self_energy}
\end{align}
where $R_1\in\C^{M\times M}$ and $R_2\in\C^{N\times N}$ denote the two diagonal blocks of $R$. The operator $\cals$ is linear and positivity preserving.

To describe the deterministic approximation of $G(z)$, define
\begin{align}
        \Om(R,z):=I_L+I_zR+\cals[R]R-\La R, \qquad R\in\C^{L\times L}. \nonumber
\end{align}
The deterministic approximation $\Pi=\Pi(z)$ is characterized by the \emph{matrix Dyson equation (MDE)}
\begin{align}
        \Om(\Pi(z),z)=0. \label{eq:intro_mde}
\end{align}
For its solution $\Pi=\Pi(z)$, we use the block decomposition
\begin{align}
        \Pi(z):=\begin{pmatrix} \Pi_{11}(z) & \Pi_{12}(z) \\
                \Pi_{21}(z) & \Pi_{22}(z)\end{pmatrix}, \nonumber
\end{align}
where $\Pi_{11}(z)\in\C^{M\times M}$ and $\Pi_{22}(z)\in\C^{N\times N}$. We further define the two normalized block traces
\begin{align}
        m(z):=\avg{\Pi_{22}(z)}_N,\qquad g(z):=\avg{\Sg\Pi_{11}(z)}_N. \label{eq:mg}
\end{align}
Then
\begin{align}
        \cals[\Pi(z)] =\begin{pmatrix} m(z)\Sg & 0       \\
                0       & g(z)I_N\end{pmatrix}. \nonumber
\end{align}
Thus the self-energy term in the MDE \eqref{eq:intro_mde} depends on $\Pi$ only through the two scalar quantities $m$ and $g$.

The next theorem collects the deterministic spectral properties needed below.

\begin{theorem}\label{thm:det_spec}
        Under Assumption~\ref{ass:basic}, the following statements hold.
        \begin{enumerate}[label=(\roman*)]
                \item For every $z\in\C^+$, the MDE \eqref{eq:intro_mde} has a unique solution $\Pi$ satisfying $\Im\Pi(z)>0$ and $\Im(z^{-1}\Pi_{11}(z))>0$. The map $z\mapsto\Pi(z)$ is holomorphic on $\C^+$.
                      There exists a unique probability measure $\bnu$ such that
                      \begin{align}
                              m(z)=\int_{[0,\infty)}\frac{\bnu(\dd x)}{x-z}, \qquad z\in\C^+. \nonumber
                      \end{align}
                      The measure $\bnu$ has support contained in $[0,C_*]$ for some constant $C_*>0$. Moreover, $\sup\supp\bnu\sim1$.
                \item The measure $\bnu$ is absolutely continuous with respect to Lebesgue measure on $(0,\infty)$ and $\bnu(\{0\})=(1-M/N)\mb{1}(M/N\leq1)$.
                      More precisely,
                      \begin{align}
                              \bnu(\dd E)=\bnu(\{0\})\delta_0(\dd E)+\rho(E)\mb 1(E>0)\dd E, \nonumber
                      \end{align}
                      where
                      \begin{align}
                              \rho(E)=\frac{1}{\pi}\lim_{\eta\downarrow0}\Im m(E+\ii\eta). \nonumber
                      \end{align}
                      The density $\rho$ is locally $1/3$-H\"older continuous on $(0,\infty)$ and real analytic on the open set $\{E>0:\rho(E)>0\}$.
                \item The functions $m$ and $g$ defined in \eqref{eq:mg} extend uniquely to locally $1/3$-H\"older continuous functions on $\{E+\ii\eta:E>0,\ \eta\geq0\}$.
        \end{enumerate}
\end{theorem}

The measure $\bnu$ will be shown to provide the deterministic approximation to the eigenvalue distribution of $W$. Since $W$ and $\calW$ have the same nonzero eigenvalues, the edges determined by $\bnu$ are used for both matrices.

By Theorem~\ref{thm:det_spec}, for $E>0$ we use the notation
\begin{align}
        m(E):=\lim_{\eta\downarrow0}m(E+\ii\eta), \qquad g(E):=\lim_{\eta\downarrow0}g(E+\ii\eta). \nonumber
\end{align}

We denote the rightmost edge by $E_+:=\sup\supp\bnu = \sup\{E>0:\rho(E)>0\}$. By Theorem~\ref{thm:det_spec}(i), $E_+\sim1$. Fix $\tau_0>0$ sufficiently small such that $E_+\in[\tau_0,\tau_0^{-1}]$ and $\tau_0^{-1}\geq C_*+1$, and define
\begin{align}
        \mb D_0:=\{E\in[\tau_0,\tau_0^{-1}]:E\in\partial\supp\bnu,\ (E,E+\tau_0)\cap\supp\bnu=\varnothing\}. \nonumber
\end{align}
Then $E_+\in\mb D_0$.

\begin{definition}[Regular right edge]\label{def:reg_r_edge}
        We call $E_*\in\mb D_0$ a regular right edge of $\bnu$ if
        \begin{align}
                s_{\min}\left(AA^\top-(E_*+g(E_*))(I_M+m(E_*)\Sg)\right) & \geq\tau_0. \label{eq:schur_nondeg}
        \end{align}
        $s_{\min}(\cdot)$ denotes the smallest singular value.
        A regular left edge is defined analogously, with the adjacent gap lying to the left of $E_*$.
\end{definition}

\begin{remark}
        Condition~\eqref{eq:schur_nondeg} provides the additional nondegeneracy needed at a regular edge. As shown in Lemma~\ref{lem:bdry_zero_dens}, it implies $|z+g(z)|\gtrsim1$ in a neighborhood of $E_*$. Together with the corresponding lower bound for the Schur complement, this yields $\|\Pi(z)\|\lesssim1$ on the same domain. The latter bound is used below in the analysis of the stability operator. In the spectral bulk, such boundedness follows directly from Lemma~\ref{lem:Pi_bd}, since $\|\Pi(z)\|\lesssim(\Im m(z))^{-1}$ and $\Im m(z)\gtrsim1$. At an edge, this argument no longer applies since $\Im m(z)$ vanishes as $z$ approaches the edge. In \cite[Assumption~(G) and Remark~2.3]{alt2020correlated}, a uniform bound on the MDE solution near an edge is imposed directly. For the present model, this bound follows as a consequence of \eqref{eq:schur_nondeg}. Together with the gap to the right in the definition of $\mb D_0$, the resulting stability analysis yields the square-root behavior of the density $\rho$ near $E_*$ (see Proposition~\ref{prop:sqrt_edge} below).  More generally, for the stability analysis below, Condition~\eqref{eq:schur_nondeg} is only used through this local boundedness of $\Pi$. Thus, it may be replaced by the assumption that $\Pi(z)$ is uniformly bounded in a fixed neighborhood of $E_*$.
\end{remark}

\begin{remark}
        Condition~\eqref{eq:schur_nondeg} recovers the standard regularity conditions in two benchmark cases.

        If $A=0$, then $m(z)=-(z+g(z))^{-1}$, and \eqref{eq:schur_nondeg} reduces to $ s_{\min}(\Sg+m(E_*)^{-1}I_M)\gtrsim 1$. By Lemma~\ref{lem:Pi_bd} and Lemma~\ref{lem:bdry_zero_dens}, $|m(E_*)|\sim1$. Hence \eqref{eq:schur_nondeg} is equivalent to $\min_i\left|m(E_*)+\sg_i^{-1}\right|\gtrsim1$, where $\sg_i$ are the eigenvalues of $\Sg$. This is precisely the  regular edge condition in \cite[Definition~2.7(i)]{knowles2017anisotropic} or \cite[Assumption~2.5]{ding2018necessary}.

        If $\Sg=I_M$, then \eqref{eq:schur_nondeg} becomes $s_{\min}(AA^\top-(E_*+g(E_*))(1+m(E_*))I_M)\geq\tau_0$. Thus $(E_*+g(E_*))(1+m(E_*))$ is required to stay uniformly separated from $\operatorname{Spec}(AA^\top)$. At the rightmost edge, this agrees with the regularity condition in \cite{zhang2024tracywidom}. Indeed, their function $w$ defined in \cite[(2.9)]{zhang2024tracywidom} satisfies $w(z)=(z+g(z))(1+m(z))$, and their parameter at the rightmost edge $\xi_r$ satisfies $w(\la_r)=\xi_r$. Hence, under the identification $E_*=\la_r$, our scalar parameter equals $\xi_r$. In that paper, $d_1$ denotes the largest singular value of the signal matrix, so $d_1^2=\la_{\max}(AA^\top)$. The parameter $\xi_r$ is defined as the largest positive critical point associated with the rightmost edge, and the signal matrix is required to satisfy $d_1^2<\xi_r$. \cite[Assumption~3]{zhang2024tracywidom} strengthens this separation to $\xi_r-d_1^2\gtrsim1$. This gives $s_{\min}(AA^\top-\xi_r I_M)=\xi_r-d_1^2\gtrsim1$, which is precisely \eqref{eq:schur_nondeg} at the rightmost edge.
\end{remark}

\begin{proposition}
        \label{prop:sqrt_edge}
        Under Assumption~\ref{ass:basic}, let $E_*\in\mb D_0$ be a regular right edge in the sense of Definition~\ref{def:reg_r_edge}. Then there exists $\eps_*>0$ such that
        \begin{align}
                \rho(E)=c_{\rm edge}(E_*-E)^{1/2}+O(E_*-E),\qquad E\in[E_*-\eps_*,E_*], \label{eq:sqrt_edge}
        \end{align}
        where $c_{\rm edge}$ is given by \eqref{eq:dens_sqrt} and satisfies $c_{\rm edge}\sim1$.
\end{proposition}

\subsection{Local laws}
To state the local laws, we introduce the notion of stochastic domination.

\begin{definition}[Stochastic domination]
        Let
        \begin{align}
                \xi=(\xi^{(N)}(u):N\in\N,\ u\in U^{(N)}),
                \qquad
                \zeta=(\zeta^{(N)}(u):N\in\N,\ u\in U^{(N)}) \nonumber
        \end{align}
        be two families of nonnegative random variables, where $U^{(N)}$ is a possibly $N$-dependent parameter set.
        We say that $\xi$ is stochastically dominated by $\zeta$, uniformly in $u$, if for all (small) $\eps>0$ and (large) $D>0$, we have
        \begin{align}
                \sup_{u\in U^{(N)}} \bb P\left(\xi^{(N)}(u)>N^\eps\zeta^{(N)}(u)\right) \leq N^{-D} \nonumber
        \end{align}
        for large $N > N_0(\eps,D)$. We write  $\xi\prec\zeta$ when $\xi$ is stochastically dominated by $\zeta$ uniformly in $u$, and $\xi=O_\prec(\zeta)$ denote $\abs{\xi}\prec\zeta$. Note that in the special case when $\xi$ and $\zeta$ are deterministic, $\xi=O_\prec(\zeta)$  means for any given $\eps>0$,  $|\xi_N(u)|\leq N^\eps\zeta_N(u)$ uniformly in $u$, for all sufficiently large $N\geq N_0(\eps)$.  Let $R$ be a family of complex square matrices and $\zeta$ a family of nonnegative random variables. We write $R=O_\prec(\zeta)$ to mean $\abs{\scalar{\bbv}{R\bbw}}\prec\zeta|\bbv||\bbw|$ uniformly for all deterministic vectors $\bbv$ and $\bbw$.
      
        In addition, we say that an event $\Xi $ holds \emph{with high probability} if
        \begin{align}
                \bb P(\Xi)\geq 1-N^{-D} \nonumber
        \end{align}
        for any large constant $D>0$ when $N$ is sufficiently large.
\end{definition}

Let $E_*$ be a regular right edge and choose $0<\eps_*<1$ as in Proposition~\ref{prop:sqrt_edge}. For any sufficiently small $\tau>0$, define the spectral domain
\begin{align}
        \mb D_*^{\eps_*}(\tau):=\{z=E+\ii\eta\in\C^+:|E-E_*|\leq\eps_*, N^{-1+\tau}\leq\eta\leq1\}. \label{eq:edge_spec_dom}
\end{align}

\begin{theorem}\label{thm:opt_ll}
        Under Assumptions~\ref{ass:basic} and \ref{ass:rand_ent}, let $E_*$ be a regular right edge in the sense of Definition~\ref{def:reg_r_edge}. Then, for any sufficiently small $\tau>0$, uniformly for $z=E+\ii\eta\in\mb D_*^{\eps_*}(\tau)$,
        \begin{align}
                G(z)-\Pi(z) = O_{\prec}\left(\sqrt{\frac{\Im m(z)}{N\eta}}+\frac1{N\eta}\right). \label{eq:opt_anis_ll}
        \end{align}
        Moreover, for every deterministic matrix $B\in\C^{L\times L}$ with $\|B\|\leq1$,
        \begin{align}
                \avg{B(G(z)-\Pi(z))}_L = O_{\prec}\left(\frac1{N\eta}\right). \label{eq:opt_avg_ll}
        \end{align}
        On the gap side, writing $\kappa:=E-E_*\geq0$, we further have
        \begin{align}
                \avg{B(G(z)-\Pi(z))}_L = O_{\prec}\left(\frac1{N(\kappa+\eta)}+\frac1{(N\eta)^2\sqrt{\kappa+\eta}}\right). \label{eq:gap_avg_ll}
        \end{align}
        If $E_*=E_+$, then \eqref{eq:opt_anis_ll}--\eqref{eq:gap_avg_ll} hold uniformly for $\{E_+-\eps_*\leq E\leq\tau^{-1}, N^{-1+\tau}\leq\eta\leq1\}$, with the gap side estimate \eqref{eq:gap_avg_ll} restricted to $E\geq E_+$.
\end{theorem}

We denote the common nontrivial eigenvalues of $W$ and $\calW$ by
\begin{align}
        \la_1\geq\la_2\geq\dots\geq\la_{M\wedge N}. \nonumber
\end{align}
At the rightmost edge $E_+$, define the classical locations by
\begin{align}
        \gamma_j:=\sup\left\{x\in\R:\int_x^\infty\rho(y)\dd y>\frac{j-1}{N}\right\},\qquad j\geq1. \nonumber
\end{align}
Thus $\gamma_1=E_+$.

\begin{corollary}\label{cor:edge_rig}
        Under the assumptions of Theorem~\ref{thm:opt_ll}, let $E_*$ be a regular right edge in the sense of Definition~\ref{def:reg_r_edge}, the following statements hold.
        \begin{enumerate}[label=(\roman*)]
                \item For every $\eps>0$ and $D>0$,
                      \begin{align}
                              \bb P\left(\spec W\cap[E_*+N^{-2/3+\eps},E_*+\eps_*]\neq\varnothing\right)\leq N^{-D}. \label{eq:no_eig_gap}
                      \end{align}
                      If $E_*=E_+$, then
                      \begin{align}
                              \bb P(\la_1>E_++N^{-2/3+\eps})\leq N^{-D}. \label{eq:lambda1_ub}
                      \end{align}
                \item If $E_*=E_+$,  then for any $j$ satisfying $E_+-\eps_*/2\leq \gamma_j\leq E_+ $,
                      \begin{align}
                              \abs{\la_j-\gamma_j}\prec N^{-2/3}j^{-1/3}. \label{eq:r_edge_rig}
                      \end{align}
                      In particular, for every fixed $K\in\N$,
                      \begin{align}
                              \max_{1\leq j\leq K}\abs{\la_j-E_+}\prec N^{-2/3}. \label{eq:lambda1_st_rig}
                      \end{align}
                \item For each nonzero eigenvalue $\la_j$, let $\bfru_j\in\R^M$ and $\bfrv_j\in\R^N$ be normalized left and right singular vectors of $Y$, respectively. If $E_*=E_+$, then for every fixed $K\in\N$, uniformly for $1\leq j\leq K$ and deterministic vectors $\bhv\in\C^M$, $\btw\in\C^N$,
                      \begin{align}
                              \abs{\scalar{\bhv}{\bfru_j}}^2\prec N^{-1}\|\bhv\|^2,\qquad \abs{\scalar{\btw}{\bfrv_j}}^2\prec N^{-1}\|\btw\|^2. \label{eq:r_edge_deloc}
                      \end{align}
        \end{enumerate}
\end{corollary}

\subsection{Edge universality}
Suppose that the rightmost edge $E_+$ is regular, and let $c_{\rm edge}$ be the square-root coefficient in \eqref{eq:sqrt_edge} at $E_+$. We define the deterministic edge scaling parameter by
\begin{align}
        \gamma_+ := (\pi c_{\rm edge})^{2/3}. \label{eq:edge_scale}
\end{align}
By Proposition~\ref{prop:sqrt_edge}, $\gamma_+\sim1$.

\begin{theorem}\label{thm:edge_univ}
        Under Assumptions~\ref{ass:basic} and~\ref{ass:rand_ent}, suppose that $E_+$ is a regular right edge in the sense of Definition~\ref{def:reg_r_edge}. Let $\mu_1^{\rr{GOE}}$ denote the largest eigenvalue of a standard $N\times N$ GOE matrix. Then, for test function $F:\R\to\R$ satisfying $ \|F\|_\infty+\|F'\|_\infty\leq C$,
        we have
        \begin{align}
                \lim_{N\to\infty}\abs{\bb E F\left(\gamma_+N^{2/3}(\la_1-E_+)\right)-\bb E F\left(N^{2/3}(\mu_1^{\rr{GOE}}-2)\right)}=0.
                \label{eq:r_edge_univ}
        \end{align}
\end{theorem}

Theorem~\ref{thm:edge_univ} immediately yields the following corollary.
\begin{corollary}\label{cor:edge_tw}
        Under Assumptions~\ref{ass:basic} and~\ref{ass:rand_ent}, suppose that $E_+$ is a regular right edge in the sense of Definition~\ref{def:reg_r_edge}. Then, for every $x\in\R$,
        \begin{align}
                \lim_{N\to\infty} \bb P( \gamma_+N^{2/3}(\la_1-E_+)\leq x ) =F_1(x). \nonumber
        \end{align}
        Here $F_1$ denotes the distribution function of the Tracy--Widom law of type~1.
\end{corollary}

\begin{remark}
        For every fixed $k\in\N$, Theorem~\ref{thm:edge_univ} extends naturally to the joint edge statistics of the top $k$ eigenvalues. The corresponding result for the Gaussian divisible model is given in \cite[Theorem~4]{ding2022tracywidom}.
\end{remark}

\begin{remark}
        Suppose that band rigidity for $\calW$ is established in the sense of \cite{alt2020correlated}. Then Theorem~\ref{thm:edge_univ} can be extended to every regular right edge $E_*$ in the sense of Definition~\ref{def:reg_r_edge}. Band rigidity determines the eigenvalue labels associated with the spectral bands. This is the mechanism used to treat internal regular edges in \cite{alt2020correlated}.
\end{remark}

\section{Deterministic properties of the matrix Dyson equation}\label{sec:det_mde}

Throughout this section, Assumption~\ref{ass:basic} holds. Theorem~\ref{thm:det_spec} follows from Proposition~\ref{prop:pi}, Theorem~\ref{thm:dens_reg}, and Corollary~\ref{cor:bdry_reg_m_g}, which are stated in Subsection~\ref{subsec:det_spec}; its proof is given at the end of that subsection. Proposition~\ref{prop:pi} is proved in Subsection~\ref{sub:st_rep}, while Theorem~\ref{thm:dens_reg} and Corollary~\ref{cor:bdry_reg_m_g} are proved in Subsection~\ref{subsec:dens_bdry} using the two-dimensional stability estimates established in Subsection~\ref{subsec:2d_stab}.

\subsection{Deterministic measure and its regularity}\label{subsec:det_spec}

Theorem~\ref{thm:det_spec} follows from the following three results, which we state first.

\begin{proposition}\label{prop:pi}
        Under Assumption~\ref{ass:basic}, the following statements hold.
        \begin{enumerate}[label=(\roman*)]
                \item There exists a unique holomorphic function $\Pi:\C^+\to\C^{L\times L}$ such that, for every $z\in\C^+$, $\Om(\Pi(z),z)=0$, $\Im\Pi(z)>0$, and $\Im(z^{-1}\Pi_{11}(z))>0$.
                \item Define $\wt T(z):=z^{-1}\Pi_{11}(z)$ and $T(z):=\Pi_{22}(z)$. There exist unique positive semidefinite matrix-valued measures $V_1$ and $V_2$ on $[0,\infty)$, with values in $\C^{M\times M}$ and $\C^{N\times N}$, respectively, such that $V_1([0,\infty))=I_M$, $V_2([0,\infty))=I_N$, and
                      \begin{align}
                              \wt T(z)=\int_{[0,\infty)}\frac{V_1(\dd\la)}{\la-z},\qquad T(z)=\int_{[0,\infty)}\frac{V_2(\dd\la)}{\la-z}. \nonumber
                      \end{align}
                      Moreover, there exists a constant $C_*<\infty$ such that $\supp V_1\cup\supp V_2\subset[0,C_*]$.
                \item Let $\bnu_1:=\avg{V_1}_M$ and $\bnu_2:=\avg{V_2}_N$. Then $\bnu_1$ and $\bnu_2$ are compactly supported probability measures on $[0,\infty)$ and, with $\wt m(z):=\avg{\wt T(z)}_M$ and $m(z):=\avg{T(z)}_N$,
                      \begin{align}
                              \wt m(z)=\int_{[0,\infty)}\frac{\bnu_1(\dd\la)}{\la-z},\qquad m(z)=\int_{[0,\infty)}\frac{\bnu_2(\dd\la)}{\la-z}. \nonumber
                      \end{align}
                      Moreover,
                      \begin{align}
                              g(z)=-\avg{\Sg}_N+\int_{[0,\infty)}\frac{\la\avg{\Sg V_1(\dd\la)}_N}{\la-z}. \nonumber
                      \end{align}
                \item Writing $\bnu:=\bnu_2$, we have $\sup\supp\bnu\sim1$. Moreover, the mass at the origin is
                      \begin{align}
                              \bnu(\{0\})=\left(1-\frac{M}{N}\right)\mb 1\left(\frac{M}{N}\leq1\right). \label{eq:origin_mass}
                      \end{align}
        \end{enumerate}
\end{proposition}

Throughout the remainder of this section, we write $\bnu:=\bnu_2$ for the probability measure associated with $m$.

\begin{theorem}\label{thm:dens_reg}
        Under Assumption~\ref{ass:basic}, the measure $\bnu$ admits the decomposition
        \begin{align}
                \bnu(\dd E) = \bnu(\{0\})\delta_0(\dd E) + \rho(E)\mb 1(E>0)\dd E, \label{eq:dens_dec}
        \end{align}
        where $\rho$ is a nonnegative function on $(0,\infty)$ given by
        \begin{align}
                \rho(E) = \frac1\pi\lim_{\eta\downarrow0}\Im m(E+\ii\eta). \nonumber
        \end{align}
        Moreover, $\rho$ is locally $1/3$-H\"older continuous on $(0,\infty)$.
        Finally, $\rho$ is real analytic on the open set $\{E\in(0,\infty):\rho(E)>0\}$.
\end{theorem}

The regularity of the density also yields boundary regularity of the scalar functions entering the MDE \eqref{eq:intro_mde}.
\begin{corollary}
        \label{cor:bdry_reg_m_g}
        Under Assumption~\ref{ass:basic}, the functions $m$ and $g$ admit unique extensions to $\{E+\ii\eta:E>0,\ \eta\geq0\}$ that are locally $1/3$-H\"older continuous.
\end{corollary}

\begin{proof}[Proof of Theorem~\ref{thm:det_spec}]
        Theorem~\ref{thm:det_spec} is a direct consequence of Proposition~\ref{prop:pi}, Theorem~\ref{thm:dens_reg}, and Corollary~\ref{cor:bdry_reg_m_g}.
\end{proof}

We now prove Proposition~\ref{prop:pi}, Theorem~\ref{thm:dens_reg}, and Corollary~\ref{cor:bdry_reg_m_g} in turn.

\subsection{Stieltjes representation and basic MDE properties}
\label{sub:st_rep}

This subsection proves Proposition~\ref{prop:pi}. Its main input is Lemma~\ref{lem:lin_mde_meas}, whose proof is given in Appendix~\ref{app:basic_mde}. The rectangular MDE \eqref{eq:intro_mde} is then recovered from the symmetric MDE \eqref{eq:sym_mde} through the change of variables $z=\zeta^2, \zeta, z\in\C^+$.

We shall repeatedly use the inverse form of the MDE,
\begin{align}
        -\Pi^{-1}=I_z-\La+\cals[\Pi]. \label{eq:inverse_mde}
\end{align}

Define the cone of matrices with nonnegative imaginary part by
\begin{align}
        \calM_+:=\{R\in\C^{L\times L}:\Im R\geq0\}. \nonumber
\end{align}
We use the following result for the symmetric linearization.

\begin{lemma}\label{lem:lin_mde_meas}
        Under Assumption~\ref{ass:basic}, for $\zeta\in\C^+$, define
        \begin{align}
                \wh\Om(R,\zeta):=I_L+\zeta R+\cals[R]R-\La R. \label{eq:sym_mde}
        \end{align}
        Then the following statements hold.
        \begin{enumerate}[label=(\roman*)]
                \item The equation $\wh\Om(R,\zeta)=0$ has a unique solution $\wh\Pi(\zeta)\in\calM_+$. The map $\zeta\mapsto\wh\Pi(\zeta)$ is holomorphic on $\C^+$ and $\Im\wh\Pi(\zeta)>0$ for every $\zeta\in\C^+$.
                \item There exists a unique positive semidefinite matrix-valued measure $\wh V$ on $\R$ such that
                      \begin{align}
                              \wh\Pi(\zeta)=\int_{\R}\frac{\wh V(\dd x)}{x-\zeta},\qquad \wh V(\R)=I_L, \nonumber
                      \end{align}
                      and, for some constant $C_1<\infty$, $\supp\wh V\subset[-C_1,C_1]$.
                \item The measure $\wh V$ satisfies
                      \begin{align}
                              \wh V(B)=J\wh V(-B)J,\qquad J=\begin{pmatrix}I_M & 0    \\
               0   & -I_N\end{pmatrix}, \nonumber
                      \end{align}
                      for every Borel set $B\subset\R$. Consequently, $\avg{\wh V}_L$ is a compactly supported symmetric probability measure on $\R$, and $\avg{\wh\Pi(\zeta)}_L$ is its Stieltjes transform.
        \end{enumerate}
\end{lemma}

The proof of Lemma~\ref{lem:lin_mde_meas} is given in Appendix~\ref{app:basic_mde}. We now deduce Proposition~\ref{prop:pi}.

\begin{proof}[Proof of Proposition~\ref{prop:pi}]
        \emph{Existence and Stieltjes representation.} Let $z\in\C^+$ and let $\zeta=\sqrt z$ be the branch satisfying $\Im\zeta>0$. By Lemma~\ref{lem:lin_mde_meas}, the equation \eqref{eq:sym_mde} has a unique solution $\wh\Pi(\zeta)\in\calM_+$. Define
        \begin{align}
                \Pi(z):=\begin{pmatrix}\zeta\wh\Pi_{11}(\zeta) & \wh\Pi_{12}(\zeta)           \\
               \wh\Pi_{21}(\zeta)      & \zeta^{-1}\wh\Pi_{22}(\zeta)\end{pmatrix}. \label{eq:rect_from_sym}
        \end{align}
        Using $z=\zeta^2$, a direct block computation gives $\Om(\Pi(z),z)=0$. Since $\zeta=\sqrt z$ is holomorphic on $\C^+$ under the branch $\Im\zeta>0$, and since $\wh\Pi$ is holomorphic on $\C^+$ by Lemma~\ref{lem:lin_mde_meas}, the map $z\mapsto\Pi(z)$ is holomorphic on $\C^+$.

        The symmetry in Lemma~\ref{lem:lin_mde_meas}(iii) implies that the diagonal measures $\wh V_{11}$ and $\wh V_{22}$ are even. Let $V_1$ and $V_2$ be their pushforwards under $x\mapsto x^2$. Then $V_1([0,\infty))=I_M$, $V_2([0,\infty))=I_N$, and $\supp V_1\cup\supp V_2\subset[0,C_1^2]$. Averaging the integrands at $x$ and $-x$ yields
        \begin{align}
                \wt T(z) & :=z^{-1}\Pi_{11}(z)=\zeta^{-1}\wh\Pi_{11}(\zeta)=\int_{[0,\infty)}\frac{V_1(\dd\la)}{\la-z},\nonumber                                  \\
                T(z)     & :=\Pi_{22}(z)=\zeta^{-1}\wh\Pi_{22}(\zeta)=\int_{[0,\infty)}\frac{V_2(\dd\la)}{\la-z}. \label{eq:rect_st_rep}
        \end{align}
        In particular, $\Im\wt T(z)>0$ and $\Im T(z)>0$. The uniqueness of $V_1$ and $V_2$ follows from the uniqueness of the scalar Stieltjes transform after testing against vectors and using polarization. Taking normalized traces in \eqref{eq:rect_st_rep} gives the stated Stieltjes representations of $\wt m$ and $m$. Since $\Pi_{11}=z\wt T$ and $V_1([0,\infty))=I_M$, we also have
        \begin{align}
                g(z)=z\avg{\Sg\wt T(z)}_N=-\avg{\Sg}_N+\int_{[0,\infty)}\frac{\la\avg{\Sg V_1(\dd\la)}_N}{\la-z}. \nonumber
        \end{align}
        Hence $\Im m(z)>0$ and $\Im g(z)\geq0$. By \eqref{eq:inverse_mde}, we have
        \begin{align}
                \Im(I_z-\La+\cals[\Pi(z)])=\begin{pmatrix}\Im m(z)\Sg & 0                  \\
               0           & (\eta+\Im g(z))I_N\end{pmatrix}>0. \label{eq:pos_inverse_mde}
        \end{align}
        Since $\Pi(z)=-(I_z-\La+\cals[\Pi(z)])^{-1}$, \eqref{eq:pos_inverse_mde} implies $\Im\Pi(z)>0$, while $\Im(z^{-1}\Pi_{11}(z))>0$ follows directly from \eqref{eq:rect_st_rep}. This proves the existence and positivity assertions in Proposition~\ref{prop:pi}(i), as well as the representations in parts (ii)--(iii). It remains to prove uniqueness and the formula in part~(iv).

        \emph{Uniqueness.} Let $\Pi^\#(z)$ be another holomorphic solution of $\Om(\Pi^\#(z),z)=0$ such that $\Im\Pi^\#(z)>0$ and $\Im(z^{-1}\Pi^\#_{11}(z))>0$. Set $\wt T^\#(z):=z^{-1}\Pi^\#_{11}(z)$ and $T^\#(z):=\Pi^\#_{22}(z)$. The inverse MDE \eqref{eq:inverse_mde} reads $-\Pi^\#(z)^{-1}=I_z-\La+\cals[\Pi^\#(z)]$, and its block form gives
        \begin{align}
                \wt T^\#(z) & =\left(\frac{1}{1+\avg{\Sg \wt T^\#(z)}_N}AA^\top-z(I_M+\avg{T^\#(z)}_N\Sg)\right)^{-1}, \label{eq:unique_whT}       \\
                T^\#(z)     & =\left(A^\top(I_M+\avg{T^\#(z)}_N\Sg)^{-1}A-z(1+\avg{\Sg \wt T^\#(z)}_N)\right)^{-1}. \label{eq:unique_T}
        \end{align}
        Since the right-hand side of the inverse MDE \eqref{eq:inverse_mde} is symmetric, $(\Pi^\#)^\top=\Pi^\#$. Taking traces of the two diagonal block equations and cancelling the identical off-diagonal trace terms yields
        \begin{align}
                z\avg{T^\#(z)}_N=c_Nz\avg{\wt T^\#(z)}_M+c_N-1. \label{eq:tr_id_uniq}
        \end{align}
        The positivity assumptions $\Im\Pi^\#(z)>0$ and $\Im(z^{-1}\Pi^\#_{11}(z))>0$ give $\Im(z\avg{\Sg\wt T^\#(z)}_N)\geq0$ and $\Im(z\avg{\wt T^\#(z)}_M)\geq0$, while \eqref{eq:tr_id_uniq} implies $\Im(z\avg{T^\#(z)}_N)\geq0$. Moreover, $\Im(1+\avg{\Sg\wt T^\#(z)}_N)^{-1}\leq0$ and $\Im(I_M+\avg{T^\#(z)}_N\Sg)^{-1}\leq0$. Hence
        \begin{align}
                \Im\left(z(I_M+\avg{T^\#(z)}_N\Sg)-(1+\avg{\Sg \wt T^\#(z)}_N)^{-1}AA^\top\right)            & \geq\eta,\nonumber \\
                \Im\left(z(1+\avg{\Sg \wt T^\#(z)}_N)-A^\top(I_M+\avg{T^\#(z)}_N\Sg)^{-1}A\right) & \geq\eta. \nonumber
        \end{align}
        Consequently, $\|\wt T^\#(z)\|+\|T^\#(z)\|\leq2\eta^{-1}$. For any unit vectors $\bbx_1\in\C^M$ and $\bbx_2\in\C^N$, the scalar functions $z\mapsto\scalar{\bbx_1}{\wt T^\#(z)\bbx_1}$ and $z\mapsto\scalar{\bbx_2}{T^\#(z)\bbx_2}$ are holomorphic on $\C^+$ and have nonnegative imaginary parts. Equations \eqref{eq:unique_whT}--\eqref{eq:unique_T} and $\|\wt T^\#(z)\|+\|T^\#(z)\|\leq2\eta^{-1}$ imply
        \begin{align}
                \ii\eta\scalar{\bbx_1}{\wt T^\#(\ii\eta)\bbx_1}\longrightarrow-1,\qquad \ii\eta\scalar{\bbx_2}{T^\#(\ii\eta)\bbx_2}\longrightarrow-1 \nonumber
        \end{align}
        as $\eta\to\infty$. Hence, by the Nevanlinna representation and polarization, there exist unique positive semidefinite matrix-valued measures $V_1^\#$ and $V_2^\#$ on $\R$, with total masses $I_M$ and $I_N$, respectively, such that
        \begin{align}
                \wt T^\#(z)=\int_{\R}\frac{V_1^\#(\dd\la)}{\la-z},\qquad T^\#(z)=\int_{\R}\frac{V_2^\#(\dd\la)}{\la-z}. \label{eq:T_st}
        \end{align}
        We next show that these measures are supported on $[0,\infty)$. Let $\bnu_1^\#:=\avg{V_1^\#}_M$ and $\bnu_2^\#:=\avg{V_2^\#}_N$.  Let $f\in C_c^\infty(\R)$ be nonnegative.
        For $z=E+\ii\eta$,
        \begin{align}
                \Im(z\avg{\wt T^\#(z)}_M) & =\int_{\R}\frac{\eta\la}{(\la-E)^2+\eta^2}\bnu_1^\#(\dd\la),\nonumber                                      \\
                \Im(z\avg{T^\#(z)}_N)     & =\int_{\R}\frac{\eta\la}{(\la-E)^2+\eta^2}\bnu_2^\#(\dd\la). \label{eq:z_st_pois}
        \end{align}
        Multiplying the first identity in \eqref{eq:z_st_pois} by $f(E)$, integrating in $E$, and using $\Im(z\avg{\wt T^\#(z)}_M)\geq0 $ and  $ \Im(z\avg{T^\#(z)}_N)\geq0$, we obtain
        \begin{align}
                0\leq\int_{\R}\la\left(\frac1\pi\int_{\R}\frac{\eta f(E)}{(\la-E)^2+\eta^2}\dd E\right)\bnu_1^\#(\dd\la). \nonumber
        \end{align}
        Letting $\eta\downarrow0$ and using the Poisson kernel approximation gives $\int_{\R}\la f(\la)\bnu_1^\#(\dd\la)\geq0$. If $\bnu_1^\#((-\infty,0))>0$, then a nonnegative $f\in C_c^\infty((-\infty,0))$ can be chosen so that the last integral is strictly negative, a contradiction. Hence $\supp\bnu_1^\#\subset[0,\infty)$. The same argument applied to the second identity in \eqref{eq:z_st_pois} gives $\supp\bnu_2^\#\subset[0,\infty)$. Since $V_1^\#$ and $V_2^\#$ are positive semidefinite matrix-valued measures, their supports are contained in those of their normalized traces. Therefore
        \begin{align}
                \supp V_1^\#\cup\supp V_2^\#\subset[0,\infty). \label{eq:T_pos_supp}
        \end{align}
        Now set $z=\zeta^2$, with $\Im\zeta>0$, and define
        \begin{align}
                \wh\Pi^\#(\zeta):=\begin{pmatrix}\zeta^{-1}\Pi^\#_{11}(\zeta^2) & \Pi^\#_{12}(\zeta^2)      \\
               \Pi^\#_{21}(\zeta^2)           & \zeta\Pi^\#_{22}(\zeta^2)\end{pmatrix}. \nonumber
        \end{align}
        A direct block computation using $\Om(\Pi^\#(\zeta^2),\zeta^2)=0$ shows that $\wh\Om(\wh\Pi^\#(\zeta),\zeta)=0$. Moreover, by \eqref{eq:T_st}--\eqref{eq:T_pos_supp},
        \begin{align}
                \Im(\zeta^{-1}\Pi^\#_{11}(\zeta^2))=\int_{[0,\infty)}\frac{\Im\zeta(\la+|\zeta|^2)}{|\la-\zeta^2|^2}V_1^\#(\dd\la)>0, \nonumber
        \end{align}
        and similarly $\Im(\zeta\Pi^\#_{22}(\zeta^2))>0$. Hence $\Im\cals[\wh\Pi^\#(\zeta)]>0$. Using the inverse form of the symmetric MDE \eqref{eq:sym_mde}, we conclude that $\Im\wh\Pi^\#(\zeta)>0$. Lemma~\ref{lem:lin_mde_meas}(i) therefore yields $\wh\Pi^\#(\zeta)=\wh\Pi(\zeta)$, and \eqref{eq:rect_from_sym} gives $\Pi^\#(z)=\Pi(z)$. This proves uniqueness.

        Now we prove \eqref{eq:origin_mass}. Taking traces in the two diagonal block equations of \eqref{eq:intro_mde} and subtracting them yields
        \begin{align}
                zm(z)-\frac1N\Tr\Pi_{11}(z)=c_N-1. \label{eq:origin_mass_id}
        \end{align}
        By \eqref{eq:rect_st_rep},
        \begin{align}
                \Re\Pi_{11}(\ii\eta)=-\int_{[0,\infty)}\frac{\eta^2}{\la^2+\eta^2}V_1(\dd\la)\leq0. \nonumber
        \end{align}
        Hence, at $z=\ii\eta$,
        \begin{align}
                \eta\Im m(\ii\eta)=1-c_N-\frac1N\Tr\Re\Pi_{11}(\ii\eta). \nonumber
        \end{align}

        Suppose first that $c_N<1$. Then $\Im m(\ii\eta)\geq(1-c_N)/\eta$. The Schur complement of \eqref{eq:inverse_mde} gives
        \begin{align}
                \Pi_{11}(\ii\eta)^{-1}=-I_M-m(\ii\eta)\Sg+\frac1{\ii\eta+g(\ii\eta)}AA^\top. \label{eq:schur_pi11}
        \end{align}
        Since $\Im(\ii\eta+g(\ii\eta))>0$, we have $-\Im\Pi_{11}(\ii\eta)^{-1}\geq(\Im m(\ii\eta))\Sg$. Hence Assumption~\ref{ass:basic} yields $\|\Pi_{11}(\ii\eta)\|\leq C\eta/(1-c_N)$, and therefore $N^{-1}\Tr\Re\Pi_{11}(\ii\eta)\to0$ as $\eta\downarrow0$. The Stieltjes representation of $m$ then gives $\bnu(\{0\})=1-c_N$.

        Suppose now that $c_N\geq1$ and $\bnu(\{0\})>0$. The Stieltjes representation of $m$ gives $\Im m(\ii\eta)\geq\bnu(\{0\})/\eta$. By \eqref{eq:schur_pi11}, we get $\|\Pi_{11}(\ii\eta)\|\leq C\eta/\bnu(\{0\})$. Letting $\eta\downarrow0$ in \eqref{eq:origin_mass_id} gives $\bnu(\{0\})=1-c_N\leq0$, a contradiction. Hence $\bnu(\{0\})=0$ for $c_N\geq1$. This proves \eqref{eq:origin_mass}.

        It remains to prove $\sup\supp\bnu\sim1$. The Stieltjes representations in parts~(ii)--(iii) extend $\Pi_{11}$, $\Pi_{22}$, $m$, and $g$ to $\C\setminus[0,\infty)$. Hence these functions are well defined at $z=-1$. The Schur complement identities from \eqref{eq:inverse_mde} also extend to this point. We have
        \begin{align}
                -g(-1)=\int_{[0,\infty)}\frac{\avg{\Sg V_1(\dd\la)}_N}{1+\la}\gtrsim1. \nonumber
        \end{align}
        Since $m(-1)>0$, the lower right Schur complement gives $m(-1)\leq(1-g(-1))^{-1}\leq1-c$. By the Stieltjes representation of $m$, we have
        \begin{align}
                1-m(-1)=\int_{[0,\infty)}\frac{\la}{1+\la}\bnu(\dd\la)\leq\sup\supp\bnu. \nonumber
        \end{align}
        Thus $\sup\supp\bnu\gtrsim1$. The reverse bound follows from $\supp\bnu\subset[0,C_*]$. This completes the proof.
\end{proof}

\subsection{Two-dimensional stability}
\label{subsec:2d_stab}

We establish the two-dimensional stability estimates needed for Theorem~\ref{thm:dens_reg} and Corollary~\ref{cor:bdry_reg_m_g}. Since the self-energy operator $\calS$ depends only on two normalized block traces, the stability problem reduces to the two-dimensional operator $\calJ$ defined in \eqref{eq:J_def}. We show that only one eigenvalue of $\calJ$ can approach $1$.

Fix sufficiently small $\tau\in(0,1)$ such that $\tau^{-1}\geq C_*+1$, and fix $\eta_0>0$ sufficiently large. We work on the spectral domain
\begin{align}
        \mb D(\tau,\eta_0):=\{E+\ii\eta:\tau\leq E\leq\tau^{-1},0<\eta\leq\eta_0\}. \label{eq:D_tau_eta}
\end{align}
Since $\tau$ can be chosen arbitrarily small, estimates uniform on $\mb D(\tau,\eta_0)$ yield the required local statements on $(0,\infty)$. The constants in this subsection may depend on $\tau$, $\eta_0$, and the constants in Assumption~\ref{ass:basic}. 

For $R,T\in\C^{L\times L}$, let
\begin{align}
        \calC_{R,T}[Q]:=RQT, \nonumber
\end{align}
and write $\calC_R:=\calC_{R,R}$. To motivate the stability operator, suppose that $G$ satisfies the perturbed MDE
\begin{align}
        \Om(G,z)=I_L+(I_z+\cals[G]-\La)G=D. \nonumber
\end{align}
By the MDE~\eqref{eq:intro_mde}, we have $\Om(\Pi,z)=0$, and hence
\begin{align}
        D=\Om(G,z)-\Om(\Pi,z)=(I_z-\La)(G-\Pi)+\cals[G]G-\cals[\Pi]\Pi. \nonumber
\end{align}
Using the linearity of $\cals$, we obtain
\begin{align}
        D=(I_z-\La+\cals[\Pi])(G-\Pi)+\cals[G-\Pi]G =-\Pi^{-1}(G-\Pi)+\cals[G-\Pi]\Pi+\cals[G-\Pi](G-\Pi), \nonumber
\end{align}
where we used the inverse form of the MDE~\eqref{eq:inverse_mde}. Multiplying by $\Pi$ from the left and rearranging gives
\begin{align}
        (\calI-\calC_\Pi\cals)[G-\Pi]=-\Pi D+\Pi\cals[G-\Pi](G-\Pi). \nonumber
\end{align}
Thus the control of the perturbation $G-\Pi$ is reduced to the stability of the linear operator $\calI-\calC_\Pi\cals$. This motivates the stability operator
\begin{align}
        \calB[R]:=R-\Pi\cals[R]\Pi=(\calI-\calC_\Pi\cals)[R]. \label{eq:B_def}
\end{align}
We further define $\calE:\C^2\to\C^{L\times L}$ and $\calP:\C^{L\times L}\to\C^2$ by
\begin{align}
        \calE[\bbh]:=\begin{pmatrix}h_1\Sg & 0 \\
               0 & h_2I_N\end{pmatrix},\qquad
        \calP[R]:=\begin{pmatrix}\avg{R_{22}}_N \\
                                  \avg{\Sg R_{11}}_N\end{pmatrix},\qquad
        \bbh=\begin{pmatrix}h_1 \\
                             h_2\end{pmatrix}\in\C^2. \nonumber
\end{align}
Then $\cals=\calE\calP$. By Assumption~\ref{ass:basic},
\begin{align}
        \|\calP\|+\|\calE\|\lesssim1,\qquad \|\calC_\Pi\|\lesssim\|\Pi\|^2. \label{eq:C_Pi_bd}
\end{align}
Define
\begin{align}
        \calJ:=\calP\calC_\Pi\calE. \label{eq:J_def}
\end{align}

We begin with the positivity identities that will be used throughout the stability analysis.

\begin{lemma}\label{lem:im_ids}
        Under Assumption~\ref{ass:basic}, for every $z=E+\ii\eta\in\C^+$, the solution $\Pi(z)$ satisfies
        \begin{align}
                \Im\Pi_{11} & =\Im m\Pi_{11}^*\Sg\Pi_{11}+(\eta+\Im g)\Pi_{21}^*\Pi_{21}, \label{eq:Im_Pi_11} \\
                \Im\Pi_{22} & =\Im m\Pi_{12}^*\Sg\Pi_{12}+(\eta+\Im g)\Pi_{22}^*\Pi_{22}. \label{eq:Im_Pi_22}
        \end{align}
        Consequently,
        \begin{align}
                \Im g & =\Im m\avg{\Pi_{11}^*\Sg\Pi_{11}\Sg}_N+(\eta+\Im g)\avg{\Pi_{21}\Sg\Pi_{21}^*}_N, \label{eq:Im_g} \\
                \Im m & =\Im m\avg{\Pi_{12}^*\Sg\Pi_{12}}_N+(\eta+\Im g)\avg{\Pi_{22}^*\Pi_{22}}_N. \label{eq:Im_m}
        \end{align}
\end{lemma}

\begin{proof}
        Taking imaginary parts in \eqref{eq:inverse_mde} gives
        \begin{align}
                \Im\Pi=\Pi^*\begin{pmatrix}\Im m\Sg & 0               \\
               0          & (\eta+\Im g)I_N\end{pmatrix}\Pi. \nonumber
        \end{align}
        The block identities \eqref{eq:Im_Pi_11}--\eqref{eq:Im_Pi_22} follow immediately. Multiplying the first identity by $\Sg$, taking normalized traces, and taking the normalized trace of the second identity yield \eqref{eq:Im_g}--\eqref{eq:Im_m}.
\end{proof}

The next lemma makes the two-dimensional reduction explicit.

\begin{lemma}\label{lem:stab_red}
        Under Assumption~\ref{ass:basic}, let $z\in\C^+$. With respect to the canonical basis of $\C^2$, the matrix representation of $\calJ(z)$ is
        \begin{align}
                \begin{pmatrix}\avg{\Pi_{21}\Sg\Pi_{12}}_N & \avg{\Pi_{22}^2}_N          \\
               \avg{(\Pi_{11}\Sg)^2}_N     & \avg{\Pi_{21}\Sg\Pi_{12}}_N\end{pmatrix}. \label{eq:J_mat}
        \end{align}
        If $\calI-\calJ$ is invertible on $\C^2$, then $\calB$ is invertible and
        \begin{align}
                \calB^{-1}=\calI+\calC_\Pi\calE(\calI-\calJ)^{-1}\calP. \label{eq:B_inverse_red}
        \end{align}
\end{lemma}

\begin{proof}
        The representation \eqref{eq:J_mat} follows by evaluating $\calJ$ on the canonical basis of $\C^2$. Assume that $\calI-\calJ$ is invertible and fix $T\in\C^{L\times L}$. Set
        \begin{align}
                \bbh:=(\calI-\calJ)^{-1}\calP[T],\qquad R:=T+\calC_\Pi\calE[\bbh]. \nonumber
        \end{align}
        Then $\calP[R]=\calP[T]+\calJ[\bbh]=\bbh$, and hence
        \begin{align}
                \calB[R]=R-\calC_\Pi\calE\calP[R]=T. \nonumber
        \end{align}
        Thus $\calB$ is surjective. If $\calB[R]=0$, then applying $\calP$ gives $(\calI-\calJ)\calP[R]=0$, so $\calP[R]=0$ and therefore $R=0$. Hence $\calB$ is invertible, and the preceding construction gives \eqref{eq:B_inverse_red}.
\end{proof}

For convenience, set $q:=\Im m/(\eta+\Im g)$. Taking normalized traces in \eqref{eq:Im_Pi_11} and \eqref{eq:Im_Pi_22}, and using $\Pi_{21}=\Pi_{12}^{\top}$ yields
\begin{align}
        qc_N\|\Sg^{1/2}\Pi_{11}\Sg^{1/2}\|_{\rr{hs}}^2 & =1-\avg{\Pi_{12}^*\Sg\Pi_{12}}_N-\frac{\eta}{\eta+\Im g}. \label{eq:q_Pi11} \\
        \|\Pi_{22}\|_{\rr{hs}}^2                       & =q(1-\avg{\Pi_{12}^*\Sg\Pi_{12}}_N), \label{eq:q_Pi22}          
\end{align}

We shall use the following a priori estimates, which are derived from the positivity identities and the block structure of the MDE \eqref{eq:intro_mde}. Their proof is given in Appendix~\ref{sub:apriori_est}.

\begin{lemma}\label{lem:unif_apriori}
        Under Assumption~\ref{ass:basic}, there exists a constant $c_1>0$ such that, uniformly for $z\in\mb D(\tau,\eta_0)$, the following estimates hold.
        \begin{enumerate}[label=(\roman*)]
                \item The off-diagonal blocks satisfy
                      \begin{align}
                              \avg{\Pi_{21}\Sg\Pi_{21}^*}_N=\avg{\Pi_{12}^*\Sg\Pi_{12}}_N\leq1-c_1. \label{eq:offdiag_gap}
                      \end{align}
                \item The diagonal blocks satisfy
                      \begin{align}
                              \|\Sg^{1/2}\Pi_{11}\Sg^{1/2}\|_{\rr{hs}}+\|\Pi_{22}\|_{\rr{hs}}\sim1. \label{eq:diag_size}
                      \end{align}
                \item The two imaginary parts are comparable:
                      \begin{align}
                              \frac{\Im m}{\eta+\Im g}\sim1,\qquad \Im m\sim\Im g. \label{eq:imag_comp}
                      \end{align}
        \end{enumerate}
\end{lemma}

We shall also need the following pointwise bounds on the solution of the MDE \eqref{eq:intro_mde}. Their proof is given in Appendix~\ref{sub:apriori_est}, together with that of Lemma~\ref{lem:unif_apriori}.

\begin{lemma}\label{lem:Pi_bd}
        Under Assumption~\ref{ass:basic}, uniformly for $z\in\mb D(\tau,\eta_0)$,
        \begin{align}
                \|\Pi^{-1}(z)\|\lesssim1,\qquad \|\Pi(z)\|\lesssim\frac1{\Im m(z)}. \nonumber
        \end{align}
        Moreover,
        \begin{align}
                \|\Pi_{11}(z)\|\lesssim\frac1{\Im m(z)+\dist(z,\supp\bnu_1)},\qquad \|\Pi_{22}(z)\|\lesssim\frac1{\Im m(z)+\dist(z,\supp\bnu_2)}. \nonumber
        \end{align}
\end{lemma}

By Lemma~\ref{lem:stab_red}, the stability problem is reduced to controlling the two eigenvalues of $\calJ$. The following proposition is the main quantitative estimate of this subsection.

\begin{proposition}\label{prop:red_stab}
        Under Assumption~\ref{ass:basic}, uniformly for $z\in\mb D(\tau,\eta_0)$, the matrix $\calI-\calJ$ is invertible and
        \begin{align}
                \|(\calI-\calJ)^{-1}\|\lesssim\frac1{\eta/\Im m+(\Im m)^2}. \label{eq:inverse_J_bd}
        \end{align}
\end{proposition}

We first identify the only potentially unstable eigenvalue of $\calJ$. By the Cauchy--Schwarz inequality and Lemma~\ref{lem:unif_apriori},
\begin{align}
        \left|\avg{\Pi_{21}\Sg\Pi_{12}}_N\right| & \leq\avg{\Pi_{12}^*\Sg\Pi_{12}}_N\leq1-c_1, \label{eq:Pi_21_bd}                                                            \\
        \left|\avg{\Pi_{22}^2}_N\right|          & \leq\avg{\Pi_{22}^*\Pi_{22}}_N=\frac{\Im m}{\eta+\Im g}(1-\avg{\Pi_{12}^*\Sg\Pi_{12}}_N), \label{eq:Pi_22_bd}   \\
        \left|\avg{(\Pi_{11}\Sg)^2}_N\right|     & \leq\avg{\Pi_{11}^*\Sg\Pi_{11}\Sg}_N=\frac{\eta+\Im g}{\Im m}(1-\avg{\Pi_{12}^*\Sg\Pi_{12}}_N)-\frac{\eta}{\Im m}. \label{eq:Pi_11_bd} 
\end{align}
In particular, all entries of $\calJ$ are uniformly bounded on $\mb D(\tau,\eta_0)$. Its eigenvalues are
\begin{align}
        \la_\pm=\avg{\Pi_{21}\Sg\Pi_{12}}_N\pm\sqrt{\avg{\Pi_{22}^2}_N\avg{(\Pi_{11}\Sg)^2}_N}. \label{eq:lambda_pm}
\end{align}
If $\avg{\Pi_{22}^2}_N\avg{(\Pi_{11}\Sg)^2}_N=0$, we take the square root in \eqref{eq:lambda_pm} to be zero. Then $|1-\la_\pm|\geq c_1$, and \eqref{eq:inverse_J_bd} follows from Lemma~\ref{lem:unif_apriori}. We may therefore assume that the product is nonzero. Choose the unique $\te\in(-\pi/2,\pi/2]$ such that
\begin{align}
        \avg{\Pi_{22}^2}_N\avg{(\Pi_{11}\Sg)^2}_N=\left|\avg{\Pi_{22}^2}_N\right|\left|\avg{(\Pi_{11}\Sg)^2}_N\right|e^{2\ii\te}, \nonumber
\end{align}
and define the square root in \eqref{eq:lambda_pm} by
\begin{align}
        \sqrt{\avg{\Pi_{22}^2}_N\avg{(\Pi_{11}\Sg)^2}_N}:=\sqrt{\left|\avg{\Pi_{22}^2}_N\right|\left|\avg{(\Pi_{11}\Sg)^2}_N\right|}e^{\ii\te}. \nonumber
\end{align}
Then $\cos\te\geq0$. Moreover,
\begin{align}
        |1-\la_-|\geq\Re(1-\la_-)\geq1-\left|\avg{\Pi_{21}\Sg\Pi_{12}}_N\right|\geq c_1. \label{eq:lambda_minus}
\end{align}
Hence only $\la_+$ can become unstable. The next two lemmas treat the case in which $\la_+$ approaches $1$. The first identifies the resulting nondegenerate structure, and the second gives the required lower bound on $1-\la_+$.

\begin{lemma}\label{lem:unst_eig_ndg}
        Under Assumption~\ref{ass:basic}, there exists a constant $c_2>0$ such that, uniformly for $z\in\mb D(\tau,\eta_0)$ satisfying $|1-\la_+|\leq c_2$,
        \begin{align}
                \left|\avg{\Pi_{22}^2}_N\right|\sim\left|\avg{(\Pi_{11}\Sg)^2}_N\right|\sim\cos\te\sim1. \nonumber
        \end{align}
\end{lemma}

The second lemma is proved by aligning the phases of the two off-diagonal entries of $\calJ$ and separating their real and imaginary parts. This converts the nondegeneracy in Lemma~\ref{lem:unst_eig_ndg} into the required lower bound on $1-\la_+$.

\begin{lemma}\label{lem:lambda_plus_lb}
        Under Assumption~\ref{ass:basic}, uniformly for $z\in\mb D(\tau,\eta_0)$ satisfying $|1-\la_+|\leq c_2$,
        \begin{align}
                |1-\la_+|\gtrsim\frac{\eta}{\Im m}+(\Im m)^2. \label{eq:lambda_plus_lb}
        \end{align}
\end{lemma}

We first deduce Proposition~\ref{prop:red_stab} from the two lemmas and then prove them.

\begin{proof}[Proof of Proposition~\ref{prop:red_stab}]
        By \eqref{eq:lambda_minus}, $|1-\la_-|\gtrsim1$. If $|1-\la_+|\leq c_2$, Lemma~\ref{lem:lambda_plus_lb} gives
        \begin{align}
                |1-\la_+|\gtrsim\frac{\eta}{\Im m}+(\Im m)^2. \nonumber
        \end{align}
        If $|1-\la_+|>c_2$, the same estimate follows from Lemma~\ref{lem:unif_apriori}, since $\eta/\Im m+(\Im m)^2\lesssim1$. Hence both eigenvalues of $\calI-\calJ$ are separated from zero on the required scale. Finally,
        \begin{align}
                (\calI-\calJ)^{-1}=\frac1{(1-\la_+)(1-\la_-)}\begin{pmatrix}1-\avg{\Pi_{21}\Sg\Pi_{12}}_N & \avg{\Pi_{22}^2}_N            \\
               \avg{(\Pi_{11}\Sg)^2}_N       & 1-\avg{\Pi_{21}\Sg\Pi_{12}}_N\end{pmatrix}, \nonumber
        \end{align}
        and the entries of the numerator matrix are uniformly bounded. This proves \eqref{eq:inverse_J_bd}.
\end{proof}

It remains to prove Lemmas~\ref{lem:unst_eig_ndg} and \ref{lem:lambda_plus_lb}, which contain the main argument behind Proposition~\ref{prop:red_stab}.

\begin{proof}[Proof of Lemma~\ref{lem:unst_eig_ndg}]
        From $|1-\la_+|\leq c_2$, \eqref{eq:Pi_21_bd}, and $\cos\te\geq0$, we obtain
        \begin{align}
                1-c_2\leq\Re\la_+\leq1-c_1+\sqrt{\left|\avg{\Pi_{22}^2}_N\right|\left|\avg{(\Pi_{11}\Sg)^2}_N\right|}\cos\te. \nonumber
        \end{align}
        Choose $c_2=c_1/2$. Then $\sqrt{\left|\avg{\Pi_{22}^2}_N\right|\left|\avg{(\Pi_{11}\Sg)^2}_N\right|}\cos\te\geq c_1/2$. The two moduli are uniformly bounded above, and $0\leq\cos\te\leq1$. It follows that each of $\left|\avg{\Pi_{22}^2}_N\right|$, $\left|\avg{(\Pi_{11}\Sg)^2}_N\right|$, and $\cos\te$ is bounded away from zero. This proves the claim.
\end{proof}

\begin{proof}[Proof of Lemma~\ref{lem:lambda_plus_lb}]
        Choose $\phi\in\R$ such that
        \begin{align}
                e^{2\ii\phi}\avg{\Pi_{22}^2}_N=\left|\avg{\Pi_{22}^2}_N\right|e^{\ii\te}, \qquad e^{-2\ii\phi}\avg{(\Pi_{11}\Sg)^2}_N=\left|\avg{(\Pi_{11}\Sg)^2}_N\right|e^{\ii\te}. \nonumber
        \end{align}
        Define the real matrices $X_{ij}$ and $Y_{ij}$ by
        \begin{align}
                e^{-\ii\phi}\Sg^{1/2}\Pi_{11}\Sg^{1/2} & =X_{11}+\ii Y_{11},                         \nonumber \\
                e^{\ii\phi}\Pi_{22}                    & =X_{22}+\ii Y_{22},\label{eq:rot_pi_22} \\
                \Sg^{1/2}\Pi_{12}                      & =X_{12}+\ii Y_{12}. \nonumber
        \end{align}
        Here $X_{11},Y_{11},X_{22},Y_{22}$ are real symmetric, while $X_{12},Y_{12}$ are real rectangular matrices. For rectangular matrices throughout this argument, $\|\cdot\|_{\rr{hs}}$ denotes the Frobenius norm normalized by $N^{-1/2}$. Since $\Pi_{21}=\Pi_{12}^\top$,
        \begin{align}
                \Re\avg{\Pi_{21}\Sg\Pi_{12}}_N=\|X_{12}\|_{\rr{hs}}^2-\|Y_{12}\|_{\rr{hs}}^2=\avg{\Pi_{12}^*\Sg\Pi_{12}}_N-2\|Y_{12}\|_{\rr{hs}}^2. \nonumber
        \end{align}
        Moreover,
        \begin{align}
                \left|\avg{\Pi_{22}^2}_N\right|\cos\te      & =\|\Pi_{22}\|_{\rr{hs}}^2-2\|Y_{22}\|_{\rr{hs}}^2,\nonumber                                    \\
                \left|\avg{(\Pi_{11}\Sg)^2}_N\right|\cos\te & =c_N(\|\Sg^{1/2}\Pi_{11}\Sg^{1/2}\|_{\rr{hs}}^2-2\|Y_{11}\|_{\rr{hs}}^2). \nonumber
        \end{align}
        It follows that
        \begin{align}
                \Re(1-\la_+) & =1-\avg{\Pi_{12}^*\Sg\Pi_{12}}_N+2\|Y_{12}\|_{\rr{hs}}^2 -\sqrt{(\|\Pi_{22}\|_{\rr{hs}}^2-2\|Y_{22}\|_{\rr{hs}}^2)c_N(\|\Sg^{1/2}\Pi_{11}\Sg^{1/2}\|_{\rr{hs}}^2-2\|Y_{11}\|_{\rr{hs}}^2)}. \nonumber
        \end{align}
        By Lemmas~\ref{lem:unif_apriori} and \ref{lem:unst_eig_ndg}, the factors appearing in the preceding square roots are uniformly nondegenerate. Hence the identity $\sqrt{x}-\sqrt{y}=(x-y)/(\sqrt{x}+\sqrt{y})$ gives
        \begin{align}
                \sqrt{c_N}\|\Pi_{22}\|_{\rr{hs}}\|\Sg^{1/2}\Pi_{11}\Sg^{1/2}\|_{\rr{hs}} & -\sqrt{(\|\Pi_{22}\|_{\rr{hs}}^2-2\|Y_{22}\|_{\rr{hs}}^2)c_N(\|\Sg^{1/2}\Pi_{11}\Sg^{1/2}\|_{\rr{hs}}^2-2\|Y_{11}\|_{\rr{hs}}^2)} \nonumber \\
                                                                                         & \gtrsim\|Y_{11}\|_{\rr{hs}}^2+\|Y_{22}\|_{\rr{hs}}^2. \nonumber
        \end{align}
        By \eqref{eq:q_Pi22} and \eqref{eq:q_Pi11}, we have
        \begin{align}
                1-\avg{\Pi_{12}^*\Sg\Pi_{12}}_N  -\sqrt{c_N}\|\Sg^{1/2}\Pi_{11}\Sg^{1/2}\|_{\rr{hs}} \|\Pi_{22}\|_{\rr{hs}}  &=\frac{\frac{\eta}{\eta+\Im g}(1-\avg{\Pi_{12}^*\Sg\Pi_{12}}_N)}{1-\avg{\Pi_{12}^*\Sg\Pi_{12}}_N+\sqrt{c_N}\|\Sg^{1/2}\Pi_{11}\Sg^{1/2}\|_{\rr{hs}} \|\Pi_{22}\|_{\rr{hs}}} \nonumber\\
                \quad& \gtrsim\frac{\eta}{\eta+\Im g}. \nonumber
        \end{align}
        Therefore,
        \begin{align}
                \Re(1-\la_+)\gtrsim\frac{\eta}{\eta+\Im g}+\|Y_{11}\|_{\rr{hs}}^2+\|Y_{12}\|_{\rr{hs}}^2+\|Y_{22}\|_{\rr{hs}}^2. \label{eq:lambda_coerc}
        \end{align}

        It remains to relate the rotated imaginary parts to $\Im m$. From \eqref{eq:rot_pi_22},
        \begin{align}
                \Im m\leq|\cos\phi|\|Y_{22}\|_{\rr{hs}}+|\sin\phi|\|X_{22}\|_{\rr{hs}}. \nonumber
        \end{align}
        The block equation $A^\top\Pi_{12}-(z+g)\Pi_{22}=I_N$ becomes
        \begin{align}
                A^\top\Sg^{-1/2}(X_{12}+\ii Y_{12})-(z+g)e^{-\ii\phi}(X_{22}+\ii Y_{22})=I_N. \nonumber
        \end{align}
        Taking imaginary parts yields
        \begin{align}
                |\Im((z+g)e^{-\ii\phi})|\|X_{22}\|_{\rr{hs}}\lesssim\|Y_{12}\|_{\rr{hs}}+\|Y_{22}\|_{\rr{hs}}.\nonumber
        \end{align}
        Furthermore,
        \begin{align}
                E\sin\phi=\eta\cos\phi+\Im(ge^{-\ii\phi})-\Im((z+g)e^{-\ii\phi}),\nonumber
        \end{align}
        and $\Im(ge^{-\ii\phi})=\avg{Y_{11}}_N$. Since $E\geq\tau/2$, $\|X_{22}\|_{\rr{hs}}\lesssim1$, and $|z+g|\lesssim1$, we conclude that
        \begin{align}
                \Im m\lesssim\eta+\|Y_{11}\|_{\rr{hs}}+\|Y_{12}\|_{\rr{hs}}+\|Y_{22}\|_{\rr{hs}}.\label{eq:rot_im_bd}
        \end{align}

        If $\Im m\lesssim\eta$, then \eqref{eq:imag_comp} implies $\eta/(\eta + \Im g)\gtrsim \Im m / (\eta + \Im g)\sim 1$, while \eqref{eq:lambda_coerc} gives $|1-\la_+|\gtrsim 1$. If $\Im m\gg\eta$, then \eqref{eq:rot_im_bd} and Cauchy--Schwarz give
        \begin{align}
                \|Y_{11}\|_{\rr{hs}}^2+\|Y_{12}\|_{\rr{hs}}^2+\|Y_{22}\|_{\rr{hs}}^2\gtrsim(\Im m)^2. \nonumber
        \end{align}
        Combining this with \eqref{eq:imag_comp} and \eqref{eq:lambda_coerc} proves \eqref{eq:lambda_plus_lb}.
\end{proof}

Proposition~\ref{prop:red_stab} and Lemma~\ref{lem:stab_red} immediately yield a bound for the stability operator $\calB$ of the MDE \eqref{eq:intro_mde}. We record this consequence separately for the derivative estimate below and the analyticity argument in Subsection~\ref{subsec:dens_bdry}.

\begin{corollary}\label{cor:full_stab}
        Under Assumption~\ref{ass:basic}, uniformly for $z\in\mb D(\tau,\eta_0)$,
        \begin{align}
                \|\calB^{-1}\|\lesssim\frac1{(\Im m)^4}. \nonumber
        \end{align}
\end{corollary}

\begin{proof}
        Lemmas~\ref{lem:stab_red}, \ref{lem:unif_apriori}, and \ref{lem:Pi_bd}, Proposition~\ref{prop:red_stab}, and \eqref{eq:C_Pi_bd} give
        \begin{align}
                \|\calB^{-1}\|\lesssim1+\frac{\|\Pi\|^2}{\eta/\Im m+(\Im m)^2}\lesssim\frac1{(\Im m)^4}, \nonumber
        \end{align}
        where we also used $\Im m\lesssim1$, which follows from Lemma~\ref{lem:unif_apriori}(ii).
\end{proof}

The final input for Theorem~\ref{thm:dens_reg} is obtained by differentiating the MDE \eqref{eq:intro_mde}. The two scalar derivatives admit a sharper bound than the full matrix derivative.

\begin{lemma}\label{lem:der_bd}
        Under Assumption~\ref{ass:basic}, uniformly for $z\in\mb D(\tau,\eta_0)$,
        \begin{align}
                \|\partial_z\Pi(z)\|              & \lesssim\frac1{(\Im m(z))^6}, \label{eq:part_Pi_bd}  \\
                |\partial_zm(z)|+|\partial_zg(z)| & \lesssim\frac1{(\Im m(z))^2}. \label{eq:part_m_g_bd}
        \end{align}
\end{lemma}

\begin{proof}
        Differentiating \eqref{eq:inverse_mde} gives
        \begin{align}
                \calB[\partial_z\Pi]=\Pi\wt I\Pi. \nonumber
        \end{align}
        Corollary~\ref{cor:full_stab} and Lemma~\ref{lem:Pi_bd} yield \eqref{eq:part_Pi_bd}. Applying $\calP$ and using $\calP[\calB[R]]=(\calI-\calJ)[\calP[R]]$, we obtain
        \begin{align}
                \begin{pmatrix}\partial_zm \\
                        \partial_zg\end{pmatrix}=(\calI-\calJ)^{-1}\begin{pmatrix}\avg{\Pi_{22}^2}_N \\
                                                                   \avg{\Pi_{21}\Sg\Pi_{12}}_N\end{pmatrix}. \label{eq:part_mg}
        \end{align}
        The vector on the right-hand side is uniformly bounded by \eqref{eq:Pi_21_bd}--\eqref{eq:Pi_22_bd}. Proposition~\ref{prop:red_stab} therefore yields \eqref{eq:part_m_g_bd}.
\end{proof}

\subsection{Regularity of the density and boundary values}
\label{subsec:dens_bdry}

We now prove Theorem~\ref{thm:dens_reg} and Corollary~\ref{cor:bdry_reg_m_g}. The derivative estimates from the previous subsection first yield H\"older continuity of the imaginary parts of $m$ and $g$ up to the real axis. Stieltjes inversion then identifies the density of $\bnu$, while the Stieltjes representations transfer the boundary regularity to the full functions $m$ and $g$.

\begin{proof}[Proof of Theorem~\ref{thm:dens_reg}]
        Fix sufficiently small $\tau>0$. For $z\in\C^+$, set
        \begin{align}
                \rho(z):=\frac1\pi\Im m(z). \nonumber
        \end{align}
        Writing $z=E+\ii\eta$ and using the Wirtinger derivative $\partial_z=\frac12(\partial_E-\ii\partial_\eta)$, the holomorphy of $m$ gives
        \begin{align}
                2\pi\ii\partial_z\rho(z)=\partial_zm(z). \nonumber
        \end{align}
        Hence Lemma~\ref{lem:der_bd} yields
        \begin{align}
                |\partial_z\rho(z)|\lesssim\rho(z)^{-2},\qquad z\in\mb D(\tau/2,\eta_0). \nonumber
        \end{align}
        This implies that the harmonic extension $\rho$ is locally $1/3$-H\"older continuous, and
       \begin{align}
        |\rho(z_1)-\rho(z_2)| \lesssim |z_1-z_2|^{1/3}, \qquad z_1,z_2\in\mb D(\tau/2,\eta_0). \label{eq:holder_rho}
        \end{align}
        Therefore $\rho$ admits a unique $1/3$-H\"older continuous extension to $\overline{\mb D(\tau/2,\eta_0)}$. In particular, for every $E\in[\tau/2,2/\tau]$,
        \begin{align}
                \rho(E):=\lim_{\eta\downarrow0}\rho(E+\ii\eta)
                =\frac1\pi\lim_{\eta\downarrow0}\Im m(E+\ii\eta)
                \nonumber
        \end{align}
        exists, and the boundary function $\rho$ is $1/3$-H\"older continuous on $[\tau/2,2/\tau]$.

        We next identify the restriction of $\bnu$ to $[\tau/2,2/\tau]$. The measure $\bnu$ has no atoms in $[\tau/2,2/\tau]$: if $\bnu(\{E\})>0$ for some $E\in[\tau/2,2/\tau]$, then Proposition~\ref{prop:pi} gives
        \begin{align}
                \Im m(E+\ii\eta)\geq\frac{\bnu(\{E\})}{\eta}, \nonumber
        \end{align}
        contradicting the boundedness of $\rho$ on $[\tau/2,2/\tau]$ from \eqref{eq:holder_rho}. Hence, for every interval $(a,b)$ with $[a,b]\subset[\tau/2,2/\tau]$, the Stieltjes inversion formula and \eqref{eq:holder_rho} yield
        \begin{align}
                \bnu((a,b))=\lim_{\eta\downarrow0}\int_a^b\rho(E+\ii\eta)\dd E=\int_a^b\rho(E)\dd E, \nonumber
        \end{align}
        where the last equality follows from the uniform convergence implied by \eqref{eq:holder_rho}. Therefore $\bnu(\dd E)=\rho(E)\dd E$ on $(\tau/2,2/\tau)$, and hence on $[\tau/2,2/\tau]$ since $\bnu$ has no atoms there. Since $\tau$ can be chosen arbitrarily small, $\bnu$ is absolutely continuous on $(0,\infty)$ and $\rho$ is locally $1/3$-H\"older continuous there. Together with $\supp\bnu\subset[0,\infty)$ from Proposition~\ref{prop:pi}, this gives
        \begin{align}
                \bnu(\dd E)=\bnu(\{0\})\delta_0(\dd E)+\rho(E)\mb 1(E>0)\dd E, \nonumber
        \end{align}
        proving \eqref{eq:dens_dec}.

        It remains to prove that $\rho$ is real analytic on $\{E\in(0,\infty):\rho(E)>0\}$. We follow the holomorphic extension argument in the proof of \cite[Proposition~2.2]{ajanki2019stability}. Fix $E_0\in(0,\infty)$ with $\rho(E_0)>0$. By \eqref{eq:holder_rho}, $\Im m$ is bounded away from zero in the intersection of a sufficiently small neighborhood of $E_0$ with the closed upper half-plane. Lemma~\ref{lem:der_bd} then gives a uniform bound on $\partial_z\Pi$ there, so $\Pi$ admits a unique continuous boundary value $ Q_0:=\lim_{\eta\downarrow0}\Pi(E_0+\ii\eta)$.
        Moreover, Corollary~\ref{cor:full_stab} and continuity imply that $\calB_0:=\calI-\calC_{Q_0}\cals $ is invertible.

        Consider the map
        \begin{align}
                F(Q):=(\calI-\calC_Q\cals)^{-1}[Q\wt I Q]. \nonumber
        \end{align}
        Since invertibility is an open condition, $F$ is holomorphic in a neighborhood of $Q_0$. Hence the initial value problem
        \begin{align}
                \partial_\omega Q(\omega)=F(Q(\omega)),\qquad Q(0)=Q_0, \nonumber
        \end{align}
        has a unique holomorphic solution for sufficiently small $|\omega|$. On the other hand, differentiating the MDE \eqref{eq:intro_mde} gives $ \partial_z\Pi(z)=(\calI-\calC_{\Pi(z)}\cals)^{-1}[\Pi(z)\wt I\Pi(z)]$.
        Therefore the continuous extension of $\Pi(E_0+\omega)$ to the upper half disc satisfies the same initial value problem. By uniqueness, $\Pi(E_0+\omega)=Q(\omega)$ there. Thus $\Pi$, and hence $m$, extends holomorphically to a complex neighborhood of $E_0$. Consequently,
        \begin{align}
                \rho(E)=\frac1\pi\Im m(E) \nonumber
        \end{align}
        is real analytic near $E_0$. Since $E_0$ was arbitrary, $\rho$ is real analytic on $\{E\in(0,\infty):\rho(E)>0\}$.
\end{proof}

\begin{proof}[Proof of Corollary~\ref{cor:bdry_reg_m_g}]
        Fix sufficiently small $\tau>0$. By \eqref{eq:holder_rho}, $\Im m$ admits a $1/3$-H\"older continuous extension to $\overline{\mb D(\tau/2,\eta_0)}$. By Lemma~\ref{lem:der_bd} and \eqref{eq:imag_comp}, we have
        \begin{align}
               |\partial_zg(z)|\lesssim(\Im m(z))^{-2}\lesssim(\Im g(z))^{-2},\qquad z\in\mb D(\tau/2,\eta_0). \nonumber
        \end{align}
        This also shows that $\Im g$ admits a $1/3$-H\"older continuous extension to $\overline{\mb D(\tau/2,\eta_0)}$.

        By Proposition~\ref{prop:pi},
        \begin{align}
                m(z)  = \int_{[0,\infty)}\frac{\bnu(\dd\la)}{\la-z}, \qquad
                g(z)  = -\avg{\Sg}_N + \int_{[0,\infty)} \frac{\la\avg{\Sg V_1(\dd\la)}_N}{\la-z}. \nonumber
        \end{align}
        The measure $\la\avg{\Sg V_1(\dd\la)}_N$ is positive and compactly supported. The bounded boundary values of its Stieltjes transform on $[\tau/2,2/\tau]$ exclude atoms there, and the Stieltjes inversion argument used in the proof of Theorem~\ref{thm:dens_reg} shows that its restriction to this interval is absolutely continuous with density $\pi^{-1}\Im g(E)$. Thus the two representing measures have bounded $1/3$-H\"older continuous densities on $[\tau/2,2/\tau]$, namely $\rho(E)$ and $\pi^{-1}\Im g(E)$.

        We may therefore apply \cite[Lemma~A.1]{alt2020dyson} on the interval $[\tau/2,2/\tau]$. Since $[\tau,\tau^{-1}]$ stays a positive distance from the boundary of $[\tau/2,2/\tau]$, for $z_1,z_2\in \overline{\mb D(\tau,\eta_0)}$ we obtain
        \begin{align}
                |m(z_1)-m(z_2)|+|g(z_1)-g(z_2)|\lesssim|z_1-z_2|^{1/3}. \nonumber
        \end{align}
        Passing to the boundary gives the same estimate for $z_1,z_2\in\overline{\mb D(\tau,\eta_0)}$. Since $\tau$ can be chosen arbitrarily small, this proves the claimed local $1/3$-H\"older continuity near the positive real axis. On compact subsets of $\C^+$ bounded away from the real axis, the Stieltjes representations give local Lipschitz continuity. Combining the two regimes yields unique locally $1/3$-H\"older continuous extensions to $\{E+\ii\eta:E>0,\ \eta\geq0\}$.
\end{proof}

\section{Analysis near a regular spectral edge}\label{sec:reg_edge}

Throughout this section, Assumption~\ref{ass:basic} holds, and $E_*$ denotes a regular right edge in the sense of Definition~\ref{def:reg_r_edge}. We prove Proposition~\ref{prop:sqrt_edge} in Subsection~\ref{sub:edge_scal_red} and Theorem~\ref{thm:sharp_stab} in Subsection~\ref{sub:sharp_stab}. Since $E_*\in\mb D_0$,
\begin{align}
        (E_*,E_*+\tau_0)\cap\supp\bnu=\varnothing. \label{eq:r_gap_ass}
\end{align}
Moreover, by Theorem~\ref{thm:det_spec},
\begin{align}
        \supp\bnu\cap(0,\infty)=\overline{\{E>0:\rho(E)>0\}}. \label{eq:supp_pos_dens}
\end{align}
We recall the definitions of the operators $\calJ$ and $\calB$ from \eqref{eq:J_def} and \eqref{eq:B_def}:
\begin{align}
        \calJ(z):=\calP\calC_{\Pi(z)}\calE, \qquad \calB(z):=\calI-\calC_{\Pi(z)}\cals. \nonumber
\end{align}
\begin{theorem}\label{thm:sharp_stab}
        Under Assumption~\ref{ass:basic}, let $E_*$ be a regular right edge in the sense of Definition~\ref{def:reg_r_edge}. Then there exists a constant $\eps_*>0$ such that, for every $z=E+\ii\eta\in\C^+$ with $|E-E_*|\leq\eps_*$ and $0<\eta\leq\eps_*$,
        \begin{align}
                \|(\calI-\calJ(z))^{-1}\|+\|\calB(z)^{-1}\|\lesssim(|E-E_*|+\eta)^{-1/2}. \label{eq:sharp_stab}
        \end{align}
        Moreover, with $\kappa:=|E-E_*|$,
        \begin{align}
                \Im m(E+\ii\eta)\sim \begin{cases} (\kappa+\eta)^{1/2},               & E\leq E_*, \\[1mm]
              \dfrac{\eta}{(\kappa+\eta)^{1/2}}, & E\geq E_*,\end{cases} \label{eq:edge_harm_scale}
        \end{align}
        and consequently
        \begin{align}
                \Im m(z)+\frac{\eta}{\Im m(z)}              & \sim (|E-E_*|+\eta)^{1/2}, \nonumber                                       \\
                \|(\calI-\calJ(z))^{-1}\|+\|\calB(z)^{-1}\| & \lesssim \frac{1}{\Im m(z)+\eta/\Im m(z)}. \label{eq:sharp_stab_im_m}
        \end{align}
        The same $\eps_*$ may be chosen so that, uniformly for $|E-E_*|\leq\eps_*$ and $\eps_*\leq\eta\leq\eta_0$,
        \begin{align}
                \|\Pi(z)\|+\|(\calI-\calJ(z))^{-1}\|+\|\calB(z)^{-1}\|\lesssim1. \label{eq:ord_edge_stab}
        \end{align}
\end{theorem}

Throughout this section, $\eps_*$ may be decreased without further comment, and we fix its value once and for all at the end of the section.
\subsection{Stability structure at the edge}

We identify the stability structure of the MDE at $E_*$. Lemma~\ref{lem:bdry_zero_dens} establishes the boundary value $\Pi(E_*)$ and local bounds for $\Pi$ near the edge, while Proposition~\ref{prop:unst_dir} identifies one unstable direction and a uniformly stable complementary direction. These results are used in Subsections~\ref{sub:edge_scal_red} and \ref{sub:sharp_stab}.

\begin{lemma}
        \label{lem:bdry_zero_dens}
        Under Assumption~\ref{ass:basic}, let $E_*$ be a regular right edge in the sense of Definition~\ref{def:reg_r_edge}. Then there exist constants $\eps_*>0$ and $c,C>0$ such that
        \begin{align}
                \inf_{\substack{|z-E_*|<\eps_*  \\ \Im z>0}}|z+g(z)|\geq c,
                \qquad
                \sup_{\substack{|z-E_*|<\eps_*  \\ \Im z>0}}\|\Pi(z)\|\leq C. \nonumber
        \end{align}
        Moreover,
        \begin{align}
                \|\Pi(z_1)-\Pi(z_2)\| \leq C|z_1-z_2|^{1/3} \nonumber
        \end{align}
        for all $z_1,z_2\in\C^+$ satisfying $|z_1-E_*|,|z_2-E_*|<\eps_*$. Consequently, $\Pi$ extends uniquely and continuously to $\{z\in\C:|z-E_*|<\eps_*,\ \Im z\geq0\}$.
        In particular,
        \begin{align}
                \Pi(E_*) := \lim_{\eta\downarrow0}\Pi(E_*+\ii\eta) \label{eq:Pi_bdry_value}
        \end{align}
        is well defined and is a real symmetric invertible matrix.
\end{lemma}

\begin{proof}
        By Corollary~\ref{cor:bdry_reg_m_g}, there exists $\eps_0>0$ such that
        \begin{align}
                |m(z)-m(E_*)|+|g(z)-g(E_*)|
                \lesssim |z-E_*|^{1/3}
                \label{eq:loc_hold_mg}
        \end{align}
        for $z\in\{z\in\C:|z-E_*|<\eps_0, \Im z\geq 0\}$.

        Choose $0<\eps_*\leq\eps_0$ sufficiently small. Then, for $z\in\{z\in\C:|z-E_*|<\eps_*, \Im z\geq 0\}$, we have
        \begin{align}
                \left\| AA^\top-(z+g(z))(I_M+m(z)\Sg)-[AA^\top-(E_*+g(E_*))(I_M+m(E_*)\Sg)]\right\|\lesssim \eps_*^{1/3}. \nonumber
        \end{align}
        Hence, by \eqref{eq:schur_nondeg},
        \begin{align}
                s_{\min}(AA^\top-(z+g(z))(I_M+m(z)\Sg) )\geq \frac{\tau_0}{2}, \label{eq:loc_schur_ndg}
        \end{align}
        for all $z\in\{z\in\C:|z-E_*|<\eps_*, \Im z\geq 0\}$.

        Since $\Im(z+g(z))>0$ for $z\in\C^+$, the lower right block in \eqref{eq:inverse_mde} is invertible. Taking the Schur complement gives
        \begin{align}
                \Pi_{11}(z) =(z+g(z)) [AA^\top-(z+g(z))(I_M+m(z)\Sg)]^{-1}. \nonumber
        \end{align}
        By $g(z)=\avg{\Sg\Pi_{11}(z)}_N$, we have
        \begin{align}
                z=(z+g(z))\left\{1-\avg{\Sg[AA^\top-(z+g(z))(I_M+m(z)\Sg)]^{-1}}_N\right\}.\nonumber
        \end{align}
        By  Assumption~\ref{ass:basic} and \eqref{eq:loc_schur_ndg}, we have
        \begin{align}
                \left|\avg{\Sg[AA^\top-(z+g(z))(I_M+m(z)\Sg)]^{-1}}_N\right|\lesssim1, \nonumber
        \end{align}
        for all $z\in\{z\in\C:|z-E_*|<\eps_*, \Im z> 0\}$. Since $E_*\in\mb D_0\subset[\tau_0,\tau_0^{-1}]$, we obtain
        \begin{align}
                \inf_{\substack{|z-E_*|<\eps_*  \\ \Im z>0}}|z+g(z)|\geq c. \label{eq:loc_zg_bd}
        \end{align}

        Combining \eqref{eq:loc_schur_ndg} and \eqref{eq:loc_zg_bd} with the Schur complement for \eqref{eq:inverse_mde} gives
        \begin{align}
                \sup_{\substack{|z-E_*|<\eps_*  \\ \Im z>0}}\|\Pi(z)\|\leq C.\nonumber
        \end{align}
        Moreover, the matrix on the right-hand side of \eqref{eq:inverse_mde} is $1/3$-H\"older continuous by \eqref{eq:loc_hold_mg}. The resolvent identity and the uniform bound on $\Pi$ therefore imply
        \begin{align}
                \|\Pi(z_1)-\Pi(z_2)\|\lesssim |z_1-z_2|^{1/3}\nonumber
        \end{align}
        for $z_j\in \{z\in\C: |z-E_*|<\eps_*, \Im z>0 \},j=1,2$.
        Hence $\Pi$ extends uniquely as a $1/3$-H\"older continuous function to $\{z\in\C:|z-E_*|<\eps_*, \Im z\geq0 \}$. Thus we can define $\Pi_*$ as in \eqref{eq:Pi_bdry_value}.

        Finally, $\rho(E_*)=0$ implies $\Im m(E_*)=0$. By \eqref{eq:imag_comp} and continuity, we also have $\Im g(E_*)=0$. At $z=E_*$, the matrix on the right-hand side of \eqref{eq:inverse_mde} is real symmetric. By \eqref{eq:loc_schur_ndg} and \eqref{eq:loc_zg_bd}, it is invertible. Hence $\Pi(E_*)$ is real symmetric and invertible.
\end{proof}

For convenience, we introduce the notation
\begin{align}
        \Pi_*:=\Pi(E_*) = \begin{pmatrix} \Pi_{*,11} & \Pi_{*,12} \\
                \Pi_{*,21} & \Pi_{*,22}\end{pmatrix}. \nonumber
\end{align}
We now describe the two eigendirections of $\calJ_*$. With respect to the canonical basis of $\C^2$,
\begin{align}
        \calJ_*:=\lim_{\eta\downarrow 0}\calJ(E_*+\ii\eta) =\begin{pmatrix} \avg{\Pi_{*,21}\Sg\Pi_{*,12}}_N    & \avg{\Pi_{*,22}^2}_N            \\[1mm]
                \avg{\Sg\Pi_{*,11}\Sg\Pi_{*,11}}_N & \avg{\Pi_{*,21}\Sg\Pi_{*,12}}_N\end{pmatrix}. \label{eq:J_star_mat}
\end{align}
By Lemma~\ref{lem:bdry_zero_dens}, all four entries in \eqref{eq:J_star_mat} are real and nonnegative. Set
\begin{align}
        \la_*:=\avg{\Pi_{*,21}\Sg\Pi_{*,12}}_N -\sqrt{\avg{\Pi_{*,22}^2}_N\avg{\Sg\Pi_{*,11}\Sg\Pi_{*,11}}_N}. \label{eq:lambda_star_def}
\end{align}
Whenever $ \avg{\Pi_{*,22}^2}_N\avg{\Sg\Pi_{*,11}\Sg\Pi_{*,11}}_N>0$, we can define
\begin{align}
        \f{r}_*         & :=\begin{pmatrix} \sqrt{\avg{\Pi_{*,22}^2}_N} \\[1mm]
                                    \sqrt{\avg{\Sg\Pi_{*,11}\Sg\Pi_{*,11}}_N}\end{pmatrix}, \qquad
        \wh{\f{r}}_*     :=\begin{pmatrix} \sqrt{\avg{\Pi_{*,22}^2}_N} \\[1mm]
                                   -\sqrt{\avg{\Sg\Pi_{*,11}\Sg\Pi_{*,11}}_N}\end{pmatrix}, \nonumber                                                                                                \\
        \f{\ell}_*      & :=\frac{1}{2\sqrt{\avg{\Pi_{*,22}^2}_N\avg{\Sg\Pi_{*,11}\Sg\Pi_{*,11}}_N}} \begin{pmatrix} \sqrt{\avg{\Sg\Pi_{*,11}\Sg\Pi_{*,11}}_N} \\[1mm]
                                                                                                             \sqrt{\avg{\Pi_{*,22}^2}_N}\end{pmatrix}, \nonumber                       \\
        \wh{\f{\ell}}_* & :=\frac{1}{2\sqrt{\avg{\Pi_{*,22}^2}_N\avg{\Sg\Pi_{*,11}\Sg\Pi_{*,11}}_N}} \begin{pmatrix} \sqrt{\avg{\Sg\Pi_{*,11}\Sg\Pi_{*,11}}_N} \\[1mm]
                                                                                                             -\sqrt{\avg{\Pi_{*,22}^2}_N}\end{pmatrix}. \label{eq:perron_vec_star}
\end{align}
A direct multiplication gives
\begin{align}
        \calJ_*\f{r}_*      & =(2\avg{\Pi_{*,21}\Sg\Pi_{*,12}}_N-\la_*)\f{r}_*, & \f{\ell}_*^*\calJ_*      & =(2\avg{\Pi_{*,21}\Sg\Pi_{*,12}}_N-\la_*)\f{\ell}_*^*, \nonumber \\
        \calJ_*\wh{\f{r}}_* & =\la_*\wh{\f{r}}_*,                                          & \wh{\f{\ell}}_*^*\calJ_* & =\la_*\wh{\f{\ell}}_*^*. \label{eq:two_eigvec_rel}
\end{align}
Moreover,
\begin{align}
        \f{\ell}_*^*\f{r}_* =\wh{\f{\ell}}_*^*\wh{\f{r}}_*=1, \qquad \f{\ell}_*^*\wh{\f{r}}_* =\wh{\f{\ell}}_*^*\f{r}_*=0. \nonumber
\end{align}
Thus $\{\f{r}_*,\wh{\f{r}}_*\}$ and $\{\f{\ell}_*,\wh{\f{\ell}}_*\}$ are dual bases, and every $\bbu\in\C^2$ admits the unique decomposition
\begin{align}
        \bbu=(\f{\ell}_*^*\bbu)\f{r}_*+(\wh{\f{\ell}}_*^*\bbu)\wh{\f{r}}_*. \label{eq:dual_basis_dec}
\end{align}
The next proposition shows that the product in the denominator of \eqref{eq:perron_vec_star} is uniformly bounded away from zero, so all these definitions are well.

\begin{proposition}
        \label{prop:unst_dir}
        Under Assumption~\ref{ass:basic}, let $E_*$ be a regular right edge in the sense of Definition~\ref{def:reg_r_edge}. Then
        \begin{align}
                \avg{\Pi_{*,21}\Sg\Pi_{*,12}}_N + & \sqrt{\avg{\Pi_{*,22}^2}_N \avg{\Sg\Pi_{*,11}\Sg\Pi_{*,11}}_N} =1, \label{eq:perron_root_one}                 \\
                \avg{\Pi_{*,22}^2}_N                                                                            \sim1, \qquad &\avg{\Sg\Pi_{*,11}\Sg\Pi_{*,11}}_N                                                              \sim1. \label{eq:unst_dir_rel}
        \end{align}
        Consequently, $1$ and $\la_*$ are the two simple eigenvalues of $\calJ_*$ and
        \begin{align}
                1-\la_* =2\sqrt{\avg{\Pi_{*,22}^2}_N\avg{\Sg\Pi_{*,11}\Sg\Pi_{*,11}}_N} \sim1. \label{eq:lambda_star_gap}
        \end{align}
        Moreover, the dual basis relations imply the spectral decomposition
        \begin{align}
                I_2 & =\f{r}_*\f{\ell}_*^*+\wh{\f{r}}_*\wh{\f{\ell}}_*^*, \qquad
                \calJ_*  =\f{r}_*\f{\ell}_*^*+\la_*\wh{\f{r}}_*\wh{\f{\ell}}_*^*. \label{eq:Jstar_spec_dec}
        \end{align}
        In particular,
        \begin{align}
                \|\f{r}_*\|+\|\wh{\f{r}}_*\|+\|\f{\ell}_*\|+\|\wh{\f{\ell}}_*\|\lesssim1. \label{eq:unst_vec_unif}
        \end{align}
\end{proposition}

\begin{proof}
        By \eqref{eq:supp_pos_dens}, choose $E_j\to E_*$ with $\rho(E_j)>0$, and set
        \begin{align}
                a_j:=\Im m(E_j), \qquad b_j:=\Im g(E_j). \nonumber
        \end{align}
        Then $a_j,b_j>0$. Taking the two normalized traces in the boundary version of Lemma~\ref{lem:im_ids} gives
        \begin{align}
                 & \begin{pmatrix} \avg{\Pi_{12}(E_j)^*\Sg\Pi_{12}(E_j)}_N    & \avg{\Pi_{22}(E_j)^*\Pi_{22}(E_j)}_N    \\[1mm]
                \avg{\Sg\Pi_{11}(E_j)^*\Sg\Pi_{11}(E_j)}_N & \avg{\Pi_{12}(E_j)^*\Sg\Pi_{12}(E_j)}_N\end{pmatrix} \begin{pmatrix}a_j \\
                                                                                                                                                                                     b_j\end{pmatrix} =\begin{pmatrix}a_j \\
                                                                                                                                                                                                       b_j\end{pmatrix}. \label{eq:dens_perron_eq}
        \end{align}
        Conjugating the matrix in \eqref{eq:dens_perron_eq} in the form $\diag(a_j,b_j)^{-1}(\cdot)\diag(a_j,b_j)$ produces a nonnegative matrix whose row sums are one. Hence its spectral radius is one. The larger eigenvalue of the matrix in \eqref{eq:dens_perron_eq} is
        \begin{align}
                 & \avg{\Pi_{12}(E_j)^*\Sg\Pi_{12}(E_j)}_N +\sqrt{\avg{\Pi_{22}(E_j)^*\Pi_{22}(E_j)}_N \avg{\Sg\Pi_{11}(E_j)^*\Sg\Pi_{11}(E_j)}_N}. \nonumber
        \end{align}
        It is therefore equal to one. Letting $E_j\to E_*$ and using Lemma~\ref{lem:bdry_zero_dens} proves \eqref{eq:perron_root_one}.

        Taking the boundary limit in \eqref{eq:Pi_21_bd} gives
        \begin{align}
                0\leq\avg{\Pi_{*,21}\Sg\Pi_{*,12}}_N\leq1-c_1. \nonumber
        \end{align}
        Together with \eqref{eq:perron_root_one}, this implies
        \begin{align}
                \sqrt{\avg{\Pi_{*,22}^2}_N\avg{\Sg\Pi_{*,11}\Sg\Pi_{*,11}}_N}\geq c_1. \nonumber
        \end{align}
        On the other hand, Lemma~\ref{lem:bdry_zero_dens} and Assumption~\ref{ass:basic} give
        \begin{align}
                \avg{\Pi_{*,22}^2}_N+\avg{\Sg\Pi_{*,11}\Sg\Pi_{*,11}}_N\lesssim1. \nonumber
        \end{align}
        Hence both factors are bounded from above and below by positive constants, proving \eqref{eq:unst_dir_rel}. In particular, the vectors in \eqref{eq:perron_vec_star} are well defined and uniformly bounded. Combining \eqref{eq:lambda_star_def} with \eqref{eq:perron_root_one} gives \eqref{eq:lambda_star_gap}, while \eqref{eq:two_eigvec_rel} identifies $1$ and $\la_*$ as the two eigenvalues.

        The spectral decomposition \eqref{eq:Jstar_spec_dec} follows from \eqref{eq:two_eigvec_rel} and \eqref{eq:dual_basis_dec}. The uniform vector bounds follow directly from \eqref{eq:unst_dir_rel}.
\end{proof}

\subsection{Scalar reduction and square-root behavior}\label{sub:edge_scal_red}

This subsection proves Proposition~\ref{prop:sqrt_edge}. We first reduce the MDE near $E_*$ to a scalar normal form for the unstable coordinate $\Theta$. Proposition~\ref{prop:scalar_nf} and Lemma~\ref{lem:nf_coeff} yield Proposition~\ref{prop:edge_sqrt_exp}, from which Proposition~\ref{prop:sqrt_edge} follows. The proofs of Proposition~\ref{prop:scalar_nf} and Lemma~\ref{lem:nf_coeff} are given after the proof of Proposition~\ref{prop:sqrt_edge}.

Let $z$ lie in the closed upper half-plane neighborhood from Lemma~\ref{lem:bdry_zero_dens}, and set
\begin{align}
        \Delta:=\Pi(z)-\Pi_*. \nonumber
\end{align}
Then
\begin{align}
        \calP[\Delta]=\begin{pmatrix}m(z)-m(E_*) \\
                              g(z)-g(E_*)\end{pmatrix}. \nonumber
\end{align}
With respect to the unstable and stable directions from Proposition~\ref{prop:unst_dir}, we decompose
\begin{align}
        \calP[\Delta]=\Theta\f{r}_*+\bbs_*, \qquad \Theta:=\f{\ell}_*^*\calP[\Delta], \qquad \bbs_*:=(\wh{\f{\ell}}_*^*\calP[\Delta])\wh{\f{r}}_*. \label{eq:unst_dec}
\end{align}
Thus $\Theta$ is the coordinate in the unstable direction, while $\bbs_*$ is the stable component.

To formulate the scalar equation, for $\bbu_1,\bbu_2,\bbu_3\in\C^2$ define
\begin{align}
        N_2[\bbu_1,\bbu_2]        & :=\calP[\Pi_*\calE[\bbu_1]\Pi_*\calE[\bbu_2]\Pi_*], \nonumber         \\
        N_3[\bbu_1,\bbu_2,\bbu_3] & :=\calP[\Pi_*\calE[\bbu_1]\Pi_*\calE[\bbu_2]\Pi_*\calE[\bbu_3]\Pi_*]. \nonumber
\end{align}
We further introduce the coefficients
\begin{align}
        \mu_1 & :=\f{\ell}_*^*\calP[\Pi_*\wt I\Pi_*], \qquad
        \mu_2  :=\f{\ell}_*^*N_2[\f{r}_*,\f{r}_*], \nonumber                                                                                                                                                                      \\
        \mu_3 & :=\f{\ell}_*^*N_3[\f{r}_*,\f{r}_*,\f{r}_*]+\frac{\wh{\f{\ell}}_*^*N_2[\f{r}_*,\f{r}_*]}{1-\la_*}\f{\ell}_*^*(N_2[\f{r}_*,\wh{\f{r}}_*]+N_2[\wh{\f{r}}_*,\f{r}_*]). \label{eq:nf_coeff}
\end{align}
The inverse form of the MDE \eqref{eq:inverse_mde} yields
\begin{align}
        \Pi(z)^{-1}=\Pi_*^{-1}-(z-E_*)\wt I-\calE[\calP[\Delta]]. \label{eq:inv_diff_edge}
\end{align}
Expanding the resolvent identity associated with \eqref{eq:inv_diff_edge} to third order and applying $\calP$ gives the following reduction.

\begin{proposition}
        \label{prop:scalar_nf}
        Under Assumption~\ref{ass:basic}, let $E_*$ be a regular right edge in the sense of Definition~\ref{def:reg_r_edge}. For $\eps_*>0$ sufficiently small, every $z$ with $|z-E_*|<\eps_*$ and $\Im z\geq0$ satisfies
        \begin{align}
                |\Theta|+\|\bbs_*\| & \lesssim|z-E_*|^{1/3}\leq\eps_*^{1/3}, \label{eq:unst_cmp_small} \\
                \|\bbs_*\|          & \lesssim|z-E_*|+|\Theta|^2. \label{eq:stab_cmp_bd}
        \end{align}
        Moreover, $\Theta$ satisfies
        \begin{align}
                0=\mu_1(z-E_*)+\mu_2\Theta^2+\mu_3\Theta^3+O(|z-E_*||\Theta|+|z-E_*|^2+|\Theta|^4), \label{eq:scalar_nf}
        \end{align}
        and the coefficients defined above are real and satisfy
        \begin{align}
                c\leq\mu_1\leq C, \qquad |\mu_2|+|\mu_3|\leq C. \label{eq:mu_bd}
        \end{align}
\end{proposition}

The following lemma gives the sign and nondegeneracy properties of the coefficients needed below.

\begin{lemma}
        \label{lem:nf_coeff}
        Under Assumption~\ref{ass:basic}, let $E_*$ be a regular right edge in the sense of Definition~\ref{def:reg_r_edge}. Then there exists a constant $c_\mu>0$ such that
        \begin{align}
                \mu_3  \geq0, \qquad
                \mu_2  \leq-c_\mu. \label{eq:mu2_neg}
        \end{align}
\end{lemma}

Proposition~\ref{prop:scalar_nf} and Lemma~\ref{lem:nf_coeff} yield the following square-root expansion.

\begin{proposition}
        \label{prop:edge_sqrt_exp}
        Under Assumption~\ref{ass:basic}, let $E_*$ be a regular right edge in the sense of Definition~\ref{def:reg_r_edge}. For $\eps_*>0$ sufficiently small,
        \begin{align}
                |\Theta(z)|^2\lesssim|z-E_*| \label{eq:Theta_sqrt_up}
        \end{align}
        for $|z-E_*|<\eps_*$ and $\Im z\geq0$, and
        \begin{align}
                \Theta^2 =-\frac{\mu_1}{\mu_2}(z-E_*)+O(|z-E_*|^{3/2}). \label{eq:Theta_sqrt}
        \end{align}
        Moreover,
        \begin{align}
                \rho(E) =\frac1\pi\sqrt{\avg{\Pi_{*,22}^2}_N} \left(\frac{\mu_1}{|\mu_2|}\right)^{1/2}(E_*-E)^{1/2} +O(E_*-E) \label{eq:dens_sqrt}
        \end{align}
        for $E\in(E_*-\eps_*,E_*)$, while $\rho(E)=0$ for $E\in(E_*,E_*+\eps_*)$.
\end{proposition}

\begin{proof}
        Since $m(E_*),g(E_*)\in\R$ and both entries of $\f{\ell}_*$ are strictly positive,
        \begin{align}
                \Im\Theta(z) =\f{\ell}_{*,1}\Im m(z)+\f{\ell}_{*,2}\Im g(z)>0,\qquad z\in\C^+. \nonumber
        \end{align}
        Hence $\Im\Theta(E)\geq0$ for every real $E$ sufficiently close to $E_*$. Moreover, \eqref{eq:unst_cmp_small} gives
        \begin{align}
                |\Theta|\lesssim|z-E_*|^{1/3}. \label{eq:Theta_hold}
        \end{align}

        Lemma~\ref{lem:nf_coeff} gives $\mu_2\leq-c_\mu$. Since we can choose $\eps_*\leq 1$, with \eqref{eq:Theta_hold}, taking absolute values in \eqref{eq:scalar_nf} gives $|\Theta|^2 \lesssim |z-E_*|$, which proves \eqref{eq:Theta_sqrt_up}. By \eqref{eq:Theta_sqrt_up}, the error terms in \eqref{eq:scalar_nf} satisfy
        \begin{align}
                |\Theta|^3+|z-E_*||\Theta|+|z-E_*|^2+|\Theta|^4
                =O(|z-E_*|^{3/2}). \nonumber
        \end{align}
        Substituting these estimates into \eqref{eq:scalar_nf} yields \eqref{eq:Theta_sqrt}.

        Let $E<E_*$ and set $\kappa:=E_*-E$. Since $\mu_2<0$, taking boundary values in \eqref{eq:Theta_sqrt} gives
        \begin{align}
                \Theta(E)^2 =-\frac{\mu_1}{|\mu_2|}\kappa+O(\kappa^{3/2}). \nonumber
        \end{align}
        Together with $\Im\Theta(E)\geq0$, this implies
        \begin{align}
                \Theta(E) =\ii\left(\frac{\mu_1}{|\mu_2|}\right)^{1/2}\kappa^{1/2}+O(\kappa). \nonumber
        \end{align}
        By \eqref{eq:stab_cmp_bd} and \eqref{eq:Theta_sqrt_up}, $\|\bbs_*(E)\|\lesssim \kappa$.
        The first component of \eqref{eq:unst_dec} gives
        \begin{align}
                m(E)-m(E_*) =\sqrt{\avg{\Pi_{*,22}^2}_N}\Theta(E)+ O(\kappa). \nonumber
        \end{align}
        Taking imaginary parts proves \eqref{eq:dens_sqrt}. The vanishing of $\rho$ on the right follows directly from \eqref{eq:r_gap_ass}.
\end{proof}

\begin{proof}[Proof of Proposition~\ref{prop:sqrt_edge}]
        Proposition~\ref{prop:edge_sqrt_exp}, together with the continuity of $\rho$, then yields the square-root behavior of $\rho$ on $[E_*-\eps_*,E_*]$ for some constant $\eps_*>0$. Finally, \eqref{eq:unst_dir_rel}, \eqref{eq:mu_bd}, and \eqref{eq:mu2_neg} show that the leading coefficient in \eqref{eq:dens_sqrt} is uniformly comparable to one, which proves \eqref{eq:sqrt_edge}.
\end{proof}

It remains to prove Proposition~\ref{prop:scalar_nf} and Lemma~\ref{lem:nf_coeff}. We do this in turn.

\begin{proof}[Proof of Proposition~\ref{prop:scalar_nf}]
        By Corollary~\ref{cor:bdry_reg_m_g},
        \begin{align}
                \|\calP[\Delta]\|\lesssim |z-E_*|^{1/3} \leq\eps_*^{1/3}. \label{eq:proj_holder}
        \end{align}
        Taking $\eps_*$ sufficiently small, we may assume that throughout $|z-E_*|<\eps_*$,
        \begin{align}
                |z-E_*|+\|\calP[\Delta]\|\leq1. \label{eq:loc_small}
        \end{align}
        The inverse identity \eqref{eq:inv_diff_edge} gives
        \begin{align}
                \Pi(z)-\Pi_* &=\Pi_*((z-E_*)\wt I+\calE[\calP[\Delta]])\Pi(z)  \nonumber\\
                             &= \sum_{k=1}^3 \Pi_*\left(((z-E_*)\wt I+\calE[\calP[\Delta]])\Pi_*\right)^k\\
                             &\quad +\Pi_*\left(((z-E_*)\wt I+\calE[\calP[\Delta]])\Pi_*\right)^3((z-E_*)\wt I+\calE[\calP[\Delta]])\Pi(z)\nonumber\\
                             &:=\sum_{k=1}^3 T_*^{(k)}+\calR_*. \label{eq:cmpd_exp}
        \end{align}
        Lemma~\ref{lem:bdry_zero_dens} and the boundedness of $\calP$ and $\calE$ show that 
        \begin{align}
                \|\calP[\calR_*]\|\lesssim (|z-E_*|+\|\calP[\Delta]\|)^4, \qquad |z-E_*|<\eps_*.\nonumber
        \end{align}
        Using the definition of $\calJ$, we have $\calP[\Pi_*\calE[\calP[\Delta]]\Pi_*]=\calJ_*[\calP[\Delta]]$.
        Moreover, by the boundedness of $\calP$ and $\calE$ and \eqref{eq:loc_small}, all terms in \eqref{eq:cmpd_exp} containing at least one factor $z-E_*$, except for the linear term in $z-E_*$, are bounded by
        \begin{align}
                O(|z-E_*|\|\calP[\Delta]\|+|z-E_*|^2). \nonumber
        \end{align}
        By the definitions of $N_2$ and $N_3$, and moving the linear term $\calP[\Pi_*\calE[\calP[\Delta]]\Pi_*]$ to the left in \eqref{eq:cmpd_exp}, we obtain
        \begin{align}
                (\calI-\calJ_*)[\calP[\Delta]] & =(z-E_*)\calP[\Pi_*\wt I\Pi_*] +N_2[\calP[\Delta],\calP[\Delta]] +N_3[\calP[\Delta],\calP[\Delta],\calP[\Delta]] \nonumber\\ 
                &\quad +O(|z-E_*|\|\calP[\Delta]\|+|z-E_*|^2+\|\calP[\Delta]\|^4). \label{eq:cmpd_eq}
        \end{align}
        A direct block computation gives
        \begin{align}
                \calP[\Pi_*\wt I\Pi_*] =\begin{pmatrix} \avg{\Pi_{*,22}^2}_N \\[1mm]
                                                \avg{\Pi_{*,21}\Sg\Pi_{*,12}}_N\end{pmatrix}. \label{eq:h_star_coord}
        \end{align}

        By \eqref{eq:unst_dec}, \eqref{eq:unst_vec_unif}, and \eqref{eq:proj_holder},
        \begin{align}
                |\Theta|+\|\bbs_*\|\lesssim\|\calP[\Delta]\|\lesssim|z-E_*|^{1/3}\leq\eps_*^{1/3}, \label{eq:edge_small}
        \end{align}
        which proves \eqref{eq:unst_cmp_small}.

        Applying $\wh{\f{\ell}}_*^*$ to \eqref{eq:cmpd_eq} and using $\wh{\f{\ell}}_*^*(\calI-\calJ_*)=(1-\la_*)\wh{\f{\ell}}_*^*$ gives
        \begin{align}
                  (1-\la_*)\wh{\f{\ell}}_*^*\calP[\Delta] &=(z-E_*)\wh{\f{\ell}}_*^*\calP[\Pi_*\wt I\Pi_*] +\wh{\f{\ell}}_*^*N_2[\calP[\Delta],\calP[\Delta]]+\wh{\f{\ell}}_*^*N_3[\calP[\Delta],\calP[\Delta],\calP[\Delta]]  \nonumber                            \\
                 & \quad+O(|z-E_*|\|\calP[\Delta]\|+|z-E_*|^2+\|\calP[\Delta]\|^4). \label{eq:stab_scal_eq}
        \end{align}
        Using $1-\la_*\sim1$, the boundedness of $N_2$ and $N_3$, and \eqref{eq:loc_small}, we obtain
        \begin{align}
                |\wh{\f{\ell}}_*^*\calP[\Delta]| \lesssim |z-E_*|+\|\calP[\Delta]\|^2. \nonumber
        \end{align}
        Since $ \|\calP[\Delta]\| \lesssim |\Theta|+|\wh{\f{\ell}}_*^*\calP[\Delta]|$ and \eqref{eq:edge_small}, it follows that
        \begin{align}
                |\wh{\f{\ell}}_*^*\calP[\Delta]| \lesssim |z-E_*|+|\Theta|^2+|\wh{\f{\ell}}_*^*\calP[\Delta]|^2\lesssim |z-E_*|+|\Theta|^2. \nonumber
        \end{align}
        The definition of $\bbs_*$ and \eqref{eq:unst_vec_unif} therefore prove \eqref{eq:stab_cmp_bd}.

        Expanding the quadratic and cubic terms in \eqref{eq:stab_scal_eq}, and using the definition of $\bbs_*$ in \eqref{eq:unst_dec} and \eqref{eq:stab_cmp_bd}, we obtain
        \begin{align}
                \bbs_* & =\frac{z-E_*}{1-\la_*}(\wh{\f{\ell}}_*^*\calP[\Pi_*\wt I\Pi_*])\wh{\f{r}}_* +\frac{\Theta^2}{1-\la_*}(\wh{\f{\ell}}_*^*N_2[\f{r}_*,\f{r}_*])\wh{\f{r}}_* + O(|z-E_*||\Theta|+|z-E_*|^2+|\Theta|^3).\label{eq:s_star_2nd}
        \end{align}
        Applying $\f{\ell}_*^*$ to \eqref{eq:cmpd_eq} eliminates the left-hand side. In the quadratic term, the first term on the right-hand side of \eqref{eq:s_star_2nd} produces only terms containing a factor $z-E_*$, while the second term produces the term proportional to $\Theta^3$. Hence
        \begin{align}
                \f{\ell}_*^*N_2[\calP[\Delta],\calP[\Delta]]               & =\Theta^2\f{\ell}_*^*N_2[\f{r}_*,\f{r}_*] +\frac{\Theta^3\wh{\f{\ell}}_*^*N_2[\f{r}_*,\f{r}_*]}{1-\la_*}\f{\ell}_*^*(N_2[\f{r}_*,\wh{\f{r}}_*]+N_2[\wh{\f{r}}_*,\f{r}_*]) \nonumber \\
                                                                           & \quad+O(|z-E_*||\Theta|+|z-E_*|^2+|\Theta|^4), \nonumber                                                                                                                            \\
                \f{\ell}_*^*N_3[\calP[\Delta],\calP[\Delta],\calP[\Delta]] & =\Theta^3\f{\ell}_*^*N_3[\f{r}_*,\f{r}_*,\f{r}_*] +O(|z-E_*||\Theta|+|z-E_*|^2+|\Theta|^4). \nonumber
        \end{align}
        Moreover, using \eqref{eq:unst_dec} and \eqref{eq:stab_cmp_bd}, the error term in \eqref{eq:stab_scal_eq} is bounded by
        \begin{align}
                O(|z-E_*||\Theta|+|z-E_*|^2+|\Theta|^4). \nonumber
        \end{align}
        Comparing the pure quadratic and cubic terms with \eqref{eq:nf_coeff} proves \eqref{eq:scalar_nf}.
        Since $\Pi_*$, $\f{r}_*$, $\wh{\f{r}}_*$, $\f{\ell}_*$, and $\wh{\f{\ell}}_*$ are real, the coefficients $\mu_1$, $\mu_2$, and $\mu_3$ are real. Furthermore, \eqref{eq:h_star_coord} and \eqref{eq:perron_vec_star} give
        \begin{align}
                \mu_1 =\frac12\left(\sqrt{\avg{\Pi_{*,22}^2}_N} +\frac{\avg{\Pi_{*,21}\Sg\Pi_{*,12}}_N}{\sqrt{\avg{\Sg\Pi_{*,11}\Sg\Pi_{*,11}}_N}} \right) \sim1. \nonumber
        \end{align}
        The upper bounds for $\mu_2$ and $\mu_3$ follow from \eqref{eq:unst_dir_rel}, \eqref{eq:lambda_star_gap}, and the uniform bounds on $\Pi_*$, $\f{r}_*$, $\wh{\f{r}}_*$, $\f{\ell}_*$, and $\wh{\f{\ell}}_*$.
\end{proof}

\begin{proof}[Proof of Lemma~\ref{lem:nf_coeff}]
        \emph{The sign of $\mu_3$.}
        Define
        \begin{align}
                \calA_* & :=\calE[\f{r}_*]^{1/2}\Pi_*\calE[\f{r}_*]^{1/2}, \qquad
                \calU_*  :=\calE[\f{r}_*]^{-1/2}\calE[\wh{\f{r}}_*]\calE[\f{r}_*]^{-1/2}. \nonumber
        \end{align}
        Since $\f{r}_*>0$, the matrix $\calE[\f{r}_*]$ is positive definite. Lemma~\ref{lem:bdry_zero_dens} implies that $\calA_*$ is self-adjoint, while the block forms of $\calE[\f{r}_*]$ and $\calE[\wh{\f{r}}_*]$ give
        \begin{align}
                \calU_*^*=\calU_*, \qquad \calU_*^2=I_L.\label{eq:D_star_invol}
        \end{align}
        For every block matrix $X$, the explicit formulas in \eqref{eq:perron_vec_star}, together with \eqref{eq:lambda_star_gap}, yield
        \begin{align}
                \f{\ell}_*^*\calP[X] & =\frac{1}{1-\la_*}\frac1N\Tr(\calE[\f{r}_*]X),\qquad
                \wh{\f{\ell}}_*^*\calP[X] =-\frac{1}{1-\la_*}\frac1N\Tr(\calE[\wh{\f{r}}_*]X).\label{eq:dual_tr}
        \end{align}

        We first record two trace identities. Since $\calJ_*\f{r}_*=\f{r}_*$ and $\f{\ell}_*^*\f{r}_*=1$, \eqref{eq:dual_tr} gives
        \begin{align}
                1 =\f{\ell}_*^*\calJ_*\f{r}_* =\f{\ell}_*^*\calP[\Pi_*\calE[\f{r}_*]\Pi_*] =\frac{1}{1-\la_*}\frac1N\Tr\calA_*^2, \nonumber
        \end{align}
        where we also used $\calJ=\calP\calC_\Pi\calE$.
        Hence
        \begin{align}
                \frac1N\Tr\calA_*^2=1-\la_*.\label{eq:Astar_sq_tr}
        \end{align}
        Similarly, using $\calJ_*\wh{\f{r}}_*=\la_*\wh{\f{r}}_*$, we obtain
        \begin{align}
                \la_* =\wh{\f{\ell}}_*^*\calJ_*\wh{\f{r}}_* =-\frac{1}{1-\la_*}\frac1N\Tr\calU_*\calA_*\calU_*\calA_*. \nonumber
        \end{align}
        Therefore
        \begin{align}
                \frac1N\Tr\calU_*\calA_*\calU_*\calA_*=-\la_*(1-\la_*).\label{eq:D_A_D_A_trace}
        \end{align}

        We next compute the scalar terms appearing in the definition of $\mu_3$. By \eqref{eq:dual_tr}, the definition of $N_2$, and cyclicity of the trace,
        \begin{align}
                \wh{\f{\ell}}_*^*N_2[\f{r}_*,\f{r}_*] =-\frac{1}{1-\la_*}\frac1N \Tr( \calE[\wh{\f{r}}_*]\Pi_* \calE[\f{r}_*]\Pi_* \calE[\f{r}_*]\Pi_* ) =-\frac{1}{1-\la_*}\frac1N \Tr(\calU_*\calA_*^3). \label{eq:stable_coord_N2}
        \end{align}
        In the same way,
        \begin{align}
                \f{\ell}_*^*N_3[\f{r}_*,\f{r}_*,\f{r}_*] =\frac{1}{1-\la_*}\frac1N\Tr\calA_*^4,\label{eq:N3_unst_coord}
        \end{align}
        while cyclicity of the trace gives
        \begin{align}
                \f{\ell}_*^*N_2[\f{r}_*,\wh{\f{r}}_*] =\f{\ell}_*^*N_2[\wh{\f{r}}_*,\f{r}_*] =\frac{1}{1-\la_*}\frac1N \Tr(\calU_*\calA_*^3). \label{eq:mixed_N2_coord}
        \end{align}
        Substituting \eqref{eq:stable_coord_N2}, \eqref{eq:N3_unst_coord}, and \eqref{eq:mixed_N2_coord} into \eqref{eq:nf_coeff}, we obtain
        \begin{align}
                \mu_3 =\frac{1}{1-\la_*}\frac1N\Tr\calA_*^4 -\frac{2}{(1-\la_*)^3} \left(\frac1N\Tr\calU_*\calA_*^3\right)^2. \label{eq:mu_3_exact_id}
        \end{align}

        By cyclicity of the trace,
        \begin{align}
                \frac1N\Tr\calU_*\calA_*^3 =\frac1N\Tr\left( \calA_*^2 \frac{\calU_*\calA_*+\calA_*\calU_*}{2} \right).\nonumber
        \end{align}
        Hence the Cauchy--Schwarz inequality yields
        \begin{align}
                \left(\frac1N\Tr\calU_*\calA_*^3\right)^2 \leq \left(\frac1N\Tr\calA_*^4\right) \frac1N\Tr\left(\frac{\calU_*\calA_*+\calA_*\calU_*}{2} \right)^2.\label{eq:cubic_trace_CS}
        \end{align}
        Using \eqref{eq:D_star_invol}, cyclicity of the trace, \eqref{eq:Astar_sq_tr}, and \eqref{eq:D_A_D_A_trace}, we obtain
        \begin{align}
                \frac1N\Tr\left(\frac{\calU_*\calA_*+\calA_*\calU_*}{2}\right)^2  &=\frac12\frac1N\Tr\calA_*^2 +\frac12\frac1N\Tr\calU_*\calA_*\calU_*\calA_*  =\frac12(1-\la_*) -\frac12\la_*(1-\la_*) \nonumber\\ 
                &=\frac{(1-\la_*)^2}{2}. \label{eq:sym_tr_bd}
        \end{align}
        Combining \eqref{eq:cubic_trace_CS} and \eqref{eq:sym_tr_bd} gives
        \begin{align}
                \left(\frac1N\Tr\calU_*\calA_*^3\right)^2 \leq \frac{(1-\la_*)^2}{2} \frac1N\Tr\calA_*^4.\nonumber
        \end{align}
        Inserting this estimate into \eqref{eq:mu_3_exact_id} proves $\mu_3\geq0$.

        \emph{The bound of $\mu_2$.}
        Let
        \begin{align}
                0<\kappa<\min\left\{\frac{\tau_0}{2},\frac{\eps_*}{2}\right\} \label{eq:gap_par_range}
        \end{align}
        be chosen below, where $\eps_*$ is the uniform radius in Proposition~\ref{prop:scalar_nf}. We evaluate the boundary value of \eqref{eq:scalar_nf} at $E_*+\kappa$.

        By \eqref{eq:imag_comp} and \eqref{eq:r_gap_ass},
        \begin{align}
                \Im g(E)=0, \qquad E\in(E_*,E_*+\tau_0). \nonumber
        \end{align}
        Since $g$ has a H\"older-continuous real boundary value on this interval, its representing measure has no mass there and its Stieltjes representation extends real analytically across the gap. Hence the Stieltjes representations from Proposition~\ref{prop:pi} may be differentiated for real $E$ in the gap, giving
        \begin{align}
                \partial_E m(E) = \int_{[0,\infty)}\frac{\bnu(\dd\la)}{(\la-E)^2}\geq0, \qquad \partial_E g(E) = \frac1N\Tr\int_{[0,\infty)}\frac{\la\Sg V_1(\dd\la)}{(\la-E)^2}\geq0. \nonumber
        \end{align}
        By the continuity of the boundary values at $E_*$, it follows that
        \begin{align}
                m(E_*+\kappa)-m(E_*)\geq0, \qquad g(E_*+\kappa)-g(E_*)\geq0. \nonumber
        \end{align}
        Since both entries of $\f{\ell}_*$ are strictly positive,
        \begin{align}
                \Theta(E_*+\kappa)\geq0. \label{eq:Theta_nonneg}
        \end{align}
        At this real spectral parameter, $m$, $g$, $\Pi$, and hence every term in the reduced scalar equation are real.

        By \eqref{eq:unst_cmp_small}, $|\Theta(E_*+\kappa)|\lesssim\kappa^{1/3}$. Hence the error term in \eqref{eq:scalar_nf} satisfies
        \begin{align}
                 & \kappa\left|\Theta(E_*+\kappa)\right| +\kappa^2 +\left|\Theta(E_*+\kappa)\right|^4 \lesssim \kappa^{4/3}.\nonumber
        \end{align}
        Taking the boundary value of \eqref{eq:scalar_nf} at $E_*+\kappa$, we obtain
        \begin{align}
                0 = \mu_1\kappa +\mu_2\Theta(E_*+\kappa)^2 +\mu_3\Theta(E_*+\kappa)^3 +O(\kappa^{4/3}). \label{eq:gap_scal_eq}
        \end{align}

        Choose $\kappa$ satisfying \eqref{eq:gap_par_range} sufficiently small, depending only on the standing constants and $\tau_0$. Then choose $c_\mu>0$ sufficiently small, depending only on $\kappa$ and the standing constants. Suppose that $\mu_2\geq-c_\mu$.  By \eqref{eq:mu_bd}, \eqref{eq:gap_par_range}, \eqref{eq:Theta_nonneg}, \eqref{eq:gap_scal_eq} and $\mu_3\geq0$, we obtain $ 0  \geq c\kappa-c_\mu|\Theta|^2-C\kappa^{4/3} \gtrsim c/2\kappa>0$, which is a contradiction. Hence $\mu_2\leq-c_\mu $. This proves \eqref{eq:mu2_neg}.

\end{proof}

\subsection{Sharp stability near the edge}\label{sub:sharp_stab}

This subsection proves Theorem~\ref{thm:sharp_stab}. Lemma~\ref{lem:red_stab_exp} gives the expansion of $\calI-\calJ(z)$ in the unstable--stable basis, and Proposition~\ref{prop:sharp_red} gives the inverse bound for $\calI-\calJ(z)$. We first prove Proposition~\ref{prop:sharp_red} and Theorem~\ref{thm:sharp_stab}, and then prove Lemma~\ref{lem:red_stab_exp}. Recall that
\begin{align}
        \calP[\Pi(z)-\Pi_*]=\Theta(z)\f{r}_*+\bbs_*(z), \qquad \|\bbs_*(z)\|\lesssim |z-E_*|+|\Theta(z)|^2. \nonumber
\end{align}

\begin{lemma}
        \label{lem:red_stab_exp}
        Under Assumption~\ref{ass:basic}, let $E_*$ be a regular right edge in the sense of Definition~\ref{def:reg_r_edge}. For $z\in\C^+$ with $|z-E_*|<\eps_*$,
        \begin{align}
                \left\|\Pi(z)-\Pi_* -\Theta(z)\Pi_*\calE[\f{r}_*]\Pi_*\right\| \lesssim |\Theta(z)|^2. \label{eq:Pi_Theta_exp}
        \end{align}
        With respect to the basis $\{\f{r}_*,\wh{\f{r}}_*\}$ and its dual basis $\{\f{\ell}_*,\wh{\f{\ell}}_*\}$, the matrix representation of $\calI-\calJ(z)$ is
        \begin{align}
                \begin{pmatrix} -2\mu_2\Theta(z)+O(|\Theta(z)|^2) & O(|\Theta(z)|)         \\
                O(|\Theta(z)|)                    & 1-\la_*+O(|\Theta(z)|)\end{pmatrix}. \label{eq:stab_mat}
        \end{align}
\end{lemma}

The lemma identifies the small entry in the matrix representation \eqref{eq:stab_mat}. We first derive the stability bounds and prove the lemma afterwards.

\begin{proposition}
        \label{prop:sharp_red}
        Under Assumption~\ref{ass:basic}, let $E_*$ be a regular right edge in the sense of Definition~\ref{def:reg_r_edge}. Then the operator $\calI-\calJ(z)$ is invertible for $z\in\C^+$ with $|z-E_*|<\eps_*$, and
        \begin{align}
                \|(\calI-\calJ(z))^{-1}\| \lesssim\frac1{|\Theta(z)|}. \label{eq:red_invTh_bd}
        \end{align}
\end{proposition}

\begin{proof}
        By Lemma~\ref{lem:nf_coeff}, $-\mu_2\sim1$.
        Moreover, \eqref{eq:lambda_star_gap} gives $1-\la_*\sim1$. Hence \eqref{eq:stab_mat} yields
        \begin{align}
                \det(\calI-\calJ(z)) =-2\mu_2(1-\la_*)\Theta(z)+O(|\Theta(z)|^2), \nonumber
        \end{align}
        and therefore
        \begin{align}
                |\det(\calI-\calJ(z))|\sim|\Theta(z)|. \label{eq:stab_det}
        \end{align}
        The entries in \eqref{eq:stab_mat} are uniformly bounded, and the bases $\{\f{r}_*,\wh{\f{r}}_*\}$ and $\{\f{\ell}_*,\wh{\f{\ell}}_*\}$ are uniformly well conditioned by \eqref{eq:unst_vec_unif}. The explicit inverse formula for a $2\times2$ matrix and \eqref{eq:stab_det} prove \eqref{eq:red_invTh_bd}.
\end{proof}

Proposition~\ref{prop:sharp_red}, Lemma~\ref{lem:stab_red}, \eqref{eq:C_Pi_bd}, and the local boundedness of $\Pi(z)$ give $\|\calB(z)^{-1}\|\lesssim|\Theta(z)|^{-1}$.

\begin{proof}[Proof of Theorem~\ref{thm:sharp_stab}]
        Choose $\eps_*$ sufficiently small so that Proposition~\ref{prop:edge_sqrt_exp} and Proposition~\ref{prop:sharp_red} apply on the rectangle $|E-E_*|\leq\eps_*$, $0<\eta\leq\eps_*$ and $[E_*-\eps_*,E_*+\eps_*]\subset[\tau,\tau^{-1}]$ for some $\tau>0$. By \eqref{eq:Theta_sqrt}, Lemma~\ref{lem:nf_coeff}, and $|z-E_*|\sim|E-E_*|+\eta$,
        \begin{align}
                |\Theta(z)|\sim(|E-E_*|+\eta)^{1/2}. \nonumber
        \end{align}
        Combining $\|\calB(z)^{-1}\|\lesssim|\Theta(z)|^{-1}$ with $|\Theta(z)|\sim(|E-E_*|+\eta)^{1/2}$ gives \eqref{eq:sharp_stab}. 
        \eqref{eq:edge_harm_scale} is a direct consequence of  Proposition~\ref{prop:sqrt_edge} and the Stieltjes representation of $m$. 
        Hence
        \begin{align}
                \Im m(z)+\frac{\eta}{\Im m(z)}\sim(|E-E_*|+\eta)^{1/2}, \nonumber
        \end{align}
        and \eqref{eq:sharp_stab_im_m} follows from \eqref{eq:sharp_stab}.

        For $\eps_*\leq\eta\leq\eta_0$, the Stieltjes representation gives $\Im m(z)\gtrsim1$, and \eqref{eq:ord_edge_stab} follows from Lemma~\ref{lem:Pi_bd}, Proposition~\ref{prop:red_stab}, and Lemma~\ref{lem:stab_red}.
\end{proof}

It remains to prove Lemma~\ref{lem:red_stab_exp}.

\begin{proof}[Proof of Lemma~\ref{lem:red_stab_exp}]
        By Lemma~\ref{lem:nf_coeff}, $-\mu_2\sim1$, while \eqref{eq:mu_bd} gives $\mu_1\sim1$. Hence \eqref{eq:Theta_sqrt} implies
        \begin{align}
                \left||\Theta(z)|^2-\frac{\mu_1}{|\mu_2|}|z-E_*|\right| \lesssim |z-E_*|^{3/2}. \nonumber
        \end{align}
        We therefore obtain
        \begin{align}
                |\Theta(z)|^2\sim|z-E_*|, \qquad \|\bbs_*(z)\|\lesssim|\Theta(z)|^2, \qquad |z-E_*|<\eps_*,\label{eq:Theta_dist_eqv}
        \end{align}
        where the second estimate follows from \eqref{eq:stab_cmp_bd}. Since $\Delta = \Pi(z)-\Pi_*$, the inverse identity \eqref{eq:inv_diff_edge} yields
        \begin{align}
                \Delta &= \Pi_*(\Pi_*^{-1}-\Pi(z)^{-1})\Pi(z)=\Pi_*((z-E_*)\wt I+\calE[\calP[\Delta]])\Pi(z)\nonumber\\
                &= \Pi_*((z-E_*)\wt I+\calE[\Theta \f{r}_*+\bbs_*])(\Pi_*+\Delta), \label{eq:Delta_res_id}
        \end{align}
        where we used $\calP[\Delta]=\Theta \f{r}_*+\bbs_*$ in the last equality.
        Since $\Pi(z)$ is uniformly bounded in the small neighborhood of $E_*$, \eqref{eq:Theta_dist_eqv} imply
        \begin{align}
                \|\Delta\|\lesssim|\Theta(z)|. \label{eq:Delta_Theta_bd}
        \end{align}
        Hence, by \eqref{eq:Delta_res_id}, we have
        \begin{align}
                \|\Pi(z)-\Pi_* -\Theta(z)\Pi_*\calE[\f{r}_*]\Pi_*\|  \leq |\Theta(z)|^2, \qquad |z-E_*|<\eps_*, \nonumber\\
        \end{align}
        this proves \eqref{eq:Pi_Theta_exp}.

        For every $\bbu\in\C^2$, the definition of $\calJ$ gives the exact identity
        \begin{align}
                (\calJ(z)-\calJ_*)[\bbu] =\calP[\Delta\calE[\bbu]\Pi_*] +\calP[\Pi_*\calE[\bbu]\Delta] +\calP[\Delta\calE[\bbu]\Delta]. \label{eq:J_diff}
        \end{align}
        Applying $\f{\ell}_*^*$ to \eqref{eq:J_diff} with $\bbu=\f{r}_*$, and using \eqref{eq:Pi_Theta_exp}, gives
        \begin{align}
                \f{\ell}_*^*(\calI-\calJ(z))\f{r}_*  =-2\Theta(z)\f{\ell}_*^*N_2[\f{r}_*,\f{r}_*]+O(|\Theta(z)|^2)  = -2\mu_2\Theta(z)+O(|\Theta(z)|^2). \nonumber
        \end{align}
        The remaining entries follow from \eqref{eq:Jstar_spec_dec}, \eqref{eq:J_diff}, and \eqref{eq:Delta_Theta_bd}:
        \begin{align}
               & \f{\ell}_*^*(\calI-\calJ(z))\wh{\f{r}}_*       =O(|\Theta(z)|), \qquad
                \wh{\f{\ell}}_*^*(\calI-\calJ(z))\f{r}_*       =O(|\Theta(z)|), \nonumber\\
               & \wh{\f{\ell}}_*^*(\calI-\calJ(z))\wh{\f{r}}_*  =1-\la_*+O(|\Theta(z)|). \nonumber
        \end{align}
        This proves \eqref{eq:stab_mat}.
\end{proof}

We henceforth fix $\eps_*$ so that all conclusions of this section hold and $0 < \eps_* < 1$.

\section{Local laws and spectral consequences}\label{sec:ll}

Throughout this section, Assumptions~\ref{ass:basic} and~\ref{ass:rand_ent} hold, and $E_*$ denotes a regular right edge in the sense of Definition~\ref{def:reg_r_edge}. We prove Theorem~\ref{thm:opt_ll} and Corollary~\ref{cor:edge_rig}. Subsection~\ref{subsec:rand_sc_eq} states the perturbed MDE for $G$ in Theorem~\ref{thm:rand_sc_eq}, whose proof is deferred to Appendix~\ref{app:rand_sc_eq}. Subsection~\ref{sub:nonlin_recon} combines the error estimates \eqref{eq:anis_rand_error} and \eqref{eq:avg_rand_error} with the deterministic edge stability to establish the nonlinear stability and reconstruction estimates. Theorem~\ref{thm:opt_ll} is proved in Subsection~\ref{subsec:ll_reg_edge}, and Corollary~\ref{cor:edge_rig} in Subsection~\ref{subsec:spec_conseq}.

The constant $\eps_*>0$ was fixed at the end of Subsection~\ref{sec:reg_edge} so that all conclusions of that subsection hold. Fix $\tau>0$ sufficiently small so that
\begin{align}
        [E_*-\eps_*,E_*+\eps_*]\subset[\tau,\tau^{-1}]. \nonumber
\end{align}
We work on the domain $\mb D_*^{\eps_*}(\tau)$ defined in \eqref{eq:edge_spec_dom}. In particular, the sharp stability bound \eqref{eq:sharp_stab} holds for all $z\in\mb D_*^{\eps_*}(\tau)$ with $\eta\le\eps_*$.
We first collect the probabilistic and resolvent tools used below. 

\begin{lemma}
        \label{lem:stoch_dom_calc}
        The following properties of stochastic domination hold.
        \begin{enumerate}[label=(\roman*)]
                \item Let $\calU$ be an index set with $\abs{\calU}\leq N^C$. If $\xi(u)\prec\zeta(u)$ uniformly for $u\in\calU$, then $\sum_{u\in\calU}\xi(u)\prec\sum_{u\in\calU}\zeta(u)$.
                \item If $\xi_1\prec\zeta_1$ and $\xi_2\prec\zeta_2$, then $\xi_1\xi_2\prec\zeta_1\zeta_2$.
                \item Suppose $0\leq\xi\leq N^C$ and $\zeta\geq N^{-C}$ for some constant $C>0$. If $\xi\prec\zeta$, then $ \bb E\xi\prec\bb E\zeta$.
        \end{enumerate}
        The conclusions remain uniform when the assumptions are uniform in additional parameters.
\end{lemma}
\begin{proof}
        The proof follows from the definitions of stochastic domination together with a union bound argument.
\end{proof}

\begin{lemma}
        \label{lem:rand_op_norm}
        Under Assumptions~\ref{ass:basic} and~\ref{ass:rand_ent},
        \begin{align}
                \|X\|+\|\calH\|+\|Y\|\prec1. \label{eq:rand_op_norm}
        \end{align}
\end{lemma}
\begin{proof}
        The proof follows from a standard application of the F\"uredi--Koml\'os argument \cite[Section~2.1.6]{anderson2009introduction}.
\end{proof}

For $R\in\C^{L\times L}$ and deterministic vectors $\bbv,\bbw\in\C^L$, we use
\begin{align}
        R_{\bbv\bbw}:=\bbv^*R\bbw, \qquad R_{a\bbw}:=\bbe_a^*R\bbw, \qquad R_{ab}:=\bbe_a^*R\bbe_b. \label{eq:anis_notation}
\end{align}
We introduce the Ward identities, which are used to control the sums of resolvent entries. The proof of the following lemma is given in Appendix~\ref{sub:ward}.
\begin{lemma}
        \label{lem:ward_ids}
        Under Assumptions~\ref{ass:basic} and~\ref{ass:rand_ent}, for every deterministic vector $\bbx\in\C^L$,
        \begin{align}
                \sum_{\mu=M+1}^{L}\abs{G_{\bbx\mu}}^2 =\frac{(\Im G)_{\bbx\bbx}}{\eta}. \label{eq:ward_2_proj}
        \end{align}
        Let $\wh\bbx:=\wh I\bbx$ and $\wt\bbx:=\wt I\bbx$. Then
        \begin{align}
                \sum_{a=1}^{L}\abs{G_{\wh\bbx a}}^2 & \prec\frac{(\Im G)_{\wh\bbx\wh\bbx}+\eta\|\wh\bbx\|^2}{\eta}, \label{eq:ward_1_cmp}  \\
                \sum_{a=1}^{L}\abs{G_{\wt\bbx a}}^2 & \prec\frac{(\Im G)_{\wt\bbx\wt\bbx}+\eta\|\wt\bbx\|^2}{\eta}. \label{eq:ward_2_cmp}
        \end{align}
\end{lemma}

We next record the resolvent derivatives used in the cumulant expansions. For $\bbq\in\R^M$, let
\begin{align}
        \wh\bbq:=\begin{pmatrix}\bbq \\
                         0\end{pmatrix}\in\R^L, \qquad \calV_{\bbq\mu}:=\wh\bbq\bbe_\mu^\top+\bbe_\mu\wh\bbq^\top, \qquad \mu\in\llbracket M+1,L\rrbracket. \nonumber
\end{align}
For a differentiable function $F=F(\calH)$, define the directional derivative
\begin{align}
        \scrD_{\bbq\mu}F(\calH) :=\left.\frac{\dd}{\dd t}F\paren{\calH+t\calV_{\bbq\mu}}\right|_{t=0}. \label{eq:dir_der}
\end{align}
This corresponds to perturbing the column of $Y$ indexed by $\mu$ in the direction $\bbq$. Since $G^{-1}=\calH+\La-I_z$, we have
\begin{align}
        \scrD_{\bbq\mu}G=-G\calV_{\bbq\mu}G, \label{eq:dir_res_der}
\end{align}
and for any $k\geq1$,
\begin{align}
        \scrD_{\bbq\mu}^kG=(-1)^k k!G\paren{\calV_{\bbq\mu}G}^k. \nonumber
\end{align}

For later use, set
\begin{align}
        \wt\Sg:=\begin{pmatrix}\Sg & 0 \\
               0   & 0\end{pmatrix}, \qquad \wh\Sg:=\begin{pmatrix}\Sg & 0   \\
               0   & I_N\end{pmatrix}. \nonumber
\end{align}
For $i\in\llbracket1,M\rrbracket$ and $\mu\in\llbracket M+1,L\rrbracket$, define
\begin{align}
        V^{i\mu} :=\calV_{\Sg^{1/2}\bbe_i,\mu} =\wh\Sg^{1/2} \paren{\bbe_i\bbe_\mu^\top+\bbe_\mu\bbe_i^\top} \wh\Sg^{1/2}. \nonumber
\end{align}
For a differentiable function of the entries of $X$, let $\partial_{i\mu}$ denote differentiation with respect to $X_{i\mu}$. Then
\begin{align}
        \partial_{i\mu}=\scrD_{\Sg^{1/2}\bbe_i,\mu}, \qquad \partial_{i\mu}\calH=V^{i\mu}, \nonumber
\end{align}
and hence
\begin{align}
        \partial_{i\mu}G=-GV^{i\mu}G, \qquad \partial_{i\mu}^kG=(-1)^k k!G\paren{V^{i\mu}G}^k. \label{eq:entry_res_der}
\end{align}

For a real random variable $h$ with moments of all orders, its $k$-th cumulant is defined by
\begin{align}
        \kappa_k(h):=(-\ii)^k \left.\frac{\dd^k}{\dd t^k}\log\bb E e^{\ii th}\right|_{t=0}. \nonumber
\end{align}

\begin{lemma}
        \label{lem:cum_bd}
        Under Assumption~\ref{ass:rand_ent},
        \begin{align}
                \kappa_1(X_{i\mu})=0, \qquad \kappa_2(X_{i\mu})=\frac1N, \qquad \kappa_k(X_{i\mu})=O(N^{-k/2}),\quad k\geq3. \nonumber
        \end{align}
\end{lemma}
\begin{proof}
        This follows from the definition of cumulants and the moment bounds in Assumption~\ref{ass:rand_ent}.
\end{proof}

We also introduce the cumulant expansion formula and a self-improving bound used in the proof of the local law. The proofs are given in \cite[Lemma~2.4]{he2018isotropic} and \cite[Lemma~2.6]{he2018isotropic}, respectively.
\begin{lemma}
        \label{lem:cum_exp}
        Let $h$ be a real random variable with moments of all orders, and let $F:\R\to\C$ be a smooth function. For every fixed integer $\ell\geq1$,
        \begin{align}
                \bb E\left(hF(h)\right) =\sum_{k=0}^{\ell}\frac{\kappa_{k+1}(h)}{k!}\bb E F^{(k)}(h)+\calR_{\ell+1}. \nonumber
        \end{align}
        For any $K>0$, the remainder term satisfies
        \begin{align}
                \abs{\calR_{\ell+1}} & \leq O(1)\left(\bb E\sup_{|t|\leq|h|}|F^{(\ell+1)}(t)|^2\bb E\left(|h|^{2\ell+4}1_{\{|h|>K\}}\right)\right)^{1/2} +O(1)\bb E|h|^{\ell+2}\sup_{|t|\leq K}|F^{(\ell+1)}(t)|. \nonumber
        \end{align}
\end{lemma}

\begin{lemma}
        \label{lem:self_imp_bd}
        Let $C>0$ be a constant,  $\xi\geq0$, and  $\zeta\in[N^{-C},N^C]$. Suppose there exists a constant $q\in[0,1)$ such that for any $\te\in[\zeta, N^C]$, we have the implication
        \begin{align}
                \xi\prec\te \quad\Longrightarrow\quad \xi\prec\te^q\zeta^{1-q}. \nonumber
        \end{align}
        If $\xi\prec N^C$, then $\xi\prec\zeta$. 
\end{lemma}

\subsection{Perturbed MDE}\label{subsec:rand_sc_eq}

This subsection states Theorem~\ref{thm:rand_sc_eq}, which provides the error estimates for the perturbed MDE used in the proof of the local law. The anisotropic and averaged estimates are proved in Appendix~\ref{sub:anis} and \ref{sub:avg}, respectively. We recall the definition of $\mb D(\tau,\eta_0)$ from \eqref{eq:D_tau_eta}, namely, 
\begin{align}
        \mb D(\tau,\eta_0):=\{z=E+\ii\eta\in\C^+: E\in[\tau,\tau^{-1}],0<\eta\leq\eta_0\}, \nonumber
\end{align}
thus $\mb D_*^{\eps_*}(\tau)\subset \mb D(\tau, \eta_0)$.

\begin{theorem}\label{thm:rand_sc_eq}
        Under Assumptions~\ref{ass:basic} and~\ref{ass:rand_ent}, fix $\tau>0$ sufficiently small. Let $z=E+\ii\eta$ satisfy
        \begin{align}
                z\in\mb D(\tau,\eta_0), \qquad \eta\geq N^{-1+\tau}, \label{eq:rand_err_dom}
        \end{align}
        and suppose that $\|\Pi(z)\|\lesssim1$. Let $\phi\in[N^{-1},N^{\tau/10}]$ and $\vphi>0$ be deterministic, with
        \begin{align}
                0<\vphi\lesssim\|\Im\Pi(z)\|+\phi, \label{eq:ward_input_size}
        \end{align}
        and assume that $ G(z)-\Pi(z)=O_\prec(\phi)$ and $ \Im G(z)=O_\prec(\vphi)$.
        Then, we have
        \begin{align}
                \Om(G(z),z) = O_\prec \left((1+\phi)^3\sqrt{\frac{\vphi+\eta}{N\eta}}\right). \label{eq:anis_rand_error}
        \end{align}
        Moreover, for every deterministic matrix $B\in\C^{L\times L}$ with $\|B\| = O(1)$, we have
        \begin{align}
                \avg{B\Om(G(z),z)}_L = O_\prec \left((1+\phi)^6\frac{\vphi+\eta}{N\eta}\right). \label{eq:avg_rand_error}
        \end{align}
\end{theorem}

The proof is given in Appendix~\ref{sub:anis} and \ref{sub:avg}.
The following bound controls the deterministic term in \eqref{eq:ward_input_size}.

\begin{lemma}\label{lem:Im_Pi_edge_bd}
        Under Assumption~\ref{ass:basic}, let $E_*$ be a regular right edge in the sense of Definition~\ref{def:reg_r_edge}. Uniformly for $z\in\mb D_*^{\eps_*}(\tau)$,
        \begin{align}
                \|\Im\Pi(z)\|\lesssim\Im m(z). \nonumber
        \end{align}
\end{lemma}

\begin{proof}
        Taking imaginary parts in the inverse MDE \eqref{eq:inverse_mde} gives
        \begin{align}
                \Im\Pi(z)=\Pi(z)^*\begin{pmatrix}\Im m(z)\Sg & 0                  \\
               0           & (\eta+\Im g(z))I_N\end{pmatrix}\Pi(z). \label{eq:ImPi_ll_id}
        \end{align}
        For $0<\eta\leq\eps_*$, the claim follows from Lemma~\ref{lem:bdry_zero_dens} and \eqref{eq:imag_comp}. On the remaining region $\eps_*\leq\eta\leq1$, Lemma~\ref{lem:Pi_bd} and \eqref{eq:imag_comp} give $\|\Pi\|\lesssim1$ and $\eta+\Im g\lesssim\Im m$. Substitution into \eqref{eq:ImPi_ll_id} proves the claim.
\end{proof}

We now record the exact identities used throughout the proof. Set
\begin{align}
        \Delta(z):=G(z)-\Pi(z), \qquad \bbu(z):=\calP[\Delta(z)]. \nonumber
\end{align}
For the random matrix $G(z)$, set
\begin{align}
        R_G(z):=(\La-I_z-\cals[G(z)])^{-1}. \label{eq:RG_def}
\end{align}
and is well defined for every $z\in\C^+$. Indeed, the resolvent identity gives $\Im G=\eta G\wt I G^*\geq0$.
Moreover, we have
\begin{align}
        -\Im(\La-I_z-\cals[G])=\eta\wt I+\cals[\Im G]=\begin{pmatrix}\avg{\Im\wt G}_N\Sg&0\\0&\left(\eta+\avg{\Im\wh G\Sg}_N\right)I_N\end{pmatrix}. \label{eq:RG_imag}
\end{align}
By \eqref{eq:G_schur}, $\wt G=(W-zI_N)^{-1}$ is invertible and $\Im\wt G=\eta\wt G^*\wt G>0$.
Thus $\avg{\Im\wt G}_N>0$, while $\avg{\Im\wh G\Sg}_N\geq0$. Since $\Sg$ is positive definite and $\eta>0$, the matrix in \eqref{eq:RG_imag} is positive definite, proving that $\La-I_z-\cals[G]$ is invertible. Whenever $z\in \mb D(\tau,\eta_0)$, by \eqref{eq:G_schur}, we have $\|\wh G\|=\|z(\calW-zI_M)^{-1}\|\lesssim 1/\eta$ and $\|\wt G\|=\|(W-zI_N)^{-1}\|\leq1/\eta$. Using the Schur complement again, the off-diagonal blocks also satisfy $\|G_{12}\|+\|G_{21}\|\lesssim 1/\eta$. Therefore, we have the trivial bound
\begin{align}
        \|G(z)\|\lesssim 1/\eta,\qquad z\in\mb D(\tau,\eta_0). \label{eq:G_triv_bd}
\end{align}

For a deterministic function $\phi$, set
\begin{align}
        \psi_\phi(z):=(1+\phi(z))^3\sqrt{\frac{\Im m(z)+\phi(z)}{N\eta}}. \nonumber
\end{align}
If $\Delta=O_\prec(\phi)$, then Lemma~\ref{lem:Im_Pi_edge_bd} gives $\Im G=O_\prec(\Im m+\phi)$. Since $\eta\lesssim\Im m$ on $\mb D_*^{\eps_*}(\tau)$, Theorem~\ref{thm:rand_sc_eq} yields
\begin{align}
        \Om(G,z)=O_\prec(\psi_\phi(z)), \qquad \avg{B\Om(G,z)}_L = O_\prec(\psi_\phi(z)^2), \label{eq:err_psi_phi}
\end{align}
whenever the control parameters lie in the range of that theorem.

By \eqref{eq:RG_def},
\begin{align}
        \Om(G,z)=I_L+(I_z-\La+\cals[G])G=I_L-R_G^{-1}G. \nonumber
\end{align}
Consequently,
\begin{align}
        \Delta=G-\Pi=R_G-\Pi-R_G\Om(G,z). \label{eq:Delta_recon}
\end{align}
Applying $\calP$ gives the exact reduced equation
\begin{align}
        \bbu-\calP[R_G-\Pi]=-\calP[R_G\Om(G,z)]. \label{eq:proj_rand_eq}
\end{align}
Moreover, since $\cals=\calE\calP$ and $-\Pi^{-1}=I_z-\La+\cals[\Pi]$,
\begin{align}
        R_G^{-1}=\Pi^{-1}-\calE[\bbu],\qquad R_G-\Pi=R_G\calE[\bbu]\Pi=\Pi\calE[\bbu]R_G. \label{eq:RG_res_id}
\end{align}

\begin{lemma}\label{lem:rand_inv_err}
        Suppose that the assumptions of Theorem~\ref{thm:rand_sc_eq} hold. Suppose that $\|R_G\|=O(1)$ with high probability. Then
        \begin{align}
                R_G\Om(G,z)=O_\prec(\psi_\phi(z)), \qquad
                \avg{BR_G\Om(G,z)}_L=O_\prec(\psi_\phi(z)^2) \nonumber
        \end{align}
        uniformly for deterministic $B$ with $\|B\|\leq1$.
\end{lemma}

\begin{proof}
        This is the two-dimensional analogue of \cite[Lemma~4.6]{he2018isotropic}. Since $\Delta=G-\Pi=O_\prec(\phi)$ and $\phi\leq N^{\tau/10}$, we have $\|\bbu\|\prec\phi$ and hence $\|\bbu\|\leq N$ with high probability. Let $\Xi$ be a high-probability event on which $\|R_G\|\lesssim1$ and $\|\bbu\|\leq N$, and define the set $U:=\{\bbu\equiv\bbu(X):X\in\Xi\}\subset\C^2$. Thus $U$ is contained in a ball of polynomial radius. Choose an $N^{-K}$-net $\widehat U\subset U$, with $K>0$ sufficiently large. Its cardinality is polynomial in $N$.

        Since $\bbh$ is attained by $\bbu$ on $\Xi$ and \eqref{eq:RG_res_id} applies, the matrix $(\Pi^{-1}-\calE[\bbh])^{-1}$ is uniformly bounded for every $\bbh\in\widehat U$. Hence \eqref{eq:anis_rand_error} and \eqref{eq:avg_rand_error} may be multiplied by this deterministic matrix, and a union bound makes the resulting anisotropic and averaged estimates simultaneous for $\bbh\in\widehat U$. For the random $\bbu$ on $\Xi$, choose $\bbh\in\widehat U$ with $\|\bbu-\bbh\|\leq N^{-K}$. The resolvent identity gives
        \begin{align}
                R_G-(\Pi^{-1}-\calE[\bbh])^{-1}=R_G\calE[\bbu-\bbh](\Pi^{-1}-\calE[\bbh])^{-1}. \nonumber
        \end{align}
        Thus the difference is $O(N^{-K})$ on $\Xi$. The definition of $\Om(G)$ and \eqref{eq:G_triv_bd} give a crude polynomial bound on $\|\Om(G)\|$ on the spectral domain. By \eqref{eq:err_psi_phi} and choosing $K$ sufficiently large, we concludes the proof.
\end{proof}

Whenever $\Delta=O_\prec(\phi)$ and $\|R_G\|=O(1)$ with high probability, combining \eqref{eq:err_psi_phi}, \eqref{eq:proj_rand_eq} and Lemma~\ref{lem:rand_inv_err} gives
\begin{align}
        \norm{\bbu-\calP[R_G-\Pi]}\prec\psi_\phi^2. \nonumber
\end{align}
Thus we obtain a two-dimensional approximate equation, whose nonlinear stability is analyzed next.

\subsection{Nonlinear stability and reconstruction}\label{sub:nonlin_recon}

This subsection proves Proposition~\ref{prop:nonlin_stab} and Lemma~\ref{lem:ll_recon}, which are used in Subsection~\ref{subsec:ll_reg_edge} to prove Theorem~\ref{thm:opt_ll}. We first analyze the two-dimensional approximate equation using the edge stability estimates from Section~\ref{sec:reg_edge} and then reconstruct $G-\Pi$ from the control of its two trace coordinates.

The relevant deterministic scale is
\begin{align}
        \gamma(z):=\Im m(z)+\frac{\eta}{\Im m(z)}. \nonumber
\end{align}
By Theorem~\ref{thm:sharp_stab}, on the part of $\mb D_*^{\eps_*}(\tau)$ with $\eta\leq\eps_*$,
\begin{align}
        \gamma(z)\sim\abs{\Theta(z)}\sim\paren{\abs{E-E_*}+\eta}^{1/2}. \label{eq:gamma_edge_comp}
\end{align}
For $\eta\geq\eps_*$, \eqref{eq:ord_edge_stab} applies. This regime is treated separately in Subsection~\ref{subsec:ll_reg_edge}.

\begin{definition}[Control function]
        Let $\calL:=\{E+\ii\eta:\eta_-\leq\eta\leq\eta_+\}$ be a vertical line segment. A deterministic function $\xi:\calL\to(0,\infty)$ is called a control function if it is continuous and $\eta\mapsto\xi(E+\ii\eta)$ is nonincreasing on $[\eta_-,\eta_+]$.
\end{definition}

\begin{lemma}
        Let $0<\xi\leq1$ and $b,d,x\in\C$. Suppose that, for some constants $c_d,C_d,C_e>0$,
        \begin{align}
                c_d\leq|d|\leq C_d,\qquad |x(b+dx)|\leq C_e\xi. \nonumber
        \end{align}
        Then
        \begin{align}
                \min\left\{|x|,\left|x+\frac bd\right|\right\}\lesssim\sqrt\xi. \label{eq:approx_qroot}
        \end{align}
        If $b\neq0$, then
        \begin{align}
                \min\left\{|x|,\left|x+\frac bd\right|\right\}\lesssim\frac{\xi}{|b|}. \label{eq:sep_quad_root}
        \end{align}
\end{lemma}

\begin{proof}
        Since
        \begin{align}
                |x|\left|x+\frac bd\right|=\frac{|x(b+dx)|}{|d|}\lesssim\xi, \nonumber
        \end{align}
        the first estimate follows. If $b\neq0$, the larger of the two factors on the left is bounded below by a constant multiple of $|b|$, which yields the second estimate.
\end{proof}

The following proposition combines the pointwise analysis of the two possible ranges of the unstable coordinate with a continuity argument; compare \cite[Lemma~4.5]{alex2014isotropic} and \cite[Lemma~3.9]{alt2020dyson}.

\begin{proposition}\label{prop:nonlin_stab}
        Under Assumption~\ref{ass:basic}, let $E_*$ be a regular right edge in the sense of Definition~\ref{def:reg_r_edge}. Fix $E\in\R$ and let
        \begin{align}
                \calL:=\{E+\ii\eta:\eta_-\leq\eta\leq\eta_+\}\subset\mb D_*^{\eps_*}(\tau),\qquad \eta_+\leq\eps_*, \nonumber
        \end{align}
        be a vertical line segment. Let $\xi:\calL\to(0,\infty)$ be a control function satisfying $N^{-C_0}\leq\xi(z)\leq N^{-c_0}$ for some fixed constants $C_0,c_0>0$. Suppose that $\bbu:\calL\to\C^2$ is continuous, $\|R_G(z)\|\leq C_R$ uniformly on $\calL$ for some fixed constant $C_R>0$, and
        \begin{align}
                \norm{\bbu(z)-\calP[R_G(z)-\Pi(z)]} \leq\xi(z) \label{eq:red_approx_eq}
        \end{align}
        uniformly for $z\in\calL$. If for some $C_{\rr{in}}>0$,
        \begin{align}
                \norm{\bbu(E+\ii\eta_+)} \leq C_{\rr{in}}\frac{\xi(E+\ii\eta_+)}{\gamma(E+\ii\eta_+)+\sqrt{\xi(E+\ii\eta_+)}}, \label{eq:red_end_bd}
        \end{align}
        then there exists a constant $C>0$, depending only on $C_{\rr{in}}$, $C_R$, and the model constants, such that
        \begin{align}
                \norm{\bbu(z)}\leq C\frac{\xi(z)}{\gamma(z)+\sqrt{\xi(z)}} \label{eq:red_vert_stab}
        \end{align}
        uniformly for $z\in\calL$.
\end{proposition}

\begin{proof}
        By \eqref{eq:RG_res_id},
        \begin{align}
                R_G-\Pi=\Pi\calE[\bbu]\Pi+\Pi\calE[\bbu]\Pi\calE[\bbu]\Pi+\Pi\calE[\bbu]\Pi\calE[\bbu]\Pi\calE[\bbu]R_G. \nonumber
        \end{align}
        Since $\Pi$ and $R_G$ are uniformly bounded, this gives the exact second order expansion
        \begin{align}
                 \bbu-\calP[R_G-\Pi] =(\calI-\calJ(z))\bbu-\calP[\Pi\calE[\bbu]\Pi\calE[\bbu]\Pi]+O(\norm{\bbu}^3). \label{eq:red_map_2nd}
        \end{align}
        By \eqref{eq:red_end_bd} and $\xi\leq N^{-c_0}$ imply $\|\bbu(E+\ii\eta_+)\|\ll 1$. Fix a sufficiently small constant $r_0>0$, and let $\wt\calL\subset\calL$ be the maximal connected subsegment containing the upper endpoint on which $\|\bbu\|\leq r_0$. We first prove \eqref{eq:red_vert_stab} on $\wt\calL$.

        Fix $z=E+\ii\eta\in\wt\calL$ and suppress the dependence on $z$. Decompose $\bbu$ in the unstable--stable basis from Proposition~\ref{prop:unst_dir} as
        \begin{align}
                \bbu=u_{\rr c}\f{r}_*+u_{\rr s}\wh{\f{r}}_*. \label{eq:rand_red_dec}
        \end{align}
        We first control the stable coordinate. By Lemma~\ref{lem:red_stab_exp},
        \begin{align}
                \wh{\f{\ell}}_*^*(\calI-\calJ(z))\bbu =O(\gamma(z))u_{\rr c}+(1-\la_*+O(\gamma(z)))u_{\rr s}. \nonumber
        \end{align}
        Projecting \eqref{eq:red_map_2nd} onto $\wh{\f{\ell}}_*$ and using \eqref{eq:red_approx_eq}, we obtain
        \begin{align}
                \left|1-\la_*+O(\gamma(z))\right||u_{\rr s}| \lesssim \gamma(z)|u_{\rr c}|+|u_{\rr c}|^2+\norm{\bbu}|u_{\rr s}|+\xi. \nonumber
        \end{align}
        Since $|1-\la_*|\sim1$, choosing $r_0$ sufficiently small allows us to absorb the term $\norm{\bbu}|u_{\rr s}|$. Hence
        \begin{align}
                |u_{\rr s}|\lesssim\gamma(z)|u_{\rr c}|+|u_{\rr c}|^2+\xi, \label{eq:rand_stab_cmp}
        \end{align}
        where we used \eqref{eq:gamma_edge_comp} to absorb the term $O(\gamma(z))$.

        We next project onto the unstable direction. By \eqref{eq:stab_mat},
        \begin{align}
                \f{\ell}_*^*(\calI-\calJ(z))\bbu =(-2\mu_2\Theta(z)+O(\gamma(z)^2))u_{\rr c}+O(\gamma(z))u_{\rr s}. \label{eq:rand_unst_lin}
        \end{align}
        Moreover, the definition of $\mu_2$ in \eqref{eq:nf_coeff}, the bound $\|\Pi(z)-\Pi_*\|\lesssim\gamma(z)$, and \eqref{eq:rand_red_dec} imply
        \begin{align}
                \f{\ell}_*^*\calP[\Pi\calE[\bbu]\Pi\calE[\bbu]\Pi] =\mu_2u_{\rr c}^2+O(\gamma(z)|u_{\rr c}|^2+\norm{\bbu}|u_{\rr s}|). \label{eq:rand_unst_quad}
        \end{align}
        Substituting \eqref{eq:rand_stab_cmp} into \eqref{eq:rand_unst_lin} and \eqref{eq:rand_unst_quad}, and using $-\mu_2\sim1$, yields
        \begin{align}
                |u_{\rr c}|\left|-2\mu_2\Theta(z)-\mu_2u_{\rr c}+O(\gamma(z)^2+\gamma(z)|u_{\rr c}|+|u_{\rr c}|^2+\xi)\right|\lesssim\xi. \nonumber
        \end{align}
        Thus there are coefficients $b$ and $d$ satisfying
        \begin{align}
                |u_{\rr c}(b+du_{\rr c})|\lesssim\xi, \qquad b=-2\mu_2\Theta(z)+O(\gamma(z)^2+\xi), \qquad d=-\mu_2+O(\gamma(z)+|u_{\rr c}|). \label{eq:rand_quad_coeff}
        \end{align}
        By the choices of  $r_0$ and \eqref{eq:gamma_edge_comp}, we have $|d|\sim1$.

         We distinguish the regimes according to whether the two approximate roots are separated on the error scale. Fix $K>0$ sufficiently large. 

        If $\paren{|E-E_*|+\eta}^{1/2}\leq K\sqrt\xi$, then \eqref{eq:gamma_edge_comp} and \eqref{eq:rand_quad_coeff} imply $|b|\lesssim\sqrt\xi$. The estimate \eqref{eq:approx_qroot} shows that $u_{\rr c}$ is within $O(\sqrt\xi)$ of either $0$ or $-b/d$. Since $|b/d|\lesssim\sqrt\xi$, both possibilities yield $|u_{\rr c}|\lesssim\sqrt\xi$. Together with \eqref{eq:rand_stab_cmp}, this gives
        \begin{align}
                \norm{\bbu}\lesssim\sqrt\xi\lesssim\frac{\xi}{\gamma(z)+\sqrt\xi}. \label{eq:rand_coal_root}
        \end{align}

        Suppose now that $\paren{|E-E_*|+\eta}^{1/2}>K\sqrt\xi$. Then \eqref{eq:gamma_edge_comp} implies $\xi\lesssim\gamma(z)^2$, and \eqref{eq:rand_quad_coeff} gives $|b|\sim\gamma(z)$. Hence the two alternatives for the unstable coordinate are separated, and \eqref{eq:sep_quad_root} yields either $|u_{\rr c}|\lesssim\xi/\gamma(z)$ or $|u_{\rr c}+b/d|\lesssim \xi/\gamma(z)$. In the first case, \eqref{eq:rand_stab_cmp} gives
        \begin{align}
                \norm{\bbu}\lesssim\frac{\xi}{\gamma(z)}\sim\frac{\xi}{\gamma(z)+\sqrt\xi}. \label{eq:rand_sep_zero}
        \end{align}
        In the second case, $|b/d|\sim\gamma(z)$, and hence $|u_{\rr c}|\gtrsim\gamma(z)$ after enlarging $K$ if necessary. Consequently, by $|u_{\rr c}|\lesssim\xi/\gamma(z)$ or $|u_{\rr c}+b/d|\lesssim \xi/\gamma(z)$, there are constants $c_1,C_1>0$ such that
        \begin{align}
                |u_{\rr c}|\leq C_1\frac{\xi}{\gamma(z)} \qquad\txt{or}\qquad |u_{\rr c}|\geq c_1\gamma(z). \label{eq:rand_unst_coord}
        \end{align}
        Enlarging $K$ once more if necessary, the two ranges in \eqref{eq:rand_unst_coord} are separated by a nonempty gap.

        Since $\xi$ is a control function, the map
        \begin{align}
                \eta\longmapsto\frac{\paren{|E-E_*|+\eta}^{1/2}}{\sqrt{\xi(E+\ii\eta)}} \nonumber
        \end{align}
        is nondecreasing. Hence the separated regime $\{\paren{|E-E_*|+\eta}^{1/2}>K\sqrt\xi\}$ is either empty or an interval containing $\eta_+$. On its complement, \eqref{eq:rand_coal_root} already proves \eqref{eq:red_vert_stab}.  After enlarging $K$ depending only on $C_{\rr{in}}$, \eqref{eq:red_end_bd} and \eqref{eq:gamma_edge_comp} select the first range in \eqref{eq:rand_unst_coord}  at $\eta_+$. Hence \eqref{eq:rand_sep_zero} holds throughout the separated regime by continuity of $\eta$. This proves \eqref{eq:red_vert_stab} on $\wt\calL$.

        Since the right-hand side of \eqref{eq:red_vert_stab} is bounded by $C\sqrt\xi\ll 1$ uniformly on $\calL$, it is smaller than $r_0/2$ for sufficiently large $N$. By continuity, $\wt\calL=\calL$, which proves the proposition.
\end{proof}

\begin{remark}\label{rem:pt_branch_stab}
        In the separated regime, if the a priori bound $\|\bbu\|\ll\gamma$ excludes the second alternative in \eqref{eq:rand_unst_coord}, then \eqref{eq:rand_sep_zero} gives
\begin{align}
        \|\bbu\|\lesssim\frac{\xi}{\gamma+\sqrt{\xi}}.
\end{align}
This estimate will be used in the final refinement on the gap side.
\end{remark}

\begin{lemma}\label{lem:ll_recon}
        Under Assumptions~\ref{ass:basic} and~\ref{ass:rand_ent}, let $E_*$ be a regular right edge in the sense of Definition~\ref{def:reg_r_edge}. Fix $z\in\mb D_*^{\eps_*}(\tau)$. Suppose that, for some fixed $\delta>0$ satisfying $10\delta<\tau$, some constants $C,c>0$, and a deterministic $\te\in[N^{-C},N^{-c}]$,
        \begin{align}
                \Delta(z)=O_\prec(N^\delta), \qquad \norm{\calP[\Delta(z)]} = O_\prec(\te).\nonumber
        \end{align}
        Then
        \begin{align}
                \Delta(z) =O_\prec\left( \te+\sqrt{\frac{\Im m(z)}{N\eta}}+\frac1{N\eta} \right). \label{eq:recon_concl}
        \end{align}
\end{lemma}

\begin{proof}
        Write $\bbu=\calP[\Delta]$. Since $\norm{\bbu}\prec\te\ll 1$ and $\|\Pi\|\lesssim1$, \eqref{eq:RG_res_id} gives
        \begin{align}
                R_G=(I_L-\Pi\calE[\bbu])^{-1}\Pi,\qquad R_G-\Pi=O_\prec(\te), \nonumber
        \end{align}
        and hence $\|R_G\|=O(1)$ with high probability. Suppose that $\Delta=O_\prec(\phi)$ for some deterministic control parameter $\phi$. Then \eqref{eq:err_psi_phi}, Lemma~\ref{lem:rand_inv_err}, and \eqref{eq:Delta_recon} yield the implication
        \begin{align}
        \Delta=O_\prec(\phi) \quad\Longrightarrow\quad \Delta=O_\prec\left(\te+(1+\phi)^3\sqrt{\frac{\Im m+\phi}{N\eta}}\right).
        \label{eq:recon_impl}
        \end{align}

        We start from $\phi=N^\delta$. Since $N\eta\geq N^\tau$ on $\mb D_*^{\eps_*}(\tau)$ and $\Im m\lesssim1$,
        \begin{align}
        (1+N^\delta)^3\sqrt{\frac{\Im m+N^\delta}{N\eta}}
        \lesssim N^{(7\delta-\tau)/2}. \nonumber
        \end{align}
        Since $10\delta<\tau$ and $\te\leq N^{-c}$, we show that the first application of \eqref{eq:recon_impl} gives a polynomially small bound on $\Delta$.

        Once $\phi \ll 1$, the factor $(1+\phi)^3$ is bounded. Moreover,
        \begin{align}
        \sqrt{\frac{\Im m+\phi}{N\eta}}
        \leq\sqrt{\frac{\Im m}{N\eta}}+\sqrt{\frac{\phi}{N\eta}}. \nonumber
        \end{align}
        Thus \eqref{eq:recon_impl} becomes the self-improving implication
        \begin{align} \Delta=O_\prec(\phi) \quad\Longrightarrow\quad \Delta=O_\prec\left( \te+\sqrt{\frac{\Im m}{N\eta}}+\sqrt{\frac{\phi}{N\eta}} \right). \nonumber
        \end{align}
        Iterating this implication gives successively better bounds on $\Delta$. The terms generated from $\te$ are absorbed using 
        \begin{align}
        \sqrt{\frac{\te}{N\eta}}\leq\te+\frac1{N\eta}. \nonumber
        \end{align}
        After a bounded number of iterations, we obtain
        \begin{align}
        \Delta =O_\prec\left( \te+\sqrt{\frac{\Im m}{N\eta}}+\frac1{N\eta} \right). \nonumber
        \end{align}
        This proves \eqref{eq:recon_concl}.
\end{proof}

\subsection{Local laws near a regular edge}\label{subsec:ll_reg_edge}

This subsection completes the proof of Theorem~\ref{thm:opt_ll}. The main inputs are Theorem~\ref{thm:rand_sc_eq}, Proposition~\ref{prop:nonlin_stab}, and Lemma~\ref{lem:ll_recon}. The continuity argument starts from Lemma~\ref{lem:ll_init_est}, whose proof is given in Appendix~\ref{sub:init_est}.

We first introduce the following monotonicity of the resolvent $G(z)$.
\begin{lemma}\label{lem:res_mono}
        Fix $E\in\R$ and let $0<\eta_1\leq\eta_2$. If $G(E+\ii\eta_2)=O_\prec(\psi)$, then $G(E+\ii\eta_1)=O_\prec\left(\eta_2/\eta_1\psi\right)$.
\end{lemma}

\begin{proof}
        For  any deterministic unit vector $\bbv\in\C^L$, we have 
        \begin{align}
                |G_{\bbv\bbv}(z+\ii h)| - |G_{\bbv\bbv}(z)|&\leq |G_{\bbv\bbv}(z+\ii h) - G_{\bbv\bbv}(z)|\leq |h| |\avg{\bbv,G(z)\wt I G(z+\ii h)\bbv}|\nonumber\\
                &\leq |h|\sqrt{\frac{\Im G(z)_{\bbv\bbv} \Im G(z+\ii h)_{\bbv\bbv}}{\eta(\eta+h)}}. \nonumber
        \end{align}
        where we used  Cauchy-Schwarz inequality and the Ward identity in \eqref{eq:ward_2_proj}.
        Letting $h\to0$, we obtain $\partial_\eta\left(\eta\abs{G_{\bbv\bbv}}\right)\geq \abs{G_{\bbv\bbv}}-\Im G_{\bbv\bbv}\geq0$.
        Hence
        \begin{align}
                \abs{G_{\bbv\bbv}(E+\ii\eta_1)}\leq\frac{\eta_2}{\eta_1}\abs{G_{\bbv\bbv}(E+\ii\eta_2)}. \nonumber
        \end{align}
        The claim follows from polarization.
\end{proof}

\begin{lemma}[Initial estimate]\label{lem:ll_init_est}
        Under the assumptions of Theorem~\ref{thm:opt_ll}, let $E_*$ be a regular right edge in the sense of Definition~\ref{def:reg_r_edge}. Uniformly for $|E-E_*|\leq\eps_*$ and $\eps_*\leq\eta\leq1$,
        \begin{align}
                \norm{\calP[G(E+\ii\eta)-\Pi(E+\ii\eta)]} = O_\prec\left(\frac1{N\eta}\right). \label{eq:init_proj_law}
        \end{align}
        On the same domain,
        \begin{align}
                G(z)-\Pi(z)=O_\prec\left(\sqrt{\frac{\Im m(z)}{N\eta}}+\frac1{N\eta}\right),\qquad \avg{B(G(z)-\Pi(z))}_L = O_\prec\left(\frac1{N\eta}\right) \label{eq:init_ll}
        \end{align}
        uniformly for deterministic $B$ with $\|B\| = O(1)$.
\end{lemma}
We now complete the proof of Theorem~\ref{thm:opt_ll}, assuming Lemma~\ref{lem:ll_init_est}.

\begin{proof}[Proof of Theorem~\ref{thm:opt_ll}]
        Lemma~\ref{lem:ll_init_est} proves all assertions on the uniformly stable part $\eps_*\leq\eta\leq1$. It remains to treat $N^{-1+\tau}\leq\eta\leq\eps_*$.

        We first work on a deterministic lattice. By \eqref{eq:G_triv_bd}, we have $\|G(z)\|\lesssim\eta^{-1}$. Moreover, differentiating the resolvent gives $\partial_zG=G\wt I G$, and hence
        \begin{align}
                \|G(z)\|\lesssim\eta^{-1},\qquad \|\partial_zG(z)\|\lesssim\eta^{-2}. \label{eq:rand_res_der}
        \end{align}
         Lemma~\ref{lem:der_bd} and Theorem~\ref{thm:sharp_stab} give a polynomial bound for $\|\partial_z\Pi(z)\|$ on the edge part of $\mb D_*^{\eps_*}(\tau)$. Hence it is enough to prove the estimates on an $N^{-K}$-net of this region, with $K>0$ sufficiently large.

        Fix a lattice value of $E$ and set $\delta:=\tau/20$. Then $N^\delta$ lies in the range of Theorem~\ref{thm:rand_sc_eq}, and $7\delta<\tau/2$. Let $K_0$ be the smallest positive integer such that $\eps_*N^{-\delta K_0}\leq N^{-1+\tau}$, and set
        \begin{align}
                \eta_k:=\eps_*N^{-\delta k},\qquad 0\leq k<K_0,\qquad \eta_{K_0}:=N^{-1+\tau}. \nonumber
        \end{align}
        For $1\leq k\leq K_0$, define
        \begin{align}
                \calL_k:=\{E+\ii\eta:\eta_k\leq\eta\leq\eta_{k-1}\}. \nonumber
        \end{align}
        We prove the local laws on $\calL_k$ by induction on $k$. For $k=1$, the estimates at the upper endpoint follow from Lemma~\ref{lem:ll_init_est}. For $k\geq2$, they follow from the induction hypothesis on $\calL_{k-1}$.

        \emph{Weak local laws.} Fix $1\leq k\leq K_0$. At the upper endpoint, $\Delta(E+\ii\eta_{k-1})=O_\prec(1)$ by the initialization for $k=1$ and by the induction hypothesis for $k\geq2$. Since $\eta_{k-1}/\eta_k\leq N^\delta$, by Lemma~\ref{lem:res_mono} we get that $G(z)=O_\prec(N^\delta)$ on $\calL_k$. Together with Lemma~\ref{lem:bdry_zero_dens}, this gives
        \begin{align}
                \Delta(z)=O_\prec(N^\delta), \qquad z\in\calL_k. \nonumber
        \end{align}
        At the upper endpoint, $\norm{\bbu}\prec(N\eta)^{-1}\ll 1$, so \eqref{eq:RG_res_id} and the boundedness of $\Pi$ imply $\|R_G\|=O(1)$ with high probability. Fix a sufficiently large constant $C_R>0$, and let $\calL\subset\calL_k$ be the maximal connected subsegment containing the upper endpoint on which $\|R_G\|\leq C_R$.  Theorem~\ref{thm:rand_sc_eq} with $\phi=N^\delta$ and $\vphi=\Im m+N^\delta$, Lemma~\ref{lem:rand_inv_err}, and \eqref{eq:proj_rand_eq} give
        \begin{align}
                \norm{\bbu-\calP[R_G-\Pi]}\prec(1+N^\delta)^6\frac{\Im m+N^\delta}{N\eta}\lesssim\frac{N^{7\delta}}{N\eta}, \qquad z\in\calL.\nonumber
        \end{align}
        Since $7\delta<\tau/2$, by a union bound over a sufficiently fine deterministic net of $\calL_k$, with high probability we have
        \begin{align}
                \norm{\bbu(z)-\calP[R_G(z)-\Pi(z)]}\leq\frac1{\sqrt{N\eta}},\qquad z\in\calL. \label{eq:weak_red_uniform}
        \end{align}
        The function $(N\eta)^{-1/2}$ is a control function and satisfies the polynomial bounds required in Proposition~\ref{prop:nonlin_stab}. The upper endpoint estimate is stronger than \eqref{eq:red_end_bd}, so Proposition~\ref{prop:nonlin_stab} gives $\norm{\bbu(z)}\prec(N\eta)^{-1/4}$ for all $ z\in\calL$.
        By \eqref{eq:RG_res_id}, chossing $C_R$ sufficiently large gives $\|R_G-\Pi\||\ll 1$ and hence $R_G \leq C_R/2 $ on $\calL$ with high probability. By continuity in $\eta$, $R_G$ cannot leave the region $\|R_G\|\leq C_R$, and therefore $\calL=\calL_k$. We conclude that
        \begin{align}
                \norm{\calP[\Delta(z)]}\prec(N\eta)^{-1/4},\qquad z\in\calL_k. \label{eq:weak_proj_bd}
        \end{align}

        \emph{Strong local laws.} Let $\te$ be a deterministic control function on $\calL_k$ satisfying $(N\eta)^{-1}\leq\te(z)\leq(N\eta)^{-1/4}$, and suppose that $\norm{\calP[\Delta(z)]}\prec\te(z)$. Since $\te\ll 1$ uniformly on $\calL_k$, \eqref{eq:RG_res_id} implies $\|R_G\|=O(1)$ with high probability. Lemma~\ref{lem:ll_recon} gives
        \begin{align}
                \Delta=O_\prec\left(\te+\sqrt{\frac{\Im m}{N\eta}}+\frac1{N\eta}\right). \nonumber
        \end{align}
        Using this bound in Theorem~\ref{thm:rand_sc_eq}, followed by Lemma~\ref{lem:rand_inv_err} and \eqref{eq:proj_rand_eq}, and using Young's inequality, $\Im m\leq\gamma$, and \eqref{eq:gamma_edge_comp}, gives
        \begin{align}
                \norm{\bbu(z)-\calP[R_G(z)-\Pi(z)]}\prec\xi(z),\qquad z\in\calL_k, \nonumber
        \end{align}
        where
        \begin{align}
                \xi(z):=\frac{\te(z)+\paren{|E-E_*|+\eta}^{1/2}}{N\eta}+\frac1{(N\eta)^2}. \nonumber
        \end{align}
        The function $\xi$ is a control function and is polynomially small uniformly on $\calL_k$. By the induction assumption and $\te\geq(N\eta)^{-1}$, we have
        \begin{align}
                \|\calP[\Delta]\|\prec\frac1{N\eta}\lesssim\frac{\xi}{\gamma+\sqrt\xi},\nonumber
        \end{align}
        at the upper endpoint  and hence satisfies the condition in Proposition~\ref{prop:nonlin_stab}. Repeating the deterministic net argument used above and applying Proposition~\ref{prop:nonlin_stab} yield the self-improving implication
        \begin{align}
                \norm{\calP[\Delta]}\prec\te\quad\Longrightarrow\quad\norm{\calP[\Delta]}\prec\frac{\xi}{\gamma+\sqrt\xi}\lesssim\frac1{N\eta}+\sqrt{\frac{\te}{N\eta}}\lesssim\te^{1/2}(N\eta)^{-1/2}. \nonumber
        \end{align}
        Starting from \eqref{eq:weak_proj_bd} and iterating this implication, by Lemma~\ref{lem:self_imp_bd}, after a bounded number of iterations, we obtain
        \begin{align}
                \norm{\calP[\Delta(z)]}\prec\frac1{N\eta},\qquad z\in\calL_k. \label{eq:proj_opt_law}
        \end{align}
        Lemma~\ref{lem:ll_recon} now gives
        \begin{align}
                \Delta(z)=O_\prec\left(\sqrt{\frac{\Im m(z)}{N\eta}}+\frac1{N\eta}\right),\qquad z\in\calL_k. \label{eq:opt_anis_lat}
        \end{align}
        Finally, taking a normalized trace in \eqref{eq:Delta_recon}, and combine \eqref{eq:RG_res_id}, \eqref{eq:proj_opt_law}, \eqref{eq:err_psi_phi}, and Lemma~\ref{lem:rand_inv_err} give
        \begin{align}
                \abs{\avg{B\Delta(z)}_L}\prec\frac1{N\eta} \nonumber
        \end{align}
        uniformly for deterministic $B$ with $\|B\| = O(1)$. The estimates at the lower endpoint of $\calL_k$ provide the upper endpoint input for $\calL_{k+1}$. Induction on $k$ proves \eqref{eq:opt_anis_ll} and \eqref{eq:opt_avg_ll}.

        \emph{Averaged estimate on the gap side.} It remains to refine the averaged estimate on the gap side. Suppose first that $E\geq E_*$ and $N^{-1+\tau}<\eta\leq\eps_*$, and set $\kappa:=E-E_*$. By the square-root behavior in Proposition~\ref{prop:sqrt_edge} and the Stieltjes transform estimate recorded in \eqref{eq:edge_harm_scale},
        \begin{align}
                \Im m(z)\sim\frac{\eta}{\gamma(z)}, \qquad \gamma(z)\sim\sqrt{\kappa+\eta}. \label{eq:gap_edge_ll}
        \end{align}
        Substituting \eqref{eq:opt_anis_lat} into \eqref{eq:err_psi_phi} , then Lemma~\ref{lem:rand_inv_err} and \eqref{eq:proj_rand_eq} give
        \begin{align}
                \|\bbu-\calP[R_G-\Pi]\|\prec \psi_\phi(z)^2 \lesssim \frac1{N\gamma(z)} +\frac1{N\eta\sqrt{N\gamma(z)}} +\frac1{(N\eta)^2}. \label{eq:gap_sq_err}
        \end{align}
        Now we take $\xi=\psi_\phi^2$ in Proposition~\ref{prop:nonlin_stab}.
        If $\gamma\lesssim\psi_\phi$, \eqref{eq:rand_coal_root} applies. If $\gamma\gg\psi_\phi$, then \eqref{eq:proj_opt_law} and $\psi_\phi\gtrsim(N\eta)^{-1}$ select the root near zero, and Remark~\ref{rem:pt_branch_stab} applies. Hence, in both cases,
        \begin{align}
                \norm{\calP[\Delta]} \prec \frac{\psi_\phi^2}{\gamma+\psi_\phi}. \nonumber
        \end{align}
        It follows from \eqref{eq:gap_sq_err} and Young's inequality that
        \begin{align}
                \norm{\calP[\Delta(z)]} \prec \frac1{N\gamma(z)^2} +\frac1{N^{3/2}\eta\gamma(z)^{3/2}} +\frac1{(N\eta)^2\gamma(z)}  \prec \frac1{N\gamma(z)^2} +\frac1{(N\eta)^2\gamma(z)}. \nonumber
        \end{align}
        Taking a normalized trace in \eqref{eq:Delta_recon} and applying \eqref{eq:err_psi_phi} and Lemma~\ref{lem:rand_inv_err} gives $\abs{\avg{B\Delta}_L}\prec N^{-1}\gamma(z)^{-2}+(N\eta)^{-2}\gamma(z)^{-1}$ and hence obtain \eqref{eq:gap_avg_ll}. For $\eps_*\leq\eta\leq1$, Lemma~\ref{lem:ll_init_est} gives \eqref{eq:gap_avg_ll}, since $\kappa+\eta\sim\eta$ on this region.

        A union bound over the edge lattice, followed by \eqref{eq:rand_res_der} and Lemma~\ref{lem:der_bd}, extends the estimates to $N^{-1+\tau}\leq\eta\leq\eps_*$. Together with Lemma~\ref{lem:ll_init_est}, this proves the local laws on the full domain $\mb D_*^{\eps_*}(\tau)$.

        \emph{Rightmost edge.} Finally, suppose that $E_*=\sup\supp\bnu$. Since the rectangle $\{E+\ii\eta:E_*+\eps_*\leq E\leq\tau^{-1},\ 0<\eta\leq1\}$ is separated from $\supp\bnu$, the Stieltjes representation gives $\Im m(z)\sim\eta$ and $\|\Im\Pi(z)\|\lesssim\eta$ uniformly on this region. Proposition~\ref{prop:red_stab}, Lemma~\ref{lem:stab_red}, and compactness yield $\|(\calI-\calJ(z))^{-1}\|+\|\calB(z)^{-1}\|\lesssim1$ uniformly throughout this rectangle. Hence the MDE is uniformly stable on this region. For $\eta\geq N^{-1+\tau}$, Theorem~\ref{thm:rand_sc_eq} and Lemma~\ref{lem:ll_recon}, together with the same argument in Lemma~\ref{lem:ll_init_est}, prove \eqref{eq:opt_anis_ll} and \eqref{eq:opt_avg_ll}. The estimate on the gap side follows from Lemma~\ref{lem:ll_init_est} when $\eta\geq\eps_*$, and from the uniform invertibility of $\calI-\calJ(z)$ and $\calB(z)$ together with $\Im m\sim\eta$ when $\eta\leq\eps_*$. This proves the final assertion and completes the proof.
\end{proof}

\subsection{Eigenvalue rigidity and delocalization}\label{subsec:spec_conseq}

We conclude this section with the proof of Corollary~\ref{cor:edge_rig}.

\begin{proof}[Proof of Corollary~\ref{cor:edge_rig}]
        \emph{No eigenvalues in the adjacent spectral gap.} Fix $\eps>0$. Since $\tau$ in Theorem~\ref{thm:opt_ll} can be chosen arbitrarily small, we choose it sufficiently small depending on $\eps$. Apply Theorem~\ref{thm:opt_ll} with this value of $\tau$, and work on the event on which its stochastic domination bounds hold uniformly with the factor $N^{\tau/2}$. Suppose that $\la$ is an eigenvalue of $W$ in the interval from \eqref{eq:no_eig_gap}, and set $\kappa:=\la-E_*,\eta:=N^{-1+\tau}/\sqrt\kappa$ and $z:=\la+\ii\eta$.
        Since $N^{-2/3+\eps}\leq\kappa\leq\eps_*<1$, the choice of $\tau$ ensures $N^{-1+\tau}\leq\eta\leq1$ for all sufficiently large $N$, and hence $z\in\mb D_*^{\eps_*}(\tau)$. By \eqref{eq:G_schur},
        \begin{align}
                \Im\avg{\wt G(z)}_N\geq\frac1{N\eta}. \label{eq:eig_res_lb}
        \end{align}
        On the other hand, \eqref{eq:edge_harm_scale} gives
        \begin{align}
                \Im m(z)\lesssim\frac{\eta}{\sqrt\kappa} \leq N^{2\tau-3\eps/2}\frac1{N\eta}. \label{eq:gap_im_bd}
        \end{align}
        Taking in \eqref{eq:gap_avg_ll} the deterministic matrix which projects onto the lower right block, and using $L\sim N$, yields on the chosen event of high probability
        \begin{align}
                \abs{\avg{\wt G(z)-\Pi_{22}(z)}_N} \leq N^{\tau/2}\left(\frac1{N(\kappa+\eta)}+\frac1{(N\eta)^2\sqrt{\kappa+\eta}}\right)  \lesssim (N^{3\tau/2-3\eps/2}+N^{-\tau/2})\frac1{N\eta}. \label{eq:gap_eig_contr}
        \end{align}
        Since $m(z)=\avg{\Pi_{22}(z)}_N$ and $\tau$ was chosen sufficiently small, the right-hand sides of \eqref{eq:gap_im_bd} and \eqref{eq:gap_eig_contr} are $N^{-c}(N\eta)^{-1}$ for some small $c>0$. This contradicts \eqref{eq:eig_res_lb} for sufficiently large $N$ and proves the first assertion of part~(i). If $E_*=E_+$, the final assertion of Theorem~\ref{thm:opt_ll} allows the contradiction argument based on \eqref{eq:eig_res_lb}--\eqref{eq:gap_eig_contr} on the whole interval from $E_++N^{-2/3+\eps}$ to $\tau^{-1}$. Together with the high probability bound $\|Y\|\leq C$, after decreasing $\tau$ if necessary, this yields \eqref{eq:lambda1_ub}.

        \emph{Rigidity and delocalization.} 
        The estimates \eqref{eq:r_edge_rig} and \eqref{eq:lambda1_st_rig} follow from the averaged local law \eqref{eq:opt_avg_ll} and its refinement on the gap side \eqref{eq:gap_avg_ll} by a standard argument using Helffer-Sj\"ostrand calculus and hence \eqref{eq:r_edge_deloc} also holds. The details are already given in \cite{erdos2013spectral,erdhos2012rigidity,pillai2014universality}.
\end{proof}

\section{Gaussian divisible model and edge universality}\label{sec:edge_univ}

Throughout this section, Assumptions~\ref{ass:basic} and~\ref{ass:rand_ent} hold, and the rightmost edge $E_+$ is regular in the sense of Definition~\ref{def:reg_r_edge}. We prove Theorem~\ref{thm:edge_univ}. In Subsection~\ref{subsec:gauss_div}, we construct the Gaussian divisible model and identify its edge parameters. We then compare its edge statistics with those of the GOE. In Subsection~\ref{subsec:rm_gauss}, we remove the Gaussian component by a short-time comparison of the edge distributions and complete the proof of Theorem~\ref{thm:edge_univ}. The auxiliary perturbation and comparison estimates used in these two steps are proved in Appendix~\ref{app:edge_univ}.

Choose $\omega,\vartheta>0$ satisfying $0<\omega<\frac1{12}$ and $0<\vartheta<\min\left\{2\omega,\frac13-\omega\right\}$. Set $\eta_*:=N^{-2/3+\vartheta}$ and $t:=N^{-1/3+\omega}$. Then, for some $\eps>0$,
\begin{align}
        N^{\eps}\eta_*\leq t^2\leq N^{-\eps},\qquad Nt^3\geq N^{\eps}. \label{eq:eta_t}
\end{align}

\subsection{Gaussian divisible model and edge parameters}
\label{subsec:gauss_div}

This subsection establishes the GOE edge comparison for the Gaussian divisible model. The regularity of the initial edge and the corresponding edge parameters are proved in  Appendix~\ref{sub:init_model}, while the deterministic and empirical edge parameters needed below are analyzed in Appendix~\ref{sub:edge_par}. We prove Proposition~\ref{prop:rect_flow_id} below and obtain the comparison for $W_t$ at the end of the subsection.

Since $\Sg\geq c_0I_M$, $\Sg_t:=\Sg-tI_M$ remains uniformly positive definite. Let $X^{\rr G}$ be an $M\times N$ real Gaussian matrix independent of $X$. Its entries are independent and have variance $N^{-1}$. Define
\begin{align}
        Y_0 & :=A+\Sg_t^{1/2}X,\qquad W_0:=Y_0^\top Y_0, \nonumber         \\
        Y_t & :=Y_0+\sqrt tX^{\rr G},\qquad W_t:=Y_t^\top Y_t. \nonumber
\end{align}
For every column index $\mu$,
\begin{align}
        \bb E(Y_t-A)_{\cdot\mu}(Y_t-A)_{\cdot\mu}^\top =\frac1N(\Sg_t+tI_M)=\frac1N\Sg. \nonumber
\end{align}
Thus the MDE associated with $Y_t$ coincides with \eqref{eq:intro_mde}. Its deterministic solution, associated measure, rightmost edge, and edge scaling parameter are therefore $\Pi$, $\bnu$, $E_+$, and $\gamma_+$, respectively.

Let $\Pi_0$ be the solution of the MDE obtained from \eqref{eq:intro_mde} by replacing $\Sg$ with $\Sg_t$, and let $\bnu_0$ and $\rho_0$ be the corresponding measure and density. For $z\in \C^+$, set
\begin{align}
        m_0(z):=\avg{\Pi_{0,22}(z)}_N,\quad h_0(z):=\avg{ \Pi_{0,11}(z)}_N, \quad g_0(z):=\avg{\Sg_t \Pi_{0,11}(z)}_N, \quad E_{+,0}:=\sup\supp\bnu_0, \nonumber
\end{align}
where $\Pi_{0,11}$ and $\Pi_{0,22}$ are the upper-left $M\times M$ and lower-right $N\times N$ blocks of $\Pi_0$, respectively. We also use $\Pi_{0,12}$ and $\Pi_{0,21}$ to denote the upper-right $M\times N$ and lower-left $N\times M$ blocks of $\Pi_0$, respectively.

\begin{proposition}\label{prop:init_edge_reg}
        Under Assumptions~\ref{ass:basic} and~\ref{ass:rand_ent}, suppose that $E_+$ is a regular right edge  in the sense of Definition~\ref{def:reg_r_edge}. Then the rightmost edge $E_{+,0}$ of $\bnu_0$ is a regular right edge  for $\Pi_0$ in the sense of Definition~\ref{def:reg_r_edge}. Moreover,
        \begin{align}
                \abs{E_{+,0}-E_+}\lesssim t, \qquad c_{{\rm edge},0}\sim1, \qquad \gamma_{+,0}:=(\pi c_{{\rm edge},0})^{2/3}\sim1, \nonumber
        \end{align}
        and the constants in Proposition~\ref{prop:sqrt_edge} and Theorem~\ref{thm:opt_ll} may be chosen uniformly for $Y_0$.
\end{proposition}

Define
\begin{align}
        m_{0,N}(z):=\frac1N\Tr(W_0-zI_N)^{-1},\qquad \la_{1,0}:=\la_1(W_0). \nonumber
\end{align}

\begin{proposition}\label{prop:eta_star_reg}
        Under assumptions of Proposition~\ref{prop:init_edge_reg}. With high probability, $Y_0$ is $\eta_*$-regular around $\la_{1,0}$ in the sense of \cite[Definition~1]{ding2022edge} or \cite[Definition~1]{ding2022tracywidom}. More explicitly, there exist constants $c_V,C_V>0$ such that
        \begin{enumerate}[label=(\roman*)]
                \item For $z=E+\ii\eta$ with $ \la_{1,0}-c_V\leq E\leq\la_{1,0}$ and $\eta_*+\sqrt{\eta_*\abs{\la_{1,0}-E}}\leq\eta\leq10$, 
                we have
                \begin{align}
                         \frac{1}{C_V}\sqrt{\la_{1,0}-E+\eta}\leq\Im m_{0,N}(E+\ii\eta) \leq C_V\sqrt{\la_{1,0}-E+\eta}, \label{eq:eta_star_reg_in} 
                \end{align}
                \item For $z=E+\ii\eta$ with $ \la_{1,0}\leq E\leq\la_{1,0}+c_V$ and $\eta_*\leq\eta\leq10$,
                we have
                \begin{align}
                        \frac{1}{C_V}\frac{\eta}{\sqrt{E-\la_{1,0}+\eta}}\leq \Im m_{0,N}(E+\ii\eta) \leq C_V\frac{\eta}{\sqrt{E-\la_{1,0}+\eta}}. \label{eq:eta_reg_out}
                \end{align}
        \end{enumerate}
       Moreover, $\la_{1,0}\sim1$ and $\|W_0\|\leq C_V$ with high probability.
\end{proposition}

The proofs of Propositions~\ref{prop:init_edge_reg} and~\ref{prop:eta_star_reg} are given in Appendix~\ref{sub:init_model}.

We next identify the deterministic law after adding the Gaussian component. For a rectangular Gaussian perturbation $V+\sqrt tX^{\rr G}$, the deterministic law of the squared singular values is described by the rectangular free convolution of the initial law with the Mar\v{c}enko--Pastur law at time $t$ with the aspect ratio $c_N$. Its Stieltjes transform is characterized by the scalar self-consistent equation in \cite[Eq.~(10)]{ding2022edge}. The corresponding subordination function and inverse map are given in \cite[Eqs.~(49)--(54)]{ding2022edge}. In our construction, the covariance $tI_M$ carried by $\sqrt tX^{\rr G}$ exactly restores the part removed in $\Sg_t=\Sg-tI_M$. We verify below that the rectangular free convolution of $\bnu_0$ with the Mar\v{c}enko--Pastur law at time $t$ coincides with $\bnu$.

For $\zeta\in\C^+$, by the definition of $m_0$ and $h_0$, we have
\begin{align}
        h_0(\zeta)=\zeta m_0(\zeta)+1-c_N. \label{eq:init_trace_id}
\end{align}
Define
\begin{align}
        a_t(\zeta):=1-tm_0(\zeta),\qquad \Phi_t(\zeta):=a_t(\zeta)(\zeta-th_0(\zeta)). \nonumber
\end{align}
Let $m_{w,t}$ be the Stieltjes transform of the rectangular free convolution in \cite[Eq.~(9)]{ding2022edge}. Following \cite[Eq.~(33)]{ding2022edge}, define
\begin{align}
        m_{\rm fc}(z):=c_Nm_{w,t}(z)-\frac{1-c_N}{z}. \label{eq:rect_mfc}
\end{align}
We also set
\begin{align}
        b_t(z):=1+c_Ntm_{w,t}(z). \label{eq:rect_bt}
\end{align}
For $z\in\C^+$, let
\begin{align}
        \zeta\equiv\zeta(z):=zb_t(z)(1+tm_{\rm fc}(z)). \label{eq:rect_zeta}
\end{align}
By \cite[Eq.~(39)]{ding2022edge}, this is precisely the subordination function $\zeta_t(z)$ considered there. Moreover, as noted in \cite[Section~5.1, after Eq.~(50)]{ding2022edge}, the functions $m_{w,t}$, $b_t$, and $\zeta$ are holomorphic on $\C^+$. The relations in \cite[Eqs.~(52)--(54)]{ding2022edge}, rewritten using \eqref{eq:init_trace_id} and \eqref{eq:rect_mfc}, give
\begin{align}
        z=\Phi_t(\zeta),\qquad m_{\rm fc}(z)=a_t(\zeta)^{-1}m_0(\zeta). \label{eq:rect_subord}
\end{align}
In particular,
\begin{align}
        a_t(\zeta)^{-1}=1+tm_{\rm fc}(z),\qquad \zeta=zb_t(z)a_t(\zeta)^{-1}. \label{eq:rect_subord_aux}
\end{align}
The following proposition identifies the relation \eqref{eq:rect_subord} at the MDE level and gives the corresponding block identity. The map $\Phi_t$ will then be used to identify the right edge $E_+$ and its square-root coefficient.

\begin{proposition}\label{prop:rect_flow_id}
        Under Assumption~\ref{ass:basic}, let $z\in\C^+$ and let $\zeta=\zeta(z)\in\C^+$ be defined by \eqref{eq:rect_zeta}. Then
        \begin{align}
                \Pi(z)=\begin{pmatrix}a_t(\zeta)\Pi_{0,11}(\zeta) & \Pi_{0,12}(\zeta) \\ \Pi_{0,21}(\zeta) & a_t(\zeta)^{-1}\Pi_{0,22}(\zeta)\end{pmatrix}. \label{eq:rect_flow_blk}
        \end{align}
        In particular,
        \begin{align}
                m(z)=a_t(\zeta)^{-1}m_0(\zeta),\qquad g(z)=a_t(\zeta)(g_0(\zeta)+th_0(\zeta)). \label{eq:rect_flow_scal}
        \end{align}
        Consequently, $m_{\rm fc}(z)=m(z)$ for all $z\in\C^+$, and hence the rectangular free convolution of $\bnu_0$ with the Mar\v{c}enko--Pastur law at time $t$ and aspect ratio $c_N$ coincides with $\bnu$.
\end{proposition}

The proof of Proposition~\ref{prop:rect_flow_id} is given in Appendix~\ref{sub:edge_par}. The following lemma identify the critical point of $\Phi_t$ .   
\begin{lemma}\label{lem:rect_crit}
        Under assumptions of Proposition~\ref{prop:init_edge_reg}. There is a unique critical point $\zeta_t$ satisfying $\zeta_t-E_{+,0}\sim t^2$. Moreover,
        \begin{align}
                \Phi_t'(\zeta_t)  & =0,\qquad m_0'(\zeta_t)\sim t^{-1},\qquad -m_0''(\zeta_t)\sim t^{-3}, \nonumber   \\
                \Phi_t''(\zeta_t) & \sim t^{-2},\qquad \Phi_t(\zeta_t)=E_+, \qquad c_{\rm edge}       =\frac1\pi\frac{m_0'(\zeta_t)}{a_t(\zeta_t)^2}\sqrt{\frac{2}{\Phi_t''(\zeta_t)}}.
                \label{eq:crit_est}
        \end{align}
\end{lemma}

\begin{lemma}\label{lem:emp_subord}
        Under assumptions of Proposition~\ref{prop:init_edge_reg}. Uniformly for real $\zeta$ satisfying $\zeta-E_{+,0}\sim t^2$,
        \begin{align}
                \abs{m_{0,N}'(\zeta)-m_0'(\zeta)}   & \prec\frac1{Nt^4}, \label{eq:emp_m0p}                            \\
                \abs{m_{0,N}''(\zeta)-m_0''(\zeta)} & \prec\frac1{Nt^6}, \label{eq:emp_m0pp}                           \\
                \abs{m_{0,N}(\zeta)-m_0(\zeta)}     & \prec\eta_*+\frac1{Nt^2}+\frac1{(N\eta_*)^2t}. \label{eq:emp_m0}
        \end{align}
\end{lemma}

For $\zeta\in\C^+$, condition on $Y_0$ and define
\begin{align}
        h_{0,N}(\zeta):=\zeta m_{0,N}(\zeta)+1-c_N, \quad a_{t,N}(\zeta):=1-tm_{0,N}(\zeta), \quad \Phi_{t,N}(\zeta):=a_{t,N}(\zeta)(\zeta-th_{0,N}(\zeta)). \nonumber
\end{align}

\begin{proposition}\label{prop:cond_edge_par}
        Under the assumptions of Proposition~\ref{prop:init_edge_reg}, with high probability, $\Phi_{t,N}$ has a unique critical point $\zeta_{t,N}$ satisfying $\zeta_{t,N}-E_{+,0}\sim t^2$.
        Define
        \begin{align}
                E_{t,N}:=\Phi_{t,N}(\zeta_{t,N}), \qquad \gamma_{t,N}:=\left(\frac{m_{0,N}'(\zeta_{t,N})}{a_{t,N}(\zeta_{t,N})^2}\sqrt{\frac{2}{\Phi_{t,N}''(\zeta_{t,N})}}\right)^{2/3}. \label{eq:cond_edge_par}
        \end{align}
        Then
        \begin{align}
                \zeta_{t,N}-\la_{1,0}\sim t^2,\qquad &\abs{\zeta_{t,N}-\zeta_t}\prec\frac1{Nt}, \nonumber\\
                \abs{E_{t,N}-E_+}\prec t\eta_*+\frac1{Nt}+\frac1{(N\eta_*)^2},\qquad &\abs{\gamma_{t,N}-\gamma_+}\prec\frac1{Nt^3}. \label{eq:cond_edge_est}
        \end{align}
\end{proposition}

The proofs of Lemmas~\ref{lem:rect_crit} and \ref{lem:emp_subord}, and Proposition~\ref{prop:cond_edge_par} are given in  Appendix~\ref{sub:edge_par}. By \eqref{eq:eta_t},
\begin{align}
        N^{2/3}\abs{E_{t,N}-E_+} & =O_\prec(N^{-c}), \qquad \abs{\gamma_{t,N}-\gamma_+}=O_\prec(N^{-c}). \label{eq:cond_tw}
\end{align}
for some small $c>0$.

We now apply \cite[Theorem~4]{ding2022tracywidom} with $p=M$ and $n=N$. The normalized Stieltjes transforms of $Y_0^\top Y_0$ and $Y_0Y_0^\top$ satisfy
\begin{align}
        \frac1N\Tr(Y_0^\top Y_0-zI_N)^{-1}=c_N\frac1M\Tr(Y_0Y_0^\top-zI_M)^{-1}-\frac{1-c_N}{z}.\nonumber
\end{align}

At a positive edge, this identity shows that the normalization in \cite[Theorem~4]{ding2022tracywidom} is exactly $\gamma_{t,N}N^{2/3}$. Let $\mu_1^{\rr{GOE}}$ denote the largest eigenvalue of a standard $N\times N$ GOE matrix. 
On the high probability event from Proposition~\ref{prop:eta_star_reg}, \cite[Theorem~4]{ding2022tracywidom} applies conditionally on $Y_0$. Since $F$ is bounded, taking expectations gives, for  test function $F:\R\to\R$ satisfying $\|F\|_\infty+\|F'\|_\infty\leq C$,
\begin{align}
        \lim_{N\to\infty}\abs{\bb E F\left(\gamma_{t,N}N^{2/3}(\la_1(W_t)-E_{t,N})\right)-\bb E F\left(N^{2/3}(\mu_1^{\rr{GOE}}-2)\right)}=0. \nonumber
\end{align}
Since $F$ is bounded and Lipschitz, \eqref{eq:cond_tw}, the estimates in \eqref{eq:cond_edge_est}, and the tightness of the rescaled largest eigenvalue implied by the GOE comparison in \cite[Theorem~4]{ding2022tracywidom} allow us to replace $E_{t,N}$ and $\gamma_{t,N}$ by $E_+$ and $\gamma_+$, respectively. Thus
\begin{align}
        \lim_{N\to\infty}\abs{\bb E F\left(\gamma_+N^{2/3}(\la_1(W_t)-E_+)\right)-\bb E F\left(N^{2/3}(\mu_1^{\rr{GOE}}-2)\right)}=0.
        \label{eq:gd_univ}
\end{align}

\subsection{Removal of the Gaussian component}
\label{subsec:rm_gauss}

In this subsection, we remove the Gaussian component by the short-time comparison in Proposition~\ref{prop:short_comp}. The proof of Theorem~\ref{thm:edge_univ} is completed at the end of the subsection. The auxiliary comparison estimates used in Proposition~\ref{prop:short_comp} are proved in Appendix~\ref{sub:comp_est}.

For $0\leq s\leq t$, define
\begin{align}
        &Y(s)             :=A+(\Sg-sI_M)^{1/2}X+\sqrt sX^{\rr G},\qquad W(s):=Y(s)^\top Y(s), \label{eq:cov_interp} \\
        &(\Sg-sI_M)+sI_M  =\Sg. \label{eq:interp_cov_id}
\end{align}
Then $W(0)=W$ and $W(t)=W_t$.
Thus the deterministic MDE \eqref{eq:intro_mde} and the edge parameters remain fixed along the interpolation.

\begin{proposition}\label{prop:short_comp}
Under Assumptions~\ref{ass:basic} and~\ref{ass:rand_ent}, suppose that $E_+$ is a regular right edge in the sense of Definition~\ref{def:reg_r_edge}. Then there exists $c>0$ such that, for every fixed $x\in\R$,
\begin{align}
        \bb P\left(\gamma_+N^{2/3}\bigl(\la_1(W_t)-E_+\bigr)\leq x-N^{-c}\right)-N^{-c}&\leq \bb P\left(\gamma_+N^{2/3}\bigl(\la_1(W)-E_+\bigr)\leq x\right) \nonumber\\
        &\leq \bb P\left(\gamma_+N^{2/3}\bigl(\la_1(W_t)-E_+\bigr)\leq x+N^{-c}\right)+N^{-c}.
        \label{eq:short_edge_comp}
\end{align}
\end{proposition}

\begin{proof}
        Fix $x\in\R$. Let $C_{\rr{cmp}}>0$ be sufficiently large. Choose $\chi>0$ sufficiently small that $\omega+C_{\rr{cmp}}\chi<1/6$. Then choose the local law parameter $\tau$ such that $5\chi<\tau<1/3-\chi$. Set $\eta_{\rr{sm}}:=N^{-2/3-\chi}$. Lemma~\ref{lem:comp_input} in Appendix~\ref{sub:comp_est} gives the local law estimates uniformly for $0\leq s\leq t$. Together with the argument used to prove Corollary~\ref{cor:edge_rig}, this yields the corresponding edge rigidity uniformly along the interpolation. We may therefore apply the smoothing argument for the eigenvalue counting function; see \cite[Section~6]{erdhos2012rigidity} and \cite[Lemmas~5.3--5.5]{ding2018necessary}.
        
        For each pair $E_1,E_2$ arising from the upper and lower smoothing bounds at the threshold $E_++\gamma_+^{-1}N^{-2/3}x$, define
        \begin{align}
                m_s(z):=\frac1N\Tr(W(s)-zI_N)^{-1}, \qquad \calX(s):=N\int_{E_1}^{E_2}\Im m_s(y+\ii\eta_{\rr{sm}})\dd y. \nonumber
        \end{align}
        The smoothing bounds approximate each probability in \eqref{eq:short_edge_comp} by expectations of bounded smooth cutoff functions of $\calX(s)$, with the rescaled threshold $x$ replaced by $x\pm N^{-c_{\rr{sm}}}$ and an additive error of size $N^{-c_{\rr{sm}}}$ for some $c_{\rr{sm}}>0$. Hence it suffices to compare $\bb Ef(\calX(t))$ and $\bb Ef(\calX(0))$ for the cutoff functions $f\in C^5(\R)$ arising from this reduction. Their derivatives up to order five are uniformly bounded. The directional derivative estimates needed below are established in the proof of Lemma~\ref{lem:comp_high_cum}.

        Set $B_s:=(\Sg-sI_M)^{1/2}$. We use the directional derivatives from \eqref{eq:dir_der}. For $0<s\leq t$,
        \begin{align}
                \frac{\dd}{\dd s}Y(s)=B_s'X+\frac1{2\sqrt s}X^{\rr G}, \qquad B_s':=-\frac12B_s^{-1}, \qquad B_s'B_s=-\frac12I_M. \label{eq:interp_der}
        \end{align}
        Hence the chain rule gives
        \begin{align}
                \frac{\dd}{\dd s}\bb Ef(\calX(s)) & =\sum_{i,\mu}\bb E X_{i\mu}\scrD_{B_s'\bbe_i,\mu}f(\calX(s)) +\frac1{2\sqrt s}\sum_{i,\mu}\bb E X^{\rr G}_{i\mu}\scrD_{\bbe_i,\mu}f(\calX(s)).\nonumber
        \end{align}

        Applying Lemma~\ref{lem:cum_exp} with $\ell=3$ to the first term and Gaussian integration by parts to the second yields
        \begin{align}
                \frac{\dd}{\dd s}\bb Ef(\calX(s)) =\frac1N\sum_{i,\mu}\bb E\scrD_{B_s'\bbe_i,\mu}\scrD_{B_s\bbe_i,\mu}f(\calX(s)) +\frac1{2N}\sum_{i,\mu}\bb E\scrD_{\bbe_i,\mu}^2f(\calX(s)) +\calR_s, \label{eq:comp_second}
        \end{align}
        where
        \begin{align}
                \calR_s:=\sum_{i,\mu}\left(\frac{\kappa_3(X_{i\mu})}{2}\bb E\scrD_{B_s'\bbe_i,\mu}\scrD_{B_s\bbe_i,\mu}^2f(\calX(s))
                +\frac{\kappa_4(X_{i\mu})}{6}\bb E\scrD_{B_s'\bbe_i,\mu}\scrD_{B_s\bbe_i,\mu}^3f(\calX(s))
                +\calR_{i\mu,N}(s)\right). \label{eq:comp_rem_def}
        \end{align}
        Here $\calR_{i\mu}(s)$ denotes the remainder in Lemma~\ref{lem:cum_exp} with $\ell=3$.
        Using $B_s'=-B_s^{-1}/2$ from \eqref{eq:interp_der} and the bilinearity of the second directional derivative, we have
        \begin{align}
                \sum_i\scrD_{B_s'\bbe_i,\mu}\scrD_{B_s\bbe_i,\mu} =-\frac12\sum_i\scrD_{\bbe_i,\mu}^2.\nonumber
        \end{align}
        Thus the two second-order terms in \eqref{eq:comp_second} cancel exactly. Lemma~\ref{lem:comp_high_cum} gives, uniformly for $0<s\leq t$, $\abs{\calR_s}\leq N^{1/6+C_{\rr{cmp}}\chi}$.
        Set $c_{\rr{cmp}}:=1/6-\omega-C_{\rr{cmp}}\chi>0$. 
        The derivative bound is uniform on $(0,t]$. The map $s\mapsto\bb Ef(\calX(s))$ is continuous at $s=0$. Hence integration over $[0,t]$ gives
        \begin{align}
                \abs{\bb Ef(\calX(t))-\bb Ef(\calX(0))}\leq tN^{1/6+C_{\rr{cmp}}\chi}=N^{-c_{\rr{cmp}}}. \nonumber
        \end{align}
        Combining this estimate with the upper and lower smoothing bounds proves \eqref{eq:short_edge_comp}, after decreasing $c>0$ if necessary.
\end{proof}

\begin{proof}[Proof of Theorem~\ref{thm:edge_univ}]
        By \eqref{eq:gd_univ} and the GOE edge limit, $\gamma_+N^{2/3}(\la_1(W_t)-E_+)$ converges weakly to the Tracy--Widom law of type~1. Proposition~\ref{prop:short_comp}, together with the continuity of the Tracy--Widom distribution function, implies that $\gamma_+N^{2/3}(\la_1(W)-E_+)$ converges weakly to the same limit. Hence \eqref{eq:r_edge_univ} follows.
\end{proof}

\section*{Acknowledgments}
The authors would like to thank Zhigang Bao and Yukun He for their insightful discussions and valuable guidance during the early stages of this research.

\bibliographystyle{plain}
\bibliography{ref}

\begin{appendix}

\section{Basic properties of the symmetric MDE}
\label{app:basic_mde}

Throughout this appendix, Assumption~\ref{ass:basic} holds. We prove the properties of the symmetric MDE \eqref{eq:sym_mde} used in Proposition~\ref{prop:pi}.

\begin{proof}[Proof of Lemma~\ref{lem:lin_mde_meas}]
        Let $\zeta=\wh E+\ii\wh\eta\in\C^+$, where $\wh\eta>0$. The existence and uniqueness of the solution to \eqref{eq:sym_mde} in $\calM_+$ were proved in \cite{helton2007operatorvalued}; see also \cite[Lemma~1.4]{he2018isotropic}. Moreover, the map $\zeta\mapsto\wh\Pi(\zeta)$ is holomorphic on $\C^+$; see also \cite{ajanki2019stability}.

        Since $\wh\Pi(\zeta)$ solves \eqref{eq:sym_mde}, and  $\Im(\zeta I_L-\La+\cals[\wh\Pi(\zeta)])>0$, we have
        \begin{align}
                \wh\Pi(\zeta) = -(\zeta I_L-\La+\cals[\wh\Pi(\zeta)])^{-1}. \label{eq:Pi_hat_inv}
        \end{align}
        Taking imaginary parts gives
        \begin{align}
                (\zeta I_L-\La+\cals[\wh\Pi(\zeta)])^*\Im\wh\Pi(\zeta)(\zeta I_L-\La+\cals[\wh\Pi(\zeta)])= \wh\eta I_L+\cals[\Im\wh\Pi(\zeta)]. \label{eq:im_part_id}
        \end{align}
        If $\Im\wh\Pi(\zeta)$ were not positive definite, then since $\Im\wh\Pi(\zeta)\ge0$, there would exist a nonzero vector $\bbx\in\C^L$ such that $\bbx^*\Im\wh\Pi(\zeta)\bbx=0$.
        Let $\bbv:=(\zeta I_L-\La+\cals[\wh\Pi(\zeta)])^{-1}\bbx$. Then $\bbv\ne0$. Testing \eqref{eq:im_part_id} against $\bbv$, we obtain
        \begin{align}
                0 = \bbx^*\Im\wh\Pi(\zeta)\bbx = \bbv^*(\wh\eta I_L+\cals[\Im\wh\Pi(\zeta)])\bbv \ge \wh\eta\|\bbv\|^2 > 0, \nonumber
        \end{align}
        which is impossible. Therefore $\Im\wh\Pi(\zeta)>0, \zeta\in\C^+$.

        We now prove the matrix-valued Stieltjes representation. First, taking imaginary parts in the inverse form of the MDE \eqref{eq:Pi_hat_inv} gives
        \begin{align}
                \Im\wh\Pi = \wh\Pi^* (\wh\eta I_L+\cals[\Im\wh\Pi])\wh\Pi\geq \wh\eta\wh\Pi^*\wh\Pi, \nonumber
        \end{align}
        which gives $\|\wh\Pi(\zeta)\|\leq\wh\eta^{-1}$.
        Let $\bbw\in\C^L$ satisfy $\|\bbw\|=1$. Since $\zeta\mapsto\wh\Pi(\zeta)$ is holomorphic on $\C^+$, the scalar function $\bbw^*\wh\Pi(\zeta)\bbw$ is also holomorphic on $\C^+$. Moreover, $\Im\wh\Pi(\zeta)>0$ implies $\Im\bbw^*\wh\Pi(\zeta)\bbw>0$. Taking $\zeta=\ii\wh\eta$,  we have $\|\cals[\wh\Pi(\ii\wh\eta)]\| \leq \wh\eta^{-1}$. Therefore
        \begin{align}
                \ii\wh\eta I_L-\La+\cals[\wh\Pi(\ii\wh\eta)] = \ii\wh\eta I_L-\La+O(\wh\eta^{-1}) = \ii\wh\eta\left(I_L+\frac{\ii}{\wh\eta}\La+O(\wh\eta^{-2})\right). \nonumber
        \end{align}
        For sufficiently large $\wh\eta$, the inverse form of the symmetric MDE \eqref{eq:Pi_hat_inv} gives $\wh\Pi(\ii\wh\eta)=-(\ii\wh\eta)^{-1}I_L+O(\wh\eta^{-2})$. Consequently,
        \begin{align}
                \ii\wh\eta \bbw^*\wh\Pi(\ii\wh\eta)\bbw \to -1, \qquad \wh\eta\to+\infty. \nonumber
        \end{align}
        Hence, there exists a probability measure $\wh V_{\bbw}$ on $\R$ such that
        \begin{align}
                \bbw^*\wh\Pi(\zeta)\bbw = \int_{\R}\frac{\wh V_{\bbw}(\dd x)}{x-\zeta}, \qquad \zeta\in\C^+. \nonumber
        \end{align}
        By polarization, there exists a unique positive semidefinite matrix-valued measure $\wh V$ on $\R$ such that
        \begin{align}
                \wh\Pi(\zeta) = \int_{\R}\frac{\wh V(\dd x)}{x-\zeta}, \qquad \zeta\in\C^+, \nonumber
        \end{align}
        and $\wh V(\R)=I_L$. The details can be found in \cite{ajanki2019stability, gesztesy2000matrix}.

        Let $\wh\mu(\dd x):=\avg{\wh V(\dd x)}_L$ be a positive measure, and $\wh\mu(\R)=\avg{\wh V(\R)}_L=1$. Thus, $\wh\mu$ is a probability measure. Taking normalized traces gives
        \begin{align}
                \wh m(\zeta)=\frac1L\Tr\wh\Pi(\zeta)=\int_{\R}\frac{\wh\mu(\dd x)}{x-\zeta}.\nonumber
        \end{align}

        It remains to prove the boundedness of the support and the symmetry. For the boundedness, we follow the argument in the proof of Proposition~2.1 in \cite{ajanki2019stability}. Then $\supp\wh V\subset [-(\|\La\|+2\|\cals\|^{1/2}), \|\La\|+2\|\cals\|^{1/2}]$, and the same inclusion holds for $\supp\wh\mu$.

        We now prove the symmetry. Let
        \begin{align}
                J := \begin{pmatrix} I_M & 0    \\
                0   & -I_N\end{pmatrix}.\nonumber
        \end{align}
        By the block off-diagonal form of $\La$ and the definition of $\cals$, we have $J\La J=-\La$ and $\cals[-JRJ]=-J\cals[R]J$. Define $\wt\Pi(\zeta) := -J\wh\Pi(-\overline{\zeta})^*J$. Since $-\overline{\zeta}\in\C^+$, the symmetric MDE \eqref{eq:sym_mde} for $\wh\Pi(-\overline{\zeta})$ gives
        \begin{align}
                -\wh\Pi(-\overline{\zeta})^{-1} = -\overline{\zeta}I_L-\La+\cals[\wh\Pi(-\overline{\zeta})]. \nonumber
        \end{align}
        Taking adjoints yields
        \begin{align}
                -\wh\Pi(-\overline{\zeta})^{*-1} = -\zeta I_L-\La+\cals[\wh\Pi(-\overline{\zeta})^*]. \nonumber
        \end{align}
        Using $J\La J=-\La$ and $\cals[-JRJ]=-J\cals[R]J$, we obtain
        \begin{align}
                -\wt\Pi(\zeta)^{-1} = \zeta I_L-\La+\cals[\wt\Pi(\zeta)].\nonumber
        \end{align}
        Moreover, $\Im\wt\Pi(\zeta) =J\Im\wh\Pi(-\overline\zeta)J\geq0$. Thus $\wt\Pi(\zeta)\in\calM_+$ solves the same symmetric MDE \eqref{eq:sym_mde} as $\wh\Pi(\zeta)$. By uniqueness, $\wh\Pi(\zeta)=-J\wh\Pi(-\overline\zeta)^*J$.
        Using the Stieltjes representation, the right-hand side can be written as
        \begin{align}
                -J\wh\Pi(-\overline{\zeta})^*J = \int_{\R}\frac{J\wh V(-\dd x)J}{x-\zeta}.\nonumber
        \end{align}
        By uniqueness of the matrix-valued Stieltjes representation, $\wh V(B)=J\wh V(-B)J$ for every Borel set $B\subset\R$. Taking normalized traces gives
        \begin{align}
                \wh\mu(B) = \frac1L\Tr\wh V(B) = \frac1L\Tr(J\wh V(-B)J) = \frac1L\Tr\wh V(-B) = \wh\mu(-B). \nonumber
        \end{align}
        Therefore $\wh\mu$ is symmetric with respect to the origin. In particular, $\supp\wh\mu=-\supp\wh\mu$.
        This completes the proof.
\end{proof}



\section{A priori estimates for the MDE solution}
\label{sub:apriori_est}

We  prove Lemmas~\ref{lem:unif_apriori} and \ref{lem:Pi_bd}. The identities \eqref{eq:Im_g}--\eqref{eq:q_Pi11}, established in Section~\ref{subsec:2d_stab}, are the starting point of the argument.

\begin{proof}[Proof of Lemma~\ref{lem:unif_apriori}]
        Since $\Pi^\top=\Pi$ and $\Sg$ is real symmetric, we have
        \begin{align}
                \avg{\Pi_{21}\Sg\Pi_{21}^*}_N=\avg{\Pi_{12}^*\Sg\Pi_{12}}_N. \nonumber
        \end{align}
        We recall the notation $q=\Im m/(\eta+\Im g)$ and the identities \eqref{eq:q_Pi22} and \eqref{eq:q_Pi11} from Section~\ref{subsec:2d_stab}.

        \emph{Proof of \eqref{eq:offdiag_gap}.} We already know from \eqref{eq:Im_m} that $\avg{\Pi_{12}^*\Sg\Pi_{12}}_N\leq1$. Suppose that $1-\avg{\Pi_{12}^*\Sg\Pi_{12}}_N\leq\delta_0$, where $\delta_0>0$ will be chosen sufficiently small. Fix $K_1>1$.

        If $q\leq K_1$, then \eqref{eq:q_Pi22} implies $\|\Pi_{22}\|_{\rr{hs}}^2\leq K_1\delta_0$, and hence $|m|\leq\sqrt{K_1\delta_0}$. For sufficiently small $\delta_0$, the matrix $I_M+m\Sg$ is uniformly invertible. The block equation $(I_M+m\Sg)\Pi_{12}=A\Pi_{22}$ therefore gives
        \begin{align}
                \avg{\Pi_{12}^*\Sg\Pi_{12}}_N=\|\Sg^{1/2}\Pi_{12}\|_{\rr{hs}}^2\lesssim\|\Pi_{22}\|_{\rr{hs}}^2\lesssim K_1\delta_0, \nonumber
        \end{align}
        which contradicts $\avg{\Pi_{12}^*\Sg\Pi_{12}}_N\geq1-\delta_0$ for sufficiently small $\delta_0$.

        If $q\geq K_1$, then \eqref{eq:q_Pi11} yields $\|\Sg^{1/2}\Pi_{11}\Sg^{1/2}\|_{\rr{hs}}^2\lesssim\frac{\delta_0}{K_1}$ and hence $|g|\lesssim\sqrt{\frac{\delta_0}{K_1}}$.
        For sufficiently small $\delta_0$, we have $|z+g|\geq\tau/4$. Using $(z+g)\Pi_{21}=A^\top\Pi_{11}$, we obtain
        \begin{align}
                \avg{\Pi_{21}\Sg\Pi_{21}^*}_N=\|\Pi_{21}\Sg^{1/2}\|_{\rr{hs}}^2\lesssim\frac{\delta_0}{K_1}, \nonumber
        \end{align}
        which again contradicts $\avg{\Pi_{12}^*\Sg\Pi_{12}}_N\geq1-\delta_0$. Choosing $\delta_0$ sufficiently small proves \eqref{eq:offdiag_gap}.

        \emph{Proof of \eqref{eq:diag_size}.} We first estimate $\Sg^{1/2}\Pi_{11}\Sg^{1/2}$. Choose $k_2>0$ sufficiently small such that $\sqrt{k_2}\|\Sg\|\leq\frac12$. If $q\geq k_2$, then \eqref{eq:q_Pi11} gives $\|\Sg^{1/2}\Pi_{11}\Sg^{1/2}\|_{\rr{hs}}\lesssim1$. If $q\leq k_2$, then \eqref{eq:q_Pi22} implies $\|\Pi_{22}\|_{\rr{hs}}^2\leq k_2$, so $|m|\leq\sqrt{k_2}$. By the choice of $k_2$, the matrix $I_M+m\Sg$ is uniformly invertible. The block equation $(I_M+m\Sg)\Pi_{11}=A\Pi_{21}-I_M$, together with $\|\Sg^{1/2}\Pi_{21}\|_{\rr{hs}}\leq 1$, yields $ \|\Sg^{1/2}\Pi_{11}\Sg^{1/2}\|_{\rr{hs}}\lesssim1$.
        For the lower bound, if $|z+g|<\tau/4$, then $|g| = |z+g-z|\geq |z|-|z+g|> E - \tau/4\geq \tau/4 $, and hence $\|\Sg^{1/2}\Pi_{11}\Sg^{1/2}\|_{\rr{hs}}\gtrsim1$. If $|z+g|\geq\tau/4$, then the block equations imply $(I_M+m\Sg-(z+g)^{-1}AA^\top)\Pi_{11}=-I_M$.
        The coefficient on the left-hand side is uniformly bounded, since $|m|\leq\|\Pi_{22}\|_{\rr{hs}}\lesssim1$ follows from $(z+g)\Pi_{22}=A^\top\Pi_{12}-I_N$. Therefore $\|\Sg^{1/2}\Pi_{11}\Sg^{1/2}\|_{\rr{hs}}\gtrsim1$. We have thus proved $\|\Sg^{1/2}\Pi_{11}\Sg^{1/2}\|_{\rr{hs}}\sim1$.

        Moreover, \eqref{eq:q_Pi11} together with $\|\Sg^{1/2}\Pi_{11}\Sg^{1/2}\|_{\rr{hs}}\sim1$ gives $q\lesssim1$, and \eqref{eq:q_Pi22} yields $\|\Pi_{22}\|_{\rr{hs}}\lesssim1$. To obtain the lower bound, choose $k_3>0$ sufficiently small such that $\sqrt{k_3}\|\Sg\|\leq\frac12$. If $|m|\geq k_3$, then $\|\Pi_{22}\|_{\rr{hs}}\geq|m|\geq k_3$. If $|m|<k_3$, then $I_M+m\Sg$ is uniformly invertible, and the second column block equations give $(A^\top(I_M+m\Sg)^{-1}A-(z+g)I_N)\Pi_{22}=I_N $. The coefficient is uniformly bounded, and hence $\|\Pi_{22}\|_{\rr{hs}}\gtrsim1$. This proves \eqref{eq:diag_size}.

        \emph{Proof of \eqref{eq:imag_comp}.} Combining  \eqref{eq:offdiag_gap}, \eqref{eq:diag_size} and \eqref{eq:q_Pi22}, we obtain $q\sim1$. Thus $\Im m\sim\eta+\Im g$. On the other hand, \eqref{eq:diag_size} and \eqref{eq:Im_g}  imply $\Im g\gtrsim\Im m$. Therefore $\Im m\sim\Im g$, which proves \eqref{eq:imag_comp}.
\end{proof}

\begin{proof}[Proof of Lemma~\ref{lem:Pi_bd}]
        Lemma~\ref{lem:unif_apriori} implies $|m|+|g|\lesssim1$ on $\mb D(\tau,\eta_0)$. Hence \eqref{eq:inverse_mde} gives $\|\Pi^{-1}\|\lesssim1$. Moreover,
        \begin{align}
                \Im(-\Pi^{-1})=\begin{pmatrix}\Im m\Sg & 0               \\
               0          & (\eta+\Im g)I_N\end{pmatrix}\gtrsim\Im mI_L, \nonumber
        \end{align}
        where we used \eqref{eq:imag_comp}. Therefore $\|\Pi\|\lesssim(\Im m)^{-1}$. The Stieltjes representations in Proposition~\ref{prop:pi} also give
        \begin{align}
                \|\Pi_{11}(z)\|\lesssim\frac1{\dist(z,\supp\bnu_1)}, \qquad \|\Pi_{22}(z)\|\lesssim\frac1{\dist(z,\supp\bnu_2)}. \nonumber
        \end{align}
        Combining the Stieltjes bounds with $\|\Pi\|\lesssim(\Im m)^{-1}$ proves the claim.
\end{proof}

\section{Error estimates for the perturbed MDE}
\label{app:rand_sc_eq}

Throughout this appendix, Assumptions~\ref{ass:basic} and~\ref{ass:rand_ent} hold. We prove Theorem~\ref{thm:rand_sc_eq}. All estimates below are first proved for a fixed $z$ satisfying \eqref{eq:rand_err_dom}. The corresponding uniform estimates follow by applying Lemma~\ref{lem:stoch_dom_calc} on an $N^{-K}$-net of the domain in \eqref{eq:rand_err_dom} and then using the resolvent identity to extend the estimates to the whole domain. 

\subsection{Ward identities}
\label{sub:ward}
We first prove Lemma~\ref{lem:ward_ids}.

\begin{proof}[Proof of Lemma~\ref{lem:ward_ids}]
        First we notice that
        \begin{align}
                G-G^*=G(I_z-I_{\bar z})G^*=2\ii\eta G\wt I G^*. \nonumber
        \end{align}
        Hence $\Im G=\eta G\wt I G^*$.
        This proves \eqref{eq:ward_2_proj}.

        For the first component, write
        \begin{align}
               \wh\bbx=\begin{pmatrix}\wh\bbx^{(1)} \\
                                   0\end{pmatrix},\qquad
                G^*\wh\bbx=\begin{pmatrix}\bbx^{(1)} \\
                                   \bbx^{(2)}\end{pmatrix}. \nonumber
        \end{align}
        The equation $(G^*)^{-1}G^*\wh\bbx=\wh\bbx$ gives $ -\bbx^{(1)}+Y\bbx^{(2)}=\wh\bbx^{(1)}$.
        By Lemma~\ref{lem:rand_op_norm},
        \begin{align}
                \|G^*\wh\bbx\|^2 =\|\bbx^{(1)}\|^2+\|\bbx^{(2)}\|^2 \prec\|\bbx^{(2)}\|^2+\|\wh\bbx^{(1)}\|^2 =\frac{(\Im G)_{\wh\bbx\wh\bbx}}{\eta}+\|\wh\bbx\|^2,\nonumber
        \end{align}
        where we used $\|\bbx^{(2)}\|^2 = \|\wt I G^*\wh\bbx\|^2$ and \eqref{eq:ward_2_proj}.
        This proves \eqref{eq:ward_1_cmp}. Similar arguments yield \eqref{eq:ward_2_cmp}.
\end{proof}

\subsection{Anisotropic estimate}\label{sub:anis}

This subsection proves the anisotropic estimate \eqref{eq:anis_rand_error} in Theorem~\ref{thm:rand_sc_eq}.

\subsubsection{Resolvent estimates and derivatives}

Write a matrix $R\in\C^{L\times L}$ in block form according to $\C^L=\C^M\oplus\C^N$. The self-energy operator in \eqref{eq:self_energy} can be written as
\begin{align}
        \cals[R] =\avg{\wt I R}_{N}\wt\Sg +\avg{\wt\Sg R}_{N}\wt I. \nonumber
\end{align}
Since $z$ is fixed throughout this subsection, we write
\begin{align}
        \Om(R):=\Om(R,z),\qquad \Om(G):=\Om(G(z),z). \nonumber
\end{align}
Set $ m_N(z):=\avg{\wt G(z)}_{N} =\avg{(W-zI_N)^{-1}}_{N}$.
By \eqref{eq:G_triv_bd}, we have
\begin{align}
        \|G(z)\|\leq\eta^{-1},\qquad z\in\mb D(\tau,\eta_0). \nonumber
\end{align}
By the definition of $G = (\calH+\La-I_z)^{-1}$,
\begin{align}
        I_L+I_zG=\calH G+\La G. \label{eq:lin_res_id}
\end{align}
For a deterministic vector $\bbv\in\C^L$, write
\begin{align}
        \bhv:=\wh I\bbv, \qquad \btv:=\wt I\bbv. \nonumber
\end{align}
Since $\bhv$ is supported on the first block, \eqref{eq:lin_res_id} and the explicit form of $\cals[G]$ give
\begin{align}
        \Om(G)_{\bhv\bbw} = I_L + (I_z+\calS[G]-\La)G =(\calH G)_{\bhv\bbw}+m_N(z)(\wt\Sg G)_{\bhv\bbw}. \label{eq:anis_err_id}
\end{align}
For a first test vector supported on the second block, set
\begin{align}
        g_N(z):=\avg{\wt\Sg G(z)}_N=\frac1N\Tr\paren{\Sg\wh G(z)}. \nonumber
\end{align}
Then
\begin{align}
        \Om(G)_{\btv\bbw} =(\calH G)_{\btv\bbw}+g_N(z)(\wt I G)_{\btv\bbw}. \label{eq:anis_err_id_2}
\end{align}
We now focus on $\Om(G)_{\bhv\bbw}$. The term $\Om(G)_{\btv\bbw}$ is treated similarly, so we postpone its discussion and return to \eqref{eq:anis_err_id_2} at the end of this subsection.

Set
\begin{align}
        \Psi(z):=(1+\phi)^3\sqrt{\frac{\vphi+\eta}{N\eta}}. \nonumber
\end{align}
Under the assumptions of Theorem~\ref{thm:rand_sc_eq}, Lemma~\ref{lem:ward_ids} gives
\begin{align}
        \sum_{a=1}^{L}\abs{G_{\wh\bbx a}}^2 +\sum_{a=1}^{L}\abs{G_{\wt\bbx a}}^2 \prec\frac{\vphi+\eta}{\eta}\|\bbx\|^2. \label{eq:ward_vec_bd}
\end{align}
Consequently, by $\bbx=\wh\bbx+\wt\bbx$, and apply \eqref{eq:ward_vec_bd} and $G^\top=G$, we have
\begin{align}
        \|G\bbx\|^2+\|G^*\bbx\|^2 \prec\frac{\vphi+\eta}{\eta}\|\bbx\|^2. \label{eq:gh_vec_bd}
\end{align}

\begin{lemma}\label{lem:ward_cons}
        Under the assumptions of Theorem~\ref{thm:rand_sc_eq}, let $B_1,B_2\in\C^{L\times L}$ be deterministic matrices with $\|B_1\|+\|B_2\|\lesssim1$. Uniformly for deterministic vectors $\bbx\in\C^L$,
        \begin{align}
                \sum_{a=1}^{L}\abs{(B_1GB_2)_{a\bbx}}^2 & \prec\frac{\vphi+\eta}{\eta}\|\bbx\|^2,\qquad \max_{1\leq a\leq L}\abs{(B_1GB_2)_{a\bbx}}\prec(1+\phi)\|\bbx\|, \nonumber \\
                \sum_{a=1}^{L}\abs{(B_1GB_2)_{a\bbx}}^4 & \prec(1+\phi)^2\frac{\vphi+\eta}{\eta}\|\bbx\|^4. \nonumber
        \end{align}
        The same three estimates hold with $G$ replaced by $G^*$.
\end{lemma}
\begin{proof}
        By \eqref{eq:gh_vec_bd},
        \begin{align}
                \sum_{a=1}^{L}\abs{(B_1GB_2)_{a\bbx}}^2 =\|B_1GB_2\bbx\|^2 \leq\|B_1\|^2\|GB_2\bbx\|^2 \prec\frac{\vphi+\eta}{\eta}\|\bbx\|^2. \nonumber
        \end{align}
        Moreover, $(B_1GB_2)_{a\bbx}=\scalar{B_1^*\bbe_a}{GB_2\bbx}$. The anisotropic hypothesis, together with $\|\Pi\|\lesssim1$, yields the second estimate uniformly in $a$. The third estimate follows from
        \begin{align}
                \sum_{a=1}^{L}\abs{(B_1GB_2)_{a\bbx}}^4 \leq\max_{1\leq a\leq L}\abs{(B_1GB_2)_{a\bbx}}^2 \sum_{a=1}^{L}\abs{(B_1GB_2)_{a\bbx}}^2. \nonumber
        \end{align}
        The estimates with $G^*$ follow by adjunction.
\end{proof}

\begin{lemma}\label{lem:entry_res_der}
        Under the assumptions of Theorem~\ref{thm:rand_sc_eq}, for every fixed integer $r\geq1$,
        \begin{align}
                \partial_{i\mu}^rG =(-1)^r r!G(V^{i\mu}G)^r. \label{eq:res_der_form}
        \end{align}
        Uniformly for deterministic vectors $\bbx,\bbw\in C^L$, 
        \begin{align}
                \abs{\partial_{i\mu}^rG_{\bbx\bbw}}                   & \prec(1+\phi)^{r+1}\|\bbx\|\|\bbw\|, \label{eq:res_entry_der}                                  \\
                \sum_{\mu=M+1}^{L} \abs{\partial_{i\mu}^rG_{\bbx\mu}} & \prec(1+\phi)^r\sqrt N \sqrt{\frac{\vphi+\eta}{\eta}}\|\bbx\|. \label{eq:res_der_ward}
        \end{align}
\end{lemma}

\begin{proof}
        Formula \eqref{eq:res_der_form} follows by induction from $\partial_{i\mu}G=-GV^{i\mu}G$. Since $V^{i\mu}$ has rank two, every entry of $\partial_{i\mu}^rG$ is a sum of at most $2^r$ monomials. Each monomial contains $r+1$ entries of $G$, with intermediate indices in $\{i,\mu\}$. The anisotropic hypothesis bounds each factor by $O_\prec(1+\phi)$, which proves \eqref{eq:res_entry_der}.
        For \eqref{eq:res_der_ward}, keep one factor carrying the free index $\mu$. Apply Cauchy--Schwarz inequality to the $\mu$-sum and use Lemma~\ref{lem:ward_ids}. The remaining factors are bounded by $O_\prec(1+\phi)$, which proves \eqref{eq:res_der_ward}. 
\end{proof}

By homogeneity, it suffices to consider deterministic unit vectors throughout the following subsection. It also suffices to prove the anisotropic estimate for real test vectors. The complex case then follows by decomposing both vectors into their real and imaginary parts.

\subsubsection{Expectation estimate}
We now display the cancellation mechanism by considering the expectation of $\Om(G)_{\bhv\bbw}$.
We first expand the first term on the right-hand side of \eqref{eq:anis_err_id}. By Lemma~\ref{lem:cum_exp}, we have
\begin{align}
        \bb E (\calH G)_{\bhv\bbw} & = \sum_{i,\mu}\bb E G_{\mu\bbw}X_{i\mu}\Sg^{1/2}_{\bhv i}  = \sum_{k=1}^l\sum_{i,\mu}\frac{\kappa_{k+1}(X_{i\mu})}{k!}\bb E\partial_{i\mu}^k(G_{\mu\bbw})\Sg^{1/2}_{\bhv i} + \calR_{l+1}  := \sum_{k=1}^l T_k + \calR_{l+1}, \nonumber
\end{align}
where the derivatives are understood in the sense of \eqref{eq:entry_res_der}.
By Lemma~\ref{lem:cum_bd}, we have
\begin{align}
        T_1 & =-\frac1N\sum_{i,\mu}\bb E( (G\wh\Sg^{1/2})_{\mu i}G_{\mu\bbw} +G_{\mu\mu}(\wh\Sg^{1/2}G)_{i\bbw} )\Sg^{1/2}_{\bhv i} =-\bb E(m_N(z)(\wt\Sg G)_{\bhv\bbw}) -\frac1N\bb E((\wt\Sg G\wt I G)_{\bhv\bbw}). \nonumber
\end{align} 
Consequently,
\begin{align}
        \bb E\Om(G)_{\bhv\bbw} =\sum_{k=2}^lT_k+\calR_{l+1} +O_\prec\left(\frac{\vphi+\eta}{N\eta}\right), \label{eq:anis_exp_cancel}
\end{align}
where Cauchy--Schwarz inequality and Lemma~\ref{lem:ward_cons} give
\begin{align}
        \abs{(\wt\Sg G\wt I G)_{\bhv\bbw}} \leq\sum_{\mu=1}^N|(\wt\Sg G)_{\bhv\mu}||G_{\mu\bbw}| \prec\frac{\vphi+\eta}{\eta}. \nonumber
\end{align}
By \eqref{eq:ward_input_size} and $\|\Pi(z)\|\lesssim1$, we have $\vphi\lesssim1+\phi\lesssim N^{\tau/10}$. Since $N\eta\geq N^\tau$,
\begin{align}
        \frac{\vphi+\eta}{N\eta}\leq N^{-9\tau/10} + N^{-1}<1. \nonumber
\end{align}
Hence the last error is $O_\prec(\Psi(z))$.
By Lemma~\ref{lem:cum_bd}, we have $\kappa_{k+1}(X_{i\mu})=O_k(N^{-(k+1)/2})$ for $k\geq 2$. Thus, we obtain for $k\geq 2$,
\begin{align}
        T_k \prec N^{-(k+1)/2} \sum_{i,\mu} \bb E |\partial_{i\mu}^k (G_{\mu\bbw})\Sg^{1/2}_{\bhv i} |\prec N^{-(k-1)/2}(1+\phi)^k\sqrt{\frac{\vphi+\eta}{\eta}}\prec\Psi(z), \label{eq:exp_high_cum}
\end{align}
where we used Lemma~\ref{lem:entry_res_der} and $\sum_i|\Sg^{1/2}_{\bhv i}|\leq \sqrt{N\sum_i|\Sg^{1/2}_{\bhv i}|^2}\lesssim \sqrt{N}$.
For the remainder term $\calR_{l+1}$, set
\begin{align}
        \calH^{(i\mu)}:=\calH-X_{i\mu}V^{i\mu}, \qquad G^{(i\mu)}(s):=(\calH^{(i\mu)}+sV^{i\mu}+\La-I_z)^{-1}.\nonumber
\end{align} Then $G=G^{(i\mu)}(X_{i\mu})$.
Similar to the \eqref{eq:G_triv_bd}, by the Schur complement formula, we have
\begin{align}
        \|G^{(i\mu)}(s)\|\lesssim\eta^{-1},\qquad z\in\mb D(\tau,\eta_0). \label{eq:interp_res}
\end{align}
 Applying Lemma~\ref{lem:cum_exp} to $F_{i\mu}(s):=G^{(i\mu)}(s)_{\mu\bbw}$, we obtain
\begin{align}
        \abs{\calR_{l+1}^{i\mu}} & \leq O(1)\left( \bb E\sup_{\abs s\leq\abs{X_{i\mu}}} \abs{\frac{\dd^{l+1}}{\dd s^{l+1}}F_{i\mu}(s)}^2 \bb E\left(\abs{X_{i\mu}}^{2l+4} 1_{\{\abs{X_{i\mu}}>K\}}\right) \right)^{1/2}\nonumber \\
                                 & \quad+ O(1)\bb E\abs{X_{i\mu}}^{l+2} \bb E\sup_{\abs s\leq K} \abs{\frac{\dd^{l+1}}{\dd s^{l+1}}F_{i\mu}(s)}, \label{eq:entry_cum_rem}
\end{align}
and hence $\calR_{l+1} = \sum_{i,\mu}\calR_{l+1}^{i\mu}\Sg_{\bhv i}^{1/2}$.
For every fixed $l$ and every $D_0>0$, the moment assumptions imply
\begin{align}
        \bb E\abs{X_{i\mu}}^{l+2}                                        & \lesssim N^{-(l+2)/2}, \label{eq:entry_moment_bd}         \\
        \bb E\left(\abs{X_{i\mu}}^{2l+4} 1_{\{\abs{X_{i\mu}}>K\}}\right) & \lesssim N^{-(l+2+D_0)}, \label{eq:trunc_mom_bd}
\end{align}
where $K=N^{-1/2+\tau/5}$ and the moment order in Markov's inequality is chosen sufficiently large.

We next record the required bounds for the interpolating resolvent.
The resolvent identity, iterated once, gives
\begin{align}
        G^{(i\mu)}(0) & =G+X_{i\mu}GV^{i\mu}G +X_{i\mu}^2G^{(i\mu)}(0)(V^{i\mu}G)^2, \label{eq:loo_exp} \\
        G^{(i\mu)}(s) & =G^{(i\mu)}(0) -sG^{(i\mu)}(0)V^{i\mu}G^{(i\mu)}(0) +s^2G^{(i\mu)}(s)(V^{i\mu}G^{(i\mu)}(0))^2. \label{eq:loo_exp_2}
\end{align}
By Assumption~\ref{ass:rand_ent} and \eqref{eq:interp_res}, we have
\begin{align}
        \|G^{(i\mu)}(0)\|\lesssim\eta^{-1}\leq N^{1-\tau}, \qquad X_{i\mu}\prec N^{-1/2}. \nonumber
\end{align}
Using \eqref{eq:loo_exp} and the anisotropic hypothesis, we obtain
\begin{align}
        \abs{G^{(i\mu)}(0)_{\bbx\bbw}} \prec1+\phi. \label{eq:loo_mat_bd}
\end{align}
 Since $\phi\leq N^{\tau/10}$ and $K=N^{-1/2+\tau/5}$, by \eqref{eq:loo_exp_2}, for every fixed integer $r\geq0$,
\begin{align}
         \sup_{\abs s\leq K} \abs{\frac{\dd^r}{\dd s^r} F_{i\mu}(s)} = \sup_{\abs s\leq K} \abs{\frac{\dd^r}{\dd s^r} G^{(i\mu)}(s)_{\bbx\bbw}} \prec(1+\phi)^{r+1}. \label{eq:interp_res_der}
\end{align}
Thus for the second term on the right-hand side of \eqref{eq:entry_cum_rem}, in which the perturbation parameter satisfies $|s|\leq K$, is $O(N^{-D})$ by \eqref{eq:entry_moment_bd} and \eqref{eq:interp_res_der}, after choosing $l=l(D)$ sufficiently large. 

For the first term on the right-hand side of \eqref{eq:entry_cum_rem}, the trival bound \eqref{eq:interp_res}  gives
\begin{align}
        \bb E\sup_{\abs s\leq\abs{X_{i\mu}}} \abs{\frac{\dd^{l+1}}{\dd s^{l+1}} F_{i\mu}(s)}^2 \leq N^{C_l}. \nonumber
\end{align}
After $l$ is fixed, we choose $D_0$ in \eqref{eq:trunc_mom_bd} sufficiently large. Hence the term also is bounded by $O(N^{-D})$. So we conclude that $\abs{\calR_{l+1}^{i\mu}}=O(N^{-D})$ uniformly in $i$ and $\mu$.
Summing over $i$ and $\mu$ yields
\begin{align}
        \abs{\calR_{l+1}} \leq\sum_{i,\mu}\abs{\calR_{l+1}^{i\mu}} \abs{\Sg_{\bhv i}^{1/2}} =O(N^{-D}). \label{eq:anis_exp_rem}
\end{align}
Combining \eqref{eq:exp_high_cum} and \eqref{eq:anis_exp_rem}, we have
\begin{align}
        \abs{\bb E\Om(G)_{\bhv\bbw}}\prec\Psi(z). \nonumber
\end{align}
\subsubsection{High moments}

We next bound $\bb E|\Om(G)_{\bhv\bbw}|^{2p}$ for any fixed integer $p \geq 1$.
Using \eqref{eq:anis_err_id} and applying Lemma~\ref{lem:cum_exp} to the terms containing $X_{i\mu}$, we obtain
\begin{align}
        \bb E\abs{\Om(G)_{\bhv\bbw}}^{2p} & =\bb E\left((\calH G)_{\bhv\bbw} \Om(G)_{\bhv\bbw}^{p-1} \overline{\Om(G)_{\bhv\bbw}}^p\right) +\bb E\left(m_N(z)(\wt\Sg G)_{\bhv\bbw} \Om(G)_{\bhv\bbw}^{p-1} \overline{\Om(G)_{\bhv\bbw}}^p\right)\nonumber                                                                                              \\
                                          & =\sum_{i,\mu}\bb E\left(G_{\mu\bbw}X_{i\mu}\Sg_{\bhv i}^{1/2} \Om(G)_{\bhv\bbw}^{p-1} \overline{\Om(G)_{\bhv\bbw}}^p\right) +\bb E\left(m_N(z)(\wt\Sg G)_{\bhv\bbw} \Om(G)_{\bhv\bbw}^{p-1} \overline{\Om(G)_{\bhv\bbw}}^p\right)\nonumber                                                                                              \\
                                          & =\sum_{k=1}^l\sum_{i,\mu} \frac{\kappa_{k+1}(X_{i\mu})}{k!} \bb E\partial_{i\mu}^k( G_{\mu\bbw}\Sg_{\bhv i}^{1/2} \Om(G)_{\bhv\bbw}^{p-1} \overline{\Om(G)_{\bhv\bbw}}^p) +\calR_{l+1}^{(1)} \nonumber \\
                                          & \quad+\bb E\left(m_N(z)(\wt\Sg G)_{\bhv\bbw} \Om(G)_{\bhv\bbw}^{p-1} \overline{\Om(G)_{\bhv\bbw}}^p\right)\nonumber                                                                                              \\
                                          & :=\sum_{k=1}^lT_k^{(1)}+\calR_{l+1}^{(1)} +\bb E\left(m_N(z)(\wt\Sg G)_{\bhv\bbw} \Om(G)_{\bhv\bbw}^{p-1} \overline{\Om(G)_{\bhv\bbw}}^p\right). \label{eq:anis_mom_exp}
\end{align}
For $T_{1}^{(1)}$, the Leibniz rule gives
\begin{align}
        T^{(1)}_1 & = \sum_{i,\mu}\frac{1}{N}\bb E \partial_{i\mu}(G_{\mu\bbw} \Sg_{\bhv i}^{1/2}\Om(G)_{\bhv\bbw}^{p-1} \overline{\Om(G)_{\bhv\bbw}}^{p}) \nonumber                        \\
                  & =\frac1N\sum_{i,\mu}\bb E\left( \partial_{i\mu}(G_{\mu\bbw})\Sg_{\bhv i}^{1/2} \Om(G)_{\bhv\bbw}^{p-1} \overline{\Om(G)_{\bhv\bbw}}^p\right)+\frac1N\sum_{i,\mu}\bb E\left( G_{\mu\bbw}\Sg_{\bhv i}^{1/2} \partial_{i\mu}( \Om(G)_{\bhv\bbw}^{p-1} \overline{\Om(G)_{\bhv\bbw}}^p)\right) \nonumber \\
                  & :=T_{11}^{(1)}+T_{12}^{(1)}. \label{eq:cum2_dec}
\end{align}
For $T_{11}^{(1)}$, repeating the computation leading to \eqref{eq:anis_exp_cancel} gives
\begin{align}
        T_{11}^{(1)}  =-\bb E\left(m_N(z)(\wt\Sg G)_{\bhv\bbw} \Om(G)_{\bhv\bbw}^{p-1}\overline{\Om(G)_{\bhv\bbw}}^p\right)-\frac1N\bb E\left((\wt\Sg G\wt I G)_{\bhv\bbw} \Om(G)_{\bhv\bbw}^{p-1}\overline{\Om(G)_{\bhv\bbw}}^p\right). \label{eq:cum2_cancel}
\end{align}
The remaining $N^{-1}$ term is bounded by
\begin{align}
          \frac1N\bb E\left( \abs{(\wt\Sg G\wt I G)_{\bhv\bbw}} \abs{\Om(G)_{\bhv\bbw}}^{2p-1} \right)  \prec \frac{\vphi+\eta}{N\eta} \bb E\abs{\Om(G)_{\bhv\bbw}}^{2p-1} \prec \Psi(z)\bb E\abs{\Om(G)_{\bhv\bbw}}^{2p-1}, \nonumber
\end{align}
where the first step follows from Lemma~\ref{lem:ward_cons} and the Cauchy-Schwarz inequality,  and the last step follows from $({\vphi+\eta})/({N\eta})<1$.
For the second term $T_{12}^{(1)}$, we have
\begin{align}
        T_{12}^{(1)} & = \frac{1}{N} \sum_{i,\mu} \bb E G_{\mu\bbw} \Sg_{\bhv i}^{1/2}\partial_{i\mu}(\Om(G)_{\bhv\bbw}^{p-1} \overline{\Om(G)_{\bhv\bbw}}^{p}) \nonumber                                      \\
                     & = \frac{p-1}{N} \sum_{i,\mu} \bb E G_{\mu\bbw} \Sg_{\bhv i}^{1/2}\partial_{i\mu}(\Om(G)_{\bhv\bbw}) \Om(G)_{\bhv\bbw}^{p-2} \overline{\Om(G)_{\bhv\bbw}}^{p} \nonumber                  \\
                     & \quad + \frac{p}{N} \sum_{i,\mu} \bb E G_{\mu\bbw} \Sg_{\bhv i}^{1/2}\partial_{i\mu}(\overline{\Om(G)_{\bhv\bbw}}) \Om(G)_{\bhv\bbw}^{p-1} \overline{\Om(G)_{\bhv\bbw}}^{p-1} \nonumber \\
                     & := T_{121}^{(1)} + T_{122}^{(1)}. \label{eq:cum2_der}
\end{align}
We return to these terms below. We first estimate $\sum_{k=2}^{l}T_{k}^{(1)}$ and $\calR_{l+1}^{(1)}$. For any fixed integer $k \geq 2$, we have
\begin{align}
        T_k^{(1)} & = \sum_{i,\mu} \frac{\kappa_{k+1}(X_{i\mu})}{k!} \bb E \partial_{i\mu}^k (G_{\mu\bbw} \Sg_{\bhv i}^{1/2}\Om(G)_{\bhv\bbw}^{p-1} \overline{\Om(G)_{\bhv\bbw}}^{p}) \nonumber \\
                  & \leq O(N^{-(k+1)/2}) \sum_{i,\mu} \bb E |\partial_{i\mu}^k (G_{\mu\bbw} \Sg_{\bhv i}^{1/2}\Om(G)_{\bhv\bbw}^{p-1} \overline{\Om(G)_{\bhv\bbw}}^{p})| \nonumber              \\
                  & = O(N^{-(k+1)/2}) \sum_{\substack{r,s,n\geq 0                                                                                                                                \\
                        r+s+n = k}} \sum_{i,\mu}\bb E |\Sg_{\bhv i}^{1/2}(\partial_{i\mu}^r G_{\mu\bbw}) ( \partial_{i\mu}^s \Om(G)_{\bhv\bbw}^{p-1} ) (\partial_{i\mu}^n \overline{\Om(G)_{\bhv\bbw}}^{p}) |. \label{eq:cum_hi_exp}
\end{align}
After expanding the derivatives of the powers of $\Om$ and $\overline{\Om}$ separately, all resulting terms have the same structure after taking absolute values. Thus it suffices to bound terms of the form
\begin{align}
         & O(N^{-(k+1)/2})\sum_{i,\mu}\bb E\left( \abs{\Sg_{\bhv i}^{1/2}\partial_{i\mu}^rG_{\mu\bbw}} \prod_{j=1}^q \abs{\partial_{i\mu}^{l_j}\Om(G)_{\bhv\bbw}} \abs{\Om(G)_{\bhv\bbw}}^{2p-1-q} \right), \label{eq:cum_hi_term}
\end{align}
where $0\leq r\leq k$, $0\leq q\leq(k-r)\wedge(2p-1)$, $l_j\geq1$, and $l_1+\dots+l_q=k-r$. Now we need to bound $\partial_{i\mu}^l\Om(G)_{\bhv\bbw}$ for fixed $l\geq1$. By \eqref{eq:anis_err_id},
\begin{align}
        \partial_{i\mu}^l\Om(G)_{\bhv\bbw} & = \partial_{i\mu}^l(\calH G)_{\bhv\bbw} +\partial_{i\mu}^l (m_N(z)(\wt\Sg G)_{\bhv\bbw}). \label{eq:err_der_dec}
\end{align}
The resolvent identity \eqref{eq:lin_res_id} gives $\calH G=I_L+(I_z-\La)G$. Since $I_z-\La$ is deterministic,
\begin{align}
        \partial_{i\mu}^l(\calH G) =(I_z-\La)\partial_{i\mu}^lG, \qquad l\geq1. \label{eq:HG_der}
\end{align}
Hence Lemma~\ref{lem:entry_res_der} yields
\begin{align}
        \abs{\partial_{i\mu}^l(\calH G)_{\bhv\bbw}} \prec(1+\phi)^{l}\max_{a\in\{i,\mu\}}|(\wt\Sg^{1/2}G)_{a\bbw}| \prec(1+\phi)^{l}\sqrt{\frac{\vphi+\eta}{\eta}},\qquad l\geq 1, \label{eq:HG_der_bd}
\end{align}
where we used $|\avg{\wt\Sg^{1/2}\bbe_a,G\bbw}|\leq\|\wt\Sg^{1/2}a\|\|G\bbw\|$ and  \eqref{eq:gh_vec_bd} in the last step. 
Morever, by \eqref{eq:ward_2_proj} and \eqref{eq:res_der_form}, we have
\begin{align}
        \abs{\partial_{i\mu}^l m_N(z)}=\abs{\frac{1}{N}\sum_\nu((GV^{i\mu})^lG )_{\nu\nu}} \prec (1+\phi)^{l-1}\frac{\vphi+\eta}{N\eta},\qquad l\geq 1. \label{eq:tr_der_bd}
\end{align}
For the self-energy term, the Leibniz rule gives
\begin{align}
        \partial_{i\mu}^l (m_N(z)(\wt\Sg G)_{\bhv\bbw}) & = \sum_{s=0}^l\binom ls \partial_{i\mu}^{l-s}(m_N(z)) \partial_{i\mu}^s(\wt\Sg G)_{\bhv\bbw}. \label{eq:self_der}
\end{align}
If $s=l$, similar to \eqref{eq:HG_der_bd}, we have
\begin{align}
|m_N(z)\partial_{i\mu}^l(\wt\Sg G)_{\bhv\bbw}|=|m_N(z)(\wt\Sg \partial_{i\mu}^lG)_{\bhv\bbw}|\prec(1+\phi)^{l+1}\sqrt{\frac{\vphi+\eta}{\eta}},\nonumber
\end{align} 
where we used the trival bound $m_N(z)\prec (1+\phi)$. 
If $s<l$, by Lemma~\ref{lem:ward_ids} and \eqref{eq:tr_der_bd}, we obtain
\begin{align}
        \abs{\partial_{i\mu}^l (m_N(z)(\wt\Sg G)_{\bhv\bbw})} \prec (1+\phi)^{l+1} \sqrt{\frac{\vphi+\eta}{\eta}},\qquad l\geq 1. \label{eq:self_der_bd}
\end{align}
Combining \eqref{eq:err_der_dec}, \eqref{eq:HG_der_bd}, and \eqref{eq:self_der_bd}, we conclude that
\begin{align}
        \abs{\partial_{i\mu}^l\Om(G)_{\bhv\bbw}} \prec (1+\phi)^{l+1} \sqrt{\frac{\vphi+\eta}{\eta}},\qquad l\geq 1. \label{eq:err_map_der}
\end{align}
Then we obtain
\begin{align}
        \eqref{eq:cum_hi_term} & \prec N^{-(k+1)/2} (1+\phi)^{k-r+q} \left(\sqrt{\frac{\vphi + \eta}{N\eta}}\right)^q N^{q/2} \bb E \sum_{i,\mu} \left( \left|\Sg_{\bhv i}^{1/2}(\partial_{i\mu}^r G_{\mu\bbw}) \Om(G)_{\bhv\bbw}^{2p-1-q}\right|\right) \nonumber \\
                                                & \prec N^{-(k+1)/2+q/2+1/2+1}(1+\phi)^{k-r+q+r}  \left(\sqrt{\frac{\vphi + \eta}{N\eta}}\right)^{q+1} \bb E|\Om(G)_{\bhv\bbw}|^{2p-1-q} \nonumber                                                                                              \\
                                                & = N^{(q+2-k)/2}(1+\phi)^{(k+q)} \left(\sqrt{\frac{\vphi + \eta}{N\eta}}\right)^{q+1} \bb E|\Om(G)_{\bhv\bbw}|^{2p-1-q} \nonumber                                                                    \\
                                                & = (N^{-1/2}(1+\phi))^{k-2-q}(1+\phi)^{2q+2} \left(\sqrt{\frac{\vphi + \eta}{N\eta}}\right)^{q+1} \bb E|\Om(G)_{\bhv\bbw}|^{2p-1-q}. \nonumber
\end{align}
For $q\leq k-2$, since $N^{-1/2}(1+\phi)<1$, we obtain
\begin{align}
        \eqref{eq:cum_hi_term}
        \prec(1+\phi)^{3(q+1)}
        \left(\sqrt{\frac{\vphi+\eta}{N\eta}}\right)^{q+1}
        \bb E\abs{\Om(G)_{\bhv\bbw}}^{2p-1-q}. \label{eq:cum_gen_bd}
\end{align}
Thus the desired estimate holds whenever $q\leq k-2$. It remains to consider $q\geq k-1$.

Since $q\geq k-1$ and $q\leq k-r$, we have $r\leq1$. Hence only the three cases $(r,q)=(0,k)$, $(1,k-1)$, and $(0,k-1)$ remain.

First suppose that $q = k$ and $r = 0$. In this case $l_1=\dots=l_q=1$, so we have
\begin{align}
        \eqref{eq:cum_hi_term} & = O(N^{-(k+1)/2}) \sum_{i,\mu}\bb E \left( \left|\Sg_{\bhv i}^{1/2} G_{\mu\bbw} (\partial_{i\mu}\Om(G)_{\bhv\bbw})^{k} \Om(G)_{\bhv\bbw}^{2p-1-k}\right|\right) \nonumber \\
                                                & \prec N^{-3/2}(1+\phi)^{2k-4}\left(\sqrt{\frac{\vphi + \eta}{N\eta}}\right)^{k-2}  \nonumber\\
                                                &\quad \sum_{i,\mu}\bb E \left( \left|\Sg_{\bhv i}^{1/2} G_{\mu\bbw} (\partial_{i\mu}\Om(G)_{\bhv\bbw})^{2} \Om(G)_{\bhv\bbw}^{2p-1-k}\right|\right) , \label{eq:cum_0k}
\end{align}
where \eqref{eq:err_map_der} with $l=1$ was used for $k-2$ of the derivative factors. We now estimate $\partial_{i\mu}\Om(G)_{\bhv\bbw}$. By \eqref{eq:anis_err_id}, we have
\begin{align}
        \partial_{i\mu}\Om(G)_{\bhv\bbw} =\partial_{i\mu}(\calH G)_{\bhv\bbw} +\partial_{i\mu}(m_N(z)(\wt\Sg G)_{\bhv\bbw}). \label{eq:err_der1}
\end{align}
Since $\calH G=I_L+(I_z-\La)G$, we have
\begin{align}
        \partial_{i\mu}(\calH G)_{\bhv\bbw} & =-((I_z-\La)G\wh\Sg^{1/2})_{\bhv i}G_{\mu\bbw}-((I_z-\La)G)_{\bhv\mu} (\wh\Sg^{1/2}G)_{i\bbw}. \label{eq:HG_der1}
\end{align}
The estimates in Lemma~\ref{lem:ward_cons} also apply to row sums by adjunction. Hence
\begin{align}
         & \sum_i\abs{\Sg_{\bhv i}^{1/2}} \abs{((I_z-\La)G\wh\Sg^{1/2})_{\bhv i}}^2 \prec(1+\phi)\sqrt{\frac{\vphi+\eta}{\eta}}, \nonumber \\
         & \sum_\mu\abs{G_{\mu\bbw}}^3 \prec(1+\phi)\frac{\vphi+\eta}{\eta}. \nonumber
\end{align}
Therefore,
\begin{align}
        \sum_{i,\mu}\abs{\Sg_{\bhv i}^{1/2}G_{\mu\bbw}}\abs{((I_z-\La)G\wh\Sg^{1/2})_{\bhv i}}^2\abs{G_{\mu\bbw}}^2\prec(1+\phi)^2\left(\frac{\vphi+\eta}{\eta}\right)^{3/2}. \nonumber
\end{align}
Similarly, for the second product in \eqref{eq:HG_der1}, we also have
\begin{align}
        \sum_{i,\mu}\abs{\Sg_{\bhv i}^{1/2}G_{\mu\bbw}}\abs{((I_z-\La)G)_{\bhv\mu}}^2\abs{(\wh\Sg^{1/2}G)_{i\bbw}}^2\prec(1+\phi)^2\left(\frac{\vphi+\eta}{\eta}\right)^{3/2}. \nonumber
\end{align}
Together with \eqref{eq:HG_der1}, these estimates yield
\begin{align}
         & \sum_{i,\mu}\bb E\left( \abs{\Sg_{\bhv i}^{1/2}G_{\mu\bbw}} \abs{\partial_{i\mu}(\calH G)_{\bhv\bbw}}^2 \abs{\Om(G)_{\bhv\bbw}}^{2p-1-k} \right)\nonumber \\
         & \prec(1+\phi)^2 \left(\frac{\vphi+\eta}{\eta}\right)^{3/2} \bb E\abs{\Om(G)_{\bhv\bbw}}^{2p-1-k}. \label{eq:err_HG_bd}
\end{align}
For the self-energy term, the product rule gives
\begin{align}
          \partial_{i\mu}(m_N(z)(\wt\Sg G)_{\bhv\bbw}) & =(\partial_{i\mu}m_N(z))(\wt\Sg G)_{\bhv\bbw}+m_N(z)(\partial_{i\mu}(\wt\Sg G)_{\bhv\bbw})  \nonumber \\
          & =(\partial_{i\mu}m_N(z))(\wt\Sg G)_{\bhv\bbw} \nonumber\\
          &\quad - m_N(z)\left( (\wt\Sg G\wh\Sg^{1/2})_{\bhv i}G_{\mu\bbw} +(\wt\Sg G)_{\bhv\mu}(\wh\Sg^{1/2}G)_{i\bbw} \right). \label{eq:mG_der1}
\end{align}
The  terms associated with $\partial_{i\mu}(\wt\Sg G)_{\bhv\bbw}$ are estimated as in \eqref{eq:err_HG_bd}, with $I_z-\La$ replaced by $\wt\Sg$. Since $\abs{m_N(z)}^2\prec (1+\phi)^2$, we obtain
\begin{align}
         & \sum_{i,\mu}\bb E\left( \abs{\Sg_{\bhv i}^{1/2}G_{\mu\bbw}} \abs{m_N(z)\partial_{i\mu}(\wt\Sg G)_{\bhv\bbw}}^2 \abs{\Om(G)_{\bhv\bbw}}^{2p-1-k} \right)  \nonumber\\
         &\quad \prec(1+\phi)^4 \left(\frac{\vphi+\eta}{\eta}\right)^{3/2} \bb E\abs{\Om(G)_{\bhv\bbw}}^{2p-1-k}. \nonumber
\end{align}
Moreover, by \eqref{eq:tr_der_bd}, we have
\begin{align}
         & \sum_{i,\mu}\bb E\left( \abs{\Sg_{\bhv i}^{1/2}G_{\mu\bbw}} \abs{\partial_{i\mu}m_N(z)}^2 \abs{(\wt\Sg G)_{\bhv\bbw}}^2 \abs{\Om(G)_{\bhv\bbw}}^{2p-1-k} \right) \nonumber \\
         & \prec(1+\phi)^3\frac{\vphi+\eta}{N\eta}\sum_{i,\mu} \bb E\abs{\Sg_{\bhv i}^{1/2}G_{\mu\bbw}}\abs{\Om(G)_{\bhv\bbw}}^{2p-1-k} \nonumber\\
         &\prec(1+\phi)^3\left(\frac{\vphi+\eta}{\eta}\right)^{3/2} \bb E\abs{\Om(G)_{\bhv\bbw}}^{2p-1-k}. \nonumber
\end{align}
Combining the last two bounds gives
\begin{align}
         & \sum_{i,\mu}\bb E\left( \abs{\Sg_{\bhv i}^{1/2}G_{\mu\bbw}} \abs{\partial_{i\mu}(m_N(z)(\wt\Sg G)_{\bhv\bbw})}^2 \abs{\Om(G)_{\bhv\bbw}}^{2p-1-k} \right)\nonumber \\
         & \prec(1+\phi)^4 \left(\frac{\vphi+\eta}{\eta}\right)^{3/2} \bb E\abs{\Om(G)_{\bhv\bbw}}^{2p-1-k}. \label{eq:self_en_der}
\end{align}
Combining \eqref{eq:cum_0k}, \eqref{eq:err_der1}, \eqref{eq:err_HG_bd}, and \eqref{eq:self_en_der}, we conclude that
\begin{align}
        \eqref{eq:cum_hi_term}\prec(1+\phi)^{2k}\left(\sqrt{\frac{\vphi+\eta}{N\eta}}\right)^{k+1}\bb E\abs{\Om(G)_{\bhv\bbw}}^{2p-1-k}\prec\Psi(z)^{k+1}\bb E\abs{\Om(G)_{\bhv\bbw}}^{2p-1-k}. \label{eq:cum_0k_bd}
\end{align}

Next suppose that $q=k-1$ and $r=1$. Then $l_1=\dots=l_{k-1}=1$, and
\begin{align}
        \eqref{eq:cum_hi_term} & \prec N^{-3/2}(1+\phi)^{2k-4} \left(\sqrt{\frac{\vphi+\eta}{N\eta}}\right)^{k-2}\sum_{i,\mu}\bb E\left( \abs{\Sg_{\bhv i}^{1/2} \partial_{i\mu}G_{\mu\bbw} \partial_{i\mu}\Om(G)_{\bhv\bbw}} \abs{\Om(G)_{\bhv\bbw}}^{2p-k} \right). \nonumber
\end{align}
The resolvent derivative is
\begin{align}
        \partial_{i\mu}G_{\mu\bbw} & =-(G\wh\Sg^{1/2})_{\mu i}G_{\mu\bbw} -G_{\mu\mu}(\wh\Sg^{1/2}G)_{i\bbw}. \label{eq:res_der1}
\end{align}
By \eqref{eq:err_der1}, we treat $\partial_{i\mu}\Om(G)_{\bhv\bbw}$ by separating the two contributions. We first consider $\partial_{i\mu}(\calH G)_{\bhv\bbw}$ and then $\partial_{i\mu}(m_N(z)(\wt\Sg G)_{\bhv\bbw})$. The first term is treated by combining \eqref{eq:res_der1} with \eqref{eq:HG_der1}. For instance,
\begin{align}
         & \sum_{i,\mu}\abs{\Sg_{\bhv i}^{1/2}} \abs{(G\wh\Sg^{1/2})_{\mu i}} \abs{((I_z-\La)G\wh\Sg^{1/2})_{\bhv i}} \abs{G_{\mu\bbw}}^2 \nonumber                                                                  \\
         & \leq \max_i\abs{((I_z-\La)G\wh\Sg^{1/2})_{\bhv i}} \left(\sum_{i,\mu}\abs{\Sg_{\bhv i}^{1/2}}^2 \abs{(G\wh\Sg^{1/2})_{\mu i}}^2\right)^{1/2} \left(\sum_{i,\mu}\abs{G_{\mu\bbw}}^4\right)^{1/2} \nonumber \\
         & \prec(1+\phi)^2\sqrt N \frac{\vphi+\eta}{\eta}. \nonumber
\end{align}
For the mixed term containing $G_{\mu\mu}(\wh\Sg^{1/2}G)_{i\bbw}$ satisfies
\begin{align}
         & \sum_{i,\mu}\abs{\Sg_{\bhv i}^{1/2}} \abs{G_{\mu\mu}} \abs{(\wh\Sg^{1/2}G)_{i\bbw}} \abs{((I_z-\La)G\wh\Sg^{1/2})_{\bhv i}} \abs{G_{\mu\bbw}} \nonumber               \\
         & \leq \left(\sum_i\abs{\Sg_{\bhv i}^{1/2}}^2\right)^{1/2} \left(\sum_i\abs{(\wh\Sg^{1/2}G)_{i\bbw}}^2 \abs{((I_z-\La)G\wh\Sg^{1/2})_{\bhv i}}^2\right)^{1/2} \max_\mu\abs{G_{\mu\mu}} \sum_\mu\abs{G_{\mu\bbw}} \nonumber                                                                                                         \\
         & \prec(1+\phi)^2\sqrt N \frac{\vphi+\eta}{\eta}. \nonumber
\end{align}
The other terms associated with $\partial_{i\mu}(\calH G)_{\bhv\bbw}$ can be treated similarly. 
For the self-energy term, the contribution from $\partial_{i\mu}(\wt\Sg G)$  in \eqref{eq:mG_der1} is treated with $I_z-\La$ replaced by $\wt\Sg$ in the terms associated with $\partial_{i\mu}(\calH G)$. If the derivative falls on $m_N$, \eqref{eq:tr_der_bd} gives the required additional factor. Therefore,
\begin{align}
         & \sum_{i,\mu}\bb E\left( \abs{\Sg_{\bhv i}^{1/2} \partial_{i\mu}G_{\mu\bbw} \partial_{i\mu}\Om(G)_{\bhv\bbw}} \abs{\Om(G)_{\bhv\bbw}}^{2p-k} \right)  \prec(1+\phi)^3\sqrt N \frac{\vphi+\eta}{\eta} \bb E\abs{\Om(G)_{\bhv\bbw}}^{2p-k}. \nonumber
\end{align}
Substituting this bound into \eqref{eq:cum_hi_term} gives
\begin{align}
        \eqref{eq:cum_hi_term}\prec(1+\phi)^{2k}\left(\sqrt{\frac{\vphi+\eta}{N\eta}}\right)^k\bb E\abs{\Om(G)_{\bhv\bbw}}^{2p-k}\prec\Psi(z)^k\bb E\abs{\Om(G)_{\bhv\bbw}}^{2p-k}. \label{eq:cum_1km1_bd}
\end{align}

Finally, suppose that $q=k-1$ and $r=0$. Then exactly one of $l_1,\dots,l_{k-1}$ equals $2$, while the remaining ones equal $1$. Without loss of generality, let $l_1=2$. By \eqref{eq:err_map_der}, we have
\begin{align}
        \eqref{eq:cum_hi_term} & \prec N^{-3/2}(1+\phi)^{2k-4} \left( \sqrt{\frac{\vphi+\eta}{N\eta}} \right)^{k-2} \sum_{i,\mu}\bb E\left( \abs{\Sg_{\bhv i}^{1/2}G_{\mu\bbw} \partial_{i\mu}^2\Om(G)_{\bhv\bbw}} \abs{\Om(G)_{\bhv\bbw}}^{2p-k} \right). \label{eq:cum_0km1}
\end{align}
By \eqref{eq:anis_err_id},
\begin{align}
        \partial_{i\mu}^2\Om(G)_{\bhv\bbw} & = \partial_{i\mu}^2(\calH G)_{\bhv\bbw} +\partial_{i\mu}^2 (m_N(z)(\wt\Sg G)_{\bhv\bbw}). \nonumber
\end{align}
Since $\calH G=I_L+(I_z-\La)G$, the exact second derivative is
\begin{align}
        \partial_{i\mu}^2(\calH G) =2(I_z-\La)GV^{i\mu}GV^{i\mu}G. \nonumber
\end{align}
This yields that 
\begin{align}
          &\sum_{i,\mu}\bb E\left( \abs{\Sg_{\bhv i}^{1/2}G_{\mu\bbw} \partial_{i\mu}^2(\calH G)_{\bhv\bbw}} \abs{\Om(G)_{\bhv\bbw}}^{2p-k} \right) \nonumber \\
          &\quad\prec \sum_{i,\mu}\bb E\left( \abs{\Sg_{\bhv i}^{1/2}G_{\mu\bbw} (GV^{i\mu}GV^{i\mu}G)_{\bhv\bbw}} \abs{\Om(G)_{\bhv\bbw}}^{2p-k} \right) \nonumber\\
         &\quad\prec (1+\phi)^2\sqrt{N}\frac{\vphi+\eta}{\eta}  \bb E\abs{\Om(G)_{\bhv\bbw}}^{2p-k}. \label{eq:HG_der2_bd}
\end{align}
For the self-energy term, the Leibniz rule gives
\begin{align}
          \partial_{i\mu}^2 (m_N(z)(\wt\Sg G)_{\bhv\bbw}) &=(\partial_{i\mu}^2m_N(z))(\wt\Sg G)_{\bhv\bbw} +2(\partial_{i\mu}m_N(z)) \partial_{i\mu}(\wt\Sg G)_{\bhv\bbw}\nonumber\\
          &+m_N(z)\partial_{i\mu}^2 (\wt\Sg G)_{\bhv\bbw}. \label{eq:self_der2}
\end{align}
By \eqref{eq:tr_der_bd}, we have
\begin{align}
         \sum_{i,\mu}\bb E\left( \abs{\Sg_{\bhv i}^{1/2}G_{\mu\bbw} (\partial_{i\mu}^2m_N(z))(\wt\Sg G)_{\bhv\bbw}} \abs{\Om(G)_{\bhv\bbw}}^{2p-k} \right) \prec (1+\phi)^3\frac{\vphi+\eta}{\eta}\bb E\abs{\Om(G)_{\bhv\bbw}}^{2p-k}, \nonumber
\end{align}
and
\begin{align}
                 \sum_{i,\mu}\bb E\left( \abs{\Sg_{\bhv i}^{1/2}G_{\mu\bbw} (\partial_{i\mu}m_N(z)) \partial_{i\mu}(\wt\Sg G)_{\bhv\bbw}} \abs{\Om(G)_{\bhv\bbw}}^{2p-k} \right)  \prec (1+\phi)^3\frac{\vphi+\eta}{\eta}\bb E\abs{\Om(G)_{\bhv\bbw}}^{2p-k}. \nonumber
\end{align}
In the last term of \eqref{eq:self_der2}, similar as in the proof of \eqref{eq:HG_der2_bd}, with $I_z-\La$ replaced by $\wt\Sg$, we have
\begin{align}
             \sum_{i,\mu}\bb E\left( \abs{\Sg_{\bhv i}^{1/2}G_{\mu\bbw} m_N(z)\partial_{i\mu}^2 (\wt\Sg G)_{\bhv\bbw}} \abs{\Om(G)_{\bhv\bbw}}^{2p-k} \right) \prec (1+\phi)^3\sqrt{N}\frac{\vphi+\eta}{\eta}  \bb E\abs{\Om(G)_{\bhv\bbw}}^{2p-k}
\end{align}
Thus,
\begin{align}
         & \quad \sum_{i,\mu}\bb E\left( \abs{\Sg_{\bhv i}^{1/2}G_{\mu\bbw} \partial_{i\mu}^2 (m_N(z)(\wt\Sg G)_{\bhv\bbw})} \abs{\Om(G)_{\bhv\bbw}}^{2p-k} \right)\nonumber \\
         & \prec (1+\phi)^3 \sqrt N \frac{\vphi+\eta}{\eta} \bb E\abs{\Om(G)_{\bhv\bbw}}^{2p-k} = (1+\phi)^3 N^{3/2} \frac{\vphi+\eta}{N\eta} \bb E\abs{\Om(G)_{\bhv\bbw}}^{2p-k}. \label{eq:self_der2_bd}
\end{align}
Combining \eqref{eq:cum_0km1}, \eqref{eq:HG_der2_bd}, and \eqref{eq:self_der2_bd}, we obtain
\begin{align}
        \eqref{eq:cum_hi_term}\prec(1+\phi)^{2k}\left(\sqrt{\frac{\vphi+\eta}{N\eta}}\right)^k\bb E\abs{\Om(G)_{\bhv\bbw}}^{2p-k}\prec\Psi(z)^k\bb E\abs{\Om(G)_{\bhv\bbw}}^{2p-k}. \label{eq:cum_0km1_bd}
\end{align}

By \eqref{eq:cum_hi_exp}, \eqref{eq:cum_hi_term}, \eqref{eq:cum_gen_bd}, \eqref{eq:cum_0k_bd}, \eqref{eq:cum_1km1_bd} and \eqref{eq:cum_0km1_bd}, for every $k\geq2$, we obtain
\begin{align}
        T_k^{(1)} \prec \sum_{n=1}^{2p}\Psi(z)^n \bb E\abs{\Om(G)_{\bhv\bbw}}^{2p-n}. \label{eq:high_cum_mom}
\end{align}

We now estimate the derivative terms $T_{121}^{(1)}$ and $T_{122}^{(1)}$ from \eqref{eq:cum2_der}. Here it is important to perform the sums over $i$ and $\mu$ before taking absolute values.
Using \eqref{eq:HG_der1}, we have
\begin{align}
         & \sum_{i,\mu}\Sg_{\bhv i}^{1/2}G_{\mu\bbw} \partial_{i\mu}(\calH G)_{\bhv\bbw} =- (\wt\Sg G)_{\bhv\bbw} ((I_z-\La)G\wt I G)_{\bhv\bbw}-(G\wt I G)_{\bbw\bbw} ((I_z-\La)G\wt\Sg)_{\bhv\bhv}. \label{eq:cum2_contr}
\end{align}
Lemma~\ref{lem:ward_cons} and Cauchy--Schwarz inequality imply
\begin{align}
        \abs{((I_z-\La)G\wt I G)_{\bhv\bbw}} +\abs{(G\wt I G)_{\bbw\bbw}} \prec\frac{\vphi+\eta}{\eta}. \nonumber
\end{align}
Consequently, the absolute value of \eqref{eq:cum2_contr} is bounded by $(1+\phi)(\vphi+\eta)/\eta$. For the self-energy term, direct differentiation gives
\begin{align}
         & \sum_{i,\mu}\Sg_{\bhv i}^{1/2}G_{\mu\bbw} m_N(z)\partial_{i\mu}(\wt\Sg G)_{\bhv\bbw}  =-m_N(z)(\wt\Sg G)_{\bhv\bbw} (\wt\Sg G\wt I G)_{\bhv\bbw} -m_N(z)(G\wt I G)_{\bbw\bbw} (\wt\Sg G\wt\Sg)_{\bhv\bhv}, \nonumber
\end{align}
and hence this expression is $O_\prec((1+\phi)^2(\vphi+\eta)/\eta)$. For the derivative of the normalized trace, summing first in $i$ and $\mu$ gives
\begin{align}
        \sum_{i,\mu}\Sg_{\bhv i}^{1/2}G_{\mu\bbw} \partial_{i\mu}m_N(z) =-\frac2N(\wt\Sg G\wt I G\wt I G)_{\bhv\bbw}. \nonumber
\end{align}
By the resolvent norm estimate \eqref{eq:G_triv_bd} and the consequence \eqref{eq:gh_vec_bd} of the Ward identities, we have
\begin{align}
        \frac2N\abs{(\wt\Sg G\wt I G\wt I G)_{\bhv\bbw}}\leq\frac2N\norm{G}\norm{\wt I G\wt\Sg\bhv}\norm{\wt I G\bbw}\prec\frac{\vphi+\eta}{N\eta^2}\leq\frac{\vphi+\eta}{\eta}, \nonumber
\end{align}
where the last inequality follows from $N\eta\geq N^\tau$.
Combining the last three identities, we obtain
\begin{align}
         & \abs{\sum_{i,\mu}\Sg_{\bhv i}^{1/2}G_{\mu\bbw} \partial_{i\mu}\Om(G)_{\bhv\bbw}} \prec (1+\phi)^2\frac{\vphi+\eta}{\eta}. \nonumber
\end{align}
Therefore,
\begin{align}
        T_{121}^{(1)}\prec(1+\phi)^2\frac{\vphi+\eta}{N\eta}\bb E\abs{\Om(G)_{\bhv\bbw}}^{2p-2}\prec\Psi(z)^2\bb E\abs{\Om(G)_{\bhv\bbw}}^{2p-2}. \nonumber
\end{align}
For $T_{122}^{(1)}$, the entries of $X$ are real and hence $\partial_{i\mu}\overline{\Om(G)_{\bhv\bbw}} =\overline{\partial_{i\mu}\Om(G)_{\bhv\bbw}}$. For example, performing the $i$- and $\mu$-sums before taking absolute values gives
\begin{align}
         & \sum_{i,\mu}\Sg_{\bhv i}^{1/2}G_{\mu\bbw} \overline{\partial_{i\mu}(\calH G)_{\bhv\bbw}}  =-\overline{((I_z-\La)G\wt\Sg)_{\bhv\bhv}} \sum_\mu\abs{G_{\mu\bbw}}^2 -\overline{(\wt\Sg G)_{\bhv\bbw}} \sum_\mu G_{\mu\bbw} \overline{((I_z-\La)G)_{\bhv\mu}}. \nonumber
\end{align}
The $G^*$-version of Lemma~\ref{lem:ward_cons} and \eqref{eq:ward_2_proj} bound this expression by $O_\prec((1+\phi)(\vphi+\eta)/\eta)$. The terms in which the derivative falls on $\wt\Sg G$ satisfy the same estimate in \eqref{eq:cum2_contr} with $I_z-\La$ replaced by $\wt\Sg$. If the derivative falls on $m_N(z)$, \eqref{eq:tr_der_bd} and the anisotropic entry bound give $O_\prec((1+\phi)^2(\vphi+\eta)/\eta)$. Hence
\begin{align}
        T_{12}^{(1)} & \prec \Psi(z)^2 \bb E\abs{\Om(G)_{\bhv\bbw}}^{2p-2}. \label{eq:cum2_der_bd}
\end{align}

We finally bound the remainder term $\calR_{l+1}^{(1)}$. For $s\in\R$, set
\begin{align}
        F_{i\mu}(s) & := G^{(i\mu)}(s)_{\mu\bbw} \Om(G^{(i\mu)}(s))_{\bhv\bbw}^{p-1} \overline{ \Om(G^{(i\mu)}(s))_{\bhv\bbw} }^{p}. \nonumber
\end{align}
The cumulant remainder satisfies
\begin{align}
        \abs{\calR_{l+1}^{(1),i\mu}} & \leq O(1)\left( \bb E\sup_{\abs s\leq\abs{X_{i\mu}}} \abs{\frac{\dd^{l+1}}{\dd s^{l+1}} F_{i\mu}(s)}^2 \bb E\left( \abs{X_{i\mu}}^{2l+4} 1_{\{\abs{X_{i\mu}}>K\}} \right)\right)^{1/2}\nonumber \\
                                     & \quad+ O(1)\bb E\abs{X_{i\mu}}^{l+2} \bb E\sup_{\abs s\leq K} \abs{\frac{\dd^{l+1}}{\dd s^{l+1}} F_{i\mu}(s)}. \label{eq:high_mom_rem}
\end{align}
By the Leibniz rule, the derivative in \eqref{eq:high_mom_rem} is a finite sum of products of derivatives of $G^{(i\mu)}(s)$, $\Om(G^{(i\mu)}(s))_{\bhv\bbw}$, and its complex conjugate. The identities \eqref{eq:HG_der} and \eqref{eq:self_der} remain valid along the interpolation. Together with \eqref{eq:interp_res_der}, for fixed $p$,we have
\begin{align}
        \sup_{\abs s\leq K} \abs{\frac{\dd^{l+1}}{\dd s^{l+1}} F_{i\mu}(s)} \prec N^{C(p)}(1+\phi)^{l+2p+1}, \label{eq:rem_small_der}
\end{align}
where $C(p)$ is independent of $l$. Since $\phi\leq N^{\tau/10}$, \eqref{eq:entry_moment_bd} and \eqref{eq:rem_small_der} show that the second term in \eqref{eq:high_mom_rem}, corresponding to $|s|\leq K$, is $O(N^{-D})$ after choosing $l=l(p,D)$ sufficiently large.

For the first term in \eqref{eq:high_mom_rem}, which contains the truncated moment from \eqref{eq:trunc_mom_bd}, the resolvent norm estimate \eqref{eq:interp_res} gives
\begin{align}
        \bb E\sup_{\abs s\leq\abs{X_{i\mu}}} \abs{\frac{\dd^{l+1}}{\dd s^{l+1}} F_{i\mu}^{(p)}(s)}^2 \leq N^{C_{p,l}}. \nonumber
\end{align}
After $p$ and $l$ are fixed, choose $D_0$ in \eqref{eq:trunc_mom_bd} sufficiently large. Summing over $i,\mu$ and using $\sum_i\abs{\Sg_{\bhv i}^{1/2}}\lesssim\sqrt N$, we obtain
\begin{align}
        \abs{\calR_{l+1}^{(1)}} =O(N^{-D}) \label{eq:high_mom_rem_bd}
\end{align}
for every fixed $D>0$.

Combining \eqref{eq:anis_mom_exp}, \eqref{eq:cum2_dec}, \eqref{eq:cum2_cancel}, \eqref{eq:high_cum_mom}, \eqref{eq:cum2_der_bd} and \eqref{eq:high_mom_rem_bd}, we conclude that
\begin{align}
        \bb E\abs{\Om(G)_{\bhv\bbw}}^{2p} & = \sum_{k=1}^l T_k^{(1)} +\calR^{(1)}_{l+1} + \bb E m_N(z) (\wt\Sg G)_{\bhv\bbw} \Om(G)_{\bhv\bbw}^{p-1} \overline{\Om(G)_{\bhv\bbw}}^{p} \nonumber\\
        &\prec \sum_{n=1}^{2p} \Psi(z)^n \bb E|\Om(G)_{\bhv\bbw}|^{2p-n}. \nonumber
\end{align}
By H\"older's inequality,
\begin{align}
        \bb E\abs{\Om(G)_{\bhv\bbw}}^{2p}\prec\sum_{n=1}^{2p}\Psi(z)^n\bb E\abs{\Om(G)_{\bhv\bbw}}^{2p-n}\leq\sum_{n=1}^{2p}\Psi(z)^n(\bb E\abs{\Om(G)_{\bhv\bbw}}^{2p})^{(2p-n)/(2p)}. \nonumber
\end{align}
For $1\leq n\leq2p-1$, Young's inequality gives, for every $\delta>0$,
\begin{align}
        \Psi(z)^n(\bb E\abs{\Om(G)_{\bhv\bbw}}^{2p})^{(2p-n)/(2p)} & \leq\frac{n}{2p}\delta^{-2p/n}\Psi(z)^{2p} +\frac{2p-n}{2p}\delta^{2p/(2p-n)} \bb E\abs{\Om(G)_{\bhv\bbw}}^{2p}. \nonumber
\end{align}
The term corresponding to $n=2p$ is simply $\Psi(z)^{2p}$. Therefore,
\begin{align}
        \bb E\abs{\Om(G)_{\bhv\bbw}}^{2p} & \prec C_\delta\Psi(z)^{2p} +\left(\sum_{n=1}^{2p-1}\frac{2p-n}{2p} \delta^{2p/(2p-n)}\right) \bb E\abs{\Om(G)_{\bhv\bbw}}^{2p} \nonumber\\
        &\leq C_\delta\Psi(z)^{2p} +\frac12\bb E\abs{\Om(G)_{\bhv\bbw}}^{2p}, \label{eq:anis_holder}
\end{align}
where $\delta>0$ is chosen sufficiently small so that the coefficient of $\bb E|\Om(G)_{\bhv\bbw}|^{2p}$ is at most $1/2$, and $C_\delta$ depends only on $p$ and $\delta$.
Absorbing the last term into the left-hand side gives
\begin{align}
        \bb E|\Om(G)_{\bhv\bbw}|^{2p} & \prec \Psi(z)^{2p}. \nonumber
\end{align}
Since $p$ is arbitrary, Markov's inequality yields
\begin{align}
        |\Om(G)_{\bhv\bbw}| & \prec\Psi(z). \nonumber
\end{align}
This proves the anisotropic estimate when the first test vector is supported on the first block. Applying the same cumulant expansion to \eqref{eq:anis_err_id_2}, with $G_{\mu\bbw}\Sg^{1/2}_{\bhv i}$ replaced by $(\btv)_\mu(\Sg^{1/2}G)_{i\bbw}$ and $m_N\wt\Sg$ by $g_N\wt I$, gives the estimate for the second block. Since $\bbv=\bhv+\btv$, the general case follows from the triangle inequality.

\subsection{Averaged estimate}\label{sub:avg}

This subsection proves the averaged estimate \eqref{eq:avg_rand_error} in Theorem~\ref{thm:rand_sc_eq}. Throughout this subsection we use Lemma~\ref{lem:ward_ids} and Lemma~\ref{lem:entry_res_der}.

We first split $B$ according to its two block columns. Since
\begin{align}
        B=B\wh I+B\wt I, \nonumber
\end{align}
it is enough to consider $B=B\wh I$ and $B=B\wt I$ separately. We first assume that
\begin{align}
        B=B\wh I =\begin{pmatrix}B_{11} & 0 \\
               B_{21} & 0\end{pmatrix}, \qquad \|B\|\leq1. \label{eq:avg_first_blk}
\end{align}
Thus  by \eqref{eq:intro_mde}  and  cyclicity of the trace, we have
\begin{align}
        \fram:=\avg{B\Om(G)}_{L}=\avg{B\calH G}_{L}+m_N\avg{B\wt\Sg G}_{L}=\avg{GB\calH}_{L}+m_N\avg{\wt\Sg GB}_{L}. \label{eq:avg_err}
\end{align}
For $B=B\wt I$, the exact starting identity is
\begin{align}
        \avg{B\Om(G)}_L =\avg{GB\calH}_L+g_N\avg{\wt I GB}_L, \qquad g_N:=\avg{\wt\Sg G}_N. \nonumber
\end{align}
We now focus on $B=B\wh I$. The term $B=B\wt I$ is treated similarly, so we postpone its discussion and return to this case at the end of this subsection.

We first collect several estimates obtained from the Ward identities in Lemma~\ref{lem:ward_ids}.
\begin{lemma}
        Under the assumptions of Theorem~\ref{thm:rand_sc_eq}, let $B_1,B_2,B_3\in\C^{L\times L}$ be deterministic matrices of bounded norm. Uniformly for $z$ satisfying \eqref{eq:rand_err_dom},
        \begin{align}
                \frac1N\abs{\avg{B_1G\wt IGB_2}_{L}}                                                      & \prec \frac{\vphi+\eta}{N\eta}, \label{eq:avg_tr_bd}                                                \\
                \frac1{N^2}\sum_{i,\mu} \abs{(B_1G\wh\Sg^{1/2})_{\mu i}}                                  & \prec \sqrt{\frac{\vphi+\eta}{N\eta}}, \label{eq:sum1_bd}                                                \\
                \frac1{NL^2}\sum_{i,\mu} \abs{(B_1G\wh\Sg^{1/2})_{\mu i}} \abs{(\wh\Sg^{1/2}GB_2)_{i\mu}} & \prec \frac1N\frac{\vphi+\eta}{N\eta} \leq\left(\frac{\vphi+\eta}{N\eta}\right)^2. \label{eq:sum2_bd}
        \end{align}
        Moreover, if $\bbv,\bbw\in \C^L$ are deterministic unit vectors, then
        \begin{align}
                \frac1{N^2}\sum_{i,\mu} \abs{(B_1G)_{\bbv i}} \abs{(B_2G)_{\mu i}} \abs{(GB_3)_{\mu\bbw}} \prec \left(\frac{\vphi+\eta}{N\eta}\right)^{3/2}. \label{eq:prod3_bd}
        \end{align}
\end{lemma}
\begin{proof}
        By Cauchy--Schwarz inequality, we have
        \begin{align}
                \abs{\avg{B_1G\wt IGB_2}_{L}} \leq\|B_1G\wt I\|_{\rr{hs}}\|\wt IGB_2\|_{\rr{hs}}. \nonumber
        \end{align}
        By Cauchy--Schwarz inequality and Lemma~\ref{lem:ward_cons},
        \begin{align}
                \|B_1G\wt I\|_{\rr{hs}}^2 +\|\wt IGB_2\|_{\rr{hs}}^2 \prec\frac{\vphi+\eta}{\eta}, \nonumber
        \end{align}
        this proves \eqref{eq:avg_tr_bd}.

        For \eqref{eq:sum1_bd}, Cauchy--Schwarz inequality  and Lemma~\ref{lem:ward_cons} give
        \begin{align}
                \sum_{i,\mu}\abs{(B_1G\wh\Sg^{1/2})_{\mu i}}\leq N\left(\sum_{i,\mu}\abs{(B_1G\wh\Sg^{1/2})_{\mu i}}^2\right)^{1/2}\prec N^{3/2}\sqrt{\frac{\vphi+\eta}{\eta}}, \nonumber
        \end{align}
        which proves \eqref{eq:sum1_bd} after division by $N^2$. Applying Cauchy--Schwarz inequality to the sum over $(i,\mu)$ and using Lemma~\ref{lem:ward_cons} for both sums, we obtain
        \begin{align}
                \sum_{i,\mu}\abs{(B_1G\wh\Sg^{1/2})_{\mu i}}\abs{(\wh\Sg^{1/2}GB_2)_{i\mu}}\prec N\frac{\vphi+\eta}{\eta}, \nonumber
        \end{align}
        which proves \eqref{eq:sum2_bd}.

        Finally, Cauchy--Schwarz inequality in the pair $(i,\mu)$ yields
        \begin{align}
                 & \sum_{i,\mu} \abs{(B_1G)_{\bbv i}} \abs{(B_2G)_{\mu i}} \abs{(GB_3)_{\mu\bbw}} \nonumber                                                                                             \\
                 & \leq \left(\sum_{i,\mu} \abs{(B_1G)_{\bbv i}}^2 \abs{(B_2G)_{\mu i}}^2\right)^{1/2} \left(\sum_{i,\mu} \abs{(GB_3)_{\mu\bbw}}^2\right)^{1/2} \nonumber                               \\
                 & \leq \left(\max_i\sum_\mu \abs{(B_2G)_{\mu i}}^2\right)^{1/2} \left(\sum_i\abs{(B_1G)_{\bbv i}}^2\right)^{1/2} \sqrt N \left(\sum_\mu\abs{(GB_3)_{\mu\bbw}}^2\right)^{1/2} \nonumber \\
                 & \prec \sqrt N\left(\frac{\vphi+\eta}{\eta}\right)^{3/2}, \nonumber
        \end{align}
        where the last line follows from Lemma~\ref{lem:ward_cons}. Dividing by $N^2$ proves \eqref{eq:prod3_bd}.
\end{proof}

We next establish the estimates for derivatives with respect to the entry $X_{i\mu}$ that will be used in the cumulant expansion. For every fixed integer $l\geq1$, similar arguments as in the proof of \eqref{eq:tr_der_bd} give
\begin{align}
        \abs{\partial_{i\mu}^l m_N} +\abs{\partial_{i\mu}^l\avg{\wt\Sg GB}_{N}} \prec(1+\phi)^{l-1}\frac{\vphi+\eta}{N\eta}. \label{eq:avg_tr_der}
\end{align}

The derivative of $\fram$ contains the following additional matrix term. Define
\begin{align}
        Q :=\frac1N\wh\Sg^{1/2} ( GB\calH G+m_NGB\wt\Sg G+\avg{\wt\Sg GB}_{N}G\wt I G )\wh\Sg^{1/2}. \label{eq:Q_def}
\end{align}
Differentiating \eqref{eq:avg_err}, we obtain
\begin{align}
        \partial_{i\mu}\fram & =\avg{GBV^{i\mu}}_{L} -\avg{V^{i\mu}GB\calH G}_{L} -m_N\avg{V^{i\mu}GB\wt\Sg G}_{L} -\avg{V^{i\mu}G\wt I G}_{L}\avg{\wt\Sg GB}_{N} \nonumber                         \\
                             & =\avg{GBV^{i\mu}}_{L} -\frac NL(Q_{i\mu}+Q_{\mu i}). \label{eq:avg_err_der}
\end{align}
We need the following bounds on $Q$ and its derivatives, which are proved after the proof of the averaged estimate.

\begin{lemma}\label{lem:avg_Q_bd}
        Under the assumptions of Theorem~\ref{thm:rand_sc_eq}, for every fixed integer $r\geq1$, uniformly for deterministic vectors $\bbv,\bbw\in\C^L$,
        \begin{align}
                \abs{\partial_{i\mu}^rQ_{\bbv\bbw}} \prec (1+\phi)^{r+1}\frac{\vphi+\eta}{N\eta}\|\bbv\|\|\bbw\|. \label{eq:Q_der_bd}
        \end{align}
        Moreover,
        \begin{align}
                \abs{Q_{\bbv\bbw}} \prec (1+\phi)^4\left(\frac{\vphi+\eta}{N\eta}\right)^{3/2}\|\bbv\|\|\bbw\|. \label{eq:Q_anis_bd}
        \end{align}
\end{lemma}

\subsubsection{High moments of the averaged error}

We now turn to the high moments of $\fram$. Fix an integer $p\geq 1$. By \eqref{eq:avg_err},
\begin{align}
        \bb E\abs{\fram}^{2p} & =\frac1L\sum_{i,j,\mu}\Sg^{1/2}_{ji} \bb E\left(X_{i\mu}(GB)_{\mu j}\fram^{p-1}\overline{\fram}^{p}\right) +\bb E\left(m_N\avg{\wt\Sg GB}_{L}\fram^{p-1}\overline{\fram}^{p}\right). \nonumber
\end{align}
Applying Lemma~\ref{lem:cum_exp} to the first term gives
\begin{align}
        \bb E\abs{\fram}^{2p}  & =\sum_{k=1}^{l}\frac1L\sum_{i,j,\mu} \frac{\kappa_{k+1}(X_{i\mu})}{k!}\Sg^{1/2}_{ji} \bb E\partial_{i\mu}^k\paren{(GB)_{\mu j}\fram^{p-1}\overline{\fram}^{p}} +\bb E\left(m_N\avg{\wt\Sg GB}_{L}\fram^{p-1}\overline{\fram}^{p}\right) +\calR_{l+1}^{(2)}, \nonumber \\
        & :=\sum_{k=1}^{l}T_k^{(2)} +\bb E\left(m_N\avg{\wt\Sg GB}_{L}\fram^{p-1}\overline{\fram}^{p}\right) +\calR_{l+1}^{(2)}. \label{eq:avg_cum} 
\end{align}
Since $X_{i\mu}$ is real, $\partial_{i\mu}^l\overline{\fram} =\overline{\partial_{i\mu}^l\fram}$,
so the corresponding absolute value estimates are the same as for derivatives of $\fram$.

We first calculate the second cumulant term. By the Leibniz rule,
\begin{align}
        T_1^{(2)} & =\frac1{NL}\sum_{i,j,\mu}\Sg^{1/2}_{ji} \bb E\partial_{i\mu}(GB)_{\mu j}\fram^{p-1}\overline{\fram}^{p}+\frac1{NL}\sum_{i,j,\mu}\Sg^{1/2}_{ji} \bb E(GB)_{\mu j}\partial_{i\mu}\paren{\fram^{p-1}\overline{\fram}^{p}}\nonumber \\
                  & :=T_{11}^{(2)}+T_{12}^{(2)}. \nonumber
\end{align}
Since $\partial_{i\mu}G=-GV^{i\mu}G$,
\begin{align}
        \partial_{i\mu}(GB)_{\mu j}=-\paren{GV^{i\mu}GB}_{\mu j}=-\paren{G\wh\Sg^{1/2}}_{\mu i}\paren{GB}_{\mu j}-G_{\mu\mu}\paren{\wh\Sg^{1/2}GB}_{ij}. \nonumber
\end{align}
Consequently,
\begin{align}
        T_{11}^{(2)} & =-\frac1{NL}\sum_{i,j,\mu}\Sg^{1/2}_{ji} \bb E\paren{G\wh\Sg^{1/2}}_{\mu i}(GB)_{\mu j} \fram^{p-1}\overline{\fram}^{p}-\frac1{NL}\sum_{i,j,\mu}\Sg^{1/2}_{ji} \bb E G_{\mu\mu}\paren{\wh\Sg^{1/2}GB}_{ij} \fram^{p-1}\overline{\fram}^{p}\nonumber \\
                     & =-\frac1N\bb E\avg{B\wt\Sg G\wt I G}_{L} \fram^{p-1}\overline{\fram}^{p}-\bb E\left(m_N\avg{\wt\Sg GB}_{L} \fram^{p-1}\overline{\fram}^{p}\right). \nonumber
\end{align}
The second term cancels the explicit self-energy term in \eqref{eq:avg_cum}. By \eqref{eq:avg_tr_bd}, the remaining $N^{-1}$ term satisfies
\begin{align}
        \frac1N\bb E\abs{\avg{B\wt\Sg G\wt I G}_{L}}\abs{\fram}^{2p-1}\prec\frac{\vphi+\eta}{N\eta}\bb E\abs{\fram}^{2p-1}\leq\Psi(z)^2\bb E\abs{\fram}^{2p-1}. \label{eq:avg_rem}
\end{align}

We next treat $T_{12}^{(2)}$. The product rule gives
\begin{align}
        \partial_{i\mu}\paren{\fram^{p-1}\overline{\fram}^{p}} & =(p-1)(\partial_{i\mu}\fram)\fram^{p-2}\overline{\fram}^{p}+p(\partial_{i\mu}\overline{\fram})\fram^{p-1}\overline{\fram}^{p-1}. \nonumber
\end{align}
The term containing $\partial_{i\mu}\overline{\fram}$ is the complex conjugate of the calculation below and has the same absolute value bound. For the term containing $\partial_{i\mu}\fram$, summing over $j$ first yields
\begin{align}
         & \frac{p-1}{NL}\sum_{i,j,\mu}\Sg^{1/2}_{ji} \bb E(GB)_{\mu j}(\partial_{i\mu}\fram) \fram^{p-2}\overline{\fram}^{p} = \frac{p-1}{NL}\sum_{i,\mu} \bb E(GB\wh\Sg^{1/2})_{\mu i}(\partial_{i\mu}\fram) \fram^{p-2}\overline{\fram}^{p}. \nonumber
\end{align}
Substituting \eqref{eq:avg_err_der}, we split this expression according to the explicit rank-two trace and the term involving $Q$. The rank-two contribution is bounded by
\begin{align}
         & \frac1{NL}\sum_{i,\mu} \bb E\abs{(GB\wh\Sg^{1/2})_{\mu i}} \abs{\avg{GBV^{i\mu}}_{L}}\abs{\fram}^{2p-2}\nonumber                                                                                      \\
         & \quad\leq\frac1{NL^2}\sum_{i,\mu} \bb E\abs{(GB\wh\Sg^{1/2})_{\mu i}} \left( \abs{(\wh\Sg^{1/2}GB\wh\Sg^{1/2})_{i\mu}} +\abs{(\wh\Sg^{1/2}GB\wh\Sg^{1/2})_{\mu i}} \right)\abs{\fram}^{2p-2}\nonumber \\
         & \quad\prec\left(\frac{\vphi+\eta}{N\eta}\right)^2 \bb E\abs{\fram}^{2p-2} \leq\Psi(z)^4\bb E\abs{\fram}^{2p-2}, \nonumber
\end{align}
where we used \eqref{eq:sum2_bd}. For the $Q$-contribution, Lemma~\ref{lem:avg_Q_bd} and \eqref{eq:sum1_bd} give
\begin{align}
         & \frac1{L^2}\sum_{i,\mu} \bb E\abs{(GB\wh\Sg^{1/2})_{\mu i}} \abs{Q_{i\mu}+Q_{\mu i}}\abs{\fram}^{2p-2}\nonumber                             \\
         & \quad\prec\sqrt{\frac{\vphi+\eta}{N\eta}} (1+\phi)^4\left(\frac{\vphi+\eta}{N\eta}\right)^{3/2} \bb E\abs{\fram}^{2p-2}\prec(1+\phi)^4\left(\frac{\vphi+\eta}{N\eta}\right)^2 \bb E\abs{\fram}^{2p-2}. \nonumber
\end{align}
Complex conjugation gives the same bound for the term containing $\partial_{i\mu}\overline{\fram}$. Hence
\begin{align}
        \abs{T_{12}^{(2)}} \prec(1+\phi)^4\left(\frac{\vphi+\eta}{N\eta}\right)^2 \bb E\abs{\fram}^{2p-2}. \label{eq:avg_T12_bd}
\end{align}

We next record the derivative estimates needed for the higher cumulants by \eqref{eq:avg_err}. 
By \eqref{eq:avg_tr_der} and Leibniz rule, we have
\begin{align}
        \abs{\partial_{i\mu}^l \left(m_N\avg{B\wt\Sg G}_{L}\right)} =\abs{\sum_{r=0}^l \binom lr (\partial_{i\mu}^{r}m_N) \partial_{i\mu}^{l-r}\avg{B\wt\Sg G}_{L}}\prec(1+\phi)^l\frac{\vphi+\eta}{N\eta},\qquad l\geq1. \label{eq:self_err_der}
\end{align}
Combining $ \calH G=I_L+(I_z-\La)G$, \eqref{eq:avg_tr_der} and \eqref{eq:self_err_der}, we conclude that
\begin{align}
        \abs{\partial_{i\mu}^l\fram} \prec(1+\phi)^l\frac{\vphi+\eta}{N\eta}, \qquad l\geq1. \label{eq:avg_err_der_bd}
\end{align}  

We now consider $T_k^{(2)}$ for fixed $k\geq 2$. By \eqref{eq:avg_cum}, we have
\begin{align}
        T_k^{(2)} & =\frac1L\sum_{i,j,\mu} \frac{\kappa_{k+1}(X_{i\mu})}{k!}\Sg^{1/2}_{ji} \bb E\partial_{i\mu}^k\paren{(GB)_{\mu j}\fram^{p-1}\overline{\fram}^{p}} \nonumber                                                                                   \\
                  & =\frac1{L}\sum_{i,\mu} \frac{\kappa_{k+1}(X_{i\mu})}{k!} \bb E\partial_{i\mu}^k\paren{(GB\wh\Sg^{1/2})_{\mu i}\fram^{p-1}\overline{\fram}^{p}}. \nonumber
\end{align}
Since $\kappa_{k+1}(X_{i\mu})\lesssim N^{-(k+1)/2}$, the absolute value of $T_k^{(2)}$ is bounded by
\begin{align}
         & N^{-(k-1)/2}N^{-2}\sum_{i,\mu}\bb E\left( \abs{\partial_{i\mu}^r(GB\wh\Sg^{1/2})_{\mu i}} \prod_{a=1}^{q}\abs{\partial_{i\mu}^{l_a}\fram} \abs{\fram}^{2p-1-q} \right), \nonumber
\end{align}
where $0\leq r\leq k$, $0\leq q\leq (k-r)\land (2p-1)$, $l_a\geq1$, and $l_1+\dots+l_q=k-r$.
By \eqref{eq:avg_err_der_bd}, this is bounded by
\begin{align}
         & N^{-(k-1)/2}(1+\phi)^{k+1} \left(\frac{\vphi+\eta}{N\eta}\right)^{q} \bb E\abs{\fram}^{2p-1-q} \leq N^{-(k-3)/2}(1+\phi)^{k+1} \left(\frac{\vphi+\eta}{N\eta}\right)^{q+1} \bb E\abs{\fram}^{2p-1-q}, \nonumber
\end{align}
where we used $N^{-1/2}\leq\sqrt{(\vphi+\eta)/(N\eta)}$. Since $\phi\leq N^{\tau/10}$ and $\tau$ is sufficiently small,
\begin{align}
        N^{-(k-3)/2}(1+\phi)^{k+1} =(1+\phi)^4[N^{-1/2}(1+\phi)]^{k-3} \lesssim(1+\phi)^4, \qquad k\geq3. \nonumber
\end{align}
Thus
\begin{align}
        \abs{T_k^{(2)}} \prec(1+\phi)^4\sum_{n=1}^{2p} \left(\frac{\vphi+\eta}{N\eta}\right)^n \bb E\abs{\fram}^{2p-n}, \qquad k\geq3. \label{eq:avg_cum_hi}
\end{align}
So it remains to consider $T_2^{(2)}$.

By \eqref{eq:avg_cum}, we have
\begin{align}
        T_2^{(2)} & =\frac1{2L}\sum_{i,j,\mu} \kappa_3(X_{i\mu})\Sg^{1/2}_{ji} \bb E\partial_{i\mu}^2\paren{(GB)_{\mu j}\fram^{p-1}\overline{\fram}^{p}} \nonumber                                                                                   \\
                  & =\frac1{2L}\sum_{i,\mu} \kappa_3(X_{i\mu}) \bb E\partial_{i\mu}^2\paren{(GB\wh\Sg^{1/2})_{\mu i}\fram^{p-1}\overline{\fram}^{p}}. \nonumber
\end{align}
Up to constants depending only on $p$, its absolute value is bounded by
\begin{align}
        N^{-5/2}\sum_{i,\mu} \bb E\abs{\partial_{i\mu}^2\left( (GB\wh\Sg^{1/2})_{\mu i}\fram^{p-1}\overline{\fram}^{p} \right)}. \nonumber
\end{align}
The terms in which derivatives hit conjugated factors have the same absolute value bounds. It therefore suffices to consider the following four structures:
\begin{align}
         & N^{-5/2}\sum_{i,\mu} \bb E\abs{\partial_{i\mu}^2(GB\wh\Sg^{1/2})_{\mu i}} \abs{\fram}^{2p-1}, \label{eq:T2_1}                         \\
         & N^{-5/2}\sum_{i,\mu} \bb E\abs{\partial_{i\mu}(GB\wh\Sg^{1/2})_{\mu i}} \abs{\partial_{i\mu}\fram}\abs{\fram}^{2p-2}, \label{eq:T2_2} \\
         & N^{-5/2}\sum_{i,\mu} \bb E\abs{(GB\wh\Sg^{1/2})_{\mu i}} \abs{\partial_{i\mu}^2\fram}\abs{\fram}^{2p-2}, \label{eq:T2_3}            \\
         & N^{-5/2}\sum_{i,\mu} \bb E\abs{(GB\wh\Sg^{1/2})_{\mu i}} \abs{\partial_{i\mu}\fram}^2\abs{\fram}^{2p-3}. \label{eq:T2_4}
\end{align}
For \eqref{eq:T2_1}, the resolvent derivative formula gives
\begin{align}
        \partial_{i\mu}^2(GB\wh\Sg^{1/2})_{\mu i} =2(GV^{i\mu}GV^{i\mu}GB\wh\Sg^{1/2})_{\mu i}. \nonumber
\end{align}
Expanding the two rank-two factors yields four products. Using $\|\wh\Sg^{1/2}\|\lesssim1$ from Assumption~\ref{ass:basic}, the four products are bounded by expressions of the forms
\begin{align}
          (G\wh\Sg^{1/2})_{\mu i}^2 (GB\wh\Sg^{1/2})_{\mu i}&, \qquad (G\wh\Sg^{1/2})_{\mu i}G_{\mu\mu} (\wh\Sg^{1/2}GB\wh\Sg^{1/2})_{ii}, \nonumber \\
         G_{\mu\mu}(\wh\Sg^{1/2}G\wh\Sg^{1/2})_{ii} (GB\wh\Sg^{1/2})_{\mu i}&, \qquad G_{\mu\mu}(\wh\Sg^{1/2}G)_{i\mu} (\wh\Sg^{1/2}GB\wh\Sg^{1/2})_{ii}. \nonumber
\end{align}
For the first product, Lemma~\ref{lem:ward_cons} gives
\begin{align}
        \sum_{i,\mu}\abs{(G\wh\Sg^{1/2})_{\mu i}}^2\abs{(GB\wh\Sg^{1/2})_{\mu i}}\prec(1+\phi)\sum_{i,\mu}\abs{(G\wh\Sg^{1/2})_{\mu i}}^2\prec(1+\phi)N\frac{\vphi+\eta}{\eta}. \nonumber
\end{align}
In each of the other three products, the diagonal matrix element is $O_\prec(1+\phi)$. Applying Cauchy--Schwarz inequality separately to the $i$- and $\mu$-sums and then Lemma~\ref{lem:ward_cons} gives the bound obtained for the first product. Consequently,
\begin{align}
        \eqref{eq:T2_1}\prec N^{-1/2}(1+\phi)\frac{\vphi+\eta}{N\eta}\bb E\abs{\fram}^{2p-1}\leq(1+\phi)\frac{\vphi+\eta}{N\eta}\bb E\abs{\fram}^{2p-1}. \nonumber
\end{align}
For \eqref{eq:T2_2}, substitute \eqref{eq:avg_err_der}. The explicit rank-two trace gives, by \eqref{eq:sum2_bd},
\begin{align}
        N^{-5/2}\sum_{i,\mu} \bb E\abs{\partial_{i\mu}(GB\wh\Sg^{1/2})_{\mu i}} \abs{\avg{GBV^{i\mu}}_{L}}\abs{\fram}^{2p-2} \prec(1+\phi)^2\left(\frac{\vphi+\eta}{N\eta}\right)^2 \bb E\abs{\fram}^{2p-2}. \nonumber
\end{align}
For the contribution containing $Q_{i\mu}+Q_{\mu i}$, by Lemma~\ref{lem:avg_Q_bd}, we obtain
\begin{align}
         & N^{-5/2}\frac NL\sum_{i,\mu} \bb E\abs{\partial_{i\mu}(GB\wh\Sg^{1/2})_{\mu i}} \abs{Q_{i\mu}+Q_{\mu i}}\abs{\fram}^{2p-2}\nonumber             \\
         & \prec N^{-5/2}N^2(1+\phi)^2  (1+\phi)^4 \left(\frac{\vphi+\eta}{N\eta}\right)^{3/2} \bb E\abs{\fram}^{2p-2}\nonumber \\
         & \prec N^{-1/2}(1+\phi)^6 \left(\frac{\vphi+\eta}{N\eta}\right)^{3/2} \bb E\abs{\fram}^{2p-2} \leq(1+\phi)^6 \left(\frac{\vphi+\eta}{N\eta}\right)^2 \bb E\abs{\fram}^{2p-2}, \nonumber
\end{align}
where we used $N^{-1/2}\leq\sqrt{(\vphi+\eta)/(N\eta)}$ in the last step.
Hence
\begin{align}
        \eqref{eq:T2_2} \prec(1+\phi)^6\left(\frac{\vphi+\eta}{N\eta}\right)^2 \bb E\abs{\fram}^{2p-2}. \nonumber
\end{align}
For \eqref{eq:T2_3}, we use the direct derivative estimate \eqref{eq:avg_err_der_bd} with $r=2$ and \eqref{eq:sum1_bd}. Namely,
\begin{align}
         & N^{-5/2}\sum_{i,\mu} \bb E\abs{(GB\wh\Sg^{1/2})_{\mu i}} \abs{\partial_{i\mu}^2\fram}\abs{\fram}^{2p-2} \nonumber                    \\
         & \qquad\prec N^{-5/2}(1+\phi)^2 \frac{\vphi+\eta}{N\eta} \sum_{i,\mu}\bb E\abs{(GB\wh\Sg^{1/2})_{\mu i}} \abs{\fram}^{2p-2} \nonumber \\
         & \qquad\prec N^{-1/2}(1+\phi)^2 \left(\frac{\vphi+\eta}{N\eta}\right)^{3/2} \bb E\abs{\fram}^{2p-2} \leq(1+\phi)^2 \left(\frac{\vphi+\eta}{N\eta}\right)^2 \bb E\abs{\fram}^{2p-2}, \nonumber
\end{align}
where the last step uses $N^{-1/2}\leq\sqrt{(\vphi+\eta)/(N\eta)}$. Thus
\begin{align}
        \eqref{eq:T2_3} \prec(1+\phi)^2\left(\frac{\vphi+\eta}{N\eta}\right)^2 \bb E\abs{\fram}^{2p-2}. \nonumber
\end{align}
Finally, \eqref{eq:avg_err_der_bd} and \eqref{eq:sum1_bd} yield
\begin{align}
        \eqref{eq:T2_4} \prec(1+\phi)^2\left(\frac{\vphi+\eta}{N\eta}\right)^3 \bb E\abs{\fram}^{2p-3}. \nonumber
\end{align}
Consequently,
\begin{align}
        \abs{T_2^{(2)}} & \prec(1+\phi)\frac{\vphi+\eta}{N\eta} \bb E\abs{\fram}^{2p-1}+(1+\phi)^6\left(\frac{\vphi+\eta}{N\eta}\right)^2 \bb E\abs{\fram}^{2p-2}\nonumber                      \\
                        & \quad+(1+\phi)^2\left(\frac{\vphi+\eta}{N\eta}\right)^3 \bb E\abs{\fram}^{2p-3}. \label{eq:T2_bd}
\end{align}

It remains to control the cumulant remainder. For $s\in\R$, let
\begin{align}
        \fram^{(i\mu)}(s) :=\avg{B\Om(G^{(i\mu)}(s))}_{L} \nonumber
\end{align}
and define
\begin{align}
        F_{i\mu}(s) :=\paren{G^{(i\mu)}(s)B\wh\Sg^{1/2}}_{\mu i} \fram^{(i\mu)}(s)^{p-1} \overline{\fram^{(i\mu)}(s)}^{p}. \nonumber
\end{align}
The interpolation bounds \eqref{eq:interp_res_der}, \eqref{eq:avg_err_der}, and \eqref{eq:avg_tr_der} give polynomial bounds for every fixed derivative of $F_{i\mu}(s)$. The same remainder argument as in the proof of \eqref{eq:high_mom_rem_bd} therefore yields
\begin{align}
        \calR_{l+1}^{(2)}=O(N^{-D}) \label{eq:avg_cum_rem}
\end{align}
for every fixed $D>0$.

Combining \eqref{eq:avg_rem}, \eqref{eq:avg_T12_bd}, \eqref{eq:T2_bd}, \eqref{eq:avg_cum_hi} and \eqref{eq:avg_cum_rem}, by H\"older's inequality followed by Young's inequality yields
\begin{align}
        \bb E\abs{\fram}^{2p} \prec(1+\phi)^{12p} \left(\frac{\vphi+\eta}{N\eta}\right)^{2p}+N^{-D}. \nonumber
\end{align}
The details already appeared in the proof of \eqref{eq:anis_holder}, so we omit them here.
Since $p$ is arbitrary, Markov's inequality implies
\begin{align}
        \abs{\avg{B\Om(G)}_{L}} =\abs{\fram} \prec\Psi(z)^2 =(1+\phi)^6\frac{\vphi+\eta}{N\eta}. \nonumber
\end{align}
The same argument for $B=B\wt I$, with $m_N,\wt\Sg$ replaced by $g_N,\wt I$, respectively, gives the same estimate. Since $B=B\wh I+B\wt I$, the averaged estimate follows for every deterministic matrix $B$.

\subsubsection{Bounds for \texorpdfstring{$Q$}{Q}}

\begin{proof}[Proof of Lemma~\ref{lem:avg_Q_bd}]
        By homogeneity, fix deterministic unit vectors $\bbv,\bbw$. We first prove the derivative estimate \eqref{eq:Q_der_bd}, which is used in the subsequent expansion of high moments for $Q$.

        We now differentiate the three terms in \eqref{eq:Q_def}. So it causes that
        \begin{align}
                \partial_{i\mu}^r Q_{\bbv\bbw} & =\frac1N\partial_{i\mu}^r (\wh\Sg^{1/2}GB\calH G\wh\Sg^{1/2})_{\bbv\bbw}-\frac1N\partial_{i\mu}^r \left(m_N(\wh\Sg^{1/2}GB\wt\Sg G\wh\Sg^{1/2})_{\bbv\bbw}\right)\nonumber \\
                                                & \quad-\frac1N\partial_{i\mu}^r \left(\avg{\wt\Sg GB}_{N} (\wh\Sg^{1/2}G\wt I G\wh\Sg^{1/2})_{\bbv\bbw}\right). \nonumber
        \end{align}       
        For the first one, Leibniz rule gives
        \begin{align}
                 \frac1N\partial_{i\mu}^r (\wh\Sg^{1/2}GB\calH G\wh\Sg^{1/2})_{\bbv\bbw}
                 & =\frac1N\sum_{r_1+r_2=r}\binom r{r_1} (\wh\Sg^{1/2}(\partial_{i\mu}^{r_1}G)B\calH (\partial_{i\mu}^{r_2}G)\wh\Sg^{1/2})_{\bbv\bbw}\nonumber                                                           \\
                 & \quad+\frac rN\sum_{r_1+r_2=r-1}\frac{(r-1)!}{r_1!r_2!} (\wh\Sg^{1/2}(\partial_{i\mu}^{r_1}G)BV^{i\mu} (\partial_{i\mu}^{r_2}G)\wh\Sg^{1/2})_{\bbv\bbw}. \nonumber
        \end{align}
       Expand each derivative by \eqref{eq:res_der_form}. Cauchy--Schwarz inequality and \eqref{eq:gh_vec_bd} control the two resolvent factors adjacent to $\bbv$ and $\bbw$, while \eqref{eq:rand_op_norm} controls $\calH$. After expanding the rank-two factors, the remaining entries with fixed indices are $O_\prec(1+\phi)$. Hence
        \begin{align}
                \frac1N\abs{\partial_{i\mu}^r (\wh\Sg^{1/2}GB\calH G\wh\Sg^{1/2})_{\bbv\bbw}} \prec(1+\phi)^{r}\frac{\vphi+\eta}{N\eta}, \nonumber
        \end{align}
        where we used $\norm{\cal H}\prec 1$ by \eqref{eq:rand_op_norm}.
        For the term containing $m_N$, we write
        \begin{align}
                 & \frac1N\partial_{i\mu}^r \left(m_N(\wh\Sg^{1/2}GB\wt\Sg G\wh\Sg^{1/2})_{\bbv\bbw}\right)\nonumber                                                                                                                         \\
                 & =\frac1N\sum_{r_0+r_1+r_2=r}\frac{r!}{r_0!r_1!r_2!} (\partial_{i\mu}^{r_0} m_N) (\wh\Sg^{1/2}(\partial_{i\mu}^{r_1}G)B\wt\Sg (\partial_{i\mu}^{r_2}G)\wh\Sg^{1/2})_{\bbv\bbw}. \label{eq:Q_mN_der}
        \end{align}
        If $r_0=0$, then $|m_N|\prec1+\phi$. In each term, Cauchy--Schwarz inequality and \eqref{eq:gh_vec_bd} control the two resolvent factors adjacent to $\bbv$ and $\bbw$, while the $r$ matrix elements introduced by the derivatives are bounded by $O_\prec(1+\phi)$. This gives
        \begin{align}
               \eqref{eq:Q_mN_der} \prec  (1+\phi)^{r+1}\frac{\vphi+\eta}{N\eta}. \nonumber
        \end{align}
        If $r_0\geq1$,  by \eqref{eq:avg_tr_der}, we also have
        \begin{align}
                \eqref{eq:Q_mN_der} \prec N^{-1} (1+\phi)^{r+1}\frac{\vphi+\eta}{N\eta}. \nonumber
        \end{align}
        Similar as the case of $m_N$, for the term containing $\avg{\wt\Sg GB}_{N}$, we have
        \begin{align}
                  \abs{\frac1N\partial_{i\mu}^r \left(\avg{\wt\Sg GB}_{N} (\wh\Sg^{1/2}G\wt I G\wh\Sg^{1/2})_{\bbv\bbw}\right)} \prec (1+\phi)^{r+1}\frac{\vphi+\eta}{N\eta}. \nonumber
        \end{align}
        Combining the three parts, for every fixed $r\geq1$,
        \begin{align}
                \abs{\partial_{i\mu}^rQ_{\bbv\bbw}} \prec(1+\phi)^{r+1}\frac{\vphi+\eta}{N\eta}. \label{eq:Q_der_est}
        \end{align}
        This proves \eqref{eq:Q_der_bd}.

        We now prove \eqref{eq:Q_anis_bd}. Since $B=B\wh I$, by \eqref{eq:Q_def}, we have
        \begin{align}
                Q_{\bbv\bbw} & =\frac1N\sum_{i,\mu} (\wh\Sg^{1/2}GB\wh\Sg^{1/2})_{\bbv i} X_{i\mu} (G\wh\Sg^{1/2})_{\mu\bbw}\nonumber\\
                &\quad+\frac1N\left(\wh\Sg^{1/2} (m_NGB\wt\Sg G+\avg{\wt\Sg GB}_{N}G\wt I G) \wh\Sg^{1/2}\right)_{\bbv\bbw}. \label{eq:Q_rand_exp}
        \end{align}
        Fix an integer $p\geq 1$. Applying Lemma~\ref{lem:cum_exp} to \eqref{eq:Q_rand_exp}, we obtain
        \begin{align}
                \bb E\abs{Q_{\bbv\bbw}}^{2p} & =\frac1N\sum_{k=1}^l\sum_{i,\mu}\frac{\kappa_{k+1}(X_{i\mu})}{k!} \bb E\partial_{i\mu}^k\left( (\wh\Sg^{1/2}GB\wh\Sg^{1/2})_{\bbv i} (G\wh\Sg^{1/2})_{\mu\bbw} Q_{\bbv\bbw}^{p-1} \overline{Q_{\bbv\bbw}}^{p}\right)\nonumber                                                                                                       \\
                                               & \quad+\frac1N\bb E\left( \left(\wh\Sg^{1/2} (m_NGB\wt\Sg G+\avg{\wt\Sg GB}_{N}G\wt I G) \wh\Sg^{1/2}\right)_{\bbv\bbw} Q_{\bbv\bbw}^{p-1} \overline{Q_{\bbv\bbw}}^{p} \right) +\calR_{l+1}^{(3)}\nonumber\\
                                               &:=\sum_{k=1}^l T_k^{(3)}+\frac1N\bb E\left( \left(\wh\Sg^{1/2} (m_NGB\wt\Sg G+\avg{\wt\Sg GB}_{N}G\wt I G) \wh\Sg^{1/2}\right)_{\bbv\bbw} Q_{\bbv\bbw}^{p-1} \overline{Q_{\bbv\bbw}}^{p} \right) \nonumber\\
                                               &\quad+\calR_{l+1}^{(3)}.\label{eq:Q_cum_exp}
        \end{align}
        For $k=1$,  we have
        \begin{align}
                T_1^{(3)} & =\frac1{N^2}\sum_{i,\mu} \bb E\partial_{i\mu}\left( (\wh\Sg^{1/2}GB\wh\Sg^{1/2})_{\bbv i} (G\wh\Sg^{1/2})_{\mu\bbw} Q_{\bbv\bbw}^{p-1} \overline{Q_{\bbv\bbw}}^{p}\right)\nonumber \\
                          & =\frac1{N^2}\sum_{i,\mu} \bb E\partial_{i\mu}\left( (\wh\Sg^{1/2}GB\wh\Sg^{1/2})_{\bbv i} (G\wh\Sg^{1/2})_{\mu\bbw}\right) Q_{\bbv\bbw}^{p-1} \overline{Q_{\bbv\bbw}}^{p}\nonumber \\
                          &\quad + \frac1{N^2}\sum_{i,\mu} \bb E (\wh\Sg^{1/2}GB\wh\Sg^{1/2})_{\bbv i} (G\wh\Sg^{1/2})_{\mu\bbw} \partial_{i\mu}\left(Q_{\bbv\bbw}^{p-1} \overline{Q_{\bbv\bbw}}^{p}\right)\nonumber\\ 
                          &:=T_{11}^{(3)}+T_{12}^{(3)}.\label{eq:Q_cum1}
        \end{align}
         From $\partial_{i\mu}G=-GV^{i\mu}G$,
        \begin{align}
                  \partial_{i\mu}(\wh\Sg^{1/2}GB\wh\Sg^{1/2})_{\bbv i}&=-(\wh\Sg^{1/2}G\wh\Sg^{1/2})_{\bbv i} (GB\wh\Sg^{1/2})_{\mu i} -(\wh\Sg^{1/2}G)_{\bbv\mu} (\wh\Sg^{1/2}GB\wh\Sg^{1/2})_{ii}, \label{eq:Q_ext1} \\
                  \partial_{i\mu}(G\wh\Sg^{1/2})_{\mu\bbw}&=-(G\wh\Sg^{1/2})_{\mu i} (G\wh\Sg^{1/2})_{\mu\bbw} -G_{\mu\mu}(\wh\Sg^{1/2}G\wh\Sg^{1/2})_{i\bbw}. \label{eq:Q_ext2}
        \end{align}
        Consequently, by $\sum_\mu G_{\mu\mu}=Nm_N$ and $\sum_i(\wh\Sg^{1/2}GB\wh\Sg^{1/2})_{ii}=N\avg{\wt\Sg GB}_{N}$, we have
        \begin{align}
                 & \frac1{N^2}\sum_{i,\mu}\partial_{i\mu}\left( (\wh\Sg^{1/2}GB\wh\Sg^{1/2})_{\bbv i} (G\wh\Sg^{1/2})_{\mu\bbw}\right)+\frac1N\left(\wh\Sg^{1/2} (m_NGB\wt\Sg G+\avg{\wt\Sg GB}_{N}G\wt I G) \wh\Sg^{1/2}\right)_{\bbv\bbw}\nonumber                                                                                         \\
                 & =-\frac1{N^2}(\wh\Sg^{1/2}G\wt\Sg B^\top G\wt I G\wh\Sg^{1/2})_{\bbv\bbw} -\frac1{N^2}(\wh\Sg^{1/2}GB\wt\Sg G\wt I G\wh\Sg^{1/2})_{\bbv\bbw}. \label{eq:Q_cum2_cancel}
        \end{align}
        Moreover, using $G^\top=G$, the first term can be expanded as
        \begin{align}
                 & (\wh\Sg^{1/2}G\wt\Sg B^\top G\wt I G\wh\Sg^{1/2})_{\bbv\bbw}=\sum_{i,\mu} (\wh\Sg^{1/2}G\wt\Sg^{1/2})_{\bbv i} (GB\wt\Sg^{1/2})_{\mu i} (G\wh\Sg^{1/2})_{\mu\bbw}. \nonumber
        \end{align}
        Hence, by Cauchy--Schwarz inequality and Lemma~\ref{lem:ward_cons}, we have
        \begin{align}
                 & \frac1{N^2}\abs{(\wh\Sg^{1/2}G\wt\Sg B^\top G\wt I G\wh\Sg^{1/2})_{\bbv\bbw}} \nonumber                                                                                                                                                            \\
                 & \leq\frac1{N^2} \left(\sum_i\abs{(\wh\Sg^{1/2}G\wt\Sg^{1/2})_{\bbv i}}^2\right)^{1/2} \left(\sum_{i,\mu}\abs{(GB\wt\Sg^{1/2})_{\mu i}}^2\right)^{1/2} \left(\sum_\mu\abs{(G\wh\Sg^{1/2})_{\mu\bbw}}^2\right)^{1/2} \nonumber \\
                 & \prec\frac1{N^2} \left(\frac{\vphi+\eta}{\eta}\right)^{1/2} \left(N\frac{\vphi+\eta}{\eta}\right)^{1/2} \left(\frac{\vphi+\eta}{\eta}\right)^{1/2} =\left(\frac{\vphi+\eta}{N\eta}\right)^{3/2}. \nonumber
        \end{align}
        For the second crossed term, expanding $\wt\Sg$ and $\wt I$ gives
        \begin{align}
                 & (\wh\Sg^{1/2}GB\wt\Sg G\wt I G\wh\Sg^{1/2})_{\bbv\bbw} =\sum_{i,\mu} (\wh\Sg^{1/2}GB\wt\Sg^{1/2})_{\bbv i} (G\wt\Sg^{1/2})_{\mu i} (G\wh\Sg^{1/2})_{\mu\bbw}, \nonumber
        \end{align}
        Applying Cauchy--Schwarz inequality and Lemma~\ref{lem:ward_cons} again gives the same bound as above. Hence, we have
        \begin{align}
                 & \frac1{N^2}\abs{(\wh\Sg^{1/2}G\wt\Sg B^\top G\wt I G\wh\Sg^{1/2})_{\bbv\bbw}} +\frac1{N^2}\abs{(\wh\Sg^{1/2}GB\wt\Sg G\wt I G\wh\Sg^{1/2})_{\bbv\bbw}}\prec\left(\frac{\vphi+\eta}{N\eta}\right)^{3/2}. \label{eq:Q_cum2_cross}
        \end{align}
        Combining \eqref{eq:Q_cum_exp}, \eqref{eq:Q_cum1}, \eqref{eq:Q_cum2_cancel} and \eqref{eq:Q_cum2_cross}, we obtain
        \begin{align}
                \bb E\abs{Q_{\bbv\bbw}}^{2p} = T_{12}^{(3)} + \sum_{k=2}^l T_k^{(3)} +\calR_{l+1}^{(3)} + O_\prec\left(\left(\frac{\vphi+\eta}{N\eta}\right)^{3/2}\bb E\abs{Q_{\bbv\bbw}}^{2p-1}\right) \label{eq:Q_cum2}.
        \end{align}

        For $k\geq2$, $T_k^{(3)}$ is bounded by an expression of the form
        \begin{align}
                 & N^{-(k-1)/2}N^{-2}\sum_{i,\mu}\bb E\left( \abs{\partial_{i\mu}^r\left( (\wh\Sg^{1/2}GB\wh\Sg^{1/2})_{\bbv i} (G\wh\Sg^{1/2})_{\mu\bbw}\right)} \prod_{a=1}^q\abs{\partial_{i\mu}^{l_a}Q_{\bbv\bbw}} \abs{Q_{\bbv\bbw}}^{2p-1-q}\right), \label{eq:Q_cum_term}
        \end{align}
        where $0\leq r\leq k$, $0\leq q\leq(k-r)\wedge(2p-1)$, $l_a\geq1$ and $l_1+\dots+l_q=k-r$.
        A differentiated factor may originate from either $Q_{\bbv\bbw}$ or its complex conjugate. This does not change the absolute value estimates.

        For any fixed $r\geq0$, \eqref{eq:res_der_form} and Lemma~\ref{lem:ward_cons} give
        \begin{align}
                 & \frac1{N^2}\sum_{i,\mu} \abs{\partial_{i\mu}^r\left( (\wh\Sg^{1/2}GB\wh\Sg^{1/2})_{\bbv i} (G\wh\Sg^{1/2})_{\mu\bbw}\right)}\prec(1+\phi)^r\frac{\vphi+\eta}{N\eta}. \label{eq:Q_ext_contr}
        \end{align}
        Combining \eqref{eq:Q_ext_contr} with \eqref{eq:Q_der_est}, the term in \eqref{eq:Q_cum_term} is bounded by
        \begin{align}
                 & N^{-(k-1)/2}(1+\phi)^{r+\sum_{a=1}^q(l_a+1)} \left(\frac{\vphi+\eta}{N\eta}\right)^{q+1} \bb E\abs{Q_{\bbv\bbw}}^{2p-1-q}\nonumber \\
                 & \qquad=N^{-(k-1)/2}(1+\phi)^{k+q} \left(\frac{\vphi+\eta}{N\eta}\right)^{q+1} \bb E\abs{Q_{\bbv\bbw}}^{2p-1-q}. \nonumber
        \end{align}
        Suppose that $q\leq k-2$. Since $N^{-1}\leq(\vphi+\eta)/(N\eta)$, we have
        \begin{align}
                N^{-(k-1)/2}\left(\frac{\vphi+\eta}{N\eta}\right)^{q+1}\leq N^{-(k-q-2)/2}\left(\frac{\vphi+\eta}{N\eta}\right)^{3(q+1)/2}. \nonumber
        \end{align}
        If $k\leq3q+4$, then $k+q\leq4(q+1)$. If $k>3q+4$, the remaining factor $N^{-(k-q-2)/2}$ absorbs $(1+\phi)^{k-3q-4}$ for sufficiently small $\tau$, since $k$ is bounded by the chosen cumulant expansion order and $1+\phi\lesssim N^{\tau/10}$. Hence
        \begin{align}
                \eqref{eq:Q_cum_term} \prec(1+\phi)^{4(q+1)}\left(\frac{\vphi+\eta}{N\eta}\right)^{3(q+1)/2} \bb E\abs{Q_{\bbv\bbw}}^{2p-1-q},\qquad q\leq k-2. \label{eq:Q_nonsat}
        \end{align}
        Thus we just need to consider the cases $q=k-1$ and $q=k$. 
        It remains to consider
        \begin{align}
                (r,q)=(0,k),\qquad(r,q)=(1,k-1),\qquad(r,q)=(0,k-1). \nonumber
        \end{align}

        \noindent\textit{Case 1: $(r,q)=(0,k)$.}
        Then $l_1=\dots=l_k=1$. Using \eqref{eq:Q_der_est} for $k-1$ derivative factors, \eqref{eq:Q_cum_term} is bounded by
        \begin{align}
                 & N^{-(k-1)/2}(1+\phi)^{2(k-1)} \left(\frac{\vphi+\eta}{N\eta}\right)^{k-1}\nonumber                                                                                                                                                                           \\
                 & \quad\times\frac1{N^2}\sum_{i,\mu}\bb E\left( \abs{(\wh\Sg^{1/2}GB\wh\Sg^{1/2})_{\bbv i}} \abs{(G\wh\Sg^{1/2})_{\mu\bbw}} \abs{\partial_{i\mu}Q_{\bbv\bbw}} \abs{Q_{\bbv\bbw}}^{2p-1-k}\right). \label{eq:Q_0k}
        \end{align}
        By \eqref{eq:Q_def}, we also have
        \begin{align}
                \partial_{i\mu}Q_{\bbv\bbw} & =\frac1N\partial_{i\mu}(\wh\Sg^{1/2}GB\calH G\wh\Sg^{1/2})_{\bbv\bbw}\nonumber\\
                &+\frac1N\partial_{i\mu}\left(\wh\Sg^{1/2} (m_NGB\wt\Sg G+\avg{\wt\Sg GB}_{N}G\wt I G) \wh\Sg^{1/2}\right)_{\bbv\bbw}. \label{eq:Q_der1}
        \end{align}
        We remark that this is difference form \eqref{eq:Q_rand_exp} because we not expand the $X_{i\mu}$ factor there. 
        For the first term,
        \begin{align}
                  \partial_{i\mu}(\wh\Sg^{1/2}GB\calH G\wh\Sg^{1/2}) &=-\wh\Sg^{1/2}GV^{i\mu}GB\calH G\wh\Sg^{1/2} \nonumber\\
                  &\quad+\wh\Sg^{1/2}GBV^{i\mu}G\wh\Sg^{1/2} -\wh\Sg^{1/2}GB\calH GV^{i\mu}G\wh\Sg^{1/2}. \label{eq:Q_GBHG_der}
        \end{align}
        Consider the middle term in \eqref{eq:Q_GBHG_der}. Expanding $V^{i\mu}$, one of its two summands gives
        \begin{align}
                 & \frac1{N^3}\sum_{i,\mu} \abs{(\wh\Sg^{1/2}GB\wh\Sg^{1/2})_{\bbv i}}^2 \abs{(G\wh\Sg^{1/2})_{\mu\bbw}} \abs{(\wh\Sg^{1/2}G\wh\Sg^{1/2})_{\mu\bbw}}\nonumber                                                              \\
                 & \quad\leq\frac1{N^3}\sum_i \abs{(\wh\Sg^{1/2}GB\wh\Sg^{1/2})_{\bbv i}}^2 \left(\sum_\mu\abs{(G\wh\Sg^{1/2})_{\mu\bbw}}^2\right)^{1/2} \left(\sum_\mu\abs{(\wh\Sg^{1/2}G\wh\Sg^{1/2})_{\mu\bbw}}^2\right)^{1/2}\nonumber \\
                 & \quad\prec\frac1{N^3}\left(\frac{\vphi+\eta}{\eta}\right)^2 \leq\left(\frac{\vphi+\eta}{N\eta}\right)^3. \label{eq:Q_0k_contr}
        \end{align}
        For the first and third terms in \eqref{eq:Q_GBHG_der}, we first use $\calH G=I_L+(I_z-\La)G$. After expanding $V^{i\mu}$, the three factors carrying summation indices are controlled by \eqref{eq:Q_0k_contr}. Any remaining entry with fixed indices is bounded by $O_\prec(1+\phi)$. Hence
        \begin{align}
                 & \frac1{N^2}\sum_{i,\mu} \abs{(\wh\Sg^{1/2}GB\wh\Sg^{1/2})_{\bbv i}} \abs{(G\wh\Sg^{1/2})_{\mu\bbw}} \frac1N\abs{\partial_{i\mu}(\wh\Sg^{1/2}GB\calH G\wh\Sg^{1/2})_{\bbv\bbw}}\prec(1+\phi)\left(\frac{\vphi+\eta}{N\eta}\right)^3. \nonumber
        \end{align}
        We next differentiate the two self-energy terms in \eqref{eq:Q_der1}.   
        If the derivative falls on $m_N$ or $\avg{\wt\Sg GB}_{N}$, we use \eqref{eq:avg_tr_der}. The two remaining resolvent factors are controlled by Cauchy--Schwarz inequality and \eqref{eq:gh_vec_bd}. If the derivative falls on a resolvent, the rank-two expansion produces the three factors in \eqref{eq:Q_0k_contr}. Hence the term with coefficient $m_N$ is bounded by $(1+\phi)^2((\vphi+\eta)/(N\eta))^3$. The last term is treated in the same way, with $m_N$ and $\wt\Sg$ replaced by $\avg{\wt\Sg GB}_{N}$ and $\wt I$. Therefore,
        \begin{align}
                 & \frac1{N^2}\sum_{i,\mu} \abs{(\wh\Sg^{1/2}GB\wh\Sg^{1/2})_{\bbv i}} \abs{(G\wh\Sg^{1/2})_{\mu\bbw}} \abs{\partial_{i\mu}Q_{\bbv\bbw}}\prec(1+\phi)^2\left(\frac{\vphi+\eta}{N\eta}\right)^3. \label{eq:Q_der_contr}
        \end{align}
        Substituting this into \eqref{eq:Q_0k} and using $N^{-1/2}\leq\sqrt{(\vphi+\eta)/(N\eta)}$, we obtain
        \begin{align}
                \eqref{eq:Q_cum_term} & \prec(1+\phi)^{2k}\left(\frac{\vphi+\eta}{N\eta}\right)^{3(k+1)/2} \bb E\abs{Q_{\bbv\bbw}}^{2p-1-k}. \label{eq:Q_0k_bd}
        \end{align}

        \noindent\textit{Case 2: $(r,q)=(1,k-1)$.}
        Here $l_1=\dots=l_{k-1}=1$. After estimating $k-2$ first derivatives of $Q_{\bbv\bbw}$ by \eqref{eq:Q_der_est}, it remains to control
        \begin{align}
                 & \frac1{N^2}\sum_{i,\mu} \abs{\partial_{i\mu}\left( (\wh\Sg^{1/2}GB\wh\Sg^{1/2})_{\bbv i} (G\wh\Sg^{1/2})_{\mu\bbw}\right)} \abs{\partial_{i\mu}Q_{\bbv\bbw}}. \label{eq:Q_1km1}
        \end{align}
        Using \eqref{eq:Q_ext1} and \eqref{eq:Q_ext2}, we can treat the term of $\partial_{i\mu}\left( (\wh\Sg^{1/2}GB\wh\Sg^{1/2})_{\bbv i} (G\wh\Sg^{1/2})_{\mu\bbw}\right)$. Expanding $\partial_{i\mu}Q_{\bbv\bbw}$ by \eqref{eq:Q_der1}, one of the resulting products is
        \begin{align}
                 & \frac1{N^3}\sum_{i,\mu} \abs{(\wh\Sg^{1/2}G\wh\Sg^{1/2})_{\bbv i}} \abs{(GB\wh\Sg^{1/2})_{\mu i}} \abs{(G\wh\Sg^{1/2})_{\mu\bbw}} \abs{\partial_{i\mu}(\wh\Sg^{1/2}GB\calH G\wh\Sg^{1/2})_{\bbv\bbw}}\nonumber \\
                 & \qquad\prec\frac{(1+\phi)^3}{N}\left(\frac{\vphi+\eta}{N\eta}\right)^{3/2} \leq(1+\phi)^3\left(\frac{\vphi+\eta}{N\eta}\right)^{5/2}, \label{eq:Q_1km1_term}
        \end{align}
        where we used  \eqref{eq:prod3_bd} in the first inequality and $N^{-1}\leq (\vphi+\eta)/\eta$. The other rank-one choices are estimated identically after placing the $i$- and $\mu$-dependent factors using the corresponding Ward estimates for the two blocks. For the self-energy terms, derivatives of the normalized traces are controlled by \eqref{eq:avg_tr_der},  while resolvent derivatives are handled by the same rank-two expansion as above. And these term will occur an extra factor of $(1+\phi)$ by the extra $G$ factor. Hence,
        \begin{align}
                \eqref{eq:Q_1km1} \prec(1+\phi)^4\left(\frac{\vphi+\eta}{N\eta}\right)^{5/2}. \nonumber
        \end{align}
        Consequently,
        \begin{align}
                \eqref{eq:Q_cum_term} & \prec N^{-(k-1)/2}(1+\phi)^{2(k-2)+4} \left(\frac{\vphi+\eta}{N\eta}\right)^{k-2+5/2} \bb E\abs{Q_{\bbv\bbw}}^{2p-k}\nonumber \\
                                                            & \prec(1+\phi)^{4k}\left(\frac{\vphi+\eta}{N\eta}\right)^{3k/2} \bb E\abs{Q_{\bbv\bbw}}^{2p-k}. \label{eq:Q_1km1_bd}
        \end{align}

        \noindent\textit{Case 3: $(r,q)=(0,k-1)$.}
        In this case $l_1+\dots+l_{k-1}=k$, hence exactly one of the $l_a$ equals two and the others equal one. After estimating the $k-2$ first derivatives by \eqref{eq:Q_der_est}, it remains to control
        \begin{align}
                 & \frac1{N^2}\sum_{i,\mu} \abs{(\wh\Sg^{1/2}GB\wh\Sg^{1/2})_{\bbv i}} \abs{(G\wh\Sg^{1/2})_{\mu\bbw}} \abs{\partial_{i\mu}^2Q_{\bbv\bbw}}. \label{eq:Q_0k_bdm1_start}
        \end{align}
        For the first term in \eqref{eq:Q_def}, differentiate \eqref{eq:Q_GBHG_der} once more. For example,
        \begin{align}
                  \partial_{i\mu}(\wh\Sg^{1/2}GV^{i\mu}GB\calH G\wh\Sg^{1/2}) &=-\wh\Sg^{1/2}GV^{i\mu}GV^{i\mu}GB\calH G\wh\Sg^{1/2} -\wh\Sg^{1/2}GV^{i\mu}GV^{i\mu}GB\calH G\wh\Sg^{1/2}\nonumber                                             \\
                 & \quad+\wh\Sg^{1/2}GV^{i\mu}GBV^{i\mu}G\wh\Sg^{1/2} -\wh\Sg^{1/2}GV^{i\mu}GB\calH GV^{i\mu}G\wh\Sg^{1/2}, \nonumber
        \end{align}
        where repeated terms only change a combinatorial coefficient depending only on the derivative order. Every term contains two rank-two factors. After using $\calH G=I_L+(I_z-\La)G$ whenever an $\calH$ remains, one of the resulting monomials is bounded by
        \begin{align}
                 & (1+\phi)^3 \abs{(\wh\Sg^{1/2}GB\wh\Sg^{1/2})_{\bbv i}} \abs{(G\wh\Sg^{1/2})_{\mu i}} \abs{(G\wh\Sg^{1/2})_{\mu\bbw}}, \nonumber
        \end{align}
        where the factor $(1+\phi)^3$ accounts for the  remaining resolvent entries between the fixed indices $i$ and $\mu$. Hence \eqref{eq:prod3_bd} gives
        \begin{align}
                 & \frac1{N^3}\sum_{i,\mu} (1+\phi)^3 \abs{(\wh\Sg^{1/2}GB\wh\Sg^{1/2})_{\bbv i}} \abs{(G\wh\Sg^{1/2})_{\mu i}} \abs{(G\wh\Sg^{1/2})_{\mu\bbw}} \nonumber \\
                 & \qquad\prec\frac{(1+\phi)^2}{N} \left(\frac{\vphi+\eta}{N\eta}\right)^{3/2} \leq(1+\phi)^2 \left(\frac{\vphi+\eta}{N\eta}\right)^{5/2}. \nonumber
        \end{align}
        For other terms in $\partial_{i\mu}^2(\wh\Sg^{1/2}GB\calH G\wh\Sg^{1/2})_{\bbv\bbw}$, we can treat those by the same ways. Specifically, expanding the two factors $V^{i\mu}$ gives finitely many monomials. In each monomial, the summation indices in the first and second blocks are estimated using the corresponding Ward identities. Cauchy--Schwarz inequality and \eqref{eq:prod3_bd} then control the normalized double sum. The remaining entries with fixed indices contribute at most $(1+\phi)^2$, while normalized traces are controlled by \eqref{eq:avg_tr_der}. Consequently,
        \begin{align}
                 & \frac1{N^3}\sum_{i,\mu} \abs{(\wh\Sg^{1/2}GB\wh\Sg^{1/2})_{\bbv i}} \abs{(G\wh\Sg^{1/2})_{\mu\bbw}} \abs{\partial_{i\mu}^2(\wh\Sg^{1/2}GB\calH G\wh\Sg^{1/2})_{\bbv\bbw}}\prec(1+\phi)^4\left(\frac{\vphi+\eta}{N\eta}\right)^{5/2}. \label{eq:Q_0k_bdm1_repr}
        \end{align}
        For the $m_N$-term, Leibniz rule produces the three classes
        \begin{align}
                 & (\partial_{i\mu}^2m_N)GB\wt\Sg G,\qquad (\partial_{i\mu}m_N)\partial_{i\mu}(GB\wt\Sg G),\qquad m_N\partial_{i\mu}^2(GB\wt\Sg G). \nonumber
        \end{align}
        For the first two class, use \eqref{eq:avg_tr_der} for $\partial_{i\mu}^l m_N$ and \eqref{eq:gh_vec_bd} for $\norm{G \bbw}$.  The third class is controlled by \eqref{eq:Q_0k_bdm1_repr} and occur an extra factor of $(1+\phi)$. The the $\avg{\wt\Sg G B}_N$-term is treated analogously. Thus
        \begin{align}
                \eqref{eq:Q_0k_bdm1_start} \prec(1+\phi)^5\left(\frac{\vphi+\eta}{N\eta}\right)^{5/2}. \nonumber
        \end{align}
        Consequently,
        \begin{align}
                \eqref{eq:Q_cum_term} & \prec N^{-(k-1)/2}(1+\phi)^{4k} \left(\frac{\vphi+\eta}{N\eta}\right)^{k-2+5/2} \bb E\abs{Q_{\bbv\bbw}}^{2p-k}\nonumber            \\
                                                            & \prec(1+\phi)^{4k}\left(\frac{\vphi+\eta}{N\eta}\right)^{3k/2} \bb E\abs{Q_{\bbv\bbw}}^{2p-k}. \label{eq:Q_0k_bdm1}
        \end{align}

       For $T_{12}^{(3)}$ in \eqref{eq:Q_cum1}, we have
        \begin{align}
                 & \partial_{i\mu}(Q_{\bbv\bbw}^{p-1}\overline{Q_{\bbv\bbw}}^{p})=(p-1)(\partial_{i\mu}Q_{\bbv\bbw})Q_{\bbv\bbw}^{p-2}\overline{Q_{\bbv\bbw}}^{p} +p(\partial_{i\mu}\overline{Q_{\bbv\bbw}})Q_{\bbv\bbw}^{p-1}\overline{Q_{\bbv\bbw}}^{p-1}. \nonumber
        \end{align}
        By complex conjugation, it suffices to estimate
        \begin{align}
                 & \frac1{N^2}\sum_{i,\mu}\bb E\left( \abs{(\wh\Sg^{1/2}GB\wh\Sg^{1/2})_{\bbv i}} \abs{(G\wh\Sg^{1/2})_{\mu\bbw}} \abs{\partial_{i\mu}Q_{\bbv\bbw}} \abs{Q_{\bbv\bbw}}^{2p-2}\right). \label{eq:Q_cum_exp2}
        \end{align}
        By \eqref{eq:Q_der_contr}, we have
        \begin{align}
                T_{12}^{(3)} \prec(1+\phi)^2\left(\frac{\vphi+\eta}{N\eta}\right)^3 \bb E\abs{Q_{\bbv\bbw}}^{2p-2}. \label{eq:Q_cum2_der}
        \end{align}

        Finally, the remainder term $\calR_{l+1}^{(3)}$ is bounded by similar arguments to \eqref{eq:high_mom_rem_bd}, thus we omit the details. So we get that, for every fixed $D>0$,
        \begin{align}
                \calR_{l+1}^{(3)}=O(N^{-D}). \label{eq:Q_rem_bd}
        \end{align}
        
        Combining \eqref{eq:Q_cum_exp}, \eqref{eq:Q_cum_exp2},  \eqref{eq:Q_nonsat}, \eqref{eq:Q_1km1_bd} ,\eqref{eq:Q_0k_bdm1}, \eqref{eq:Q_cum2_der} and \eqref{eq:Q_rem_bd}, we obtain
        \begin{align}
                \bb E\abs{Q_{\bbv\bbw}}^{2p} \prec\sum_{n=1}^{2p}(1+\phi)^{4n}\left(\frac{\vphi+\eta}{N\eta}\right)^{3n/2} \bb E\abs{Q_{\bbv\bbw}}^{2p-n}+N^{-D}. \nonumber
        \end{align}
        By Young's inequality in \eqref{eq:anis_holder} and Markov's inequality imply
        \begin{align}
                \abs{Q_{\bbv\bbw}} \prec(1+\phi)^4\left(\frac{\vphi+\eta}{N\eta}\right)^{3/2}, \nonumber
        \end{align}
        which proves \eqref{eq:Q_anis_bd}.
\end{proof}

\subsection{Initial estimate for the local laws}\label{sub:init_est}

\begin{proof}[Proof of Lemma~\ref{lem:ll_init_est}]
        For $\eps_*\leq\eta\leq \eta_0$, by Theorem~\ref{thm:sharp_stab} and \eqref{eq:G_triv_bd}, we have $\|G\|+\norm{\Pi(z)}+\norm{(\calI-\calJ)^{-1}}\lesssim1$ and therefore $G-\Pi=O_\prec(1)$. Moreover, by \eqref{eq:RG_imag}, we have $\Im R_G(z)^{-1}\gtrsim 1$ and hence $\|R_G(z)\|\lesssim1$.  Theorem~\ref{thm:rand_sc_eq}, with $\phi=\vphi=1$, together with Lemma~\ref{lem:rand_inv_err} and \eqref{eq:proj_rand_eq}, yields
        \begin{align}
                \norm{\bbu-\calP[R_G-\Pi]}\prec N^{-1}. \nonumber
        \end{align}
        By \eqref{eq:RG_res_id},
        \begin{align}
                R_G-\Pi=\Pi\calE[\bbu]\Pi+\Pi\calE[\bbu]\Pi\calE[\bbu]R_G, \nonumber
        \end{align}
        and hence
        \begin{align}
                \calP[R_G-\Pi]=\calJ(z)\bbu+O(\norm{\bbu}^2). \nonumber
        \end{align}
        Therefore, taking a union bound over a deterministic net of the spectral domain,  we obtain
        \begin{align}
                \norm{\bbu}\prec  N^{-1}+\norm{\bbu}^2, \qquad z\in\{z=E+\ii\eta: |E-E_*|\leq \eps_*,\eps_*\leq\eta\leq\eta_0.\} \label{eq:init_proj_quad}
        \end{align}

        Moreover, by $G^{-1}=\calH+\La-I_z$, we have $\wt G(z)=(W-zI_N)^{-1}$ and $\wh G(z)+I_M=Y\wt G(z)Y^\top $.
        Hence, on the high probability event $\{\|Y\|\leq C\}$,  for $z=E+\ii\eta_0$, we have $\|\wt G(z)\|+\|\wh G(z)+I_M\|\lesssim\eta_0^{-1}$.
        On the other hand, the Stieltjes representations of the diagonal blocks of $\Pi$ give $\|\Pi_{22}(z)\|+\|\Pi_{11}(z)+I_M\|\lesssim\eta_0^{-1}$.
        Therefore, with high probability, we have
        \begin{align}
                \|\bbu(E+\ii\eta_0)\|\lesssim\eta_0^{-1}. \label{eq:init_proj_small}
        \end{align}
        Choosing $\eta_0$ sufficiently large, we have $\|\bbu(E+\ii\eta_0)\|\leq c_*/4$ with high probability. Hence \eqref{eq:init_proj_quad} implies $\|\bbu\|\prec N^{-1}$ at $\eta=\eta_0$. By a continuity argument, we obtain \eqref{eq:init_proj_law}.

        Applying Lemma~\ref{lem:ll_recon} with \eqref{eq:init_proj_law} proves the anisotropic estimate in \eqref{eq:init_ll}. Finally, taking a normalized trace in \eqref{eq:Delta_recon} and using \eqref{eq:RG_res_id}, \eqref{eq:err_psi_phi}, and Lemma~\ref{lem:rand_inv_err} gives the averaged estimate. This proves the lemma.
\end{proof}

\section{Auxiliary estimates for edge universality}
\label{app:edge_univ}

Throughout this appendix, Assumptions~\ref{ass:basic} and~\ref{ass:rand_ent} hold, and the rightmost edge $E_+$ is regular in the sense of Definition~\ref{def:reg_r_edge}. We prove the perturbative edge estimates and comparison bounds used in Section~\ref{sec:edge_univ}.

\subsection{Initial model}
\label{sub:init_model}

This subsection proves Propositions~\ref{prop:init_edge_reg} and \ref{prop:eta_star_reg}, stated in Subsection~\ref{subsec:gauss_div}.

\begin{proof}[Proof of Proposition~\ref{prop:init_edge_reg}]
        For $z_1,z_2\in\C^+$, since $\Sg_t=\Sg-tI_M$, the inverse MDEs for $\Pi_0(z_1)$ and $\Pi(z_2)$ give
        \begin{align}
                \Pi_0(z_1)^{-1} =\Pi(z_2)^{-1}-(z_1-z_2)\wt I -\calS[\Pi_0(z_1)-\Pi(z_2)] +t(m_0(z_1)\wh I+h_0(z_1)\wt I). \label{eq:Pi0_Pi_diff}
        \end{align}
        Since $\Sg_t$ satisfies Assumption~\ref{ass:basic} uniformly for all sufficiently large $N$, Lemma~\ref{lem:unif_apriori}(ii) and \eqref{eq:init_trace_id} give
        \begin{align}
                |m_0(z)|+|h_0(z)|\lesssim1 \label{eq:ord_edge_mh}
        \end{align}
        uniformly for $z\in \mb D(\tau,\eta_0)$.

        We first compare $\Pi_0$ and $\Pi$ near $E_+$.
        Set $\Delta_0(z) = \Pi_0(z) - \Pi(z)$ for $z\in\C^+$.
        Taking $z_1=z_2=z$ in \eqref{eq:Pi0_Pi_diff} and using the resolvent identity, we obtain
        \begin{align}
                \Delta_0(z) =\Pi(z)\calS[\Delta_0(z)]\Pi_0(z) -t\Pi(z)(m_0(z)\wh I+h_0(z)\wt I)\Pi_0(z). \nonumber
        \end{align}
        Substituting $\Pi_0(z)=\Pi(z)+\Delta_0(z)$ yields
        \begin{align}
                \calB(z)[\Delta_0(z)]
                = & \,\Pi(z)\calS[\Delta_0(z)]\Delta_0(z) -t\Pi(z)(m_0(z)\wh I+h_0(z)\wt I)\Pi(z) \nonumber \\
                  & -t\Pi(z)(m_0(z)\wh I+h_0(z)\wt I)\Delta_0(z).
                \label{eq:Delta0_eq_stab}
        \end{align}
        By Lemma~\ref{lem:bdry_zero_dens} and Theorem~\ref{thm:sharp_stab}, we have $ \|\Pi(z)\|\lesssim 1$ and $\|\calB(z)^{-1}\| \lesssim 1+(|E-E_+|+\eta)^{-1/2}$
        uniformly for $|E-E_+|\leq\eps_*, 0<\eta\leq\eta_0$.
        Combining this with \eqref{eq:ord_edge_mh}, we obtain
        \begin{align}
                \|\Delta_0(z)\| \lesssim  \left(1+(|E-E_+|+\eta)^{-1/2}\right) \left(t+\|\Delta_0(z)\|^2 +t\|\Delta_0(z)\|\right). \label{eq:ord_edge_stab1}
        \end{align}
        On $\eta=\eta_0$, we have $|E-E_+|+\eta_0\sim1$. By \eqref{eq:ord_edge_stab}, $\|\Pi(E+\ii\eta_0)\|+\|\calB(E+\ii\eta_0)^{-1}\|\lesssim1$. Since $\Sg_t$ satisfies Assumption~\ref{ass:basic} uniformly, Lemmas~\ref{lem:unif_apriori} and~\ref{lem:Pi_bd}, applied to $\Pi_0$, also give $\|\Pi_0(E+\ii\eta_0)\|\lesssim1$. Hence \eqref{eq:ord_edge_stab1} yields $\|\Delta_0(E+\ii\eta_0)\|\lesssim t+\|\Delta_0(E+\ii\eta_0)\|^2$.
        By continuity from $t=0$, where $\Delta_0=0$, we obtain $\|\Delta_0(E+\ii\eta_0)\|\lesssim t$ uniformly for $|E-E_+|\leq\eps_*$.
        Starting from $\eta=\eta_0$ and decreasing $\eta$, by a standard continuity argument with \eqref{eq:ord_edge_stab1}, for any sufficiently large fixed $C$,
        \begin{align}
                \|\Delta_0(z)\| \lesssim t\left(1+(|E-E_+|+\eta)^{-1/2}\right)
                \label{eq:ord_edge_stab2}
        \end{align}
        uniformly whenever $|E-E_+|\leq\eps_*, 0<\eta\leq\eta_0$ and $|E-E_+|+\eta\geq Ct$.
        
        We now locate the left side of the perturbed edge $E_{+,0}$. Since $t\ll 1$, we have $Ct<\eps_*$ for all sufficiently large $N$. Proposition~\ref{prop:sqrt_edge}, applied at the original edge $E_+$, gives $\Im m(E_+-Ct)\sim\sqrt{Ct}$.
        On the other hand, \eqref{eq:ord_edge_stab2} and Corollary~\ref{cor:bdry_reg_m_g}, after taking $\eta\downarrow0$, give $|m_0(E_+-Ct)-m(E_+-Ct)|  \lesssim t+C^{-1/2}\sqrt t$.
        Therefore, after choosing $C$ sufficiently large, we have $\Im m_0(E_+-Ct)>0$ for all sufficiently large $N$. This proves $E_+-Ct\in\supp\bnu_0$.
        Hence $E_{+,0}\geq E_+-Ct$.
        
        We next exclude spectrum to the right of $E_++Ct$. Since $E_+$ is the rightmost edge, when $E_++\eps_*\leq E\leq C_*+1, 0<\eta\leq\eta_0$, the Stieltjes representations gives $\Im m(E+\ii\eta)\sim\eta$ and $\norm{\Pi(z)}\lesssim 1$. Thus, by Lemma~\ref{lem:stab_red} and Proposition~\ref{prop:red_stab}, we have $\norm{\calB(z)^{-1}}+\norm{(\calI-\calJ(z))^{-1}}\lesssim 1$. By the same argument as in \eqref{eq:ord_edge_stab2}, the estimate \eqref{eq:ord_edge_stab2} holds uniformly for $E_++Ct\leq E\leq C_*+1, 0<\eta\leq\eta_0$.
        Moreover, for $E_++Ct\leq E\leq C_*+1$, by \eqref{eq:ord_edge_stab2} and let $\eta\downarrow 0$, we have
        \begin{align}
                \left(1+(E-E_+)^{-1/2}\right)\norm{\Delta_0(E)}^2\lesssim t\left(1+(E-E_+)^{-1/2}\right)^2\lesssim t+C^{-1}. \nonumber
        \end{align}
        We apply $\calB(E)^{-1}$ to \eqref{eq:Delta0_eq_stab} at $\eta=0$ and regard the resulting equation as a fixed-point equation for $\Delta_0(E)$. Let $K$ be a sufficiently large fixed constant. By \eqref{eq:ord_edge_mh}, the constant term in the fixed-point equation is bounded by $O\left(t\left(1+(E-E_+)^{-1/2}\right)\right)$.

        On the ball $ \left\{\Delta_0:\|\Delta_0\|\leq Kt\left(1+(E-E_+)^{-1/2}\right)\right\}$, the terms containing $\calS[\Delta_0]\Delta_0$ and $t(m_0\wh I+h_0\wt I)\Delta_0$ are bounded by $$O\left(K^2t^2\left(1+(E-E_+)^{-1/2}\right)^3+Kt^2\left(1+(E-E_+)^{-1/2}\right)^2\right).$$
        Hence, after fixing $K$ sufficiently large, we choose $C$ sufficiently large and then take $N$ sufficiently large so that the fixed-point map sends this ball into itself. Since $m_0$ and $h_0$ are linear functionals of $\Pi_0$, the difference of the fixed-point equations corresponding to any two points in this ball is bounded by $O\left(Kt\left(1+(E-E_+)^{-1/2}\right)^2\right)$
        times the distance between the two points.
        By the contraction mapping theorem, for every $E_++Ct\leq E\leq C_*+1$, the MDE associated with $\Sg_t$ has a unique solution satisfying $$\|\Delta_0(E)\|\leq Kt\left(1+(E-E_+)^{-1/2}\right).$$ Since $E>E_+$ and $E_+$ is the rightmost edge, $\Pi(E)$ is real. The MDE associated with $\Sg_t$ also has real coefficients at the real spectral parameter $E$. Hence the complex conjugate of the local solution is another solution in the same ball, and the local uniqueness implies that this solution is real. On the other hand, \eqref{eq:ord_edge_stab2} shows that, for $\eta>0$, $$\|\Delta_0(E+\ii\eta)\|\lesssim t\left(1+(E-E_++\eta)^{-1/2}\right).$$
        Letting $\eta\downarrow0$ and using the boundary continuity of $\Pi_0$ and $\Pi$, we obtain $$\|\Delta_0(E)\|\leq Kt\left(1+(E-E_+)^{-1/2}\right).$$  Hence $\Pi_0(E)$ belongs to the above ball. This coincides with the unique local solution $\Pi_0(E)$ and is real. Consequently, $\Im m_0(E)=0$ for $E_++Ct<E\leq C_*+1$. Consequently, $E_{+,0}\leq E_++Ct$. Together with the lower bound above, this proves $|E_{+,0}-E_+|\lesssim t$.

        It remains to verify the regularity of $E_{+,0}$. \eqref{eq:ord_edge_stab2} gives $|m_0(E_+-Ct)-m(E_+-Ct)|+|g_0(E_+-Ct)-g(E_+-Ct)|\lesssim \sqrt t$.
        Since $\Sg_t$ satisfies Assumption~\ref{ass:basic} uniformly, Corollary~\ref{cor:bdry_reg_m_g} applies to both $\Pi_0$ and $\Pi$ with uniform constants. Using $|E_{+,0}-E_+|\lesssim t$, we obtain
        \begin{align}
                |m_0(E_{+,0})-m(E_+)|
                +|g_0(E_{+,0})-g(E_+)|
                \lesssim t^{1/3}+\sqrt t\ll 1. \nonumber
        \end{align}
        Since $\Sg_t-\Sg=-tI_M$, it follows that
        \begin{align}
                \norm{(E_{+,0}+g_0(E_{+,0}))(I_M+m_0(E_{+,0})\Sg_t)-(E_++g(E_+))(I_M+m(E_+)\Sg)}\ll 1. \nonumber
        \end{align}
        Hence \eqref{eq:schur_nondeg}, we have
        \begin{align}
                s_{\min}\left(
                AA^\top-(E_{+,0}+g_0(E_{+,0}))
                (I_M+m_0(E_{+,0})\Sg_t)
                \right)\geq\frac{\tau_0}{2} \nonumber
        \end{align}
        for all sufficiently large $N$. Moreover, since $E_{+,0}=\sup\supp\bnu_0$, this shows that $E_{+,0}$ is a regular right edge for the MDE associated with $\Sg_t$ in the sense of Definition~\ref{def:reg_r_edge}. Proposition~\ref{prop:sqrt_edge}, applied to $\Pi_0$, gives $c_{{\rm edge},0}\sim1$ and $\gamma_{+,0}=(\pi c_{{\rm edge},0})^{2/3}\sim1$. Finally, since $\Sg_t$ satisfies Assumption~\ref{ass:basic} uniformly and the regularity constants at $E_{+,0}$ are uniform, the constants in Proposition~\ref{prop:sqrt_edge} and Theorem~\ref{thm:opt_ll} may also be chosen uniformly for $Y_0$.
\end{proof}

\begin{proof}[Proof of Proposition~\ref{prop:eta_star_reg}]
         By Proposition~\ref{prop:init_edge_reg}, $E_{+,0}$ is a regular right edge for $\Pi_0$ in the sense of Definition~\ref{def:reg_r_edge}. Hence Theorem~\ref{thm:sharp_stab}, applied to $\Pi_0$, gives
        \begin{align}
                \Im m_0(E+\ii\eta)\sim \begin{cases} \sqrt{E_{+,0}-E+\eta}, & E\leq E_{+,0}, \\[1mm] \dfrac{\eta}{\sqrt{E-E_{+,0}+\eta}}, & E\geq E_{+,0}, \end{cases} \quad |E-E_{+,0}|\leq\eps_*, \quad 0<\eta\leq\eta_0. \label{eq:det_edge_est}
        \end{align}
        Appling Corollary~\ref{cor:edge_rig} to $W_0$, yields
        \begin{align}
                \abs{\la_{1,0}-E_{+,0}}\prec N^{-2/3}\ll \eta_*. \label{eq:init_rig}
        \end{align}
        Hence, for $\eta\geq\eta_*$, $\abs{E-E_{+,0}}+\eta\sim\abs{E-\la_{1,0}}+\eta$ with high probability.
        Since $E_+\sim1$ and $\abs{E_{+,0}-E_+}\lesssim t$ by Proposition~\ref{prop:init_edge_reg}, we have $\la_{1,0}\sim1$ with high probability, and hence $\|W_0\|\lesssim1$ with high probability.

        Notice that $Y_0$ satisfying the assumptions of Theorem~\ref{thm:opt_ll}. Hence, by \eqref{eq:opt_avg_ll} and \eqref{eq:gap_avg_ll}, we have that for $E_{+,0}-\eps_*\leq E\leq E_{+,0}$ and $\eta_*\leq \eta\leq \eta_0$, 
        \begin{align}
                \abs{m_{0,N}(E+\ii\eta)-m_0(E+\ii\eta)}\prec\frac1{N\eta}, \label{eq:avg_star_in}
        \end{align}
        and for $E_{+,0}\leq E\leq E_{+,0}+\eps_*$ and $\eta_*\leq \eta\leq \eta_0$,
        \begin{align}
                \abs{m_{0,N}(E+\ii\eta)-m_0(E+\ii\eta)}\prec\frac1{N(\kappa+\eta)}+\frac1{(N\eta)^2\sqrt{\kappa+\eta}}, \qquad \kappa:=|E-\la_{1,0}|. \label{eq:avg_star_out}
        \end{align}
        Combining \eqref{eq:r_edge_rig}, \eqref{eq:det_edge_est}, \eqref{eq:init_rig}, \eqref{eq:avg_star_in} and \eqref{eq:avg_star_out}, we obtain that there exist constants $0<c_V<\eps_*/2$ and $C_V>0$ such that \eqref{eq:eta_star_reg_in} holds for $\la_{1,0}-c_V\leq E\leq\la_{1,0}$ and $\eta_*+\sqrt{\eta_*\abs{\la_{1,0}-E}}\leq\eta\leq 10$, while \eqref{eq:eta_reg_out} holds for $\la_{1,0}\leq E\leq\la_{1,0}+c_V$ and $\eta_*\leq\eta\leq 10$. The details can be found in \cite[Eqs. (A.26)--(A.29)]{ding2022tracywidom}.
        
        Thus, we have proved Proposition~\ref{prop:eta_star_reg}.
\end{proof}

\subsection{Deterministic and empirical edge parameters}
\label{sub:edge_par}

This subsection proves Proposition~\ref{prop:rect_flow_id}, Lemmas~\ref{lem:rect_crit} and \ref{lem:emp_subord} and Proposition~\ref{prop:cond_edge_par}, stated in Subsection~\ref{subsec:gauss_div}.

\begin{proof}[Proof of Proposition~\ref{prop:rect_flow_id}]
        The inverse MDE associated with $Y_0$ is
        \begin{align}
                -\Pi_0(\zeta)^{-1}=\begin{pmatrix} I_M+m_0(\zeta)\Sg_t & -A \\ -A^\top & (\zeta+g_0(\zeta))I_N\end{pmatrix}. \nonumber
        \end{align}
        Write $a=a_t(\zeta)$ and denote the right-hand side of \eqref{eq:rect_flow_blk} by $\Pi^{\rr{fc}}$. Then
        \begin{align}
                -(\Pi^{\rr{fc}})^{-1}=\begin{pmatrix} a_t^{-1}(I_M+m_0(\zeta)\Sg_t) & -A \\ -A^\top & a_t(\zeta+g_0(\zeta))I_N\end{pmatrix}. \label{eq:rect_flow_inv}
        \end{align}
        Since $a_t=1-tm_0(\zeta)$ and $\Sg_t=\Sg-tI_M$,
        \begin{align}
                a_t^{-1}(I_M+m_0(\zeta)\Sg_t)=I_M+a_t^{-1}m_0(\zeta)\Sg. \nonumber
        \end{align}
        Moreover, by \eqref{eq:rect_subord},
        \begin{align}
                z+a_t(g_0(\zeta)+th_0(\zeta))=a_t(\zeta+g_0(\zeta)). \nonumber
        \end{align}
        The diagonal traces of $\Pi^{\rr{fc}}$ are
        \begin{align}
                \avg{\Pi^{\rr{fc}}_{22}}_N=a_t^{-1}m_0(\zeta),\qquad \avg{\Sg\Pi^{\rr{fc}}_{11}}_N=a_t(g_0(\zeta)+th_0(\zeta)). \label{eq:rect_traces}
        \end{align}
        Hence \eqref{eq:rect_flow_inv} is the inverse MDE \eqref{eq:inverse_mde} at $z$.

        It remains to verify the positivity conditions in Proposition~\ref{prop:pi}. By \eqref{eq:rect_mfc}, \eqref{eq:rect_bt}, and \cite[Lemma~2(i)]{ding2022edge},
        \begin{align}
                \Im b_t(z)>0,\qquad \Im(zb_t(z))>0,\qquad \Im m_{\rm fc}(z)>0,\qquad \Im(zm_{\rm fc}(z))>0. \label{eq:rect_pos_scal}
        \end{align}
        By \eqref{eq:rect_subord_aux}, we have $\Im a_t^{-1}=t\Im m_{\rm fc}(z)>0$ and $\Im(za_t^{-1})=\Im z+t\Im(zm_{\rm fc}(z))>0$.
        By Proposition~\ref{prop:pi} applied to $\Pi_0$, there exists a positive semidefinite matrix-valued measure $V_{0,1}$ on $\R_+$ such that
        \begin{align}
                \zeta^{-1}\Pi_{0,11}(\zeta)=\int_{\R_+}\frac{V_{0,1}(\dd\la)}{\la-\zeta}. \nonumber
        \end{align}
        Since $\zeta=zb_t(z)a_t^{-1}$ by \eqref{eq:rect_subord_aux}, we obtain
        \begin{align}
                z^{-1}\Pi^{\rr{fc}}_{11}=b_t(z)\int_{\R_+}\frac{V_{0,1}(\dd\la)}{\la-\zeta},\qquad \Pi^{\rr{fc}}_{11}=zb_t(z)\int_{\R_+}\frac{V_{0,1}(\dd\la)}{\la-\zeta}. \nonumber
        \end{align}
        For every $\la\geq0$,
        \begin{align}
                \Im\frac{b_t(z)}{\la-\zeta}=\frac{\la\Im b_t(z)+|b_t(z)|^2\Im(za_t^{-1})}{|\la-\zeta|^2}>0,\qquad   \Im\frac{zb_t(z)}{\la-\zeta}=\frac{\la\Im(zb_t(z))+|zb_t(z)|^2\Im a_t^{-1}}{|\la-\zeta|^2}>0. \nonumber
        \end{align}
        Hence we have $\Im(z^{-1}\Pi^{\rr{fc}}_{11})>0$ and $\Im\Pi^{\rr{fc}}_{11}>0$.
        Together with \eqref{eq:rect_traces}, this gives
        \begin{align}
                \Im(-(\Pi^{\rr{fc}})^{-1})=\begin{pmatrix}\Im m_{\rm fc}(z)\Sg & 0 \\ 0 & \left(\Im z+\Im\avg{\Sg\Pi^{\rr{fc}}_{11}}_N\right)I_N\end{pmatrix}>0. \nonumber
        \end{align}
        Therefore $\Im\Pi^{\rr{fc}}=(\Pi^{\rr{fc}})^*\Im(-(\Pi^{\rr{fc}})^{-1})\Pi^{\rr{fc}}>0$. 
        Thus $\Pi^{\rr{fc}}$ satisfies both positivity conditions in Proposition~\ref{prop:pi}. Since $\Pi_0$ and the subordination function $\zeta$ are holomorphic on $\C^+$, the relations above show that $\Pi^{\rr{fc}}$ is holomorphic on $\C^+$. Hence, uniqueness in Proposition~\ref{prop:pi} yields $\Pi^{\rr{fc}}=\Pi$. This proves \eqref{eq:rect_flow_blk}, and \eqref{eq:rect_flow_scal} follows from \eqref{eq:rect_traces}.
\end{proof}

\begin{proof}[Proof of Lemma~\ref{lem:rect_crit}]
        Using \eqref{eq:init_trace_id} and $a_t'=-tm_0'$, we obtain
        \begin{align}
                \Phi_t'(\zeta) =a_t(\zeta)^2-tm_0'(\zeta)[2a_t(\zeta)\zeta-t(1-c_N)]. \label{eq:subord_der}
        \end{align}
        Appling Proposition~\ref{prop:sqrt_edge} to $\Pi_0$, and the Stieltjes representation of $m_0$ give for $s>0$ sufficiently small,
        \begin{align}
                m_0'(E_{+,0}+s)\sim s^{-1/2}, \qquad -m_0''(E_{+,0}+s)\sim s^{-3/2}. \label{eq:init_st_der}
        \end{align}
        Moreover, $a_t(\zeta)\sim1$ and $2a_t(\zeta)\zeta-t(1-c_N)$ is uniformly positive for $\zeta$ in a fixed neighborhood of $E_{+,0}$. Hence $\Phi_t'$ is negative at $\zeta-E_{+,0}=c_1t^2$ for sufficiently small $c_1>0$, and positive at $\zeta-E_{+,0}=c_2t^2$ for sufficiently large $c_2$. Differentiating \eqref{eq:subord_der} and using \eqref{eq:init_st_der} gives $\Phi_t''(\zeta)\sim t^{-2}>0$ uniformly for $c_1t^2\leq\zeta-E_{+,0}\leq c_2t^2$. Thus there is a unique critical point satisfying $\zeta_t-E_{+,0}\sim t^2$. The derivative estimates and the estimate for $\Phi_t''(\zeta_t)$ in \eqref{eq:crit_est} follow from \eqref{eq:init_st_der}. Proposition~\ref{prop:rect_flow_id} gives the identity $\Phi_t(\zeta_t)=E_+$ and the formula for $c_{\rm edge}$ in \eqref{eq:crit_est}.
\end{proof}

\begin{proof}[Proof of Lemma~\ref{lem:emp_subord}]
        Let $\la_{1,0}\geq\la_{2,0}\geq\dots\geq\la_{N,0}$ be the eigenvalues of $W_0$, and let $\gamma_{j,0}$ be the classical locations of $\rho_0$. Appling Corollary~\ref{cor:edge_rig} to $W_0$ in view of Proposition~\ref{prop:init_edge_reg} and Proposition~\ref{prop:sqrt_edge} give, uniformly for $\gamma_{j,0}\geq E_{+,0}-\eps_*/8$,
        \begin{align}
                \abs{\la_{j,0}-\gamma_{j,0}}\prec N^{-2/3}j^{-1/3},\qquad E_{+,0}-\gamma_{j,0}\sim\left(\frac jN\right)^{2/3},\qquad \gamma_{j,0}-\gamma_{j+1,0}\sim N^{-2/3}j^{-1/3}. \nonumber
        \end{align}
        For $x\in[\gamma_{j+1,0},\gamma_{j,0}]$, we have $\abs{\la_{j,0}-x}\prec N^{-2/3}j^{-1/3}$.  Moreover, since $\zeta-E_{+,0}\sim t^2$ and $t^2\gg N^{-2/3}$, we get $ \abs{\la_{j,0}-\zeta}\sim\abs{x-\zeta}\sim t^2+N^{-2/3}j^{2/3}$.
        
        Since $\int_{\gamma_{j+1,0}}^{\gamma_{j,0}}\rho_0(x)\dd x=N^{-1}$, we compare each eigenvalue with the integral over the corresponding quantile interval. For $\gamma_{j,0}<E_{+,0}-\eps_*/8$, both $\la_{j,0}$ and the corresponding quantile interval stay at a fixed positive distance from $\zeta$, so their total contribution is $O(1)$. For the remaining indices, the rigidity estimate gives, for every fixed $k\geq1$,
        \begin{align}
                \abs{\partial_\zeta^k(m_{0,N}(\zeta)-m_0(\zeta))}\prec\frac1N\sum_{\gamma_{j,0}\geq E_{+,0}-\eps_*/8}\frac{N^{-2/3}j^{-1/3}}{(t^2+(j/N)^{2/3})^{k+2}}+1. \nonumber
        \end{align}
        Splitting the sum at $j=Nt^3$ and using $\sum_{j\leq L}j^{-1/3}\lesssim L^{2/3}$ and $\sum_{j>L}j^{-(2k+5)/3}\lesssim L^{-(2k+2)/3}$, we obtain
        \begin{align}
                \frac1N\sum_{\gamma_{j,0}\geq E_{+,0}-\eps_*/8}\frac{N^{-2/3}j^{-1/3}}{(t^2+(j/N)^{2/3})^{k+2}}
                \lesssim\frac1{Nt^{2k+2}}. \nonumber
        \end{align}
        Therefore,
        \begin{align}
                \abs{\partial_\zeta^k(m_{0,N}(\zeta)-m_0(\zeta))}\prec\frac1{Nt^{2k+2}}+1,\qquad k\geq1. \nonumber
        \end{align}
        Taking $k=1,2$ proves \eqref{eq:emp_m0p} and \eqref{eq:emp_m0pp}.

        For the zeroth derivative, the averaged law on the gap side \eqref{eq:gap_avg_ll} at $\zeta+\ii\eta_*$ and $t^2\gg\eta_*$ gives, 
        \begin{align}
                \abs{m_{0,N}(\zeta+\ii\eta_*)-m_0(\zeta+\ii\eta_*)}\prec\frac1{Nt^2}+\frac1{(N\eta_*)^2t}. \nonumber
        \end{align}
                For $x\in[\gamma_{j+1,0},\gamma_{j,0}]$ with $\gamma_{j,0}\geq E_{+,0}-\eps_*/8$, direct subtraction gives
        \begin{align}
                \abs{\left(\frac1{\la_{j,0}-\zeta}-\frac1{\la_{j,0}-\zeta-\ii\eta_*}\right)-\left(\frac1{x-\zeta}-\frac1{x-\zeta-\ii\eta_*}\right)}
                \lesssim\eta_*\frac{\abs{\la_{j,0}-x}}{(t^2+(j/N)^{2/3})^3}. \nonumber
        \end{align}
        For $\gamma_{j,0}<E_{+,0}-\eps_*/8$, both $\la_{j,0}$ and the corresponding quantile interval stay at a fixed positive distance from $\zeta$, so their total contribution is $O(\eta_*)$. Hence, by rigidity and the same summation estimate as above with $k\geq1$,
        \begin{align}
                \abs{[m_{0,N}(\zeta)-m_0(\zeta)]-[m_{0,N}(\zeta+\ii\eta_*)-m_0(\zeta+\ii\eta_*)]}
                \prec\frac{\eta_*}{Nt^4}+O(\eta_*)\ll\frac1{Nt^2}+O(\eta_*). \nonumber
        \end{align}
       Combining the last two estimates proves \eqref{eq:emp_m0}.
\end{proof}

\begin{proof}[Proof of Proposition~\ref{prop:cond_edge_par}]
        Direct differentiation gives
        \begin{align}
                \Phi_{t,N}'(\zeta)
                 & =a_{t,N}(\zeta)^2-tm_{0,N}'(\zeta)[2a_{t,N}(\zeta)\zeta-t(1-c_N)], \nonumber\\
                \Phi_{t,N}''(\zeta)
                 & =-tm_{0,N}''(\zeta)[2a_{t,N}(\zeta)\zeta-t(1-c_N)]\nonumber\\
                 &\quad -4ta_{t,N}(\zeta)m_{0,N}'(\zeta)+2t^2\zeta m_{0,N}'(\zeta)^2. \nonumber
        \end{align}
        Uniformly for $\zeta-E_{+,0}\sim t^2$, Lemma~\ref{lem:emp_subord} and \eqref{eq:init_st_der} give $a_{t,N}(\zeta)\sim1$, $m_{0,N}'(\zeta)\sim t^{-1}$, and $-m_{0,N}''(\zeta)\sim t^{-3}$. Hence the second derivative above is positive and comparable to $t^{-2}$. Moreover, \eqref{eq:emp_m0p}, \eqref{eq:emp_m0}, and \eqref{eq:eta_t} give
        \begin{align}
              \abs{\Phi_{t,N}'(\zeta)-\Phi_t'(\zeta)}\prec\frac1{Nt^3}. \label{eq:emp_subord_der}
        \end{align}
        At $\zeta=E_{+,0}+c_1t^2$ and $\zeta=E_{+,0}+c_2t^2$, as shown in the proof of Lemma~\ref{lem:rect_crit}, $\Phi_t'$ has opposite signs and is bounded away from zero. By \eqref{eq:emp_subord_der}, the function $\Phi_{t,N}'(\zeta)$
        has the same signs at these two points with high probability. Moreover, \eqref{eq:emp_subord_der} and the lower bound for $\Phi_t''$ in the proof of Lemma~\ref{lem:rect_crit} imply $ \Phi_{t,N}''(\zeta)\gtrsim t^{-2}$ uniformly between these two points with high probability. Hence $\Phi_{t,N}'(\zeta)$ has a unique zero $\zeta_{t,N}$ satisfying $\zeta_{t,N}-E_{+,0}\sim t^2$. Corollary~\ref{cor:edge_rig} gives $\abs{\la_{1,0}-E_{+,0}}\prec N^{-2/3}\ll t^2$, proving the first estimate in \eqref{eq:cond_edge_est}.
        
        Since $\Phi_t'(\zeta_t)=0$ and $\Phi_{t,N}'(\zeta_{t,N})=0$, we have
        \begin{align}
                \Phi_{t,N}'(\zeta_t)= \Phi_{t,N}'(\zeta_t)-\Phi_{t,N}'(\zeta_{t,N}) =-\int_{\zeta_t}^{\zeta_{t,N}}\Phi_{t,N}''(s)\dd s. \nonumber
        \end{align}
        Therefore,  using  $\zeta_t-E_{+,0}\sim t^2$ , $\zeta_{t,N}-E_{+,0}\sim t^2$ and hence $\Phi_{t,N}''(s)\gtrsim t^{-2}$ on the interval between $\zeta_t$ and $\zeta_{t,N}$, we have $\abs{\Phi_{t,N}'(\zeta_t)}\gtrsim t^{-2}\abs{\zeta_{t,N}-\zeta_t}$. Combining this with \eqref{eq:emp_subord_der} gives $\abs{\zeta_{t,N}-\zeta_t}\prec(Nt)^{-1}$, proving the second estimate in \eqref{eq:cond_edge_est}.

        Using $h_{0,N}(\zeta)=\zeta m_{0,N}(\zeta)+1-c_N$, \eqref{eq:init_trace_id}, and \eqref{eq:emp_m0}, we obtain
        \begin{align}
                \abs{\Phi_{t,N}(\zeta_{t,N})-\Phi_t(\zeta_{t,N})}
                \prec t\eta_*+\frac1{Nt}+\frac1{(N\eta_*)^2}. \nonumber
        \end{align}
        Since $\Phi_t'(\zeta_t)=0$, \eqref{eq:crit_est} and the second estimate in \eqref{eq:cond_edge_est} yield $\abs{\Phi_t(\zeta_{t,N})-\Phi_t(\zeta_t)}\prec N^{-2}t^{-4}\ll(Nt)^{-1}$. Together with $\Phi_t(\zeta_t)=E_+$, this proves the third estimate in \eqref{eq:cond_edge_est}.

        It remains to compare the edge scales. From \eqref{eq:emp_m0p}, \eqref{eq:emp_m0}, \eqref{eq:init_st_der}, and the second estimate in \eqref{eq:cond_edge_est},
        \begin{align}
                \frac{m_{0,N}'(\zeta_{t,N})}{m_0'(\zeta_t)}
                =1+O_\prec\left(\frac1{Nt^3}\right),\qquad
                \frac{a_{t,N}(\zeta_{t,N})}{a_t(\zeta_t)}
                =1+O_\prec\left(\frac1{Nt^3}\right). \nonumber
        \end{align}
        The Stieltjes representation and the square-root behavior at $E_{+,0}$ give $\abs{m_0'''(\zeta)}\lesssim t^{-5}$ for $\zeta-E_{+,0}\sim t^2$, hence $\abs{\Phi_t'''(\zeta)}\lesssim t^{-4}$. Using $\abs{\Phi_t'''(\zeta)}\lesssim t^{-4}$ together with \eqref{eq:emp_m0p}, \eqref{eq:emp_m0pp}, and the estimates for $m_{0,N}'(\zeta_{t,N})/m_0'(\zeta_t)$ and $a_{t,N}(\zeta_{t,N})/a_t(\zeta_t)$ gives
        \begin{align}
               \Phi_{t,N}''(\zeta_{t,N})
                =\Phi_t''(\zeta_t)\left(1+O_\prec\left(\frac1{Nt^3}\right)\right). \nonumber
        \end{align}
        Combining these three comparisons with \eqref{eq:cond_edge_par}, the formula for $c_{\rm edge}$ in \eqref{eq:crit_est}, and \eqref{eq:edge_scale} gives
        \begin{align}
                \gamma_{t,N}=\gamma_+\left(1+O_\prec\left(\frac1{Nt^3}\right)\right). \nonumber
        \end{align}
        Since $\gamma_+\sim1$, this proves the last estimate in \eqref{eq:cond_edge_est}.
\end{proof}

\subsection{Short-time comparison estimates}
\label{sub:comp_est}

This subsection proves Lemmas~\ref{lem:comp_input} and \ref{lem:comp_high_cum}, which are used in Subsection~\ref{subsec:rm_gauss} to prove Proposition~\ref{prop:short_comp}.

\begin{lemma}\label{lem:comp_input}
        Along the interpolation \eqref{eq:cov_interp}, the estimates of Theorems~\ref{thm:rand_sc_eq} and~\ref{thm:opt_ll} hold uniformly for $0\leq s\leq t$ with the same deterministic solution $\Pi$. With $z_j=E_j+\ii\eta_{\rr{sm}}$ as in the proof of Proposition~\ref{prop:short_comp}, uniformly for deterministic unit vectors $\bbv,\bbw\in\C^L$,
        \begin{align}
                \abs{\scalar{\bbv}{(G(s,z_j)-\Pi(z_j))\bbw}} &\prec N^{-1/3+C_{\rr{cmp}}\chi},\qquad \abs{\scalar{\bbv}{G(s,z_j)\bbw}}\prec1, \nonumber\\
                \abs{\scalar{\bbv}{(G(s,z_2)-G(s,z_1))\bbw}} &\prec N^{-1/3+C_{\rr{cmp}}\chi}. \nonumber
        \end{align}
\end{lemma}

\begin{proof}
        Along the interpolation \eqref{eq:cov_interp}, the random part of $Y(s)$ is $(\Sg-sI_M)^{1/2}X+\sqrt{s}X^{\rr G}$. Since $t\ll 1$ and $\Sg\geq c_0I_M$, Assumption~\ref{ass:basic} gives $\|(\Sg-sI_M)^{1/2}\|+\|(\Sg-sI_M)^{-1/2}\|\lesssim1$ uniformly for $0\leq s\leq t$.

        Since $X$ and $X^{\rr G}$ are independent, the cumulant expansions in Appendix~\ref{app:rand_sc_eq} may be applied separately to their entries. Differentiation with respect to $X_{i\mu}$ acts in the direction $(\Sg-sI_M)^{1/2}\bbe_i$, while differentiation with respect to $X^{\rr G}_{i\mu}$ acts in the direction $\sqrt{s}\bbe_i$. By \eqref{eq:interp_cov_id}, their second-order contributions combine to $\sum_{i=1}^M(\Sg-sI_M)^{1/2}\bbe_i\bbe_i^\top(\Sg-sI_M)^{1/2}+s\sum_{i=1}^M\bbe_i\bbe_i^\top=\Sg$. Hence the second-order terms are the same as in the proof of Theorem~\ref{thm:rand_sc_eq}, and the self-energy operator remains $\cals$. For the higher-order terms, Lemma~\ref{lem:cum_bd} gives $\kappa_k(X_{i\mu})=O(N^{-k/2})$ for every fixed $k\geq3$, while all cumulants of $X^{\rr G}_{i\mu}$ of order at least three vanish. The uniform bounds above show that the derivative estimates and the remainder estimates used in Appendix~\ref{app:rand_sc_eq} remain uniform for $0\leq s\leq t$. Therefore Theorem~\ref{thm:rand_sc_eq} holds uniformly along the interpolation.

        The deterministic MDE is thus \eqref{eq:intro_mde} for every $0\leq s\leq t$, with the same solution $\Pi$, and the deterministic stability estimates used in the proof of Theorem~\ref{thm:opt_ll} are unchanged. Repeating that proof gives its estimates uniformly for $0\leq s\leq t$, which proves the first two estimates. Finally, Lemma~\ref{lem:bdry_zero_dens} and Theorem~\ref{thm:sharp_stab} imply $\|\Pi(z_2)-\Pi(z_1)\|\lesssim|z_2-z_1|^{1/2}\lesssim N^{-1/3+C_{\rr{cmp}}\chi}$. Combining this with the first estimate at $z_1$ and $z_2$ proves the third estimate.
\end{proof}

\begin{lemma}\label{lem:comp_high_cum}
        Under the assumptions of Proposition~\ref{prop:short_comp}, the remainder in \eqref{eq:comp_second} satisfies
        \begin{align}
                \sup_{0<s\leq t}\abs{\calR_s}\leq N^{1/6+C_{\rr{cmp}}\chi}. \label{eq:comp_cum_rem}
        \end{align}
\end{lemma}

\begin{proof}
        Let $G(s,z)$ be the linearized resolvent associated with $Y(s)$. By \eqref{eq:dir_res_der},
        \begin{align}
                \scrD_{\bbq\mu}G(s,z)=-G(s,z)\calV_{\bbq\mu}G(s,z),\qquad \partial_zG(s,z)=G(s,z)\wt I G(s,z). \nonumber
        \end{align}
        Repeated differentiation gives, for every fixed $k\geq1$,
        \begin{align}
                \scrD_{\bbq_1\mu_1}\dots\scrD_{\bbq_k\mu_k}G(s,z)=(-1)^k\sum_{\pi\in S_k}G(s,z)\calV_{\bbq_{\pi(1)}\mu_{\pi(1)}}G(s,z)\dots\calV_{\bbq_{\pi(k)}\mu_{\pi(k)}}G(s,z), \nonumber
        \end{align}
        where $S_k$ denotes the permutation group of $\{1,\ldots,k\}$, and $\pi(j)$ denotes the image of $j$ under $\pi$. 
        Set $z_j:=E_j+\ii\eta_{\rr{sm}}$, $j=1,2$. Since the lower right block of $G(s,z)$ is $(W(s)-zI_N)^{-1}$,
        \begin{align}
                \calX(s)=\int_{E_1}^{E_2}\Im\Tr\wt I G(s,x+\ii\eta_{\rr{sm}})\dd x. \nonumber
        \end{align}
        Cyclicity of the trace and $\partial_zG=G\wt I G$ give
        \begin{align}
                \scrD_{\bbq\mu}\calX(s)=-\Im\bigl(\Tr G(s,z_2)\calV_{\bbq\mu}-\Tr G(s,z_1)\calV_{\bbq\mu}\bigr). \label{eq:comp_end1}
        \end{align}
        For $k\geq2$, repeated differentiation of \eqref{eq:comp_end1} shows that $\scrD_{\bbq_1\mu_1}\dots\scrD_{\bbq_k\mu_k}\calX(s)$ is $(-1)^k$ times the imaginary part of the difference between the values at $z_2$ and $z_1$ of
        \begin{align}
                \sum_{\pi}\Tr G(s,z)\calV_{\bbq_{\pi(2)}\mu_{\pi(2)}}G(s,z)\dots\calV_{\bbq_{\pi(k)}\mu_{\pi(k)}}G(s,z)\calV_{\bbq_1\mu_1}, \label{eq:comp_end_mix}
        \end{align}
        where the sum runs over all permutations $\pi$ of $\{2,\dots,k\}$. Since each $\calV_{\bbq_j\mu_j}$ has rank at most two, every trace in \eqref{eq:comp_end_mix} is a finite sum of products of anisotropic resolvent entries, with an overall factor $\prod_j\|\bbq_j\|$.  Expanding the difference between the values at $z_2$ and $z_1$ by telescoping gives one factor $G(s,z_2)-G(s,z_1)$ in every term. Lemma~\ref{lem:comp_input} therefore yields, for $1\leq k\leq5$,
        \begin{align}
                \abs{\scrD_{\bbq_1\mu_1}\dots\scrD_{\bbq_k\mu_k}\calX(s)}\prec N^{-1/3+C_{\rr{cmp}}\chi}\prod_{j=1}^k\|\bbq_j\|. \label{eq:sm_count_der}
        \end{align}
        The chain rule and the boundedness of the derivatives of $f$ then give, for $1\leq k\leq5$,
        \begin{align}
                \abs{\scrD_{\bbq_1\mu_1}\dots\scrD_{\bbq_k\mu_k}f(\calX(s))}\prec N^{-1/3+C_{\rr{cmp}}\chi}\prod_{j=1}^k\|\bbq_j\|. \label{eq:test_der}
        \end{align}

        By \eqref{eq:comp_rem_def}, $\calR_s$ is the sum of the third and fourth cumulant contributions and the cumulant-expansion remainder. By Lemma~\ref{lem:cum_bd}, $|\kappa_r(X_{i\mu})|\lesssim N^{-r/2}$. Together with $MN\sim N^2$, \eqref{eq:test_der}, and Lemma~\ref{lem:stoch_dom_calc}, this gives $N^{1/6+C_{\rr{cmp}}\chi}$ and $N^{-1/3+C_{\rr{cmp}}\chi}$ for the third and fourth cumulant contributions, respectively, after absorbing the arbitrarily small loss from stochastic domination into $N^{C_{\rr{cmp}}\chi}$.
        For the cumulant-expansion remainder, apply Lemma~\ref{lem:cum_exp} with $\ell=3$ and cutoff $K=N^{-1/2+\chi}$. Since $5\chi<\tau$, we have $K\leq N^{-1/2+\tau/5}$. Hence \eqref{eq:interp_res} and \eqref{eq:interp_res_der} apply uniformly in the interpolation used in Lemma~\ref{lem:cum_exp}. These estimates remain valid with $\Sg^{1/2}$ replaced by $(\Sg-sI_M)^{1/2}$, since $\|(\Sg-sI_M)^{1/2}\|+\|(\Sg-sI_M)^{-1/2}\|\lesssim1$ uniformly for $0\leq s\leq t$. Thus \eqref{eq:test_der} also holds uniformly for the derivatives in the remainder term. Lemma~\ref{lem:stoch_dom_calc} then gives
        \begin{align}
                \sum_{i,\mu}\abs{\calR_{i\mu,N}(s)}\leq N^{-5/6+C_{\rr{cmp}}\chi}. \nonumber
        \end{align}
        Combining these three bounds and using \eqref{eq:comp_rem_def}, we obtain
        \begin{align}
                \abs{\calR_s}\leq N^{1/6+C_{\rr{cmp}}\chi}, \nonumber
        \end{align}
        which proves \eqref{eq:comp_cum_rem}.
\end{proof}

\end{appendix}

\end{document}